\font \z=msbm10 
\font \fraktur=eufm10 
\def \QQ{\hbox{\z \char "51}}
\def \RR{\hbox{\z \char "52}}       
\def \ZZ{\hbox{\z \char "5A}}
\def \CC{\hbox{\z \char "43}}      
\def \cal{\mathcal}
\def \beginaxioms{\begin{quotation}\begin{trivlist}\labelwidth 15pt\itemindent 0pt}
\def \endaxioms{\end{trivlist}\end{quotation}}
\def \theoremproven{\vspace{-0.2in} \qed}
\def \flexskip{\par}
\font\vsmall cmr10 scaled 700

\documentclass{memo-l}

\newtheorem{theorem}{Theorem}[chapter]
\newtheorem{lemma}{Lemma}[chapter]
\newtheorem{proposition}{Proposition}[chapter]

\theoremstyle{definition}

\theoremstyle{remark}
\newtheorem*{remark}{Remark}
\theoremstyle{corollary}
\newtheorem*{corollary}{Corollary}
\theoremstyle{remark}
\newtheorem*{example}{Example}
\theoremstyle{claim}
\newtheorem*{claim}{Claim}

\numberwithin{section}{chapter}
\numberwithin{equation}{chapter}
\numberwithin{figure}{chapter}
\numberwithin{table}{chapter}

\begin{document}
\frontmatter
\title{Some Generalized Kac-Moody Algebras \\ With Known Root Multiplicities}

\author{Peter Niemann}
\address{Department of Pure 
Mathematics and Mathematical Statistics, University of Cambridge, 
16 Mill Lane, CB2 1SB, UK}
\curraddr{Logica UK Ltd., 75 Hampstead Road, London NW1 2PL, UK}
\email{niemannp@logica.com}
\thanks{The author was supported by the Science and Engineering 
Research Council (UK), and Peterhouse, Cambridge.}


\date{1 July, 1997, revised 20 July 1999}

\begin{abstract}
Starting from Borcherds' fake monster Lie algebra we construct
a sequence of six generalized Kac-Moody algebras whose denominator formulas,
root systems and all root multiplicities can be described explicitly. 
The root systems decompose space into convex holes, of finite and affine
type, similar to the situation in the case of the Leech lattice. As a
corollary, we obtain strong upper bounds for the root multiplicities of
a number of hyperbolic Lie algebras, including $AE_3$. 
\end{abstract}

\maketitle

\tableofcontents

\setcounter{page}{5}

\mainmatter
\chapter*{Introduction}

In recent years the area of infinite-dimensional Lie algebras has attracted
considerable attention because of its numerous connections with other topics
in mathematics and, not least, its importance in theoretical
physics. The state space of physical theories will sometimes be a
representation space of an infinite-dimensional Kac-Moody algebra.
\flexskip

Surprisingly little is known about many obvious problems associated with such
Lie algebras. One of these questions regards their root multiplicities, that
is the dimensions of their root spaces. The cases of finite and affine Lie
algebras are fully understood. Let us turn to more general Lie algebras of
indefinite type which allow roots of negative norm. The problem now becomes far
more complicated and only very partial answers are known. As \cite {Kac90} remarks,
the multiplicities of all roots of an indefinite-type Kac-Moody algebra are not
known explicitly in any single case. At the same time, numerical calculations
yield intriguing results, in particular for the simplest class of such Lie
algebras which are called hyperbolic. They are defined by the condition on
their Dynkin diagram that every subdiagram be of finite or affine type. \flexskip

There are various kinds of results known about root multiplicities. There are
global upper bounds for all Lie algebras which work well in some cases but are
useless in others. 
There are recursive formulas which are useful for numerical calculations,
notably the results of Kac and Peterson \cite {KP83}, and of Berman and Moody
\cite {BM79}. Furthermore, there are explicit results for some low level roots 
of selected algebras such as the treatments of Feingold and Frenkel \cite {FF83}
(on the algebra denoted $AE_3$ in the notation of \cite {Kac90}), and
Kac, Moody, and Wakimoto \cite {KMW88} (on $E_{10} = T_{7,3,2}$). Recently, 
Gebert and Nicolai \cite{GN97} produced some intriguing numerical results 
on the simple roots of $E_{10}$ for level 2 and 3.\flexskip

Borcherds introduces in \cite {Bor92} the notion of generalized Kac-Moody algebras.
These differ from ordinary Kac-Moody algebras in that they may possess
imaginary simple roots. Generalized Kac-Moody algebras sometimes form the space
of states for quantized chiral strings.
It was a remarkable achievement that Borcherds then managed to construct a
generalized Kac-Moody algebra, called the fake monster Lie algebra, and to
give explicit root multiplicities for all its roots. The set of its
roots can be identified with the 26-dimensional even unimodular lattice
$II_{25,1} = \Lambda \oplus II_{1,1}$. Here, $\Lambda$ stands for the
24-dimensional Leech lattice and $II_{1,1}$ denotes the unique even
2-dimensional unimodular Lorentzian lattice. \flexskip

Borcherds suggested in \cite {Bor92} that the fake monster Lie algebra might only
be one of a whole class of generalized Kac-Moody algebras whose root
multiplicities can be described explicitly. Let $N$ be such that $N+1$ divides
24, let $M=\frac{24}{N+1}$. Thus $N$ will be one of 2, 3, 5, 7, 11, or 23.
Let $\sigma$ be an automorphism of order $N$, cycle
shape $1^MN^M$, of the Leech lattice.
The aim of this thesis is to prove some of Borcherds' conjectures by
constructing a series of generalized Kac-Moody
algebras ${\cal G}_N$ whose systems of simple roots can be identified with
the fixed point lattices $\Lambda^\sigma$ and certain elements of the dual
lattices ${\Lambda^\sigma}^*$. (Note that, unlike $\Lambda$, $\Lambda^\sigma$ is
no longer self-dual.) The Weyl denominator formula then becomes particularly
simple and allows us to calculate the root multiplicities for these generalized
Kac-Moody algebras explicitly.  \flexskip

We will now give a brief survey of the contents of the 6 chapters of
this thesis. 
Chapter 1 recalls the basic definition and construction of the fake monster Lie
algebra. Following \cite {Bor92} we introduce the notion of twisted denominator
formulas. They are obtained from the Weyl denominator formula of the fake
monster Lie algebra by the action of the Leech lattice automorphism $\sigma$.
We then proceed to formulate the main theorem of this work which
claims the existence of a series of generalized Kac-Moody algebras and states
their root multiplicities explicitly. Following ideas of \cite {Bor92} the proof is
reduced to an equality between two modular forms.
One of these is the $\theta$-function of the lattice ${\Lambda^\sigma}^*$. The
other function is derived from the Dedekind $\eta$-function.  \flexskip

Both the above functions are modular forms with respect to some subgroup of the
modular group, which has finite index.
The strategy of the proof must now be as follows. We check that
the modularity properties of both functions are equal and that a sufficient
number of leading coefficients in a Laurant expansion around zero coincide.
As the space of modular forms with respect to such a modular subgroup is
finite-dimensional this shows the claimed equality. It is easy to see that,
the larger the transformation group is, the smaller will be the dimension of
the space of modular forms, thus requiring fewer leading coefficients in the
remainder of the argument. The most suitable group for our purposes turns out
to be $\Gamma_0(N)$.\flexskip

Chapter 2 recalls the basic notions of modular group and modular form, as well
as the Dedekind $\eta$-function. We then proceed to define the specific modular
forms we need in this work and establish their precise modularity properties
under the elements of $\Gamma_0(N)$.  \flexskip

Chapter 3 recalls some basic results about the lattices in question and then
determines the exact transformation properties, including characters,
of their $\theta$-functions under the elements of $\Gamma_0(N)$.  \flexskip

In chapter 4 we determine the leading coefficients of the above modular forms.
This is straightforward for those modular forms which were
derived from the $\eta$-function. For $\theta$, however, this corresponds to the
enumeration of the short vectors of ${\Lambda^\sigma}^*$. For $N=2$ and $N=3$
it is well known that the fixed point lattices are the Barnes-Wall lattice
$\Lambda_{16}$ and the Coxeter-Todd lattice $K_{12}$. They have been considered
in the literature in their own right. For the remaining $N$ we devise a
strategy to count the short vectors. This makes use of the fact that the
lattices in question are all induced from the Leech lattice which in turn is
built upon the 24-dimensional Golay code. It is possible to reduce the
required enumeration in ${\Lambda^\sigma}^*$ to one of Golay code elements.
\flexskip

We then put everything together and obtain the desired equality, thus
completing the proof of the main theorem of this work. We thus have proven the
existence of a series of generalized Kac-Moody algebras ${\cal G}_N$ and have
determined their root multiplicities explicitly. At the same time, their Weyl
denominator formulas can be interpreted as new combinatorial identities,
similar to the Macdonald identities.
\flexskip

In the second part of this work we investigate the generalized
Kac-Moody algebras ${\cal G}_N$ in more detail. In chapter 5 we determine their
simple roots. We then identify this set with a set ${\cal R}$ consisting
of the elements of the fixed point lattice $\Lambda^\sigma$
and some elements of its dual lattice ${\Lambda^\sigma}^*$. It is well known
(see e.g. \cite {CS88}, chapter 25) that the elements of the Leech lattice generate a
decomposition of space into convex holes of radius less or equal to $\sqrt 2$.
Now the elements of the Leech lattice can be identified with the simple roots
of the fake monster Lie algebra. It turns out that holes of radius $\sqrt 2$
correspond to affine subalgebras, whereas holes of radius less than $\sqrt 2$
correspond to subalgebras of finite type. We generalize this to the sets ${\cal
R}$. The main difference is that now the algebras have simple roots of
two different lengths. We therefore have to generalize the notion of radius and
centre accordingly. If we do so it remains true that the decomposition
produces holes of (generalized) radius less or equal $\sqrt 2$, corresponding
to subalgebras of finite and affine type respectively. We further show that all
subalgebras of ${\cal G}_N$ of finite or affine type can be identified as
(generalized) holes in ${\cal R}$. As a corollary, we can relate the covering
radius of the fixed point lattices $\Lambda^\sigma$ to the covering radius of
the Leech lattice. \flexskip

One of our aims in chapter 6 is the complete classification of all
(generalized) holes of ${\cal R}$. The remainder of chapter 5 therefore develops
a number of techniques which will enable us to carry out the decomposition and
check its correctness. One such check is the volume formula which simply states
that the sum of the volumes of the individual holes must be the total volume of
space. We therefore determine the volumes of any finite and affine holes.
Another test concerns the automorphism groups of the individual holes, related
to the automorphism group of the fixed point lattice as a whole. We derive a
number of technical results, relating the generalized holes of ${\cal R}$ to
the (known) holes in the decomposition of the Leech lattice. \flexskip

Chapter 6 carries out the classification. This presents no theoretical
problems but involves long and repetitive calculations
which may only be carried out by computer.
The main purpose of the first two
sections of chapter 6 is to demonstrate how the various results of chapter 5
come together to achieve the classification and why the output of a computer
program constitutes a mathematical proof. To this end, we give the explicit
calculations of the two simplest cases, that is $N=23$ and $N=11$. The case
$N=23$ is almost trivial but very useful as it is 2-dimensional and so helps to
visualize the problem. The 4-dimensional case $N=11$ is already sufficiently
general to demonstrate how the program acts. We can therefore restrict
ourselves to giving the input data for the remaining cases. The complete
classification is given as appendix A. For the most complicated case $N=2$
there are 475 different types of (generalized) holes in ${\cal R}$.  \flexskip

As a corollary to this decomposition we can now identify all hyperbolic
subalgebras of the ${\cal G}_N$. The known root multiplicities of ${\cal G}_N$
then form upper bounds for the root multiplicities of these hyperbolic Lie
algebras. There are many interesting
examples for which the upper bounds obtained by the above technique improve on
the upper bounds given in the existing literature. Of particular interest is
the hyperbolic Lie algebra $AE_3$ which was also
investigated in \cite {FF83}. The algebra $AE_3$ is defined by the following Cartan
matrix: 
\[
\begin{pmatrix} 2&-2&0 \\ -2&2&-1 \\ 0&-1&2 \end{pmatrix}.
\]
We obtain strong upper bounds for the root multiplicities of this algebra.
These bounds are close to some existing conjectures of \cite {Kac90} (see exercise
13.37). Our results explain why the multiplicities of the roots in the
algebra $AE_3$ are often equal to the values of a partition function
$p_n (1-r^2/2)$.\flexskip

However, for other hyperbolic Lie algebras our new upper bounds are not
always sharp or at least an improvement on existing ones. In section 6.2 we
list the successful cases. In appendix B we provide some numerical data for
all hyperbolic Lie algebras of rank 7 to 10. The results for the algebra 
$T_{4,3,3}$ are of particular interest as they show that there are roots
$r$ of multiplicities both larger and smaller than the partition function
$p_6 (1-r^2/2)$. In the concluding section 6.3
we derive some conditions which are necessary if we want to create
useful upper bounds. This will then enable us to conjecture how far the
strategy of this work may be extended to other generalized Kac-Moody algebras.
\flexskip\flexskip

{\bf Acknowledgements:} I wish to thank Richard Borcherds, my research 
supervisor, for suggesting many of the problems discussed in this work
and for the advice and encouragement he provided throughout my research.
Furthermore, I would like to thank Simon Norton for his patient
explanations of the ATLAS, and Elizabeth Jurisich for clarifications
regarding the nature of specializations. Finally, I would like to thank
the referee for valuable comments on the first version of this paper.

\chapter{Generalized Kac-Moody Algebras} \flexskip

The main object of study in this work is a series of generalized Kac-Moody
algebras whose root multiplicities will be determined explicitly. This chapter
introduces the concept and describes the construction of these generalized
Kac-Moody algebras. Chapters 2-4 will then complete the proof of their
existence and their root multiplicity formulas.
It must be emphasised in this context that chapters 2-4 are independent
of the results of the present chapter 1. Chapter 1 has been placed in its
position ahead of chapters 2-4 only to provide the framework identifying
which auxiliary calculations are required in chapters 2-4.
Sections 1.1 and 1.2 recall the definition and elementary properties of
generalized Kac-Moody algebras, and specifically their Weyl denominator formula.
The series of generalized Kac-Moody algebras in this work will be derived
from the fake monster Lie algebra as introduced in \cite {Bor90b}. This, in turn,
is constructed from subspaces of a certain vertex algebra. In section 1.3 we
therefore briefly recall the construction of the vertex algebra in question.
Section 1.4 outlines the construction and some properties of the fake
monster Lie algebra, following \cite {Bor92}. Its root lattice can be identified
with the Leech lattice and all root multiplicities are known explicitly.
Section 1.5 quotes \cite {Bor92} to recall how automorphisms of the Leech lattice
can be applied to the denominator formula of the fake monster Lie algebra,
resulting in new identities which may be regarded as `twisted'
denominator formulas.
In sections 1.6 and 1.7 we restrict our attention to a specific series of
six Leech lattice automorphisms. For these we show how the new `twisted'
identities can be interpreted as
denominator formulas of six new generalized Kac-Moody algebras whose root
multiplicities can again be calculated explicitly.
For general automorphisms, the new algebras will not necessarily be
generalized Kac-Moody algebras but may be Lie superalgebras where
negative root multiplicities may occur. Section 1.7, in conjunction with
chapters 2-4, establishes an explicit form of the denominator formulas
(equation 1.26), for the series of six new generalized Kac-Moody
algebras, corresponding to the series of Leech lattice automorphisms. This
yields an explicit (non-recursive) formula for the root multiplicities of
the generalized Kac-Moody algebras (corollary to theorem 1.7), and carries
out a programme suggested in \cite {Bor92} (14. Examples 1 and 2).

\section{Definition and Fundamental Properties}

This section establishes the definition of generalized Kac-Moody algebras
which we will use within this work. Numerous variants have been proposed
in recent publications. We will follow the approach of \cite {Jur98}, but
specialize to the variant considered in \cite {Bor92} when appropriate. Throughout
this work we will use the initials GKM for generalized Kac-Moody algebra.

\subsection {The Definition}

\cite {Jur98} defines GKMs through generalized Cartan matrices.
A real matrix $C=(c_{ij})$,
$i,j$ in some index set $I$ (possibly countably infinite),
shall be called a generalized Cartan matrix if it
satisfies the following conditions:

\beginaxioms

\item[\bf(C1)] $C$ is symmetric;
\item[\bf(C2)] $c_{ij}\le0$ if $i\not=j$;
\item[\bf(C3)] if $c_{ii}>0$ then $2c_{ij}/ c_{ii} \in \ZZ$
                for all $j\in I$.
\endaxioms

\cite {Jur98} defines the generalized Kac-Moody algebra (GKM) $G=G(C)$ 
associated with this generalized Cartan matrix
to be the Lie algebra generated by elements $e_i$,
$f_i$, $h_i$, for $i \in I$, where the generators
satisfy the following relations for $i,j,k \in I$:

\beginaxioms
\item[\bf(R1)] $[h_i, h_j] = 0$
\item[\bf(R2)] $[e_i, f_j] = \delta_i^j h_i$ 
\item[\bf(R3)] $[h_i, e_k] = c_{ik}e_k,
    [h_i, f_k] = - c_{ik}f_k$
\item[\bf(R4)] If $c_{ii}>0$ and $i\not= j$ then ${\rm ad}(e_i)^ne_j= 
    {\rm ad}(f_i)^nf_j =0$, where $n=1-2c_{ij}/c_{ii}$.
\item[\bf(R5)] If $c_{ii}\le 0$, $c_{jj}\le 0$, and $c_{ij}=0$ then 
    $[e_i,e_j] = [f_i, f_j]= 0$. 
\endaxioms

We follow \cite {Jur98} in restricting the definition to symmetric Cartan
matrices as we will not encounter any non-symmetrisable Cartan matrices 
in the context of the present work. \flexskip

As in the finite dimensional case,
the elements $h_i$, $i\in I$, span an abelian subalgebra $H$ of $G(C)$,
called its Cartan subalgebra. Let $E$ be the subalgebra
generated by the $e_i$, $i\in I$, and $F$ be the subalgebra generated by
the $f_i$, $i\in I$. Then the GKM $G(C)$ has the triangular decomposition
\cite {Jur96}
\begin{equation}
G(C) = E \oplus H \oplus F.
\end{equation}
Every non-zero ideal of $G(C)$ has non-zero intersection with $H$.
The centre of $G(C)$ is contained in $H$. \cite {Jur96}

\subsection{Roots, Central Ideals, Central Extensions} 

Roots of finite dimensional Lie algebras are commonly defined as elements
of the dual space of the Cartan subalgebra. Hence, quotienting out some
ideal of the Cartan subalgebra, or extending the Cartan subalgebra centrally
will also affect the space where roots are defined. For the definition of
roots for GKMs, there consequently exist various options which relate
to the choice of central extension for the GKM defined above. 

Alternatively, the approach of \cite {Bor92} defines roots as a free
abelian group (of abstract symbols), and identifies a natural homomorphism
to the elements $h_i$ of the Cartan subalgebra $H$. 

If we consider roots within the dual space $H^*$ with $G(C)$
constructed as above, the simple roots
will not necessarily be linearly independent. \cite {Jur96} therefore extends the
GKM $G(C)$ by an algebra of `degree derivations'. This increases the
dimension of the Cartan subalgebra and ensures that the simple roots will be
defined linearly independent. Note that, besides the approaches of \cite {Jur96}
and \cite {Bor92}, others have been proposed, such as the approach described
in \cite {Kac90} whose `realization' of the Cartan matrix again guarantees the
linear independence of the resulting simple roots. \flexskip

As is discussed in detail in \cite {Jur98}, the denominator formula for an
arbitrary Cartan matrix with linearly dependent simple roots will not
necessarily be well defined as some terms may be infinite.
A general theory can therefore only be formulated in a framework which
guarantees the linear independence of the simple roots. This will either
be achieved as in \cite {Bor92}, or through a suitable extension of the Cartan
subalgebra as in \cite {Jur96}. Nevertheless, it may be possible to define
denominator formulas for certain other Lie algebras, through a process
called `specialization' in \cite {Jur96}. Such specializations will quotient out
subalgebras of the (extended) Cartan subalgebra. This operation will also
quotient out the corresponding subspaces of the dual space where roots are
defined. As is pointed out in the remark following definition 3 of \cite {Jur98},
specializations are valid, as long as they are well defined.
For the GKMs constructed in the present work we will find that the
denominator formula remains well defined under specialization. \flexskip

Where simple roots are linearly independent they obviously have multiplicity
one. The process of specialization may map multiple,
distinct simple roots to the same image. Therefore, following
specialization, some simple roots may have multiplicity greater than one.

\subsection*{The approach of Jurisich} 

\cite {Jur96} extends the GKM by all `degree derivations'.
Let ${\rm deg}(e_i) = - {\rm deg}(f_i) = (0,\dots,0,1,0,\dots)$ where 1
appears in the $i^{\rm th}$ position, and let ${\rm deg}(h_i) = (0,\dots)$.
Degree derivations $d_i$ are defined by letting $d_i$ act on the degree
$(n_1,n_2,\dots)$ subspace of the GKM as multiplication by the scalar $n_i$.
\cite {Jur96} defines the extended Lie algebra $G^e=G^e(C)$ as the semidirect
product of the GKM $G(C)$ with the space of all degree derivations. The
extension is central. Let $H^e$ denote the Cartan subalgebra of $G^e(C)$.
The extended Lie algebra has the decomposition
\begin{equation}
G^e(C) = E \oplus H^e \oplus F.
\end{equation}
Roots of $G^e(C)$ defined
in the space ${H^e}^*$ are necessarily linearly independent. \flexskip

The extended Lie algebra $G^e(C)$ is clearly not identical to the GKM $G(C)$
constructed from the original Cartan matrix $C$, but may be considered
naturally associated with it. More generally, if a Lie algebra can be
mapped to $G(C)$ modulo some central ideals, some central extensions, or
outer derivations (as is the case in specializations), we will consider
that Lie algebra naturally associated with the Cartan matrix $C$ and the
GKM $G(C)$. \flexskip

Having selected the Cartan algebra $H^e$, roots may be defined as in the
finite dimensional case. For $r \in {H^e}^*$, let 
\[
G^r =\{\  x\in G \ \vert \ [h,x] = r(h)x{\rm\ for\ all\ }h\in H^e\}.
\]
The roots of $G$ are the non-zero elements $r$ of ${H^e}^*$ such that
$G^r \not= 0$. The elements $r_i\in {H^e}^*$ such that the
generators $e_i$ are in $G^{r_i}$ are called simple roots. $G^r$ is the
root space of $r\in {H^e}^*$. A root $r$ is called positive if it is the
sum of simple roots, and negative otherwise.

\subsection*{The approach of Borcherds}

\cite {Bor92} defines the root lattice of $G(C)$ as the free abelian group
generated by elements $r_i$, for $i\in I$, with the bilinear form given by 
$(r_i, r_j) = c_{ij}$. Here, the $c_{ij}$ are the elements of the symmetric
generalized Cartan matrix $C$. The elements $r_i$ correspond to the simple
roots. The GKM is graded by the root lattice if we let $e_i$ have degree
$r_i$ and $f_i$ have degree $-r_i$. If $r$ is in the root lattice then the
vector space of elements of the Lie algebra of that degree is called the
root space of $r$. There is a natural homomorphism of
abelian groups from the root lattice to the Cartan subalgebra $H$ taking
$r_i$ to $h_{i}$ which preserves the bilinear forms. This homomorphism will
not necessarily be injective. If, for example, the null space of the
bilinear form has been quotiented out, this may have introduced relations
among the elements $h_i$ of $H$.

\subsection{Norm, Weyl Vector, Weyl Chamber and Cartan Involution}

We define the norm of an element $r$ of the
root lattice as the scalar product 
\begin{equation}
{\rm norm}(r)=(r,r).
\end{equation}
Note that we use the square of the standard norm, because the scalar 
product is not positive definite. A root of a GKM is called real if 
it has positive norm $(r,r)$, and imaginary otherwise. \flexskip

In the framework of the extended root space ${H^e}^*$, the Weyl vector
$\rho$ is any element of ${H^e}^*$ which satisfies $(\rho,r_i)=
-(r_i,r_i)/2$. Note that this implies $\rho(h_i)=-c_{ii}/2$. The
(fundamental) Weyl chamber is the set of vectors $v$ of the root space
${H^e}^*$ which satisfy $(v, r_i)\le 0$ for all real simple roots $r_i$.
Equivalently, \cite {Bor92} defines the Weyl vector as the additive map from the
free abelian group of roots to $\RR$, taking $r_i$ to $-(r_i,r_i)/2$ for
all $i\in I$. The (fundamental) Weyl chamber is the set of all vectors
$v\in H$ with $(v,h_i)\le 0$ for all $h_i \in H$ that correspond
to real simple roots. Following \cite {Bor90b}, we define the height of a root
$r$ as $-(\rho,r)$. \flexskip

Note that \cite {Bor92} uses non-standard sign conventions in this place, and the 
present paper will follow his conventions. In the case of a specialization,
where simple roots may be linearly dependent,
there is no reason why a Weyl vector should exist in general.
In the example of the affine algebra $A_2$, the relation $r_2 = -r_1$ of
the two simple roots (which holds in unextended 2-dimensional dual space)
shows that no Weyl vector can exist. \flexskip

$G(C)$ has an involution $\omega$ with $\omega(e_i) = -f_i$,
$\omega(f_i)= -e_i$, called the Cartan involution. There is a unique invariant
bilinear form $( {^.}, {^.})$ on $G(C)$ such that $(e_i,f_i)=1 $ for all $i$, 
and it also has the property that $-(g,\omega(g)) >0$ whenever $g$ is a 
homogeneous element of non-zero degree.

\subsection{Universal Central Extension}

Let $C$ be a generalized Cartan matrix, satisfying conditions (C1) to (C3).
For an alternative characterisation of some GKMs,
\cite {Bor92} defines the universal generalized Kac-Moody algebra (UGKM) $U(C)$
of this matrix to be the Lie algebra generated by elements $e_i$, 
$f_i$, $h_{ij}$, for $i,j \in I$ satisfying the following relations
($i,j,k,l\in I$):\flexskip

\beginaxioms
\item[\bf(U1)] $[h_{ij}, h_{kl}] = 0$ 
\item[\bf(U2)] $[e_i, f_j] = h_{ij}$
\item[\bf(U3)] $[h_{ij}, e_k] = \delta_i^j c_{ik}e_k,
    [h_{ij}, f_k] = -\delta_i^j c_{ik}f_k$
\item[\bf(U4)] If $c_{ii}>0$ and $i\not= j$ then ${\rm ad}(e_i)^ne_j= 
    {\rm ad}(f_i)^nf_j =0$, where $n=1-2c_{ij}/c_{ii}$.
\item[\bf(U5)] If $c_{ii}\le 0$, $c_{jj}\le 0$, and $c_{ij}=0$ then 
    $[e_i,e_j] = [f_i, f_j]= 0$.
\endaxioms

The UGKM is therefore an extension of the GKM defined in section 1.1.1
above. It is extended precisely by the additional central generators
$h_{ij}, i \ne j$. \cite {Bor92} observes that
the element $h_{ij}$ is 0 unless the $i$'th and $j$'th column of the Cartan
matrix $C$ are equal. The elements $h_{ij}$ for which the $i$'th and $j$'th 
column are equal form a basis of the Cartan subalgebra $H$ of $U(C)$.
In the case of ordinary Kac-Moody algebras, the $i$'th and
$j$'th column of $C$ cannot be equal unless $i=j$, so the only non-zero 
elements $h_{ij}$ are those of the form $h_{ii}$, which are usually denoted by
$h_i$. The centre of $U(C)$ contains all elements $h_{ij}$ for $i\not= j$.
\flexskip

Using the definition of the universal extension, \cite {Bor92} provides an
alternative characterisation of some GKMs (see also \cite {Jur98}).

\begin{theorem}
Suppose that $G$ is a Lie algebra satisfying the
following three properties: \flexskip

(1) $G$ can be $\ZZ$-graded as $G = \bigoplus_{i\in\ZZ} G_i$, and
$G_i$ is finite dimensional if $i\not=0$. \flexskip

(2) $G$ has an involution $\omega$ which maps $G_i$ into $G_{-i}$
and acts as $-1$ on $G_0$. \flexskip

(3) $G$ has a Lie algebra invariant bilinear form $( {^.}, {^.})$, which is
also invariant under $\omega$ such that $G_i$ and $G_j$ are orthogonal if
$i\not= -j$, and such that $-(g,\omega(g)) > 0$
if $g$ is a nonzero homogeneous element of $G$ of nonzero degree. 
\flexskip

Then there is a unique UGKM, graded
by putting ${\rm deg}(e_i) = -{\rm deg}(f_i) = n_i$ for some positive integers
$n_i$, with a homomorphism $f$ (not necessarily unique) to $G$ such that \flexskip

(a) f preserves the gradings, involutions and bilinear forms 
(as defined above). \flexskip

(b) The kernel of $f$ is in the centre of the UGKM (which is contained
in the abelian subalgebra spanned by the elements $h_{ij}$). \flexskip

(c) The image of $f$ is an ideal of $G$, and $G$ is the semidirect product of 
this subalgebra and a subalgebra of the abelian subalgebra $G_0$. Moreover, 
the images of all the generators $e_i$ and $f_i$ are eigenvectors of $G_0$.
\end{theorem}

\cite {Bor92} defines GKMs through properties (1) to (3) of theorem 1.1.
\cite {Jur98} points out that the converse of theorem 1.1 is not true:
GKMs constructed from generalized Cartan matrices cannot necessarily
be graded satisfying both conditions (1) and (3). For this reason,
the present work follows \cite {Jur98} and adopts the wider definition
of GKMs directly from generalized Cartan matrices. All GKMs considered 
in this work do, however, permit a grading of the type described in 
theorem 1.1. \flexskip

The above exposition shows that the only
major difference between GKMs and ordinary Kac-Moody algebras is that
GKMs may have imaginary simple roots. 
A further generalization of GKMs are Lie superalgebras. Here, we allow
the imaginary simple roots to have negative multiplicity. These are then 
called superroots. The Cartan matrix of a Lie superalgebra may depend on the
$\ZZ$-grading chosen. We will not encounter Lie superalgebras in the course 
of this work. 

\section{The Denominator Formula}

We recall from equations (1.1) and (1.2) that any GKM can be written as the 
direct sum $E \oplus H \oplus F$ where $H$ is the
Cartan subalgebra and $E$ and $F$ are the subalgebras corresponding to the 
positive and negative roots. The homology groups of a Lie algebra are defined
as the homology groups of the standard sequence of exterior powers,
\begin{equation}
  \textstyle
   \dots \rightarrow \bigwedge^2(E) \rightarrow \bigwedge^1(E) 
         \rightarrow \bigwedge^0(E) \rightarrow 0,
\end{equation}
(see Cartan and Eilenberg's introduction to Lie algebra homology, \cite {CE56}.)
We consider the following two virtual vector spaces.
\[
\textstyle \bigwedge(E)=   \bigwedge^0(E) 
\ominus \bigwedge^1(E) \oplus \bigwedge^2(E) \dots,
\]
which is the alternating sum of the exterior powers of $E$, and 
\[
H_*(E) =  H_0(E) \ominus H_1(E) \oplus H_2(E) \dots,
\]
which is the alternating sum of the homology groups $H_i(E)$. If $L$ is the
root lattice of the Lie algebra then both spaces are $L$-graded virtual vector
spaces whose homogeneous pieces are finite dimensional, so the infinite sums
are meaningful. From the definitions it follows that
\begin{equation}
{\textstyle \bigwedge}(E) = H_*(E),
\end{equation}
as virtual $L$-graded vector spaces (Euler-Poincare principle).
This formula can be used as a starting point to calculate the denominator 
formula for the GKM. To do this, \cite {Bor92} identifies the spaces $H_i(E)$
and then calculates the formal character
\begin{equation}
\chi(V) = \sum_{\lambda \in L} ({\rm dim} V_\lambda) e^\lambda,
\end{equation}
on both sides. 
The left hand term of (1.5), $\bigwedge(E)$, can be dealt with by a standard 
combinatorial argument, counting occurrences. The argument is analogue to that
in the case of ordinary Kac-Moody algebras, which is well documented, see, for 
example, chapter 10 of \cite {Kac90}.
Let us consider the right hand term of (1.5), $H_*(E)$. \cite {Bor92} reports
that the techniques developed by Garland and Lepowsky in \cite {GL76} can be 
adapted. A detailed discussion can be found in \cite {Jur96} (theorem 3.13).
\flexskip

\begin{theorem} 
Let $G$ be a GKM with Weyl group $W,$ Weyl vector $\rho$, and root lattice 
$L$. \flexskip

a) $H_i(E)$ is the subspace of $\bigwedge^i(E)$ spanned
 by the homogeneous vectors of $\bigwedge^i(E)$ whose degrees $r\in L$ satisfy
$(r+\rho)^2 = \rho^2$. \flexskip

b) Let $S$ denote the subspace of $H(E)$ of elements whose degree $r$ has the 
property that $r+\rho$ is in the fundamental Weyl chamber. Then $S$ is 
isomorphic to the subspace of $\bigwedge(E)$ of all elements that can be 
written in the form $e_1 \wedge e_2 \wedge \dots$ where the $e_i$'s are
vectors in the root spaces of pairwise orthogonal imaginary roots.
\end{theorem}

It must be recalled in this context that the fundamental Weyl chamber for 
GKMs is still determined by the {\it real} simple roots.
Using theorem 1.2, \cite {Bor92} calculates the formal character of $H_*(E)$.
Again, a detailed discussion can be found in \cite {Jur96} (theorem 3.16).
\flexskip

\begin{theorem}
Let $G$ be a GKM with Weyl group $W,$
Weyl vector $\rho$, root lattice $L$, and denote the positive roots by $L^+$.
If $w\in W$ then ${\rm det}(w)$ is defined to be $+1$ or $-1$, depending on 
whether $w$ is the product of an even or odd number of reflections. (If the 
root lattice is finite dimensional this is just the usual determinant of $w$.)
We define $\epsilon(\alpha)$ for $\alpha \in L$ to be $(-1)^n$ if $\alpha$ is
the sum of a set of $n$ pairwise orthogonal imaginary simple roots, and $0$
otherwise. \flexskip

Then the denominator formula of $G$ is
\[
e^\rho \prod_{\alpha \in L^+}\left(1-e^\alpha\right)^{{\rm mult}(\alpha)} =
   \sum_{w\in W}{\rm det}(w)w\left(e^\rho
	 \sum_{\alpha \in L^+} \epsilon(\alpha)e^\alpha\right).
\]
\end{theorem}

Thus, the formula is very similar to the well known denominator formula for
ordinary Kac-Moody algebras and reduces to it if there are no imaginary simple
roots. In fact, the sum over $\alpha\in L^+$ is precisely the character of the
subspace $S$ described in theorem 1.2b. If there are no imaginary simple roots 
this sum collapses to 1. Note that the definition of $\epsilon(\alpha)$
in theorem 1.3 assumes that the simple roots are linearly independent. In the
case of a specialization, $\epsilon(\alpha)$ will denote the sum over all
relevant terms $(-1)^{n_i}$, each corresponding to a representation of
$\alpha$ as the sum of $n_i$ simple roots. \flexskip

We can recover the full character formula for GKMs in the same way,
starting, in the place of (1.4), with the generalized chain complex whose 
vector spaces are spaces $\bigwedge^j(E,V)$, that is spaces $\bigwedge^j(E)$, 
tensored with any lowest weight Lie algebra module $V$. We will, 
however, not make use of this in the remainder of this paper.
\flexskip

For more information on the significance of especially the term related
to imaginary simple roots, we refer to the detailed discussion in \cite {Jur96}.
The individual terms of the denominator formula will depend on the selection
of the Cartan subalgebra, and thus the root space.
The identity is valid as an identity in the free abelian group, as
in the approach of \cite {Bor92}. Equally, the formula is meaningful and valid in
general if, as in the approach of \cite {Jur96}, we work with a suitably
extended Lie algebra $G^e$ in the place of the GKM $G(C)$ so that all simple
roots are linearly independent. \flexskip

As we have seen, the roots of the unextended GKM $G(C)$ defined in the
space $H^*$ will not necessarily be linearly independent. \cite {Jur98} discusses
that this may lead to denominator identities with infinite terms and other
problems. The affine algebra $A_2$ provides a simple example of such
problems: Here, the two real simple roots satisfy $r_2 = -r_1$ in the
unextended root space. \flexskip

There are, however, circumstances where it is possible to formulate results
in $H^*$ through specialization of the results in the extended root space
${H^e}^*$. This specialization is achieved through projection from ${H^e}^*$
to the space $H^*$. As an abstract identity, the specialization of the
denominator formula remains valid provided no multiplicities become
infinite. Note that we may encounter further problems related to the
interpretation of that abstract identity as a denominator formula. Examples
of potential ambiguities are the identification of simple and non-simple
roots, or of positive and negative roots. \cite {Jur98} discusses why
specialization works for Borcherds' monster Lie algebra. Similarly, we will
find that this is also the case for all GKMs constructed in this present
work. 

\section{Vertex Algebras}

Vertex algebras had already been used extensively in theoretical physics, when
Borcherds \cite {Bor86} formalized the definition and showed how to construct
vertex algebras for any even lattice. We follow \cite {Bor92} to give a brief survey 
of the definitions and results which will be relevant to the present work. 
\flexskip

A vertex algebra over the real numbers is a vector space $V$ over $\RR$ with an
infinite number of bilinear products, written $u_nv$, 
where $u,v,u_nv\in V$ and $n\in\ZZ$, such that

\beginaxioms

\item[(1)] $u_nv =0$ for $n$ sufficiently large (depending on $u$ and $v$).
\item[(2)] 
\[
            \sum_{i\in\ZZ} {m\choose i} \bigl(u_{q+i}v \bigr)_{m+n-i}w
\]
\[
            = \sum_{i\in\ZZ} (-1)^i {q\choose i} \biggl(u_{m+q-i}
            (v_{n+i}w) - (-1)^q \bigl(v_{n+q-i}(u_{m+i}w)\bigr)\biggr)
\]
\item[]     for all $u$, $v$, and $w$ in $V$ and all integers $m$, $n$,
            and $q$.
\item[(3)] There is an element $1\in V$ such that $v_n1=0$ if $n\ge0$ and
            $v_{-1}1=v$.
\endaxioms
The operators $u_n$ may then be combined into the vertex operator 
\[
Q(u,z) = \sum_{n\in\ZZ} u_nz^{-n-1}
\]
which is an operator valued formal Laurent series in the formal variable $z$.
\flexskip

\cite {Bor86} defines an operator $D$ on the vertex algebra by $D(v)=v_{-2}1$.
The vector space $V/DV$ is a Lie algebra, where the bracket is defined by 
$[u,v] = u_0 v$. \flexskip

A conformal vector of dimension or central charge $c\in \RR$ of a vertex 
algebra $V$ is defined to be an element $\omega$ of $V$ such that

\beginaxioms
\item[(1)] $\omega_0 v = D(v)$ for any $v\in V$,
\item[(2)] $\omega_1 \omega = 2\omega$,
\item[(3)] $\omega_3 \omega = c/2$,
\item[(4)] $\omega_i \omega = 0$, if $i=2$ or $i>3$,
\item[(5)] any element of $V$ is a sum of eigenvectors of $\omega_1$ with
            integral eigenvalues.
\endaxioms

\cite {Bor86} defines the operators $L_i$ on $V$ for $i\in \ZZ$ by 
$L_i = \omega_{i+1}$.
It can then be shown that the operators $L_i$ satisfy the relations
\[
[L_i, L_j] = (i-j)L_{i+j} + {i+1\choose 3} {\frac{c}{2}} \delta_{-j}^i.
\]
They make $V$ into a module over the Virasoro algebra. 
For $n \in \ZZ$, the physical space 
$P^n$ is defined to be the space of vectors $w\in V$ such that
\begin{subequations}
\begin{align}
L_0 (w) &= \omega_1 (w) = nw\\
L_i (w) &= 0 \text{, if } i>0. 
\end{align}
\end{subequations}
The space $P^1/(DV\cap P^1)$ is a subalgebra of the Lie algebra $V/DV$. In the 
cases considered in \cite {Bor92} (and thus in the cases considered in this work) 
this is equal to $P^1/DP^0$. \flexskip

In \cite {Bor86} examples of vertex algebras are constructed from lattices.
We present the explicit spaces, operators, and elements, both as an
illustration of the above definitions and for later reference. The results
are well-known but not necessarily in the context of Borcherds' framework.
Let $L$ be an even lattice. There exists a central extension
by a group of order 2, $\hat L$, which is uniquely characterized by the
following properties. The elements of $\hat L$ will be written
$\epsilon^n e^r$, $r\in L$, where $n=0$  or $n=1$. The
commutator is 
\begin{equation}
e^{r_1} e^{r_2}= \epsilon^{(r_1,r_2)}e^{r_2}e^{r_1} 
\mskip 20mu {\rm where} \mskip 20mu \epsilon^2=1.
\end{equation}

The underlying vector space $V(L)$ of the vertex algebra associated with 
$L$ is defined as follows; 
\begin{equation}
V(L) = \RR (\hat L)\bigotimes S\bigl(\bigoplus_{i>0}
(L_{(i)}\otimes\RR)\bigr).
\end{equation}

In \cite {Bor86}, this space is referred to as Fock space. This term has since 
been used to denote a slightly different space so that we will not use the 
term. $\RR (\hat L)$ is the twisted group ring of $\hat L$ and
$S=S\bigl(\bigoplus_{i>0}(L_{(i)}\otimes\RR)\bigr)$ is the symmetric
algebra on the sum of a countable number of copies of the lattice $L$. 
Thus, a general element of $V(L)$ will be a linear combination of elements
of the form
\begin{equation}
v = e^r \prod_{i=1}^q t_i(p_i)
\end{equation}
where $r, t_i\in L$, $q \geq 0$, $p_i\geq 1$, and 
$t_i(p_i)$ is an element of the copy $L_{(p_i)}\otimes\RR$
within the symmetric algebra, and the integers $p_i$ are not necessarily
distinct. \flexskip

In order to construct vertex operators for all elements, we begin by  
defining `annihilation' and `creation' operators. For $t\in L$, and $j\in\ZZ$,
define $t(j)$ as a linear map on $V(L)$. It is fully characterized by its 
action on elements of $V(L)$ of the form (1.10). \flexskip

if $j>0$ then $t(j)$ is multiplication by $t(j)$; 

if $j=0$ then $t(0) v = (t,r) v$;

if $j<0$ then $t(j)$ acts as a derivation so that $t(j) e^r = 0$, and
\vspace{-.35in}
\begin{multline*} \\
\shoveright{t(j) t_i(p_i) = -j(t,t_i) \delta^{-j}_{p_i}.} \\
\end{multline*}
\vspace{-0.35in} 

We define the vertex operator of a general element of $V(L)$ in three steps.
Let $t\in L$, $p>0$, let $z\in \CC$ be a complex number, and let $v$ be as 
in (1.10).
\[
Q(t,z) = \sum_{j \ne 0} t(j) {\frac{z^j}{j}} + t(0){\rm log} z + t
\]
\[
Q\left(t(p),z\right) = 
   {\frac{1}{(p-1)!}} \left( {\frac{\rm{d}}{\rm{d}z}} \right) Q(t,z)
\]
\[
Q(v,z) = Q\left( e^r \prod t_i(p_i), z \right) = \mskip 10mu
   :e^{Q(r,z)}\prod Q\left(t_i(p_i),z\right):
\]
Here the `:' is the standard notation for normal ordering, such that 
all `creation' operators ($e^r$, $t(p), p\geq 1$) occur to the left of all 
`annihilation' operators ($t(p), p\leq 0$). 
The product $v_n(w)$ for $n\in \ZZ$ and $w \in V(L)$ is then defined as the 
coefficient of $z^{-n-1}$ in $:Q(v,z): (w)$. \flexskip

Choose a basis $s_i, i\in I$, of the 
lattice $L$, and a dual basis $s_i'$. Then the vector
\[
\omega = {\frac{1}{2}} \sum_{i\in I} s_i(1)s_i'(1)
\]
is a conformal vector in $V(L)$. 
Its central charge can be identified as the dimension
of the lattice $L$. Applying the general vertex operator construction to 
$\omega$, we identify the operators $L_n = \omega_{n+1}$ as
\[
L_n = \sum_{j\in \ZZ} \sum_{i\in I} :s_i(j)s_i'(-n-j):
\]
In particular, we note the action of $L_0$
on homogeneous elements of the form (1.10)
\begin{equation}
L_0 v = \bigl({\frac{r^2}{2}} + \sum_i p_i\bigr) v. 
\end{equation}
Thus, $L_0$ defines a $\ZZ$-grading which we will refer to below as 
${\rm deg}_{\ZZ}$, or ${\rm deg}_{\ZZ}^{(L)}$ if we require to specify the  
referenced underlying lattice $L$. 
$L_{-1}$ is the derivation $D$. For $r,t\in L$, 
and $p\geq 1$ we obtain $L_{-1}e^r = r(1)e^r$, and $L_{-1}t(p) = pt(p+1)$.
\flexskip

Suppose $A$ is an ordinary Kac-Moody algebra with simple roots $a_i$ of norm
$2$. The derived algebra $A'$ is the algebra generated by generators
$e_i, f_i, h_i$ for each simple root $a_i$.
If the lattice $L$ contains the root lattice of $A$ then  
we can map the derived algebra $A'$ of $A$ to $P^1/DP^0$, as follows;
$e_i \mapsto e^{a_i}$, $f_i \mapsto e^{-a_i}$, $h_i \mapsto a_i(1)$.
As an example of the definitions and constructions above, we calculate 
the bracket $[e_i,h_j] = (e_i)_0 (h_j)$ in $P^1/DP^0$.
\[
Q\left(e^{a_i},z\right) a_j(1) = 
   e^{a_i}
   {\rm exp} \Bigl( {\sum_{n>0} {\frac{a_i(n)}{n}}z^n} \Bigr)
   {\rm exp} \Bigl( {\sum_{n<0} {\frac{a_i(n)}{n}}z^n} \Bigr)
   z^{a_i(0)}a_j(1)
\]
Here, $z^{a_i(0)}$ acts on $e^r\in V(L)$ as multiplication by 
$z^{(a_i,r)}$. It acts as 
identity on the symmetric algebra. Taylor expand the exponentials to obtain
\[
Q\left(e^{a_i},z\right) a_j(1)= 
     e^{a_i} \Bigl(1 + a_i(1)z^1 + \dots \Bigr) 
     \Bigl(1 - a_i(-1)z^{-1} + \dots \Bigr) a_j(1)
\]
The bracket $[e_i,h_j]$ is the coefficient of $z^{-1}$ of the expansion. 
The terms $a_i(-n)a_j(1)$ vanish for $n \neq 1$. Hence, when applied to 
$a_j(1)$, only the term $e^{a_i} \times (1) \times ( - a_i(-1)z^{-1})$ gives 
a non-zero contribution to the Lie-algebra bracket. Thus, 
\[
[e^{a_i}, a_j(1)] = - (a_i, a_j) e^{a_i}.
\]
This provides an illustration of how the notions of the roots and root spaces 
of an ordinary Kac-Moody algebra are embedded in the framework of vertex
algebras. If we define the natural $L$-grading of $V(L)$ as 
\begin{equation}
{\rm deg}_L(e^r) = r, \mskip 5mu {\rm deg}_L(t(p)) = 0,
\end{equation}
then we find that the Lie algebra
$P^1/DP^0$ has a natural decomposition into root spaces. 
Furthermore, for any root $r\in L$, the root space of $r$ is the 
homogeneous subspace of the Lie algebra of degree $r$. \flexskip
 
We can define a unique bilinear form on $V(L)$ through the following two 
requirements \cite {Bor86};
\begin{subequations}
\begin{align}
&(e^{r_1},e^{r_2}) = \begin{cases}
				1&  \text{ if } r_1=-r_2,\\
				0&  \text{ otherwise;}
                           \end{cases} \\
&\text{ the adjoint of } t(p) \text{ is }-t(-p).
\end{align}
\end{subequations}
This bilinear form induces an invariant bilinear form on the Lie algebra
$P^1/DP^0$ (\cite {Bor92}, proof of theorem 6.1). 
Evaluating the bilinear product on $V(L)$ we find that $L_i$ is the 
adjoint of $L_{-i}$. Let $v\in P^0$, and $w\in P^1$. Then $Dv$ is in $P^1$, 
and $(Dv,w) = (L_{-1}v,w) = (v,L_1w) = 0$. The kernel of the bilinear form
$( {^.}, {^.})$ on $P^1$ thus contains $DP^0$. \flexskip

Using this explicit construction of the Lie algebra elements, \cite {Bor86}
announces the following results about the Lie algebra $P^1/DP^0$. \flexskip

\begin{theorem} Let $L$ be a non-singular, even lattice.
Let the physical space $P^1 \subset V(L)$ be defined as in (1.7). \flexskip

a) Let $A$ be a Lie algebra, $A$ connected, 
simply laced, and not 
affine. If the lattice $L$ contains the root lattice of $A$ (possibly 
quotiented out by some null lattice) then $A$ can be mapped to $P^1/DP^0$ 
such that the kernel is in the centre of $A$. \flexskip

b) Let $d$ be the dimension of $L$, and let $r\in L$ be a root such that 
$r^2 \le0$. Let $p_d(n)$ denote the number of partitions of $n$ into $d$ 
colours. Then the dimension of the degree $r$ subspace of $P^1/DP^0$ is equal
to 
\[ 
p_{d-1}\left(1-{\frac{r^2}{2}}\right) - p_{d-1}\left(-{\frac{r^2}{2}}\right).
\]
This forms an upper limit of the multiplicities of the roots of any Lie algebra
$A$ which satisfies the assumptions in a). 
\end{theorem}

To provide an indication of the proof and the type of elements in the 
space $P^1/DP^0$, we recall from (1.10) that 
the degree $r$ part of $P^1$ is spanned by vectors of the form 
$v=e^r \prod t_i(p_i)$. Being an element of $P^1$, $v$ must be 
eigenvector of $L_0$ with eigenvalue 1, hence we require 
$1 = {\frac{r^2}{2}} + \sum p_i$. The sum over $p_i$ 
defines a partition (with colours) of the total of $1-{\frac{r^2}{2}}$. 
The colours of the partition correspond to the dimension $d$ of the 
underlying lattice $L$ - note that the $t_i(p_i)$ are not necessarily 
linearly independent.
The conditions $L_iv=0, i>0$, of formula (1.7), introduce further
restrictions which reduce the dimension $d$ by 1.
Thus, $p_{d-1}(1-r^2/2)$ is the dimension of the space 
$P^1$. Similarly, $p_{d-1}(-r^2/2)$ is the dimension of $P^0$, and 
we note that $D$ maps $P^0$ injectively into $P^1$.
\cite {Bor86} announces that the above formulas 
can be generalized for the case that $A$ is not simply laced. We will
not require the generalized formulas in this work, though. 
We are now in the position to quote Borcherds' version of the no-ghost theorem.
(The original form of the no-ghost theorem was proven by Goddard and Thorn.)
\flexskip

\begin{theorem} Suppose that $V$ is a vector space with a non-singular
bilinear form $( {^.}, {^.})$ and suppose that V is acted on by the Virasoro algebra in
such a way that the adjoint of $L_i$ is $L_{-i}$, the central element of the 
Virasoro algebra acts as multiplication by $24$, any vector of $V$ is a sum of 
eigenvectors of $L_0$ with non-negative integral eigenvalues, and all the 
eigenspaces of $L_0$ are finite dimensional. We let $V_n$ be the subspace of 
$V$ on which $L_0$ has eigenvalue $n$. Assume that $V$ is acted on by a group
$G$ which preserves all this structure. We let $V(II_{1,1})$ be the vertex
algebra of the double cover of the two dimensional even unimodular 
Lorentzian lattice $II_{1,1}$ (so that $V(II_{1,1})$ is $II_{1,1}$-graded,
has a bilinear form $( {^.}, {^.})$ and is acted on by the Virasoro algebra as above).
We let $P^1$ be the subspace of the vertex algebra $V\otimes V(II_{1,1})$
of vectors $v$ with $L_0(v) = v$, $L_i(v) = 0$ for $i>0$, and we let $P^1_r$
be the subspace of $P^1$ of degree $r\in II_{1,1}$. All these spaces inherit
an action of $G$ from the action of $G$ on $V$ and the trivial action of $G$
on $V(II_{1,1})$ and $\RR^2$. Then the quotient of $P^1_r$ by the nullspace
of its bilinear form is naturally isomorphic, as a $G$ module with an invariant
bilinear form, to $V_{1-r^2/2}$ if $r\not=0$ and to $V_1\oplus\RR^2$ if $r=0$.
\end{theorem}

In the next section we will describe how \cite {Bor90b} used the no-ghost 
theorem to identify the root spaces of certain GKMs explicitly. Note that
the theorem only applies in the case of central charge 24, which is one reason 
why Borcherds' construction cannot be generalized straightforwardly to lattices
of dimension other than 24.

\section{The Fake Monster Lie Algebra}

We will now recall  
Borcherds' construction of the fake monster Lie algebra. The original
calculations were published in \cite {Bor90b}, where the algebra is called the 
monster Lie algebra. The construction uses the Leech lattice, which is an even
lattice of 24 dimensions, such that all the results of the previous section,
including the no-ghost theorem, apply. \flexskip

For this construction we consider the vertex algebras of three different
lattices, each with their two respective natural gradings, as defined above
in equations (1.11) and (1.12). The three lattices are \flexskip

$\Lambda$ - the Leech lattice, \flexskip

$II_{1,1}$ - the unique even 2-dimensional unimodular Lorentzian lattice, \flexskip

$II_{25,1}$ - the 26-dimensional unimodular Lorentzian lattice 
	  $\Lambda \oplus II_{1,1}$.\flexskip

If $L$ is any of the above three lattices and $V(L)$ the vertex algebra 
associated with $L$, we use the $L$-grading (1.12) and $\ZZ$-grading (1.11) 
to define $V(L)_r$, $V(L)_n$ and $V(L)_{(r,n)}$ as the parts of $V(L)$ of 
${\rm deg}_L = r$ or ${\rm deg}_{\ZZ}^{(L)} = n$, respectively. We will also 
use $S_n$, where $S$ is the symmetric algebra, introduced in formula (1.9), 
and the $\ZZ$-grading on $S$ is the restriction of the $\ZZ$-grading on 
$V(L)$. \flexskip

We define the norm of an element $(m,n) \in II_{1,1}$ as $-2mn$. Therefore, 
the norm of an element $(\lambda,m,n) \in II_{25,1}$ is 
$\lambda^2 - 2mn$. As introduced in equation (1.3), the norm is
the square of the usual one.  \flexskip

We consider the physical space 
$P^1{\scriptstyle(V(\Lambda)\otimes V(II_{1,1}))}$ 
of $V(\Lambda)\otimes V(II_{1,1})$, and a general element $(m,n)$ of 
$II_{1,1}$. Then 
$\bigl(P^1{\scriptstyle(V(\Lambda)\otimes V(II_{1,1}))}\bigr)_{(m,n)}$ is 
the part of \\ $P^1{\scriptstyle(V(\Lambda)\otimes V(II_{1,1}))}$ 
whose $II_{1,1}$-grading is $(m,n)$. On this space, a bilinear form is
defined as in (1.13).
Let $K$ denote the null space of this bilinear form. 
The no-ghost theorem 1.5 identifies 
\[
\bigl(P^1{\scriptstyle(V(\Lambda)\otimes V(II_{1,1}))}\bigr)_{(m,n)}/K 
 \mskip 5mu \tilde =  \mskip  5mu V(\Lambda)_{1-{\frac{(m,n)^2}{2}}} 
 \mskip 5mu = \mskip  5mu V(\Lambda)_{1+mn}.
\]
Using the $\Lambda$-grading of $V(\Lambda)$ on both sides, Borcherds 
refines the isomorphism to 
\[
\bigl( P^1{\scriptstyle(V(\Lambda)\otimes V(II_{1,1}))}\bigr)_{(\lambda,m,n)}
 /K \mskip 5mu \tilde = \mskip  5mu V(\Lambda)_{(\lambda,1+mn)}.
\]
Next, consider the vertex algebra of the 26-dimensional Lorentzian
lattice $V(II_{25,1})$ which is isomorphic to 
$V(\Lambda)\otimes V(II_{1,1})$ \cite {Bor92}. Consider the physical space 
$P^1{\scriptstyle (V(II_{25,1}))}$. 
Let $K'$ denote the kernel of the bilinear form 
$( {^.}, {^.})$ on $P^1{\scriptstyle (V(II_{25,1}))}$. 
Then \cite {Bor90b} defines the fake monster Lie algebra as 
\[
M_\Lambda \mskip 5mu = \mskip  5mu 
 \bigl( P^1{\scriptstyle(V(II_{25,1}))} \bigr) / K' 
 \mskip 5mu \tilde = \mskip 5mu
 \bigl( P^1{\scriptstyle(V(\Lambda)\otimes V(II_{1,1}))}\bigr) / K.
\]
(The explicit construction first quotients out $DP^0$ but we know that
$DP^0 \subset K'$.) Combining the two isomorphisms above, we obtain that
\[
(M_\Lambda)_{(\lambda,m,n)} \mskip 5mu \tilde= \mskip 5mu 
 V(\Lambda)_{(\lambda,1+mn)}.
\]
Note that the left hand side is derived from the vertex
algebra of the Lorentzian lattice $II_{25,1}$, and is graded in $II_{25,1}$.
The right hand side is a piece of the vertex algebra of the Leech lattice
with $\Lambda$ and $\ZZ$ gradings. \flexskip

\cite {Bor90b} shows that the Lie algebra $M_\Lambda$ satisfies all conditions
of theorem 1.1 above. $M_\Lambda$ is therefore associated with a GKM,
through central extensions, and quotienting out the null space of the
bilinear form. Section 5 of \cite {Bor90b} identifies the universal central
extension $\hat M$ of $M_\Lambda$. The no-ghost theorem provides us with an
explicit description of the building blocks of $M_\Lambda$, in terms of
subspaces of vertex algebras. \cite {Bor90b} uses this, the denominator formula
for GKMs (theorem 1.3), and the theory of modular forms to determine the
simple roots of $M_\Lambda$.\flexskip

The root lattice of the Lie algebra $P^1/DP^0$, as a subset of the dual of
the (unspecialized) Cartan subalgebra, is infinite-dimensional and singular
with regard to the scalar product. Quotienting out the null space $K'$
corresponds to a specialization of the root space in the terminology of
\cite {Jur96} as in section 1.1.2, above. Theorem 1.3 provides a denominator
formula for a suitably extended Lie algebra ${\hat M}^e$. Before we can
formulate the denominator formula of $M_\Lambda$ itself (that is, without any
extension of the Cartan subalgebra) we need to verify the validity of the
specialization. Let $\pi_0$ denote the projection from $P^1/DP^0$ to
$M_\Lambda$, and equally the projection of the dual spaces of the Cartan
subalgebras (that is, the spaces within which the roots are defined). 
The subscript 0 is used to indicate that the projection is defined
through the nullspace of the scalar product. \cite {Bor92} shows that, under the
projection $\pi_0$, the set of roots is mapped into the lattice $II_{25,1}$.
Defining the norm zero vector $\rho = (0,0,1)\in II_{25,1}$, \cite {Bor92}
establishes that the real simple roots are projected
to $(\lambda,1,{\frac{\lambda^2}{2}}-1)$ with multiplicity 1 and the imaginary
simple roots project to $n\rho = (0,0,n)$ with multiplicity 24 for all $n>0$.
We have $(\rho,\pi_0 r)=-{\frac{r^2}{2}}$ for all simple roots $r$, which
makes $\rho$ a Weyl vector. The set of positive roots projects to
the set of vectors in $II_{25,1}$ of norm at most 2 which are either
positive multiples of $\rho$ or have negative inner product with $\rho$.
\flexskip

An explicit description of the images of all positive roots under the
projection $\pi_0$ can be given as follows. A root $r$ of $M_\Lambda$ is
positive if for $\pi_0 r=(\lambda,m,n)$ either $m>0$, or $m=0$ and 
$n>0$ holds. A root is negative if either $m<0$, or $m=0$ and $n<0$.
Note that no element $r$ with $\pi_0 r=(\lambda,0,0)$ can be a root as, in
that case, $r^2 \ge 4$. Thus, the above description covers all roots of
$M_\Lambda$. \flexskip

Using the explicit description, it is straightforward to see that the
projection $\pi_0$ does not create infinite terms, nor does it identify
positive and negative roots. Evaluating the actual root
multiplicities provided by Borcherds' denominator formula, we also find that
the projection does not identify simple and non-simple roots. \flexskip

Furthermore, the explicit identification of the roots then allows to build
the matrix of scalar products which satisfies conditions (C1) to (C3) of a
generalized Cartan matrix. The inner product with the Weyl vector can be
used as a $\ZZ$-grading of $M_\Lambda$. \flexskip

Having verified that the specialization of the roots to $II_{25,1}$ is both
valid as an abstract identity, and meaningful as the denominator formula
of $M_\Lambda$, we can apply the full theory of denominator formulas to
$M_\Lambda$, with roots in $II_{25,1}$. We will write $M_{(\lambda,m,n)}$
for $(M_\Lambda)_{(\lambda,m,n)}$. Let furthermore $E$ denote the part of
$M_\Lambda$ corresponding to the positive roots. \flexskip

We can summarize the above results in the (specialized) denominator formula
of $M_\Lambda$. The remark following theorem 1.3 identified
$\epsilon( {^.} )$ in the case that the simple roots are not linearly
independent, as occurs in this specialization. All positive multiples of
$\rho$ are simple roots of multiplicity 24 and are perpendicular to each
other, so $\vert \epsilon(n\rho) \vert$ is $p_{24}(n)$, the number of
partitions of $n$ into parts of $24$ = dim$(\Lambda)$ colours (notation
introduced in theorem 1.4). Thus $\epsilon(n\rho)$ equals the coefficient
of $q^n$ in $\prod_n (1-q^n)^{24} = q^{-1}\eta^{24}(q)$ where $\eta$ is the
Dedekind eta-function. We will give details on $\eta$ in chapter 2.
\cite {Bor90b} obtains
\begin{equation}
e^\rho \prod_{r\in (II_{25,1})^+}\bigl(1-e^r\bigr)^{p_{24}(1-r^2/2)} = 
\sum_{w\in W} {\rm det}(w) w\bigl(\eta^{24}(e^\rho)\bigr).
\end{equation}
Here, $(II_{25,1})^+$ is the set of positive roots identified above,
and $W$ is the Weyl group, which is the group of isometries of the root
lattice generated by the reflections corresponding to the real simple
roots. \cite {Bor90b} gives a full description of the Cartan subalgebra of the  
universal central extension $\hat M$ of $M_\Lambda$, as introduced in
section 1.1.4. An explicit basis for the Cartan subalgebra is the sum of
a one dimensional space for each vector of the Leech lattice and a space
of dimension $24^2=576$ for each positive integer $n$.
The significance of the latter is that, for each $n$, the imaginary
simple root $n\rho$ has multiplicity $24$. Correspondingly, there are
$24^2 = 576$ generators $h_{ij}, i,j=1,\dots24$. Because of the nature of
its construction, it would be tedious to describe the Cartan subalgebra of
$M_\Lambda$ itself. For details of the constructions, see \cite {Bor90b}.\flexskip

We conclude this section with a brief remark about the monster Lie algebra.
The name was in \cite {Bor90b} used for the GKM which we now call fake monster Lie
algebra. In \cite {Bor92} the (proper) monster Lie algebra was constructed from
the vertex algebra of the monster Lie group, as introduced in \cite {FLM88}. The
steps from the vertex algebra to the Lie algebra are parallel to the case of
the fake monster Lie algebra, above. The vertex algebra is tensored with
$V(II_{1,1})$, the physical spaces $P^1_r$ and the kernel $K$ of the
bilinear form are defined as above. Again, the no-ghost theorem 1.5 applies
to the pieces $P^1_r/K$ and provides an explicit description of the building
blocks. \cite {Bor92} uses this to show that the Lie algebra $P^1/K$ satisfies all
conditions of theorem 1.1, above, and hence is associated with a GKM.
\cite {Bor92} calculates an explicit denominator formula, analogue to
formula (1.14). Again, this is a specialization in the terminology of
\cite {Jur96}. Theorem 6.1 of \cite {Jur98} discusses the details of this 
specialization and identifies a GKM {\fraktur g}$(M)$ and an ideal 
{\fraktur c} so that the  Monster Lie algebra can be recovered as the quotient 
{\fraktur g}$(M)/${\fraktur c}. In the specialization, the simple
roots of the monster Lie algebra are identified as the set $\{(1,i) \ \vert
\ i=-1,1,2,3,\dots\}$. The steps to prove that the specialization is well
defined are similar to the case of the fake monster Lie algebra
because the set of projected simple roots of the monster Lie algebra
displays characteristics very similar to the set of projected simple
roots of the fake monster Lie algebra. 

\section{The Twisted Denominator Formula}

\cite {Bor92} uses the concept of `twisted' denominator formulas, which are 
a generalization of the ordinary denominator formulas. 
Starting from the equality of graded virtual spaces (1.5),
${\textstyle \bigwedge}(E)=H_*(E),$
we obtained the ordinary denominator formulas by calculating the formal
character
\[
\chi = \sum_\lambda ({\rm dim}\  V_\lambda) e^\lambda
         = \sum_\lambda ({\rm Tr}\ {\rm id} |\ V_\lambda) e^\lambda, 
\]
on both sides, as in (1.6). As before, $V_\lambda$ corresponds to a lattice
grading of $V$. Let $\sigma$ be an automorphism on $V$. Then we may,
more generally, calculate the trace of $\sigma$ 
\begin{equation}
\sum_\lambda ({\rm Tr}\ \sigma |\ V_\lambda) e^\lambda 
\end{equation}
on both sides of (1.5). 
We will arrive at some generalized (`twisted') denominator formula.
In the case of the monster Lie algebra introduced in the concluding 
paragraph of section 1.4, we may choose 
$\sigma$ to be an automorphism of the monster group.
The resulting twisted denominator formula is closely related to the Thompson 
series $T_\sigma (q)$ (see \cite {Bor92} for more details). \flexskip

Returning to the fake monster Lie algebra, we consider automorphisms of
the vertex algebra $V(\Lambda \oplus II_{1,1})$ which are induced by
automorphisms of $\hat \Lambda$, that is the Leech lattice, centrally
extended by a group of order $2$, as defined in equation (1.8). If
$\sigma\in \hbox{Aut} (\hat \Lambda)$ then we define, using the notation of
equation (1.10) for a general element of $V(\Lambda\oplus II_{1,1})$,
\begin{equation}
\sigma(v) := \sigma(e^r)
\bar \sigma t_1(p_1) \dots \bar \sigma t_q(p_q).
\end{equation}
Here $\bar \sigma$ is the projected action on elements of $\Lambda$
which is well defined because
$\hbox{Aut}(\hat \Lambda)$ is an extension of $\hbox{Aut}(\Lambda)$ of the form
$2^{24} {^.}\hbox{Aut} (\Lambda)$, such that
$\sigma(e^r)=\epsilon_\sigma(r)e^{\bar\sigma r}$
where $\epsilon_\sigma$ satisfies equation (1.8).
The `$^.$' indicates that the extension is nonsplit, cf. the Atlas \cite {Con85}.
Note in particular that $\hbox{deg}_\Lambda(\sigma(\hbox{v}))=
\bar\sigma(\hbox{deg}_\Lambda(\hbox{v}))$.
Now let $\sigma$ be an automorphism of $\hat\Lambda$ of finite order $N$. Thus
we can regard $\bar\sigma$ as an element of $SL_{24}(\ZZ)$. Its Jordan normal
form is diagonal and the eigenvalues $\epsilon_j,j=1,\dots,24$ are
$n^{\hbox{th}}$ roots of unity where $n\vert N$. Equivalently, we can describe
$\sigma$ through its cycle shape $a_1^{b_1}\dots a_s^{b_s}$. We associate
$\sigma$ to the following modular form:
\begin{equation}
\eta_\sigma(q):=\eta(\epsilon_1q) \dots \eta(\epsilon_{24}q) =
\eta(q^{a_1})^{b_1} \dots \eta(q^{a_s})^{b_s}.
\end{equation}
Here $\eta$ stands for the Dedekind eta-function. We will give details on
$\eta$ in chapter 2. \flexskip

Restricting further, let $\sigma\in \hbox{Aut}(\hat \Lambda)$ 
be of prime order $N$, let $\Lambda^\sigma$ denote the fixed
point lattice and $L={\Lambda^\sigma}\oplus II_{1,1}$ the corresponding
Lorentzian lattice. Following \cite {Bor92} (section 13, p. 438) we assume for
simplicity that any power $\sigma^n$ of $\sigma$ fixes all elements of
$\hat\Lambda$ which are in the inverse image (with respect to the
projection $\hat\Lambda \to \Lambda$) of any vector of $\Lambda$
fixed by $\sigma^n$. This final assumption will be satisfied by all
automorphisms $\sigma$ of interest in this work. \flexskip

For the dual lattices we find $L^* = {\Lambda^\sigma}^*\oplus II_{1,1}$.
The $(II_{25,1})$-grading of $M_\Lambda$ induces an $L^*$-grading as follows.
It is well known that the projection $\pi_\sigma:\Lambda \to
{\Lambda^\sigma}^*$ maps onto the dual. (See theorem 3.1 below for a proof.)
For $r = (\lambda^*,m,n) \in L^*$ let
\begin{equation}
\tilde M_r = \tilde M_{(\lambda^*,m,n)} =
\bigoplus_{ \lambda \in \Lambda : \pi_\sigma (\lambda) =\lambda^*}
M_{(\lambda,m,n)}.
\end{equation}
Also, we will write $\tilde E_r$ and $E_r$ for $\tilde M_r$ and $M_r$ if we
want to emphasize that we consider positive roots. \flexskip

If we consider both sides of equation (1.5), $\bigwedge(E)=H_*(E)$, as 
$L^*$-graded Aut$(\hat\Lambda)$ modules then we can calculate the 
trace of $\sigma$. For any $r\in L^*$ define the numbers mult$(r)$ as
\begin{equation}
{\rm mult}(r)=   \sum_{d,s>0, ds\vert ((r,L),N)}
{\frac{\mu(s)}{ds}}\hbox{Tr}(\sigma^d\vert \tilde E_{\frac{r}{ds}})
\end{equation}
\cite {Bor92}, (13.2). Here $\mu(s)$ is the M\"obius function, and $(r,L)$ denotes
the greatest common divisor of the numbers $(r,a)$ for $a\in L$. Then the 
trace (1.15) on the equation (1.5) of graded virtual vector spaces is
\begin{equation}
e^\rho \prod_{r\in {L^*}^+}(1-e^r)^{\rm{mult}(r)} 
 = \sum_{w\in W^\sigma} \hbox{det}(w)w(\eta_\sigma(e^\rho))
\end{equation}
cf. \cite {Bor92}, (13.3), with $\eta_\sigma$ as in (1.17) and $\rho= (0,0,1)
\in L$. The proof is quite
similar to the proof of theorem 1.3, that is the untwisted character formula.
For the left hand side, again, the combinatorial argument applies. For the 
right hand side, again, the significant contribution is the 
term $w(\eta_\sigma(e^\rho))$, which is derived directly from the imaginary 
roots of $M_\Lambda$, as described in theorem 1.2b. \flexskip

Equation (1.20) contains two terms, $W^\sigma$ and ${L^*}^+$, which need to
be discussed in more detail. \cite {Bor92} introduces the term $W^\sigma$ to denote
the subgroup of the Weyl group $W$ of $M_\Lambda$ consisting of reflections
which commute with $\sigma$. This group is a subgroup of the reflection group
of $L$. Theorem 2.2 of \cite {Bor90a} provides a number of useful equivalent
characterisations of the group $W^\sigma$. Let $\pi_\sigma$ denote the
projection from the span of $II_{25,1}=\Lambda\oplus II_{1,1}$ to the span
of $L=\Lambda^\sigma\oplus II_{1,1}$ induced by the automorphism $\sigma$:
\flexskip

\begin{lemma} Let $\sigma\in \hbox{Aut} (\hat \Lambda)$ be as above,
let $W$ be the Weyl group of $M_\Lambda$ as introduced in equation (1.14).
Let $W^\sigma$ be the group of automorphisms introduced in equation (1.20).
Then the following are equivalent characterizations of $W^\sigma$.

\beginaxioms
\item[($\alpha$)] $W^\sigma$ is the group of elements of $W$ commuting with
$\sigma$. \flexskip

\item[($\beta$)] $W^\sigma$ is the group of elements of $W$ fixing the
subspace $L$. \flexskip

\item[($\gamma$)] $W^\sigma$ is generated by  the reflections of the vectors
$\pi_\sigma r$ as $r$ runs through the simple roots of $W$ whose projections
$\pi_\sigma r$ have positive norm. \flexskip

\item[($\delta$)] Same as $(\gamma)$, with `simple roots' replaced by `roots'.
\endaxioms
\end {lemma}
\begin{proof} Theorem 2.2 of \cite {Bor90a}. 
\end{proof}

\cite {Bor92} identifies ${L^*}^+$ as the set of all
$r=(\lambda^*,m,n) \in L^* $ such that $m>0$, or $m=0$ and $n>0$, which can
be derived directly form the description of the positive roots of 
$M_\Lambda$ in section 1.4, above.
Hence, the projection $\pi_\sigma$ does not identify positive and non-positive
roots of $M_\Lambda$. \cite {Bor92} observes that $\pi_\sigma \rho$ is a norm 0
Weyl vector for $W^\sigma$ which we will again denote $\rho$ for simplicity.
The description $(\gamma)$ of lemma 1.1, above, shows that the simple roots
of $W^\sigma$ can be represented through vectors in ${L^*}^+$.

\section{Construction of the GKMs}

We continue to use the notation of the previous section. That is, let 
$\sigma\in \hbox{Aut}(\hat \Lambda)$ be of prime order $N$.
Now assume furthermore that $N$ is any one of $2$, $3$, $5$, $7$, $11$, or
$23$, and that $\sigma$ is of cycle
shape $1^MN^M$ where $M=24/(N+1)$. Lemma 12.1 of \cite {Bor92} confirms that
these automorphisms satisfy all assumptions placed on $\sigma$ in the
previous section 1.5. Let $\Lambda^\sigma$ denote the fixed
point lattice and $L={\Lambda^\sigma}\oplus II_{1,1}$ the corresponding
Lorentzian lattice as above. Starting from equation (1.20),
the aim of this section will be to construct a new generalized Cartan
matrix, and thus a new GKM for each $N$. \cite {Bor92} identifies a
subset of ${L^*}^+$ which constitutes the set of `prospective simple roots'.
Clearly, they are not linearly independent. Therefore, we expect that
they will be roots resulting from a specialization in the sense of
\cite {Jur98}, and that we will recover (1.20) as the denominator formula of the
new (specialized) GKM. \flexskip

Before we can verify the properties of a generalized Cartan matrix we need
to recall a few results of \cite {Bor90a} about the fixed point lattice $L$ and
the Weyl group $W^\sigma$. As seen in lemma 1.1, the group $W^\sigma$ is a
subgroup of the reflection group of the fixed point lattice $L$. In general,
it will not be the full reflection group of $L$. One obvious necessary
condition for it to be the full group is that the lattice $\Lambda^\sigma$
has no roots. Any root $\lambda\in\Lambda^\sigma$ induces a root
$(\lambda,0,0)$ of $L$. However, in the cases considered in this work, the
condition is also sufficient.
The reflection induced by a root $r\in L$ is the same as that induced by $nr$
where $n$ is any non-zero integer. Hence we may restrict ourselves to roots
$r$ which are primitive in the sense that, if $n\in \ZZ$, $\vert n \vert>1$,
then $\frac{r}{n}$ is not in $L$. \flexskip

\begin{lemma} Let $\sigma$ be an automorphism of the Leech lattice
$\Lambda$ such that the sublattice $\Lambda^\sigma$ fixed by $\sigma$ has
no roots, and let $L=\Lambda^\sigma \oplus II_{1,1}$. Let $\rho=(0,0,1) \in
II_{25,1}=\Lambda \oplus II_{1,1}$ denote the Weyl vector of $II_{25,1}$.
\flexskip

a) The group $W^\sigma$ is the full reflection group of the lattice $L$. \flexskip

b) The simple roots of $W^\sigma$ are exactly the roots $r \in L$ such that
$(r,\rho)$ is negative and divides $(r,v)$ for all vectors $v$ of $L$. \flexskip

c) The vector $\rho$ is also norm 0 Weyl vector for the lattice $L$. That
is, the primitive simple roots of $W^\sigma$ are exactly the primitive roots
$r$ of $W^\sigma$ satisfying $(r,\rho)= -r^2/2$.
\end{lemma}

\begin{proof} The introduction of \cite {Bor90a} asserts that all roots considered
are implicitly understood to be primitive in that paper. \cite {Bor92} uses the
results of \cite {Bor90a} so that the implicit assumption also applies to \cite {Bor92}.
Given this implicit understanding, claim b) is part of theorem 3.3 of
\cite {Bor90a}. Claim a) is an auxiliary result from the proof of the same
theorem. Claim c) is quoted from section 13 of \cite {Bor92}.  
\end{proof}

In order to apply lemma 1.2 we still need to verify the assumptions about
$\Lambda^\sigma$. These depend on specific properties of the Leech lattice
and the cycle shape of $\sigma$. We will establish these properties in
chapter 4.1, below. At this point, we merely state \flexskip

\begin{lemma} Let $\sigma$ be an automorphism of the Leech lattice
$\Lambda$ of cycle shape $1^MN^M$ where $N$ is any of the primes $2$, $3$,
$5$, $7$, $11$, $23$, and $M=24/(N+1)$. Then $\Lambda^\sigma$ has no roots.
\end{lemma}

\begin{remark} The proof will be provided in section 4.1, below. 
\end{remark}

Let us now return to the twisted denominator formulas.
If we compare equation (1.20) to the denominator formula of a GKM (theorem
1.3 of section 1.2 above) we observe that the right hand side of (1.20) is
precisely the right hand side of (the specialization of) a denominator
formula for a GKM with Weyl group $W^\sigma$.
Until the definition of the GKM and the verification of the specialization
have been completed, we will need to distinguish carefully between the
linearly independent roots of the suitably extended Lie algebra, and those
of the specialization to $L^*$. We will refer to the linearly independent
roots as `roots', and to the (specialized) roots in $L^*$ as `images of
roots'. \flexskip

Assuming that there exists a suitable GKM with well-defined specialization,
the images of its simple roots under the specialization must be selected
in ${L^*}^+$ as follows:

\beginaxioms
\item[(G1)] As images of the real simple roots for the GKM, select the
simple roots of the reflection group $W^\sigma$. From lemmas 1.1 and 1.2 it
follows that each such simple root of $W^\sigma$ may be represented through
a vector $r$ in ${L^*}^+$ within the fundamental Weyl chamber of $W^\sigma$.
The scaling factor of the primitive representative $r$ is such that
$(r,\rho) = - r^2/2$.\flexskip

\item[(G2)] As images of the imaginary simple roots, select the positive
multiples $n\rho$ of the Weyl vector $\rho$, with multiplicities equal to
\begin{equation}
\hbox{mult}(n\rho)= 
  \sum_{a_k\ {\rm divides}\ n} b_k, \ \ \ n\in\ZZ, \ n>0, 
\end {equation}
if $\sigma$ has generalized cycle shape $a_1^{b_1}a_2^{b_2}\dots$
\endaxioms

The value for the multiplicity of imaginary simple roots specified in
equation (1.21) can be obtained by carefully counting occurrences, using
the general description of the space $S$ of theorem 1.2b in section 1.2 -
note that all multiples of $\rho$ are mutually orthogonal to each other.
Thus, all `prospective simple roots' selected under (G1) and (G2) have been
represented in the lattice ${L^*}^+$. Let $\{ r_i\}$ be the set of all such
`prospective simple roots'. A scalar product
for this set is naturally inherited from the lattice. We form the matrix
of scalar products $C=\Bigl((r_i,r_j)\Bigr)$. \flexskip

\begin{lemma} The matrix $C$ satisfies the defining conditions (C1)
to (C3) of a generalized Cartan matrix.
\end{lemma}

\begin{proof}
(C1) is clearly satisfied as the scalar product is symmetric. \flexskip

(C2) We begin by considering matrix elements $c_{ij}$ where either $r_i$ or
$r_j$ are candidates for imaginary simple roots selected under (G2).
All candidates for imaginary simple roots are positive multiples of
$\rho = (0,0,1)$. For any element $r=(\lambda,m,n)$ of ${L^*}^+$, the
scalar product satisfies $(r,\rho)=-m\le0$.
Now consider a matrix element $c_{ij}$ where both $r_i$ and $r_j$ are
candidates for real simple roots selected under (G1).
By selection, the candidates are all part of the same fundamental Weyl
chamber of $W^\sigma$. Hence, $(r_i, r_j) \le 0$ for all $i,j$. \flexskip

(C3) Let $r_i$ be a candidate selected under (G1), that is $r_i^2>0$,
and let $r_j$ be any candidate. By selection $(r_i,r_i) =-2(r_i,\rho)$.
Having verified all assumptions about $\Lambda^\sigma$ in lemma 1.3, we may
apply lemma 1.2. Therefore, the simple roots of $W^\sigma$ are all in $L$
(considered as a subset of $L^*$), rather than the whole of $L^*$.
In particular, this holds for $r_j$. Furthermore, lemma 1.2 shows that the
simple roots of $W^\sigma$ are exactly the roots $r\in L$ such that
$(r,\rho)$ is negative and divides $(r,v)$ for all $v\in L$. \flexskip

Hence, we have for $v=r_j$,
\[   
   2 \frac{c_{ij}}{c_{ii}}
          =  2 \frac{(r_i , r_j)}{(r_i, r_i)}
          = -\frac{(r_i , r_j)}{(r_i, \rho)} \in \ZZ. 
\]
\end{proof}

The explicit calculation of all Cartan matrix elements will be carried
out in section 5.1 (theorem 5.3). This will be simplified by a deeper
understanding of the fixed point lattices involved, which we will gain in
chapter 4.1 (which also contains the proof of lemma 1.3). \flexskip

For the present argument, however, the existence of the generalized
Cartan matrix suffices. Starting from the Cartan matrix constructed above,
we may now construct the GKM $G_N=G_N(C)$ with generators $e_i$,
$h_i$, and $f_i$ and relations (R1) to (R5) where $N$ is the order of the
automorphism $\sigma$. Quotient out the null space of the bilinear form, and
denote the resulting Lie algebra ${\cal G}_N$. We may extend the GKMs $G_N$
by the standard set of outer derivations, to obtain a new suite of Lie
algebras ${\cal G}_N^e$. The results of \cite {Jur96} apply and from equation
(1.2) it follows that the Lie algebra ${\cal G}_N^e$ has the decomposition
${\cal G}_N^e = E_N \oplus H_N^e \oplus F_N$. Note that the outer
derivations are part of the null space of the bilinear form. \flexskip

Consider the linearly independent simple roots of ${\cal G}_N^e$ in the
extended root space ${H_N^e}^*$. From theorem 1.3, we obtain a denominator
formula in the extended root space. Recall that we began the construction by
selecting candidates for simple roots in $L^*$. In order to establish that
a specialization to $L^*$ is meaningful we need to check that all
terms remain finite. Let $\pi_0$ again denote the projection
$\pi_0 : {\cal G}_N^e = E_N \oplus H_N^e \oplus F_N \to
         {\cal G}_N   = E_N \oplus H_N   \oplus F_N$, 
defined by quotienting out the null space of the bilinear form.
We will use the term $\pi_0$ also for
the corresponding map of the dual spaces of the Cartan subalgebras, that
is, $\pi_0 : {H_N^e}^* \to {H_N}^*$. By construction, we have
${H_N}^* \tilde= [L^*]$ where $[L^*]=\RR \otimes
L^*$ denotes the span of $L^*$. Recall that, for the duration of the
verification of the specialization, we refer to the roots of ${\cal G}_N^e$
in ${H_N^e}^*$ as `roots'. We denote the (specialized) roots of ${\cal G}_N$
in $L^*$ as `images of roots'.
Let $\pi_\sigma$ denote the projection
$\pi_\sigma : \Lambda \oplus II_{1,1} \to L=\Lambda^\sigma \oplus II_{1,1}$
induced by $\sigma$. We observe that the $II_{1,1}$ components are not
affected by the projection $\pi_\sigma$. \flexskip

The selection rules (G1) and (G2) imply that only finitely many simple
roots are identified by the projection $\pi_0$. For the image of a general
root, recall that ${L^*}^+$ is the set of all $(\lambda,m,n) \in L^* $ such
that $m>0$, or $m=0$ and $n>0$. Consider a general positive root $r$ of
${\cal G}_N^e$ and its decomposition into simple roots $r_i$
\[
\pi_0 r = (\lambda,m,n)= \sum_{i=1}^p \pi_0r_i =
      \sum_{i=1}^p (\lambda_i, m_i, n_i)
\]
where the $r_i$ are $p>0$ not necessarily distinct simple roots
such that $\pi_0 r_i = (\lambda_i, m_i, n_i) \in {L^*}^+$.
The explicit description of the set ${L^*}^+$ yields that $m_i \ge 0$ for
the images of all simple roots under $\pi_0$. If $m_i=0$ then we find
$n_i>0$ directly from the description of ${L^*}^+$. If $m_i>0$ then,
for the images of all real simple roots, the relation
\[
\frac{1}{2} (\lambda_i^2 - 2m_in_i) = \frac{1}{2} (\pi_0 r_i)^2
   = -(\pi_0 r_i, \rho) = m_i
\] 
implies $n_i = \lambda_i^2/2m_i -1 \ge -1$.
These constraints imply that $r$ must be the sum of at most $m$ simple roots
such that $m_i>0$ and at most $m+n$ simple roots such that $m_i=0$.
Conversely, the constraints $m_i \ge 0$ and $n_i \ge -1$ imply that
each simple root $r_i$ must satisfy the conditions that
$m_i \le m$, $n_i < m+n$ and $\lambda_i^2 \le 2m(m+n+1)$.
Hence, for the images of any roots of ${\cal G}_N^e$ in $L^*$ we have shown;

\beginaxioms
\item[($\alpha$)] Only finitely many simple roots are identified by
    the projection $\pi_0$. \flexskip

\item[($\beta$)] For all positive roots $r$, the maximum number of
    summands $p$, within any decomposition, is finite. \flexskip

\item[($\gamma$)] For all positive roots $r$, there exists a finite set
    of simple roots $S_r=\{r_i\}$ so that no $r_j \not\in S_r$  may be
    part of any decomposition of the image of the root $r$.
\endaxioms

Hence, all multiplicities in the specialized denominator formula will
remain finite, and the formula is well defined. Having confirmed that the
specialization is valid we conclude that the specialization of the
denominator formula of the Lie algebra ${\cal G}_N^e$ to $L^*$ is the
denominator formula of the Lie algebra ${\cal G}_N$. By construction, the
specialized denominator formula is the twisted denominator formula (1.20) we
started with. As before, we extend the term GKM to Lie algebras which differ
from a GKM in the strict sense of defining relations (R1) to (R5) in outer
derivations and central extensions and quotients only. For the remainder
of this work, we will make the distinction implicit, and refer to the
{\it specialized} Lie algebra ${\cal G}_N$ with roots in $L^*$ as the GKM
${\cal G}_N$. We have thus proven the following theorem announced in \cite {Bor92}
(the outstanding proof of lemma 1.3 will be given in chapter 4.1). \flexskip

\begin{theorem} Suppose $N$ is any of the primes $2$, $3$, $5$, $7$,
$11$, or $23$, such that $(N+1)$ divides $24$. Suppose that $\sigma\in
\hbox{Aut}(\hat \Lambda)$ is of order $N$ and cycle shape $1^MN^M$ where
$M=\frac{24}{N+1}$. Let $\rho = (0,0,1)$ denote the Weyl vector of the
reflection group $W^\sigma$.
Then there exists a GKM ${\cal G}_N$ with simple roots as follows;

\beginaxioms
\item[(1)]  The real simple roots are the simple roots of the reflection
            group $W^\sigma$ in ${L^*}$, which are the roots $r$ with
            $(r,\rho) = -(r,r)/2$. That is, $\rho$ is also Weyl vector
            for the GKM ${\cal G}_N$. \flexskip

\item[(2)]  The imaginary simple roots are the positive multiples $n\rho$,
            with multiplicities equal to
            $$\hbox{\rm mult}(n\rho)=
              \sum_{a_k\ {\rm divides}\ n} b_k, \ \ \ n\in\ZZ, \ n>0,$$
            if $\sigma$ has generalized cycle shape
            $a_1^{b_1}a_2^{b_2}\dots$.
\endaxioms

The denominator formula of ${\cal G}_N$ is given by (1.20). 
\end{theorem} \theoremproven

This construction of Cartan matrices and GKMs may be performed for more
general automorphisms than the suite of $\sigma$ considered above.
The left hand side of equation (1.20) thus gives the multiplicities of
the positive roots of this algebra. As these multiplicities could possibly be
negative the algebra will, in general, be a Lie superalgebra.

\section{Root Multiplicities}

We continue to use the notation of the previous section. That is, let
$N$ be any one of the primes $2$, $3$, $5$, $7$, $11$, $23$. Let
$\sigma\in \hbox{Aut}(\hat \Lambda)$ be of prime order $N$ and cycle shape
$1^MN^M$ where $M=24/(N+1)$. Let $\Lambda^\sigma$ denote the fixed point
lattice and $L={\Lambda^\sigma} \oplus II_{1,1}$ the corresponding Lorentzian
lattice as above. We will evaluate formula (1.19)
\[
{\rm mult}(r)=
  \sum_{ds\vert ((r,L),N)} \frac{\mu(s)}{ds} \hbox{Tr}(\sigma^d\vert
 \tilde E_{\frac{r}{ds}})
\]
to find explicit formulas for the root multiplicities.
We recall that in equation (1.19) the brackets $( {^.}, {^.})$ denote the
greatest common divisor.
$N$ divides $(r,L)$ if and only if $r\in NL^*$, otherwise we conclude
$\bigl((r,L),N\bigr)=1$. Thus the pair $(d,s)$ can only be one of the following:
$(1,1), (1,N), (N,1)$. The values of the M\"obius function are $\mu(1)=1$ and
$\mu(N)=-1$ because $N$ is prime. Using that the automorphism $\sigma$
satisfies $\sigma^N=\hbox{id}$ we obtain from (1.19)
\begin{equation}
\hbox{mult}(r) = \hbox{Tr}(\sigma\vert\tilde E_r) -
 \frac{1}{N} \hbox{Tr}(\sigma\vert\tilde E_{r/N})+
 \frac{1}{N} \hbox{Tr}(\hbox{id}\vert\tilde E_{r/N})
\end{equation}
We will determine $\hbox{Tr}(\sigma\vert\tilde E_{r})$ and
$\hbox{Tr}(\hbox{id}\vert\tilde E_{r}) =$ dim$(\tilde E_{r})$ for $r\in L^*$
in lemmas 1.6 and 1.7, however we first require\flexskip

\begin{lemma} Let $V,W$ be $\ZZ$-graded vector spaces, $\sigma,\tau$
automorphisms of $V,W$ respectively. Let $V\otimes W$ be graded by the tensor
grading and define the generating function $$\phi_{\sigma,V}(q) = \sum_n 
\hbox{Tr}(\sigma\vert V_n) q^n $$ and analogously for $\tau,W$ and 
${\sigma \otimes\tau,V\otimes W}$. Then \\ $\phi_{\sigma,V}(q)\ \phi_{\tau,W}(q)
=  \phi_{\sigma\otimes\tau,V\otimes W}(q)$.
\end{lemma}

\begin{proof} We observe that $\hbox{deg}(v\otimes w)= \hbox{deg}(v) +
\hbox{deg}(w)$. The proof now is a straightforward exercise. 
\end{proof}

\begin{lemma} $\hbox{Tr}(\sigma\vert \tilde E_r)$ is equal to the
coefficient of $q^{(1-r^2/2)}$ in $q\eta_\sigma(q)^{-1}$ if $r\in L$ and zero
otherwise.
\end{lemma}

\begin{proof} If $r \not\in L$ then it follows from the defining equation
(1.18) that $\sigma$ does not fix the individual roots $(\lambda,m,n)\in L$
which project to $r=(\pi_\sigma\lambda,m,n)$, hence the trace must be zero.
The same argument
proves for any $r\in L$ that $\hbox{Tr}(\sigma\vert \tilde E_r) =
\hbox{Tr}(\sigma\vert E_r)$. If $r=(\lambda,m,n)$ then $E_r=e^\lambda
S_{1+mn-\frac{\lambda^2}{2}}$. (Recall the $\ZZ$-grading of $S$ from formula
(1.11) and the construction of $M_\Lambda$ in section 1.4.)
Let $x_j$ denote the eigenvector of $\sigma$ of eigenvalue $\epsilon_j$
(see equation (1.17)) and let $(i)$ indicate the specific copy of $\Lambda$. 
Define $S^{(j)}=S(\bigoplus_{i>0}{x_j}_{(i)})$. For the
complexified symmetric algebras we obtain 
\[
S=S(\bigoplus_{i>0} \Lambda_{(i)}) = \bigotimes_{j=1}^{24} 
S(\bigoplus_{i>0}{x_j}_{(i)})=\bigotimes_{j=1}^{24} S^{(j)}.
\]
A vector space basis for the factor $S^{(j)}$
in the tensor product is given by all elements of the form
$\prod_k x_{j(m_k)}^{n_k}$. Such an element has $\ZZ$-degree
$\sum_k n_km_k$ and eigenvalue $\epsilon_j^{\sum n_k}$. Thus the generating
function for $S^{(j)}$ is the product over the collection of copies
$\Lambda_{(i)}$:
\[
\phi_{\sigma,S^{(j)}}=\prod_{i>0}\bigl(\sum_{n\geq 0} \epsilon_j^n q^{in}\bigr)
= \prod_{i>0} \sum_{n\geq 0}(\epsilon_j q^i)^n =
    \prod_{i>0}  \frac{1}{1-\epsilon_j q^i}.
\]
Thus the generating function of the trace is
\begin{equation}
\phi_{\sigma,S}(q) = \prod_j \phi_{\sigma,S^{(j)}}(q) =
\prod_j \prod_{i>0} (1-\epsilon_jq^i)^{-1} = q\eta_\sigma(q)^{-1}
\end{equation}
As $1-\frac{r^2}{2}=1+mn-\frac{\lambda^2}{2}$
we can read off the relevant trace as claimed. 
\end{proof}

\begin{lemma} Let $r=(\lambda^*,m,n)\in L^*$ and let
${\lambda^\perp}^*\in {{\Lambda^\sigma}^\perp}^*$ such that
$\lambda^*+{\lambda^\perp}^*\in\Lambda$. The dimension of $\tilde E_r$ is
equal to the coefficient of $q^{1-r^2/2}$ in $\frac{q}{\Delta(q)}
\theta_{{\Lambda^\sigma}^\perp + {\lambda^\perp}^*}(q)$ where
$\theta_{{\Lambda^\sigma}^\perp + {\lambda^\perp}^*}(q)$ is the theta-function
of the translated lattice ${\Lambda^\sigma}^\perp+{\lambda^\perp}^*$.
\end{lemma}

\begin{remark} For a precise definition of $\theta$-functions, see chapter 3
below. We will use the shorter notation $\theta_{{\lambda^\perp}^*}$ in the
sequel. 
\end {remark}

\begin{proof} 
\[
\tilde E_r
= \bigoplus_{\pi_\sigma \lambda = \lambda^*} E_{(\lambda,m,n)}
= \bigoplus_{\pi_\sigma \lambda
= \lambda^*} V_{(\lambda,1+mn)}
\]
\[
= \bigoplus_{\pi_\sigma \lambda = \lambda^*} \RR(\hat \Lambda)_
{\lambda} \otimes S_{1+mn-\frac{\lambda^2}{2}} 
\]
Now $1-\frac{r^2}{2}=1+mn-\frac{(\lambda^*)^2}{2}$. Hence we will find 
$\hbox{dim}(\tilde E_r)$ as the coefficient of $q^{1-\frac{r^2}{2}}$ in 
the series
\[
\sum_p \sum_{\pi_\sigma\lambda=\lambda^*}
\hbox{dim}\bigl(S_{p+\frac{(\lambda^*)^2}{2}-\frac{\lambda^2}{2}}\bigr)q^p=
\sum_{\pi_\sigma\lambda=\lambda^*} \sum_{p'}\hbox{dim}
(S_{p'})q^{p'+\frac{\lambda^2}{2}-\frac{(\lambda^*)^2}{2}}=
\]
\[
 = \sum_{\pi_\sigma\lambda=
     \lambda^*} q^{{\lambda^2\over2}-{(\lambda^*)^2\over2}}\
     \sum_{p'}\hbox{dim}(S_{p'})q^{p'}
=\theta_{{\lambda^\perp}^*} \ {q\over\Delta(q)}.
\]
Here, $\Delta(q)=\eta(q)^{24}$ and $q\over\Delta(q)$ is well known to be the
generating function for the symmetric algebra (it describes the number of
partitions into 24 colours, see \cite {Bor92}, section 12).
This completes the proof of lemma 1.7. 
\end{proof}

Using the modular form $\eta_\sigma$ as defined in equation (1.17) we define
\begin{equation}
\sum_{j>0} p_\sigma(1+j)q^{1+j} = q/\eta_\sigma(q).
\end{equation}
This is a
generalized partition function. Using the Dedekind $\eta$-function we define
for all elements $\tau$ of the upper half plane ${\cal H}$ the function
\[
\psi_j(\tau)=\eta({\tau+j\over N}+j).
\]
Setting $q=e^{2\pi i\tau}$, we also define $\psi_j(q)$ for $|q|<1$. (For 
more details see chapter 2.2, especially equation (2.5) below.) Let
$\delta(r\in L)$ be defined to be of value $1$ if $r\in L$ and of value $0$
otherwise. We are now in the position to state the central theorem 

\begin{theorem} Suppose $N$ is any of the primes 
$2$, $3$, $5$, $7$, $11$, or $23$, such that $(N+1)$ divides $24$. 
Suppose that $\sigma\in \hbox{Aut}(\hat \Lambda)$ is
of order $N$ and cycle shape $1^MN^M$ where $M={24\over N+1}$. Then
\begin{equation}
\theta_{{\lambda^\perp}^*}(q) = \eta(q)^{NM}
\left(\eta(q^N)^{-M}\delta({\lambda^\perp}^*\in {\Lambda^\sigma}^\perp)+
\sum_{0\leq j<N} e^{-j({\lambda^\perp}^*)^2\pi i}\psi_j(q)^{-M}\right)
\end{equation}
and the twisted denominator formula (1.20) has the following explicit form:
\begin{equation}\begin{split}
e^\rho &\prod_{r\in L^+} (1-e^r)^{p_\sigma(1-r^2/2)} \prod_{r\in
N{L^*}^+}(1-e^r)^{p_\sigma(1-r^2/2N)} 
\\ &=\sum_{w\in W^\sigma} \hbox{\rm det}(w)
w\left(e^\rho\prod_{i>0}(1-e^{i\rho})^M(1-e^{Ni\rho})^M\right)
\end{split}\end{equation}
\end{theorem}

\begin{proof} The proof of (1.25) will be given in chapter 4. Assuming (1.25),
we proceed to prove the denominator formulas (1.26). The right hand side
follows straightforwardly from equation (1.20) and the definition of
$\eta_\sigma$. For the left hand side we compare the multiplicities given by
(1.22) with those claimed in the theorem.
We distinguish 4 cases\flexskip

Case 1: $r\not\in L$. Obviously the multiplicities are 0. (Lemma 1.6) \flexskip

Case 2: $r\in L, r\not\in NL^*$. In (1.22) only the first summand is nonzero,
lemma 1.6 proves the claimed formula.\flexskip

Case 3: $r\in NL^*,r\not\in NL$. Suppose $r=(N\lambda^*,m,n)$. We choose
${\lambda^\perp}^*$ as in lemma 1.7.
Under the assumptions of case 3 the first and third term of (1.22) are
nonzero. The first accounts for the exponent $p_\sigma(1-r^2/2)$. By lemma
1.7, dim $\tilde E_{r/N}$ is the coefficient of
$q^{1-{r^2\over2N^2}}$ in ${q\over\Delta(q)}\theta_{{\lambda^\perp}^*}(q)$ or
equivalently the coefficient of $q^{1-{r^2\over2N}}$ in
${q\over\Delta(q^N)}\theta_{{\lambda^\perp}^*}(q^N)$. \flexskip

In case 3, $r/N \not\in L$ and we obtain from (1.25)
\begin{equation}
  {1\over N} {q\over\Delta(q^N)}\theta_{{\lambda^\perp}^*}(q^N) =
  {1\over N} q \eta(q^N)^{-M}\sum_{0\le j<N} e^{-j(r/N)^2\pi i}
\eta(\epsilon^jq)^{-M}.
\end{equation}
Now 
\[
\sum_{0\le j<N} e^{-j({r\over N})^2\pi i} \epsilon^{j{r^2\over2N}} =
\sum_{0\le j<N} e^{-j({r\over N})^2\pi i}
   \left(e^{2\pi i\over N}\right)^{j{r^2\over2N}}= N.
\]
Thus the coefficient of $q^{1-r^2/2N}$ in (1.27) is equal to that of
$q\eta(q^N)^{-M}\eta(q)^{-M} = q\eta_\sigma^{-1}(q)$. This completes case 3.
\flexskip

Case 4: $r\in NL$. Here all three terms of (1.22) are nonzero.
As ${r\over N}\in L$, we have an additional term in ${q\over\Delta}\theta$,
namely $q\eta(q^N)^{-M} \eta(q^{N^2})^{-M}=q\eta_\sigma^{-1}(q^N)$. Thus it
cancels exactly with the contribution of $\hbox{Tr}(\sigma\vert \tilde E_{r/N})$.
Now the argument is the same as in case 3.
\end{proof}

\begin{corollary} Let $N$ be any of $2$, $3$, $5$, $7$, $11$, $23$. Then 
the generalized Kac-Moody algebras ${\cal G}_N$ constructed in theorem 1.6
have root lattice $L=\Lambda^\sigma \oplus II_{1,1}$. Their
root multiplicities are functions of the norm of the roots as follows:
\[
{\rm mult}(r)=p_\sigma(1-{r^2\over2})\mskip 40mu r\in L,\ r\not\in NL^*
\]
 and
\[
{\rm mult}(r)=p_\sigma(1-{r^2\over2})+p_\sigma(1-{r^2\over2N})\mskip 40mu
r\in NL^*.
\]
\end{corollary}

\begin{proof} The explicit multiplicities for any root follow by
comparison of theorem 1.7, equation (1.26) with the general denominator
formula theorem 1.3. 
\end{proof}

We conclude the chapter with a number of remarks:\flexskip

1.) In chapter 5 below we will identify the real simple roots of the 
    GKM ${\cal G}_N$, as the following:
    \[
    \bigl(\lambda,1,{\lambda^2\over2}-1\bigr) \mskip 30mu {\rm for}  
    \mskip 30mu \lambda\in \Lambda^\sigma
    \]
    and
    \[
    \bigl(\lambda,N,{\lambda^2\over2N}-1\bigr)
    \mskip 30mu {\rm for} \mskip 30mu \lambda\in N{\Lambda^\sigma}^* 
    \mskip 15mu {\rm such}\ \ {\rm that}\mskip 15mu N|\ \bigl(
    {\lambda^2\over 2N}-1\bigr).
    \]
    The former real simple roots have height $1$ and norm $2$ whereas the
    latter have height $N$ and norm $2N$.
    The imaginary simple roots have already been identified explicitly in
    theorem 1.6.    \flexskip

2.) Equation (1.26) can be read as a new combinatorial identity, independent
    of the theory of generalized Kac-Moody algebras, in the same way that the
    denominator formulas of affine Kac-Moody algebras give the Macdonald
    identities. \flexskip

3.) The explicit root multiplicities form upper bounds for the root
    multiplicities of any subalgebras of the ${\cal G}_N$. This will be used
    in chapter 6. \flexskip

4.) Finally, we will briefly discuss the limitations 
of the main methods of constructing GKMs outlined in this chapter.
This will emphasize the pivotal role of the no-ghost theorem in the
argument. 
The method of construction may be summarised in the following diagram
where $V$ symbolizes the construction of vertex algebra and physical
spaces, $K$ symbolises the operation of quotienting out the kernel of the
bilinear form, and $\sigma$ the application of the projection defined by the
automorphism of order $N$:
\[
  \Lambda \mskip 30mu 
  {V \atop \longmapsto} 
  \mskip 30mu P^1/DP^0 \mskip 30mu 
  {K \atop \longmapsto} 
  \mskip 30mu M_\Lambda \mskip 30mu 
  {\sigma \atop \longmapsto} 
  \mskip 30mu {\cal G}_N.
\]
For the fake monster Lie algebra, \cite {Bor92} presents the following results
\flexskip

a) the Cartan matrix (identifying generators and relations), \flexskip

b) all root multiplicities, \flexskip

c) an explicit construction (as subspace $V_{1+mn}$ of a vertex algebra). 
\flexskip

We have seen in section 1.3 above that the result c) considerably simplifies
the identification of Lie algebra elements, allowing to calculate explicit
root multiplicities (theorem 1.4). This work provides results analogous to
a) (all elements of the generalized Cartan matrix will be calculated in
theorem 5.3) and b) (corollary to theorem 1.7) for the new family of GKMs
${\cal G}_N$. We do not know a natural construction for these GKMs. \flexskip

It is natural to try and apply our construction method to other lattices, 
in particular, we may consider the sequence of operations symbolised in the 
following diagram;
\[
  \Lambda \mskip 30mu 
  {\sigma \atop \longmapsto} 
  \mskip 30mu \Lambda^\sigma \mskip 30mu 
  {V \atop \longmapsto} 
  \mskip 30mu P^1(N)/DP^0(N) \mskip 30mu 
  {K \atop \longmapsto} 
  \mskip 30mu M_{\Lambda^\sigma}(N).
\]
Consider the Lie algebras $P^1(N)/DP^0(N)$ first. Theorem 1.4, above, applies 
and provides us with the root multiplicities of these Lie algebras. 
Turning to $M_{\Lambda^\sigma}(N)$, \cite {Bor92} reports that the bilinear form 
can be shown to be almost positive definite 
(the kernel $K$ is empty if the dimension of the 
lattice $\Lambda^\sigma$ is less than 24), which makes this Lie algebra a 
GKM. However, the no-ghost theorem 1.5 does not apply. Therefore, in this 
case we have neither an explicit construction of the simple roots, nor can 
we express any multiplicity formulas for the simple roots. \flexskip

Neither  ${\cal G}_N \subset M_{\Lambda^\sigma}(N)$  nor
$M_{\Lambda^\sigma}(N) \subset {\cal G}_N$ hold. The first relation can 
be disproven by considering that ${\cal G}_N$ contains real roots 
of norm $2N$. Hence, their $\ZZ$-grading is $N$, and they will not 
be elements of $P^1(N)$. Consequently, they will not be element of 
$M_{\Lambda^\sigma}(N)$ either. \flexskip

For the second relation, we consider the case
$N=23$. The hyperbolic algebra $AE_3$, as defined in the introduction,
can be embedded into both ${\cal G}_{23}$ and 
$M_{\Lambda^\sigma}(23)$. (For details see chapter 6.2.1, below.) 
We consider its unique root of norm $-2$ (up to Weyl automorphisms), 
which may be expressed in terms of simple roots as $r=(2,2,1)$. Its  
multiplicity in $AE_3$ is $2$. The corollary to theorem 1.7 shows that the
multiplicity of $r$ in ${\cal G}_{23}$ is also $2$. However, 
${\rm dim}(P^1_r) = p_2(2) = 5$, and ${\rm dim}(P^0_r )= p_2(1) = 2$. 
Section 1.3 identifies the explicit vertex algebra bases for the space 
$P^1_r$, and its bilinear form. Using these, it is straightforward to 
calculate that the kernel of the bilinear form (1.13) on $P^1_r$ 
(paired with $P^1_{-r}$) is equal to the 
2-dimensional $DP^0_r$. Hence the multiplicity of $r$ in 
$M_{\Lambda^\sigma}(23)$ is 3. \flexskip

    The considerations of the previous paragraph do not only 
    disprove any inclusions but also indicate that upper bounds for the 
    root multiplicites of $AE_3$ provided by $M_{\Lambda^\sigma}(23)$
    will be inferior to those of ${\cal G}_{23}$.


\chapter{Modular Forms}

This chapter is the first in a series of three chapters which give the proof of
the central theorem 1.7. In section 2.1 we recall the definitions of the
modular group, its subgroups, and of modular forms. We proceed to recall some
properties of modular forms. The material can be found in any standard
treatment of modular forms and it has been included here solely to set up the
notation. Section 2.2 defines a number of modular forms related to
the $\eta$-function. Their modularity properties are fully
understood in principle. We will, however, need to know their exact
characters explicitly in chapter 4. Therefore, section 2.2 will establish the
transformation behaviour under the generators $S$ and $T$ of the modular group.

\section{Review of Modular Group and Modular Forms}

As all the material included in this section is standard textbook material
we will not always give specific references.
Let $\Gamma = SL_2(\ZZ)$ denote the modular group, $\Gamma(N)$ its subgroup
consisting of all elements which are identical to the identity
matrix modulo $N$ and
\[
\Gamma_0(N) = \{ 
\begin{pmatrix} a&b \\ c&d \end{pmatrix} 
\in \Gamma \mskip 10mu \Bigl\arrowvert \mskip 10mu c\equiv0\ (\hbox{mod}\ N) \}.
\]
Recall that $\Gamma(N)$ and $\Gamma_0(N)$ are subgroups of the modular group of
finite index. Furthermore let
\[ 
S = 
\begin{pmatrix} 0&-1 \\ 1&0 \end{pmatrix}, \mskip 50mu T = 
\begin{pmatrix} 1&1 \\ 0&1 \end{pmatrix}, \mskip 50mu F = 
\begin{pmatrix} 0&-1 \\ N&0 \end{pmatrix}.
\]
$F$ is called the Fricke involution.
Then $\Gamma$ is generated by $T$ and $S$, $\Gamma_0(2)$ is generated by $T$
and $ST^2S$, and for $N\neq2$ it is known that $\Gamma_0(N)$ is generated by
$T$ and the collection of
\begin{equation}
V_k = S T^k S T^{k'} S = 
\begin{pmatrix} -k'&1 \\ kk'-1&-k \end{pmatrix}.
\end{equation}
Here $k=1,...,N-1$ and $k'$ is any integer such that $kk'\equiv1\ (\hbox{mod}\
N)$ (cf. \cite {Apo76}, theorem 4.3 or \cite {Rad29}). $F^{-1}={-1\over N}F$ and
\[
F\begin{pmatrix} a&b \\ Nc&d \end{pmatrix}
F^{-1} \ =\ \begin{pmatrix} d&-c \\ -Nb&a \end{pmatrix}
\]
In particular, if $kk' = N+1$ then $FV_kF^{-1} = - V_{-k'}$. Let ${\cal H}$
denote the upper half complex plane, that is the set of all $\tau\in \CC$ with
positive imaginary part.
$\Gamma$ acts on ${\cal H}$ by ${{ab} \choose {cd}}(\tau) =
{{a\tau+b}\over{c\tau+d}}$. The Fricke involution acts in the same way. Note
that $SF(\tau) = N\tau$ and $FS(\tau) = {\tau\over N}$. \flexskip

The metaplectic group $Mp$ is a double cover of $\Gamma$ consisting of elements
$(A,j)$ where $A\ =\ {{ab} \choose {cd}} \in \Gamma$ and $j$ is a holomorphic
function in $\tau \in {\cal H}$ such that $j(\tau)^2 = c\tau +d$. The group
multiplication is defined to be $\bigl(A,j(\tau))(B,k(\tau)\bigr) =
\bigl(AB,j(B\tau)k(\tau)\bigr)$.
For any $k \in {1\over2}\ZZ$ the metaplectic group then
acts on the space of functions on ${\cal H}$ by
\begin{equation}
f\vert_{(A,j,k)} (\tau) = j(\tau)^{-2k} f(A(\tau)).
\end{equation}
In particular, we will write $(\tau)^{1/2}_+$ to denote the branch of the
square root that is positive for positive $\tau$. $(S,+)$ will be the
corresponding element of the metaplectic group using this branch.
$(T,+)$ uses the constant square root of $+1$. Note that $(S,+)^2=(-1,i)$. We
will sometimes write $(G,*)$ when it is not required to specify the branch.
Note further that $f^{-1}\vert_{(A,j,k)} = \bigl( f\vert_{(A,j,-k)} \bigr)^{-1}$.
Let $\Gamma'$ be a subgroup of $\Gamma$ of finite
index. An orbit in $\QQ \cup i\infty$ under $\Gamma'$ is called a cusp.
Let $\Gamma'$ be a subgroup of finite index of $\Gamma$ of level $N$, that is,
$\Gamma'\ \supset \ \Gamma(N)$.
A function $f:\ {\cal H}\ \to\ \CC$ is called a modular
form of weight $k$ and multiplier system $\chi$ for $\Gamma'$ if $f$ is
holomorphic on ${\cal H}$ and $f\vert_{(A,j,k)} = \chi(A,j) f$ for all
$A={{ab}\choose{cd}} \in \Gamma'$ and $\tau\in{\cal H}$.
We can substitute $q=e^{2\pi i\tau}$, such that $|q|<1$.
$f$ is called holomorphic modular form if in addition at each cusp of
$\Gamma'$ $f$ has a $q$-expansion without negative powers. \flexskip

We can now define the Dedekind eta-function. For $\tau\in{\cal H}$, let
\begin{equation}
\eta(\tau) =
e^{2\pi i\tau\over24} \prod^\infty_1 (1-e^{2\pi in\tau}).
\end{equation}
In terms of $q$ we obtain the simpler form
\begin{equation*} \tag{2.3$'$}
\eta(q) = q^{1\over24} \prod^\infty_1 (1-q^n). 
\end{equation*}
Note that (2.3$'$) is well defined only
by the additional prescription $q^{1\over24} = e^{2\pi i\tau/24}$. This will be
implicitly understood whenever either variable $q$ or $\tau$ is used in this
work. We will also consider $\Delta(q) = \eta^{24}(q)$.
It is well known that $\eta(S\tau) = e^{-\pi i/4} (\tau)^{1/2}_+\eta(\tau)$ and
$\eta(T\tau) = e^{\pi i/12} \eta(\tau)$. Thus $\eta$ is a modular form of
weight $1\over2$ and $\Delta$ is a
holomorphic modular form of $\Gamma$ of weight 12 and trivial multiplier
system. Furthermore, $\Delta$ has a zero of first order at $q=0$ and has no
other zeros for $\vert q\vert < 1$.   \flexskip

Modular forms are uniquely determined by a certain number of initial
coefficients. We will in particular use the following theorem.\flexskip

\begin{theorem} Let $\Gamma'$ be a congruence subgroup of the
modular group $\Gamma$ of finite index. Let $f$ be a holomorphic modular form
of $\Gamma'$ of integer weight $k$ and a character $\chi$ such that
$\chi^{12}(g)=1$ for all $g\in\Gamma'$. Suppose that
the $q$-expansions of $f(\tau)$ at all cusps do not contain any terms
of $q$ of exponent $\leq {k\over12}$. Then $f\equiv0$. 
\end{theorem}

\begin{proof} $f^{12}$ is holomorphic in ${\cal H}$ and has trivial multiplier.
By assumption its expansion at any cusp
does not contain terms of exponent less or equal $k$.
The $q$-expansion of $\Delta^k$ at any cusp has a $k^{\hbox{th}}$ order zero,
hence ${f^{12}\over\Delta^k}$ is a holomorphic modular form with trivial
multiplier of $\Gamma'$.
Its weight is 0. However, the only such forms are constant,
cf. \cite {Kob84}, III.3.Prop.18. As the constant term in the $q$-expansion was
assumed to be zero we conclude $f\equiv0$. 
\end{proof}

\begin{remark} The proof shows that a similar theorem holds for forms of 
any character of finite order.
\end{remark}

\begin{corollary} Let $N$ be a prime and $\Gamma' = \Gamma_0(N)$. Let $f$
be as above. If the $q$-expansions of both $f$ and $f\vert_S$ at zero do 
not contain any terms of exponent $\leq {k\over12}$ then $f\equiv0$.
\end{corollary}

\begin{proof} $\Gamma_0(N)$ has exactly two cusps, represented by 0 and $\infty$,
cf \cite {Shi71}, p.26. As $S$ interchanges these two cusps the given
criterion is equivalent to that of the theorem. 
\end{proof}

\section{Some Modular Forms Related to Eta}

For the remainder of this chapter let $N$ be an integer such that $N+1$ divides 
24, and let $M={24\over N+1}$.
We will need to know the transformation properties under the modular group for
a number of functions related to the $\eta$-function. We begin by considering
$\eta(N\tau)$, which corresponds to $\eta(q^N)$ in terms of the argument $q$.
Then 
\begin{subequations}
\begin{equation}
\eta(N\tau)\vert_{\scriptstyle(S,+,{1\over2})} =
\left( (\tau)^{1/2}_+ \right)^{-1} \eta(NS\tau) =
\left( (\tau)^{1/2}_+ \right)^{-1} \eta(S({\tau \over N}))
 = {e^{-\pi i/4} \over \sqrt{N}} \eta({\tau \over N}) 
\end{equation} 
\begin{equation}
\eta(N\tau)\vert_{\scriptstyle(T,+,{1\over2})} = \eta(NT\tau)
= \eta(T^N(N\tau)) = e^{N\pi i\over 12} \eta(N\tau)
\end{equation} 
\end{subequations}
Now let 
\begin{equation}
\psi_j(\tau) = \eta({{\tau+j} \over N}+j)
\end{equation}
which corresponds to $\psi_j(q) = \eta(\epsilon^j q^{1/N})$ in somewhat 
ambiguous notation where $\epsilon$ is an $N^{\hbox{th}}$ root of unity. Note 
that these are the $\psi_j$ as defined in chapter 1, and used in theorem 1.7.
\begin{equation}
\psi_j\vert_{\scriptstyle(T,+,{1\over2})}(\tau) = \psi_j(T\tau) =
\eta({{\tau+(1+j)}\over N}+(1+j)-1) = e^{-\pi i/12} \psi_{j+1}(\tau).
\end{equation}
To determine $\psi_j^M\vert_{\scriptstyle(S,+,{1\over2})}$ we first notice that
${\tau + (j+N)\over N} + (j+N) = \left( {\tau + j \over N} + j \right) +(N+1)$.
Now $\left( e^{\pi i/12} \right)^{(N+1)M} = 1$ and thus
$\psi_{j+N}^M(\tau) =  \psi_j^M(\tau).$
Given $j$, choose $j'$ such that $jj'\equiv1(\hbox{mod}N)$ and define
\begin{equation}
G=-FV_jF^{-1}=-FST^{j}ST^{j'}SF^{-1}
= \begin{pmatrix} j& \frac{jj'-1}{N} \\ N&j' \end{pmatrix} \in \Gamma_0(N).
\end{equation}
The action of $G$ on $\tau$ equals that of $-G$, hence we obtain that
\[
{S\tau+j \over N }+j = T^jFST^jS\tau =
 T^jGT^{j'}T^{-j'}FST^{-j'}\tau = T^jGT^{j'}\left({\tau-j'\over N}-j'\right)
\]
Note furthermore that the holomorphic function $j$ associated to the element
$G$ of equation (2.7) satisfies $j^2({\tau-j'\over N}) = \tau$.
Hence 
\begin{equation}\begin{split}
  \psi_j\vert_{\scriptstyle(S,+,{1\over2})}(\tau) 
  &= \left( (\tau)^{1/2}_+\right)^{-1} \psi_j(S\tau) =
   \left( (\tau)^{1/2}_+\right)^{-1} \eta\left(T^jG({\tau-j'\over N})\right)=\\
  &= \left( (\tau)^{1/2}_+\right)^{-1} e^{j\pi i/12}
   \eta\left(G({\tau-j'\over N})\right) =
   e^{j\pi i/12} \chi(G) \eta({\tau-j'\over N}) = \\
  &= e^{(j+j')\pi i/12} \chi(G) \psi_{-j'}(\tau).
\end{split}\end{equation}
The character $\chi(G)$ will have to be calculated case by case using the
following two lemmata. \flexskip

\begin{lemma} Suppose that, for the element $G$ of equation (2.7), \\
$G=ST^aST^bST^cS$. Then 
\[ 
\psi_j\vert_{\scriptstyle (S,+,\frac{1}{2})}(\tau) =
 e^{-\pi i} e^{(j+j'+a+b+c)\pi i/12} i^J \psi_{-j'}(\tau)
\]
Here $J=J(a,b,c)$ is an integer that will be specified in the proof.
\end{lemma}
\begin{proof} Let $y={{\tau-j'}\over{N}}$. Following the line of
argument in formula (2.8) we find
\begin{equation*} \begin{split}
  \psi_j \vert_{\scriptstyle(S,+,{1\over2})}(\tau) 
  = &\left( (\tau)^{1/2}_+\right)^{-1} e^{j\pi i/12}
     \eta(ST^aST^bST^cS(y)) = \\
  = &\left( (\tau)^{1/2}_+\right)^{-1} e^{j\pi i/12} 
     e^{-4\pi i/4} e^{(a+b+c)\pi i/12} * 
     \Bigl(y\Bigr)^{1/2}_+ \Bigl( {{-1}\over y}+c\Bigr)^{1/2}_+ \\
    &\left({{-1}\over {{{-1}\over y}+c}}+b\right)^{1/2}_+
     \left({{-1}\over{{{-1}\over {{{-1}\over y}+c}}+b}}+a\right)^{1/2}_+
     e^{j'\pi i/12} \psi_{-j'}(\tau)
\end{split}\end{equation*}
The factors containing $\tau$ or $y$ other than $\eta$ itself
multiply to a constant as already seen in equation (2.8) above.
To determine the correct branch of the square root it is sufficient to
evaluate for $\tau=R+{i\over R}$ where $R$ is a large positive real. Then
$1\over y$ becomes arbitrarily small and it suffices to evaluate
\[
 \Bigl(c\Bigr)^{1/2}_+ \Bigl({{-1}\over {c}}+b\Bigr)^{1/2}_+
	\Bigl({{-1}\over{{{-1}\over c}+b}}+a\Bigr)^{1/2}_+
   = \Bigl(c\Bigr)^{1/2}_+ \Bigl({{bc-1}\over c}\Bigr)^{1/2}_+
	\Bigl({{abc-a-c}\over {bc-1}}\Bigr)^{1/2}_+.
\]
Note that $a,b,c$ are integers and in particular real.
If $J$ is the number of changes of sign in the sequence $1,c,bc-1,abc-a-c=N$
then the correct branch is $i^J$. This completes and proves the
claim. 
\end{proof}

\begin{lemma} Suppose that, in the above notation, $jj'=N+1$. Then
$G=ST^{-j'}ST^{-j}S$ and 
\[
\psi_j\vert_{\scriptstyle(S,+,{1\over2})}(\tau) =
\hbox{\rm sign}(-j') e^{-3\pi i/4} \psi_{-j'}(\tau).
\]
\end{lemma}
\begin{proof} The argument is a simplification of the proof of lemma 2.1. Here we
have to count the number of sign changes in the sequence $1,-j',jj'-1=N$.
\end{proof}

With these tools we can now describe the transformation of ${\psi_j}^M$ for
any $j$ and $N,M$ as above. \flexskip

\begin{theorem} Let $N$ be any of $2$, $3$, $5$, $11$. Then
\[
   \psi_j^M\vert_{\scriptstyle(S,+,{M\over2})}(\tau) =
    e^{{-3\pi i\over4}M} \psi_{-j'}^M(\tau).
\]
\end{theorem}
\begin{proof} Claim: We can use lemma 2.2 for all $j$ and $N$.
The cases $j=1$ ($G=-ST^{-N}S$) and $j=-1$ ($G=ST^NS$) are trivial.
We recall that we only need to consider $j$ modulo $N$.
Furthermore, as $jj'=(-j)(-j')$ it suffices to consider $j=2,...,{N-1\over2}$.
However, if $N=5$, $2*3=5+1$, and if $N=11$ we have the decompositions
$jj'=2*6=3*4=(-6)*(-2)=11+1$ (note that $-6 \equiv 5$ mod 11).
Application of lemma 2.2 now gives the result as the factor
sign$(-j')$ is cancelled by the even power $M$. 
\end{proof}

\begin{theorem} Let $N=7$ or $N=23$. Then
\[  
    \psi_j^M\vert_{\scriptstyle(S,+,{M\over2})}(\tau) =
    e^{{-3\pi i\over4}M} \left({-j'\over N}\right)
    \psi_{-j'}^M(\tau).
\]
Here $\left({x\over y}\right)$ denotes the Legendre symbol.
\end{theorem}

\begin{proof} If $N=7$ we use the decompositions $jj'=2*4=(-4)*(-2)=7+1$ where
$3\equiv -4$ mod 7. Squares modulo 7 are 1,2,4.
Hence the claim follows from Lemma 2.2. \flexskip

Now let $N=23$. The squares modulo $23$ are $1,2,3,4,6,8,9,12,13,16,18$.
The decompositions $jj'=1*24=2*12=3*8=4*6=23+1$ show which $j$ are covered by
Lemma 2.2. For the cases of the remaining $j$ first observe that if
\[
   \psi_j^M\vert_{\scriptstyle(S,+,{M\over2})} =
   e^{{-3\pi i\over4}M} \left({-j'\over N}\right) \psi_{-j'}^M 
\]
then
\begin{equation*}\begin{split}
   \psi_{-j'}^M\vert_{\scriptstyle(S,+,{M\over2})} 
   &= e^{{3\pi i\over4}M}
      \left({-j'\over N}\right) \psi_j^M\vert_{\scriptstyle(-1,i,{M\over2})}
      = e^{{3\pi i\over4}M} (-1)^M \left({j\over N}\right) i^{-M}\psi_j^M= \\
   &= e^{{-3\pi i\over4}M} \left({j\over N}\right) \psi_j^M. 
\end{split}\end{equation*}
(Note that $\left({j\over N}\right) = \left({j'\over N}\right)$ as the Legendre
symbol is multiplicative.) Thus
we are only left to consider separately $j=5,-5,7,-7$ as we obtain
$j'=-9,9,10,-10$.\flexskip

Now $G=-FST^5ST^{-9}SF^{-1} = ST^{-5}ST^{-2}ST^{2}S$ hence by Lemma 2.1 we
obtain
\begin{equation*}\begin{split}
   \psi_5\vert_{\scriptstyle(S,+,{1\over2})}(\tau) 
   &= e^{-\pi i} e^{(5+(-9)-5-2+2)\pi i/12} i^{J(-5,-2,2)} \psi_9(\tau)\\
   &= e^{-\pi i} e^{-9\pi i/12} (-1)\psi_9(\tau)
      = e^{-3\pi i/4} \psi_9(\tau)
\end{split}\end{equation*}
As $9$ is a square modulo $23$ this is the claimed transformation. \flexskip

The other cases are analogous, and it suffices to give the factor
\[
e^{-\pi i} e^{(j+j'+a+b+c)\pi i/12} i^{J(a,b,c)}.
\]
For $(j,j')=(-5,9)$ we obtain $(a,b,c)=(4,-2,-3)$ and thus a factor
of 
\[e^{-\pi i}e^{(-5+9+4-2-3)\pi i/12}(-1) = e^{1\pi i/4} =
e^{-3\pi i/4} \left({-9\over 23}\right). 
\]
For $(j,j')=(7,10)$ we obtain $(a,b,c)=(-3,3,-2)$ and thus a factor
of 
\[
e^{-\pi i}e^{(7+10-3+3-2)\pi i/12}(-1) = e^{5\pi i/4} =
e^{-3\pi i/4} \left({-10\over 23}\right). 
\]
For $(j,j')=(-7,13)$ we obtain $(a,b,c)=(3,-4,-2)$ and thus a factor
of 
\[
e^{-\pi i}e^{(-7+13+3-4-2)\pi i/12}(-1) = e^{1\pi i/4} =
e^{-3\pi i/4} \left({-13\over 23}\right). 
\]
This concludes the proof of theorem 2.3.
\end{proof}

\chapter{Lattices and their Theta-Functions}

This chapter establishes a number of properties of integral lattices which will
be used in chapter 4. The lattices in question are those which are part of
formula (1.25). Section 3.1 quotes a number of results concerning the dual of
a sublattice of a self-dual lattice, and modularity properties of
$\theta$-functions under the generators $S$ and $T$ of the modular group.
We then restrict our attention to the particular type of lattice that arises
in (1.25). Again, the modularity properties of their $\theta$-functions are
fully understood in principle. Section 3.2 establishes the explicit
modularity properties (including characters) of these $\theta$-functions.
Please note that in the notation of chapter 3, the letter $L$ denotes a
more general object than in chapter 1.

\section{Review of Results about Lattices} \flexskip

For any lattice $L$, let $[L]=L\otimes \RR$ denote the vector space spanned
by the elements of $L$. Let $\Lambda$ again denote the 24-dimensional
unimodular self-dual Leech lattice. For any sublattice $L$ of $\Lambda$,
we define the orthogonal projection $\pi_L : [\Lambda] \to [L]$. We begin our
review with a well known description for the dual of some lattices. \flexskip

\begin{theorem} Let $L$ be a sublattice of $\Lambda$ such that
$\Lambda \cap [L] = L$. Define the projection $\pi_L : [\Lambda] \to [L]$ as
above. Then $L^* = \pi_L \Lambda$. 
\end{theorem}

\begin{proof} The projection $\pi_L(\Lambda)$ is contained in the dual $L^*$
because, for any $\lambda \in \Lambda$ and any $\mu \in L$, the
inner products satisfy $(\pi_L\lambda,\mu)=(\lambda,\mu)$.
The other inclusion will be proved by induction. \flexskip

For any proper sublattice $L_r$ of $\Lambda$ such that $L_r=\Lambda \cap [L_r]$
we can choose $\lambda_{r+1} \in \Lambda$ such that
$L_{r+1} := L_r \oplus \ZZ\lambda_{r+1}$ as lattices and that
$L_{r+1} = \Lambda \cap [L_{r+1}]$. Below, we will show how, given
$\lambda_r^* \in L_r^*$, to find $\lambda_{r+1}^* \in L_{r+1}^*$ such that
$\pi_{L_r}(\lambda_{r+1}^*) = \lambda_r^*$.
The proof can then be completed as follows: Given $\lambda^* \in L^*$ we choose
$r_0=$ dim $L$, $L_{r_0}=L$, $\lambda_{r_0}^* = \lambda^*$ and continue
inductively until $L_{r_1}=\Lambda$. Then $\lambda_{r_1}^*\in \Lambda^* =
\Lambda$ and hence $\lambda^* = \pi_L(\lambda_{r_1}^*) \in \pi_L(\Lambda)$
which finishes the proof. \flexskip

We now construct $\lambda_{r+1}^*$ for given $\lambda_r^* \in L_r^*$.
Suppose $L_r$ as a lattice has a basis ${\lambda_1,...,\lambda_r}$. We choose
rational $a_i, i=1,...,r$ such that
$\lambda_{r+1}-\sum_{i=1}^r a_i\lambda_i \in [L_r^\perp]$.
(This is possible because $L$ is an integer lattice.) We use the ansatz
$$\lambda_{r+1}^* = \lambda_r^* + a\left( \lambda_{r+1}-\sum_{i=1}^r a_i
\lambda_i \right)$$ with $a\in\QQ$. We then have to satisfy the condition that
the inner product of this expression with any element of $L_{r+1}$ is integer.
However, the inner product with any element of $L_r$ is trivially integer. Thus
we complete the proof by choosing $a$ such that the inner product with
$\lambda_{r+1}$ is integer.
\end{proof}

We also quote the following fact from \cite {CS88} (chapter 4, theorem 1):
\flexskip

\begin{theorem} Let $L$ be a sublattice of $\Lambda$.
Let $L^\perp$ denote the orthogonal complement of $L$ within $\Lambda$.
Then the abelian group $L^*/L$ is isomorphic to the group
$(L^\perp)^*/(L^\perp)$.
\end{theorem}

Obviously theorems 3.1 and 3.2 remain true if we replace $\Lambda$ by
any self-dual lattice. Let $\theta_{L+\lambda^*}$ denote the theta-function of
the lattice $L$ translated by $\lambda^*$, that is
\begin{equation}
\theta_{L+\lambda^*} (q) = \sum_{\lambda \in L+\lambda^*}
				   q^{(\lambda,\lambda)/2}.
\end{equation}
As before, non-integer powers of $q=e^{2\pi i\tau}$ will be understood to be
defined by $q^x=e^{2\pi i\tau x}$. We can now quote the following theorem which
goes back to Jacobi and is a special case of theorem 13.5 of \cite {Kac90}:

\begin{theorem} Let $L$ be an $l$-dimensional integral lattice,
let $\lambda^*$ be an element of $L^*$. Then
$$\theta_{L+\lambda^*}\vert_{\scriptstyle (S,+,{l\over2})} (q) =
  {e^{-l\pi i/4} \over \sqrt{\vert L^*/L \vert}}
  \sum_{\mu^* \in L^*/L} e^{-2\pi i(\lambda^*,\mu^*)} \theta_{L+\mu^*}(q)$$
and $$\theta_{L+\lambda^*}\vert_{\scriptstyle (T,1,{l\over2})} (q) =
e^{\pi i(\lambda^*,\lambda^*)} \theta_{L+\lambda^*}(q).$$
Furthermore, the action of the modular group on the vector space spanned by the
functions $\theta_{L+\lambda^*}, \lambda^* \in L^*/L$ is unitary.
\end{theorem}

\begin{corollary} For $\lambda^*=0$, 
\[
\theta_L\vert_{\scriptstyle(S,+,{l\over2})}
(q) ={e^{-l\pi i/4} \over \sqrt{\vert L^*/L \vert}} \theta_{L^*}(q).
\]
\end{corollary}

\section{The Character of Theta}

Let $L$ be an even $l$-dimensional integer lattice with dual $L^*$.
Let $N$ be the least positive integer such that $NL^*\subseteq L$. Suppose
that $N$ is prime but not equal to 2. We calculate the
transformation of $\theta_L$ under the generators of $\Gamma_0(N)$. Clearly
$\theta_L$ is invariant under $(T,1)$. For $(V_k,*)=(ST^kST^{k'}S,*)$
(cf equation (2.1)) we obtain
\begin{equation*} \begin{split}
  \theta_L&\vert_{\scriptstyle
     \bigl((S,+)(T^k,1)(S,+)(T^{k'},1)(S,+),{l\over2}\bigr)} \\
  &= {e^{-l\pi i/4} \over \sqrt{\vert L^*/L \vert}} \sum_{\lambda \in L^*/L}
     \theta_{L+\lambda}\vert_{\scriptstyle
     \bigl((T^k,1)(S,+)(T^{k'},1)(S,+),{l\over2}\bigr)} =  \\
  &= \dots ={e^{-3l\pi i/4} \over {\sqrt{\vert L^*/L \vert}}^3}
     \sum_{\lambda,\mu,\nu \in L^*/L}
     e^{\pi ik(\lambda,\lambda)-2\pi i(\lambda,\mu)+\pi ik'(\mu,\mu)
     -2\pi i(\mu,\nu)} \theta_{L+\nu} \\
  &=:  \sum_{\nu \in L^*/L} c_\nu(k) \theta_{L+\nu},
\end{split} \end{equation*}
where this defines coefficients $c_\nu(k)$. We first evaluate the coefficient
$c_0(k)$. Note that because $kk' \equiv 1$ (mod $N$)
\[
k(\lambda-k'\mu)^2 = k\lambda^2 -2kk'(\lambda,\mu)+k(k')^2\mu^2 \equiv
k\lambda^2 -2(\lambda,\mu)+k'\mu^2 (\hbox{mod\ }2\ZZ).
\]
Hence 
\begin{equation*} \begin{split}
  c_0(k) 
  &= {e^{-3l\pi i/4} \over {\sqrt{\vert L^*/L \vert}}^3} 
     \sum_{\lambda,\mu \in L^*/L} e^{\pi ik(\lambda,\lambda)-2\pi 
     i(\lambda,\mu)+\pi ik'(\mu,\mu)} = \\
  &= {e^{-3l\pi i/4} \over {\sqrt{\vert L^*/L \vert}}^3} 
     \sum_{\lambda,\mu \in L^*/L} e^{\pi ik(\lambda-k'\mu)^2}
\end{split} \end{equation*}

Now the sum runs through all elements of $L^*/L$ exactly
$\vert L^*/L \vert$ times. Thus
\begin{equation}
  c_0(k) = {e^{-3l\pi i/4} \over \sqrt{\vert L^*/L \vert}} \sum_{\lambda
  \in L^*/L} e^{\pi ik\lambda^2}
\end{equation}
$L^*/L$ is an abelian group of, say, $n$ generators $x_i$ of order $r_i$,
$i=1,\dots n$. There cannot be other relations. As we
assume $N$ is prime all $r_i$ must equal $N$.
Furthermore we note that $\vert L^*/L \vert = N^n$.\flexskip

Hence any $\lambda\in L^*/L$ is of the form $\sum_{i=1}^n a_ix_i$
with $0\leq a_i \leq N-1$. Consider
\[
\sum_{a_n=0}^{N-1} e^{\pi ik(y+a_nx_n)^2}
\]
where $y$ is a
fixed element of $L^*/L$ of the form $\sum_{i=1}^{n-1} a_ix_i$. Choose
integers $\alpha,\beta$ such that $x_n^2={2\alpha\over N}$, and
$(x_n,y)\equiv{2\beta\over N}(\ZZ)$. Note that $\alpha\neq 0$. 
Choose $\alpha'$ such that $\alpha\alpha'\equiv1(N)$.
We observe 
\[
(y+a_nx_n)^2 \equiv (y-\alpha'\beta x_n)^2 + {2\alpha'\over N}
(\alpha a_n+\beta)^2 \ (\hbox{mod \ }2\ZZ)
\]
Thus
\begin {equation}
  \sum_{a_n=0}^{N-1} e^{\pi ik(y+a_nx_n)^2} = 
  e^{\pi ik(y-\alpha'\beta x_n)^2} \sum_{a_n=0}^{N-1} 
  e^{\pi ik{2\alpha'\over N} (\alpha a_n+\beta)^2}.
\end{equation}
The second factor on the right hand side of (3.3) can be calculated explicitly.
As $\alpha'\neq 0$, we find for $N\equiv3(4)$
\begin{subequations}
\begin{equation}
= \sum_{r=0}^{N-1} e^{{2\pi i\over N} k\alpha' r^2} = \left({k\alpha'\over
N}\right) {\sqrt N}i  = \left({k\alpha\over N}\right) {\sqrt N}i
\end{equation} 
and for $N\equiv1(4)$
\begin{equation}
= \sum_{r=0}^{N-1} e^{{2\pi i\over N} k\alpha' r^2} = \left({k\alpha'\over
N}\right){\sqrt N}=\left({k\alpha\over N}\right) {\sqrt N}.
\end{equation} \end{subequations}
Now consider the first factor. Let $\phi(y)=y-(y,x_n)N\alpha'x_n$. The
set 
\[
\left\{\phi(y)\ \vert\ y=\sum_{i=1}^{n-1} a_i x_i, 0\leq a_i \leq N-1
\right\}
\]
forms an abelian group on $n-1$ generators.
It has exactly the $n-1$ relations inherited from $L^*/L$.
Hence we evaluate (3.2) using equations (3.4) and induction on the
number of generators as in equation (3.3). We obtain for $N\equiv3(4)$
\[
c_0(k) = {e^{-3l\pi i/4} \over \sqrt{\vert L^*/L \vert}} \prod_{r=1}^n \left(
i\left({k\alpha_r\over N}\right) \sqrt N \right) = e^{n\pi i/2-3l\pi i/4}
\left({k^n\over N}\right) \left({\prod_{r=1}^n\alpha_r\over N}\right)
\]
For $N\equiv1(4)$ we obtain
\[
c_0(k) = {e^{-3l\pi i/4} \over \sqrt{\vert L^*/L \vert}} \prod_{r=1}^n \left(
\left({k\alpha_r\over N}\right) \sqrt N \right) = e^{-3l\pi i/4}
\left({k^n\over N}\right) \left({\prod_{r=1}^n\alpha_r\over N}\right)
\]
Next we observe that $(V_1,*) = (S,+)(T,1)(S,+)(T,1)(S,+) = (T^{-1},1)$ acts
trivially, i.e. $c_0(1)=1$.
Hence, for general $k$, and any prime $N\neq2$,
 $c_0(k)=\left({k\over N}\right)^n$.\flexskip

Because the coefficient $c_0(k)$ has modulus 1 and because the
transformation induced by the action of $V_k$ is unitary (see the final part of
theorem 3.3), the remaining coefficients must be zero. The Legendre symbol is
multiplicative and thus we have proven

\begin{theorem} Let $L$ be even $l$-dimensional integral lattice, let
$L^*/L$ have $n$ generators, let $N\ne2$ be prime such that $NL^*\subseteq L$.
Then $\theta_L$ is a modular form for $\Gamma_0(N)$ of weight $l\over2$ with
character $\chi{ab\choose cd} = \left({d\over N}\right)^n$.
(The character of course depends on the branch and is as stated when the branch
of any $(G,*)$ is the product of the branches of the $(V_k,*)$ above.)
\end{theorem} \theoremproven

We observe two interesting consequences of the above argument.
The fact that the action of $(V_1,1)$ is trivial provides us with some relations
between $l,n$ and $N$. In the case $N\equiv 3(4)$ we conclude that $n\equiv
{l\over2}$ mod 2. Thus ${l\over2}$ can replace $n$ in theorem 3.4. In the
case $N\equiv 1(4)$ we conclude that ${l\over2}$ must be even. \flexskip

For $N=2$ we recall that $\Gamma_0(2)$ was generated by $T$ and $ST^2S$. The
following theorem can be proved along the same lines as theorem 3.4.

\begin{theorem} Let $L$, $l$, $N$, $n$ be as in theorem 3.4, except
that we now assume $N=2$. Then 
\[
\theta_L\vert_{\scriptstyle \bigl(
(S,+)(T^2,1)(S,+),{l\over2} \bigr)} (q) = e^{-l{\pi i\over2}} \theta_L(q). 
\]
\end{theorem} \theoremproven



\chapter{The Proof of Theorem 1.7} 

In this chapter we prove formula (1.25). This will then complete the proof of
the central result of chapter 1. In section 4.1, we describe the lattice
$L={\Lambda^\sigma}^\perp$ and its dual $L^*$. Note that in this chapter, as
in chapter 3, $L$ does not denote the same object as in chapter 1.
The explicit identification of the fixed point lattices can then be used to
complete the proof of lemma 1.3. We identify the modularity properties of
the theta-function of $L$ and some lattice symmetries. In section 4.2, we
proceed to count the number of elements of $L^*/L$ according to norm. This
amounts to counting sums of quadratic residues. Section 4.3 classifies the
vectors of $L$ of norm less than 2. Section 4.4 puts everything together to
conclude the proof of the central theorem 1.7. \flexskip

We start with the observation that for $N=2$ and $N=3$ the fixed point lattices
$\Lambda^\sigma$ are well known. For $N=2$, \cite{CS88} identify the fixed point
lattice in chapter 4, section 10, as the Barnes-Wall lattice $\Lambda_{16}$. Its
orthogonal complement is the lattice $\sqrt2 E_8$. For $N=3$, \cite{CS88} identify
the fixed point lattice in chapter 4, section 9, as the Coxeter-Todd lattice
$K_{12}$ whose orthogonal complement is $K_{12}$. Thus in these cases the
theta-functions are well known.

\section{The Theta-Function of L*}

\cite{CS88}, chapter 10, theorem 25 provides the following description
of the Leech lattice $\Lambda$ in $\RR^{24}$. Let ${\cal C}$ be the set of the
elements of the 24-dimensional Golay code. The vector
$(x_\infty,x_0,...,x_{22}) \in \ZZ^{24}$ is in $\sqrt8\Lambda$ if and only if
\begin{subequations}
\begin{align}
& \text{the co-ordinates $x_i$ are all congruent modulo 2, to $m$, say;} \\
& \text{the set of $i$ for which $x_i$ takes any given value modulo 4    
is a ${\cal C}$-set;} \\
& \text{the co-ordinate-sum is congruent to $4m$ modulo 8.}
\end{align}
\end{subequations}

Equivalently we recall from \cite{CS88}, chapter 4 that any element of
$\sqrt8\Lambda$ is of one of the following two shapes:
\begin{equation}
2c+4x \mskip 30mu \hbox{or} \mskip 30mu (1^{(24)})+2c+4y
\end{equation}
Here $c \in {\cal C}$ is understood as an element of $\RR^{24}$,
$(x_n)\in \ZZ^{24}$ is
such that $\sum x_n$ even, $y=(y_n)\in \ZZ^{24}$ is such that $\sum y_n$ odd.
In line with equation (1.3), the norm of an actual Leech lattice element in
the above descriptions is $1\over8$ times the sum of the squares of the
(integer) co-ordinates. \flexskip

We consider simultaneously the cases $N=2,3,5,7,11,23$.
Let $\sigma$ be an automorphism of the Leech lattice of order $N$ and cycle
shape $1^M N^M$ where $M={24 \over N+1}$. Thus $\sigma$ acts by permuting the
co-ordinates of $\RR^{24}$. The character of the 24-dimensional representation
(where $\sigma$ acts as permutation of the co-ordinates) is registered in the
Atlas \cite{Con85} as $\chi_{102}$. From the cycle shape of $\sigma$ we conclude
that $\chi_{102}(\sigma) = M$. Hence we can identify the conjugacy class
of $\sigma$ in the Atlas \cite{Con85} as 2A, 3B, 5B, 7B, 11A, and 23A respectively.
\flexskip

In the case $N=5$ we observe that the five elements of a cycle of $\sigma$
exactly determine an octad. The action of $\sigma$ maps this into another octad
as $\sigma$ preserves the Leech lattice as a whole. This new octad also
contains the 5-cycle, hence the octad is preserved by $\sigma$. Thus the 3
remaining points of the octad must be among the 4 fixed points of $\sigma$.
A direct check of the examples
given in \cite{CS88}, chapter 10, section 2.1, yields that every 7-cycle
together with one of the fixed points forms an octad of ${\cal C}$ and that
every 11-cycle together with one of the fixed points forms a dodecad of ${\cal
C}$.\flexskip

To determine the elements of $\Lambda^\sigma$ we observe that
$\lambda\in\sqrt8\Lambda$ is invariant under $\sigma$ if and only if
$\lambda=(a_1,...,a_M,$ $b_1^{(N)},...,b_M^{(N)})$. ($b^{(N)}$ denotes $N$
entries of $b$ in the positions of a cycle of $\sigma$. We may order the
co-ordinates according to the cycles of $\sigma$.) Thus the fixed point lattice
$\Lambda^\sigma$ is a $2M$-dimensional lattice. We define
$L = (\Lambda^\sigma)^\perp$, the lattice orthogonal to the fixed point
lattice of dimension $24-2M$. By theorem 3.2,  $L^*/L$ and
$(\Lambda^\sigma)^*/\Lambda^\sigma$ are isomorphic. By theorem 3.1,
$(\Lambda^\sigma)^* = \pi \Lambda$, where $\pi=\pi_{\Lambda^\sigma}$ is the
projection of the Leech lattice into the subspace spanned by the fixed point
lattice.

\begin{lemma} If an element $\lambda=(a_1,...,b_M^{(N)})$ of
$\sqrt8(\Lambda^\sigma)^*$ has integer co-ordi-nates, it is in
$\sqrt8\Lambda^\sigma$. 
\end{lemma}
\begin{proof} For $N=2,3$ this follows from the explicit description of the
lattices. For the remaining $N$ we observe that $\pi$ acts by averaging over
the cycles. Hence it follows from the description in (4.2) that
the co-ordinates of $\lambda$ are either all even or all odd.
We subtract a suitable multiple of
$(-3^{(1)}, 1^{(23)}) \in \sqrt8\Lambda^\sigma$ to obtain even co-ordinates
everywhere. From the description of the cycles of $\sigma$ with respect to the
Golay code we conclude that for every cycle there exists an octad (or dodecad
respectively) with entries precisely in one cycle of $\sigma$
and in a number of co-ordinates fixed by $\sigma$.
Subtracting suitable multiples of these we obtain from
$\lambda$ an element $\lambda' \in \sqrt8(\Lambda^\sigma)^*$ with non-zero
entries only in the $M$ co-ordinates fixed by $\sigma$. This, in turn, implies
that $\lambda'$ must be an element of $\sqrt8\Lambda$ by the following argument.
The product of $\lambda'$ with any element of $\sqrt8 L$ must be in $8\ZZ$.
$\lambda'$ has non-zero entries only in up to four places. Now there are
elements of $\sqrt8L$ having entries $(2,0,0,0,*,...,*)$. Thus all the entries
of $\lambda'$ must be divisible by 4. Furthermore, there are elements of shape
$(1,1,1,1,*,...,*)$ of $\sqrt8L$. Hence the sum of the entries of $\lambda'$
must be divisible by 8. Thus $\lambda'$ satisfies the conditions (4.1)
and hence is an element of $\sqrt8\Lambda$. 
\end{proof}

The projection $\pi$ acts by averaging the co-ordinates of the cycles of the
permutation. Thus the entries must be in ${1\over N} \ZZ$.
Hence we have shown that $N(\Lambda^\sigma)^* \subset \Lambda^\sigma$ and
more precisely $(\Lambda^\sigma)^*/(\Lambda^\sigma)$ is an
$M$-dimensional vector space over $\ZZ/N\ZZ$ and has $N^M$ elements.
Moreover, we can give an explicit basis for $(\Lambda^\sigma)^*/\Lambda^\sigma$
as follows:
\[
{\lambda_j^\sigma}^*={1\over\sqrt8}(4,0,\dots,0,{4\over N}^{(N)},0,\dots)
\]
where $j=1,\dots,M$, one for each $N$-cycle of $\sigma$.
Note that the ${\lambda_j^\sigma}^*$ are mutually orthogonal modulo $2\ZZ$ and
of norm ${2\over N}(2\ZZ)$. \par

By theorem 3.2, the same holds for $L^*/L$.
For the transition we observe that any $\lambda\in\Lambda$ can be written as
\[
\lambda=\pi_L(\lambda)+\pi_{\Lambda^\sigma}(\lambda).
\]
In particular, this implies that $\pi_L(\lambda)\equiv\pi_L(\mu)
(\text{mod }L)$ if and only if
$\pi_{\Lambda^\sigma}(\lambda)\equiv\pi_{\Lambda^\sigma}(\mu)\
(\text{mod }\Lambda^\sigma)$ and further
\[
 \left(\pi_L(\lambda),\pi_L(\mu)\right)\equiv
-\left(\pi_{\Lambda^\sigma}(\lambda),\pi_{\Lambda^\sigma}(\mu)\right) \
\hbox{mod}\ \ZZ
\]
$L$ is even because $\Lambda$ is and $L$ has dimension $24-2M$.
We can now apply theorems 3.4 and 3.5 to obtain

\begin{theorem} For all $N=2,3,5,7,11,23$, the theta-function of $L$,
$\theta_L$, is a modular form of $\Gamma_0(N)$ of weight ${24-2M\over2}=
12{N-1\over N+1}$. The character is $\chi{ab\choose cd} =
\left({d\over N}\right)$ for $N=7$ and $N=23$. The character is trivial in the
other cases.
\end{theorem} \theoremproven

Using the explicit identification of the dual lattice elements, we are now
able to complete the

\begin{proof} [Proof of Lemma 1.3] 
The cases $N=2$ and $N=3$ are clear as the fixed point lattices are
$\Lambda_{16}$ and $K_{12}$ respectively. So we assume $N\ge5$ and $M\le4$.
Now suppose there were a root $r\in \Lambda^\sigma$. The reflection induced
by $r$ is the same as that of $nr\in \Lambda^\sigma$ where $n$ is any
nonzero integer. Hence we can assume that $r$ is primitive in the sense that,
if $n\in\ZZ, |n|>1$, then $r\over n$ is not in $\Lambda^\sigma$. For any
$v\in\Lambda^\sigma$, the inner product $(r,v)$ is integer, hence the
greatest common divisor $d=(r,\Lambda^\sigma)$ of all the inner products
$(r,v), v\in\Lambda^\sigma$, is defined. It follows that ${r\over d}\in
{\Lambda^\sigma}^*$. We have seen above that $N{\Lambda^\sigma}^* \subset
\Lambda^\sigma$. As $N$ is prime it follows that the only possible cases are
$d=1$ or $d=N$. 

The reflection through $r$ takes a vector $v$ to the vector
$v-{2(r,v)\over(r,r)}r$. If $d=1$ we conclude that $2\over(r,r)$ must be
integer, which is impossible as the Leech lattice does not contain elements
of norm 1 or norm 2. If $d=N$ we conclude that $2N\over(r,r)$ must be
integer. This is equivalent to $2\over N$ being an integer multiple of
$({r\over N}, {r\over N})$. Now ${r\over N}\in {\Lambda^\sigma}^*=
\pi_{\Lambda^\sigma}(\Lambda)$. As identified above, elements of
${\Lambda^\sigma}^*$ are of shape $${1\over\sqrt8}
\bigl(a_1,\dots,a_M,{b_1\over N}^{(N)},\dots, {b_M\over N}^{(N)}\bigr).$$
Here $a_i, b_i\in\ZZ$, and the raised $^{(N)}$ indicates $N$ equal entries.
The norm of the above element is $${1\over8}
\bigl(a_1^2+\dots+a_M^2+{b_1^2\over N}+\dots+{b_M^2\over N}\bigr).$$ It
remains to prove that this cannot be smaller or equal $2\over N$ in the
relevant cases.

The projection $\pi$ from $\Lambda$ to $\Lambda^\sigma$ acts by averaging the
entries of the $N$-cycles. Equation (4.2) stated that the vectors in
$\Lambda$ are of two different shapes. For vectors of shape
$(1^{(24)}+2c+4y)$ we note that the vector $(1^{(24)})$ is preserved under
$\pi$. Hence
\[
\|\pi(1^{(24)}+2c+4y)\|\ge{1\over8}(M\times1+M\times{1\over N})=
{1\over8}{M(N+1)\over N} = {3\over N}.
\]
For vectors of the remaining shape, $(2c+4x)$ where $c$ is an element of
the Golay-code and the sum of the entries of $x$ is even we observe that the
norm will be less or equal ${2\over N}$ only if all $a_i$ equal 0 and
at most two of the $b_i$ have modulus at most 1. However, if the sum of
coordinates of a cycle is nonzero, it is at least 2, because all entries are
even. Thus there cannot be
any vectors of norm less or equal $2\over N$.
\end{proof}

The group of all automorphisms of the Leech lattice is identified as $2.Co_1$
in the Atlas \cite{Con85}. If $\sigma$ is an automorphism as above we consider
the group $\langle\sigma\rangle$ of order $N$, generated by $\sigma$. The
normalizer of $\sigma$ in $2.Co_1$ is the subgroup of $2.Co_1$ of all
elements $\psi$ such that $\psi \langle\sigma\rangle = \langle\sigma\rangle
\psi$. If $\psi$ is an element of the normalizer of $\sigma$ it maps the
fixed point lattice of $\sigma$ to itself and hence induces an automorphism
of $L$ and moreover $L^*$. For the same reason it induces an automorphism on
$L^*/L$. For $N$ prime, we write $\ZZ_N$ for the finite field $\ZZ/N\ZZ$. We
follow the notation of the Atlas \cite{Con85} for orthogonal groups over finite
fields. For even dimension $2m$, the `$+$' type of an orthogonal group
corresponds to maximum Witt index $m$, and the `$-$' type corresponds to
Witt index $m-1$. 
\begin{claim} For $N=11$, $7$, $5$, and $M=24/(N+1)$, consider $\sigma$
of order $N$ as above. For $N=11$, the centralizer of $\sigma$ in $2.Co_1$
acts on $L^*/L$ as the orthogonal group $O_2^-(11)$ (up to a conjugacy class
of order $2$). For $N=7$ and $N=5$, the normalizers of $\sigma$ in $2.Co_1$
act on $L^*/L$ as the orthogonal groups $SO_3(7)$ and $GO^+_4(5)$,
respectively.
\end{claim}
\begin{proof} All calculations required to prove this claim were carried out
during the composition of the Atlases \cite{Con85}, \cite{Jan95}, but not all have
been documented. Theorem 3.2 applies such that we obtain $L^*/L =$
$({{\Lambda^\sigma}^\perp}^*)/ ({{\Lambda^\sigma}^\perp}) \tilde =$
${\Lambda^\sigma}^*       /  {\Lambda^\sigma}$. Hence, we may
consider the fixed point lattice, rather than its orthogonal complement.
This reduces the dimension of the lattices under consideration
from $24-2M$ to $2M$, over $\CC$.
Furthermore, explicit descriptions of the lattice elements are far more
readily available, and a basis has been identified on p. 41 above. \flexskip

For the case $N=11$, the element 11A has centralizer $11\times D_{12}$ in
$2.Co_1$. In accordance with the notations of the Atlas \cite{Con85}, $D_{12}$
here denotes the dihedral group of order $12$. We note the isomorphism
$D_{12} \tilde= $ $O_2^-(11)$. The action of the centralizer on the quotient
space ${\Lambda^\sigma}^*/ \Lambda^\sigma$ induces a representation of the
orthogonal group over $\ZZ_{11}$ of degree 2. 

Using the Suzuki construction of the Leech lattice we find that the order 3
automorphism within the centralizer corresponds to a rotation through 120
degrees in that construction. Hence its trace must be $-12$.
The character of the 24-dimensional representation of $2.Co_1$ is recorded
in the Atlas \cite{Con85} as $\chi_{102}$. This identifies the
automorphism class as 3A.

When restricted to the action of $D_{12}$, the 24-dimensional representation
of $2.Co_1$ decomposes into a number of irreducible constituents. These are
sufficiently characterised by the trace of the automorphism of order 3.
Counting multiplicities, there are 12 irreducible constituents whose
characters all satisfy
\[
\chi(1A) = 2, \ \ \ \chi(2A)=0,\ \ \ \chi(3A)=-1,\ ...
\]
This, in turn, enables us to identify the explicit matrix form
of the order 3 automorphism on the 4-dimensional vector space
$\Lambda^\sigma$ over $\CC$, which has trace $-2$.
It is given in formula (6.14) below.

We now consider representations over the finite field $\ZZ_{11}$.
These are in one-to-one correspondence with the representations over $\CC$
because the characteristic 11 of $\ZZ_{11}$ does not divide the group order
12 of $D_{12}$. Using the explicit basis for ${\Lambda^\sigma}^*/
\Lambda^\sigma$, we find that the automorphism of class 3A continues to act
non-trivially on the 2-dimensional space ${\Lambda^\sigma}^*/\Lambda^\sigma$.
Its matrix form over the finite field $\ZZ_{11}$ is
\[
\begin{pmatrix} 5& -3 \\ 3& 5 \end{pmatrix} \ \ {\rm mod}\ 11.
\]

The character table of $D_{12}$ now sufficiently characterises the induced
action of $D_{12}$ for our purposes. Note that we have not determined the
action of one conjugacy class of order $2$ which may, according to character
table, act as either $+{\rm id}$ or $-{\rm id}$. Up to this class however, we
may nevertheless conclude that the induced action of the centralizer is the
orthogonal matrix action of the orthogonal group $D_{12} \tilde= O_2^-(11)$.
\flexskip

For the case $N=7$, the element 7B in $Co_1$ has normalizer
$(7.3 \times L_2(7)).2$, which contains $L_2(7).2 \tilde=$ $SO_3(7)$, and
centralizer $7\times L_2(7)$, which contains $L_2(7)\tilde=$ $O_3(7)$.
Note that, in the case $N=7$, we may consider the normalizer and centralizer
in $Co_1$ because the double cover of $L_2(7)$ in $2.Co_1$ is not proper.
The action of the normalizer on the quotient space ${\Lambda^\sigma}^*/
\Lambda^\sigma$ induces a representation of the orthogonal group
of degree 3 over $\ZZ_7$.

The centralizer contains an element $\psi$ of order 3 which cyclically
permutes the 3 octads underlying the automorphism $\sigma$. Hence,
tr$(\psi) = 0$. Using the explicit basis of
${\Lambda^\sigma}^*/ \Lambda^\sigma$, we find that $\psi$ induces an
automorphism of the 3-dimensional quotient space (over $\ZZ_7$),
again with trace 0. 

The Brauer character table of the (unextended) group $L_3(2)\tilde=L_2(7)$
in the Atlas \cite{Jan95} shows that this trace will only be matched by the
character of the irreducible representation $\phi_2$, which, at the same
time, is the character of the orthogonal representation of $SO_3(7)$. Hence,
we have uniquely identified the induced action of the normalizer as the
orthogonal matrix action of the orthogonal group $SO_3(7)$. 
\flexskip

We turn to the case $N=5$. The normalizer of an element of class 5B in $Co_1$
is identified in the Atlas \cite{Con85} as $(D_{10} \times (A_5 \times
A_5).2).2$. In $2.Co_1$, the normalizer of an element of class 5B contains
the orthogonal group $GO_4^+(5) \tilde=$ $2.O_4^+(5).2^2$ where $O_4^+(5)
\tilde=$ $A_5 \times A_5$. The centralizer of class 5B in $2.Co_1$ is
$5 \times 2.(A_5 \times A_5).2$. Unlike $N=7$, in the case $N=5$ we must
consider the normalizer and centralizer in the double cover $2.Co_1$.
The representations involved do not restrict to representations of subgroups
of $Co_1$.
The action of the normalizer on the quotient space ${\Lambda^\sigma}^*/
\Lambda^\sigma$ induces a representation of the orthogonal group $GO_4^+(5)$
of degree 4 over $\ZZ_5$.

The leading `$5\times$' corresponds to the element
5B itself. In terms of representations, the trailing `$.2$' relates to split
or fused representations. For our purposes, it is therefore sufficient
to consider the irreducible representations of the
group $2.(A_5 \times A_5)$. Over both $\CC$ and the finite field $\ZZ_5$,
these can be constructed from the representations of the group
$(2.A_5) \times (2.A_5)$, by quotienting out the kernel of the
representation (which corresponds to quotienting out the central 2,
identifying the two leading `$2.$'s). The characters of all representations
of $2.A_5$ are listed in the Atlases \cite{Con85} and \cite{Jan95}.

Again, it is easy to identify explicitly a number of elements of the
centralizer. There are permutations of the 4 octads underlying the
automorphism $\sigma$ of order 2, 3, 4. There is also the order 3
automorphism which corresponds to a rotation through 120 degrees in the
Suzuki construction. The explicit form of these elements of the centralizer
is sufficient to conclude that neither of the two subgroups $A_5$ acts
trivially on any 1-dimensional subspace of the 8-dimensional space spanned
by the fixed point lattice. Equally, this holds for the 4-dimensional
quotient space $L^*/L$ (over $\ZZ_5$). 

This sufficiently characterises the representation in the Atlas character
tables of $2.A_5$. Over $\CC$, we obtain the 8-dimensional real
representation of $2.A_5 \times 2.A_5$ by tensoring the 2-dimensional
(conjugate) representations $\chi_6$ or $\chi_7$ and fusing with their
algebraic conjugate: $8 = 2 \times 2 + 2 \times 2$. Over $\ZZ_5$, we obtain
the 4-dimensional representation of $2.A_5 \times 2.A_5$ uniquely by
tensoring the 2-dimensional representation $\phi_4$ such that $4=2\times 2$.
At the same time, this is the character of the orthogonal representation of
$2.O^+_4(5) \tilde= 2.(A_5 \times A_5)$. (Note that the `generically
simple' group  $O^+_4(5) \tilde= A_5 \times A_5$ has a faithful
{\it projective} 4-dimensional representation only.) Hence, we have
uniquely identified the induced action of the normalizer as the orthogonal
matrix action of the orthogonal group $GO^+_4(5)$. 
\end{proof}

\begin{theorem} Let $N$ be any of $23$, $11$, $7$, $5$, $3$, or $2$.
Let $\sigma$ of order $N$ and $L={\Lambda^\sigma}^\perp$ be as above.
Let $\lambda^*,\mu^*$ be elements of $L^*$, not in $L$,
such that ${(\lambda^*)^2\over2} \equiv {(\mu^*)^2\over2}$ modulo $\ZZ$.
Then $\theta_{L+\lambda^*} = \theta_{L+\mu^*}$.
\end{theorem}
\begin{proof}
Let $[ {^.} ]$ denote the residue class in $L^*/L$.
It suffices to construct an automorphism of $L^*$ so that its
restriction to $L^*/L$ takes the residue class $[\lambda^*]$ of
$\lambda^*$ to the residue class $[\mu^*]$ of $\mu^*$.
Now $L^*/L$ is a vector space over the finite field $\ZZ_N$ of dimension
$M$. We recall the inclusion $N(\Lambda^\sigma)^* \subset \Lambda^\sigma$.
Hence, for $N\ne2$ we obtain a non-degenerate inner product
with values in $\ZZ_N$ if we define
\[
([\lambda^*],[\mu^*])= N (\lambda^*,\mu^*)\ \hbox{mod}\ N
\]
where the bilinear form on the right hand side is the standard form.
Calculating the Witt index shows that this form is of `$-$' type for
$N=11$, and of `$+$' type for $N=5$. (For $N=23$, the form is trivial,
for $N=7$ it is unique. The cases $N=2$ and $N=3$ are handled separately.) 
We will now consider the various $N$
individually. The claim is trivial for $N=23$ as $M=1$. \flexskip

Case $N=11$:
The 121 elements of the 2-dimensional vector space over $\ZZ_{11}$ fall
into 10 classes (of norms 1 to 10, respectively) with 12 elements each,
plus the zero vector. Co-ordinate interchange and application of 
$-{\rm id}$ are obvious symmetries. It is straightforward to check that
the automorphism ${5 \ -3 \choose 3 \ \ 5}$ mod 11 of order 3
explicitly identified above completes the required symmetries. \flexskip

Cases $N=7$ and $N=5$: We found above that the group $SO_3(7)$, and
hence $GO_3(7)$, describes symmetries of the quotient lattice $L^*/L =$
$({{\Lambda^\sigma}^\perp}^*)/ ({{\Lambda^\sigma}^\perp})$ for $N=7$.
Equally we found that $GO^+_4(5)$ describes symmetries of the quotient
lattice for $N=5$.
Witt's Extension Theorem (\cite{Jac85}, Chapter 6.5) asserts that, given any
non-singular bilinear form (i.e. `$+$' and `$-$' type for even dimensions,
unique for odd dimensions), any isometry of subspaces (i.e. the map taking
$[\lambda^*]$ to $[\mu^*]$ ) can be extended to an isometry of the whole
vector space, i.e. an element of the group $GO^\epsilon_M(N)$. \flexskip

Case $N=3$: \cite{CS83} prove that the group of isometries of the
Coxeter-Todd lattice acts transitively on vectors of norms 4, 6, and 8
respectively, which is sufficient to prove the claim as
$L^*={1\over\sqrt3}K_{12}$.\flexskip

Case $N=2$: In \cite{CS88}, chapter 4, section 10 (p.131) the centralizer
of $\sigma$ is identified as $O_8^+(2)$. $O_8^+(2)$ acts transitively on
vectors of norm 1 and 2 in $L^*={1\over\sqrt2}E_8$. 
\end{proof}

We summarize: The translated theta-function $\theta_{L+\lambda^*}$, $\lambda^*
\not\in L$, does only depend on ${(\lambda^*)^2\over2}\ \hbox{mod}\ \ZZ$.
Hence we can pick $\lambda_r^*$ such that ${(\lambda_r^*)^2\over2} \equiv
{-r\over N}\ \hbox{mod}\ \ZZ$ for $r=0,\dots,N-1$ and write
\begin{equation}
\theta_{L^*} = \theta_L +\sum_{r=0}^{N-1} \rho_M(r,N) \theta_{L+\lambda_r^*}.
\end{equation}
Here
$\rho_M(r,N)$ is the number of non-zero elements with half-norm identical to
${-r\over N}$ in $L^*/L$. It will be determined in the next section. 

\section{Sums of Quadratic Residues}

Let $N$ and $M={24\over N+1}$ be as above. For $\lambda\in \Lambda$ we
calculate the norm of $\pi_L \lambda \in L^*$ as
\begin{equation}
(\pi_L \lambda)^2= \Bigl((1-\pi_{\Lambda^\sigma})
\lambda \Bigr)^2= \lambda^2 - \bigl(\pi_{\Lambda^\sigma} \lambda \bigr)^2.
\end{equation}

Hence the number $\rho_M(r,N)$ of residue classes
of a certain half-norm ${-r\over N}\ \hbox{mod}\ \ZZ$ in $L^*/L$ is equal to the
number of residue classes of half-norm ${r\over N}\ \hbox{mod}\ \ZZ$ in
${\Lambda^\sigma}^*/\Lambda^\sigma$.
We recall that we found an explicit basis $\{{\lambda_j^\sigma}^*\}$ of
${\Lambda^\sigma}^*/\Lambda^\sigma$ such that
${({\lambda_j^\sigma}^*)^2\over2}\equiv{1\over N}\ \hbox{mod}\ \ZZ$ and
${({\lambda_j^\sigma}^*,{\lambda_k^\sigma}^*)\over2}\equiv0\ \hbox{mod}\ \ZZ$.
Thus $\rho_M(r,N)$ is equal to the number
of non-zero solutions to the following equation:
\begin{equation}
\sum_{i=1}^M {x_i}^2 \equiv r \;\; \hbox{mod} \;\; N 
\end{equation}
for $x_i \in \ZZ_N$.
We now define $\tilde\rho_M(r,N)$ to be the number of all solutions to the
above equation. Hence for $r\not\equiv 0(N)$, $\tilde\rho_M(r,N)=\rho_M(r,N)$,
and $\tilde\rho_M(0,N)=\rho_M(0,N)+1$. We observe that $\tilde\rho_M(m^2r,N) =
\tilde\rho_M(r,N)$ for any non-zero $m$. Note that we always consider ordered
pairs, or $n$-tuples, of residues.

\begin{lemma} Let $N\not= 2$ be a prime. If $N\equiv1(4)$ (that is, $-1$ is
a square modulo $N$) then $\tilde\rho_2(0,N) = 2N-1$ and $\tilde\rho_2(r,N) =
N-1$ for all $r\not\equiv 0(N)$. If $N\equiv3(4)$ (that is, $-1$ is a non-square
modulo $N$) then $\tilde\rho_2(0,N) = 1$ and $\tilde\rho_2(r,N)= N+1$ for all
$r\not\equiv 0(N)$. 
\end{lemma}
\begin{proof} If $-1$ is a square, say $z^2\equiv -1$ then $x^2+y^2 \equiv 0$ has
the following $2N-1$ solutions $(x, zx), (x,-zx)$ for non-zero $x$ and $(0,0)$.
(note that $z \not \equiv -z$ as $N\not= 2$.) $x^2+y^2 \equiv r$ is equivalent
to $x^2-(zy)^2 = (x+zy)(x-zy) \equiv r$. This clearly has the same number of
solutions for any non-zero $r$. The result follows as there are $N^2-2N+1$
possible choices of $(x,y)$ and $N-1$ different non-zero values for $r$. \par
If $-1$ is a non-square modulo $N$, $\tilde\rho_2(0,N)=1$ obviously. Let us
first consider $\tilde\rho_2(r,N)$ for a square $r=z^2$.
The equation $x^2+y^2\equiv z^2 \;\;\hbox{mod}\;\; N$
is equivalent to $x^2 \equiv (z+y)(z-y) \;\;\hbox{mod} \;\;N$. We
substitute $a=z+y$ and $b=z-y$. Now if we choose $a=0$ then $b$ is arbitrary
and $x=0$. (Note that we consider free $z$.) This gives $N$ solutions.
If we choose $a\not= 0$ ($N-1$ choices) then we can furthermore choose $x$
arbitrarily which fixes $b$. This gives $(N-1)N$ solutions. We exclude the
solution $(x,y,z)=(0,0,0)$ and are left with $N^2-1$ solutions such that
$z\not= 0$. Of course, the number of solutions for any of the $N-1$ non-zero $z$
must be the same, hence there are $N+1$ such. This proves $\tilde\rho_2(r,N) =
N+1$ for any square $r$. $\tilde\rho_2(r,N)=N+1$ for $r$ a non-square then
follows from the fact that there are ${1\over2}(N-1)(N+1)$ pairs $(x,y)$ left
and that there are ${1\over2}(N-1)$ non-squares.
\end{proof}

\begin{corollary} For $M$ even, $\tilde\rho_M(r,N)$ is equal for any
non-zero $r$.
\end{corollary}
\begin{proof} We obtain the result by induction using $$\tilde\rho_M(r,N)=
\sum_{s=0}^{N-1} \tilde\rho_2(s,N) \; \tilde\rho_{M-2}(r-s,N).$$
\end{proof}
\begin{theorem} For even $M$, that is $N=2,3,5,11,$ we have
\[
\tilde\rho_M(r\neq0,N) =
  {1\over N}\left(N^M-(-N)^{M\over2}\right)  \mskip 10mu \hbox{\it and}
\]
\[
\tilde\rho_M(r=0,N) =
  N^{M-1}+(-1)^{M\over2}\left(N^{M\over2}-N^{{M\over2}-1}\right).
\]
For odd $M$, that is $N=7,23,$ we have
\[
\tilde\rho_M(r\neq0,N) =
{1\over N}\left(N^M+\left(r\over N\right)N^{M+1\over2}\right)
\mskip 10mu \hbox{\it and} \mskip 10mu \tilde\rho_M(r=0,N) = N^{M-1}.
\]
Here, $\left({r\over N}\right)$ stands for the Legendre symbol as before.
\end{theorem}
\begin{proof}
We use lemma 4.2 to proceed inductively to the cases $\tilde\rho_M(r,N)$ of
interest. The case $(N,M)=(23,1)$ is trivial as any square can be represented
in 2 ways and any non-square cannot be represented.
The case $(N,M)=(11,2)$ is covered by lemma 4.2.

Case $(N,M)=(7,3)$: 
\begin{align*}
\tilde\rho_3(0,7) 
 &= \tilde\rho_1(0,7)\;\tilde\rho_2(0,7)
    + \sum_{s=1}^{N-1\over2} \tilde\rho_1(s^2,7)\;\tilde\rho_2(-s^2,7) \\
 &= 1\times1 +{(N-1)\over2}\times2\times(N+1) = N^2 \\
\tilde\rho_3(-1,7) 
 &= \tilde\rho_1(0,7)\;\tilde\rho_2(-1,7) +
    \sum_{s=1}^{N-1\over2} \tilde\rho_1(s^2,7)\;\tilde\rho_2(-1-s^2,7) \\
 &= 1\times(N+1) + {(N-1)\over2}\times2\times(N+1) = N(N+1) \\
\intertext{Note that
$-1$ is not a square, hence $-1-s^2$ is never zero. The value for $r$ a square
now follows as the difference to the total $N^3$. }
\tilde\rho_3(1,7) 
 &= {{N^3 - \tilde\rho_3(0,7) - {(N-1)\over2} \tilde\rho_3 
    (-1,7)}\over {N-1\over2}} = N(N-1) \\
\intertext{Case $(N,M)=(5,4)$:}
\tilde\rho_4(0,5) 
 &= \tilde\rho_2(0,5)\;\tilde\rho_2(0,5) +
\sum_{s=1}^{N-1} \tilde\rho_2(s,5)\;\tilde\rho_2(-s,5) \\
 &= (2N-1)\times(2N-1) + (N-1)\times(N-1)\times(N-1) \\
 &= N^3+N^2-N \\
\intertext{As $M=4$ is even, for non-zero $r$ we obtain}
\tilde\rho_4(r,N) 
 &= {{N^4-(N^3+N^2-N)}\over{N-1}} = N^3-N \\
\end{align*}

Case $(N,M)=(3,6)$: The squares modulo 3 are 0 and 1. The sum of six squares
will be zero modulo 3 if and only if 0, 3, or 6 summands equal 1.
Hence there are
${6 \choose 0} 2^0 + {6 \choose 3} 2^3 + {6 \choose 6}2^6= 225$ solutions.
As $M=6$ is even, for non-zero $r$ we obtain 
\[
\tilde\rho_6(r,3) = {{3^6-225} \over2}=252.
\]
Case $(N,M)=(2,8)$: The squares modulo 2 are 0 and 1. The sum of eight squares
will be zero modulo 2 if and only if an even number of summands equal 1.
of 8 squares by either 0,2,4,6,8 summands 1. Hence
$ \tilde\rho_8(0,2) = {8 \choose 0} + {8 \choose 2} + {8 \choose 4} +
{8 \choose 6} + {8 \choose 8} = 136$ and $\tilde\rho_8(1,2) = 2^8 - 135 = 120$
which completes the proof.
\end{proof}

\section{The Short Vectors of L*}

In order to identify $\theta_L$ we will determine the leading
coefficients of both $\theta_L$ and $\theta_{L^*}$. As $L$ is a lattice of
dimension $24{N-1\over N+1}$, by theorems 3.4 and 3.5 its theta-function has
weight $12{N-1\over N+1}$. We will want to apply theorem 2.1.
Hence we need to determine the terms of a $q$-expansion
of exponent up to and including exponent ${N-1\over N+1}$. $L$ is a sublattice of
the Leech lattice. Hence any non-zero vector has norm greater or equal to 4.
Thus the expansion of $\theta_L - 1$ has a leading term of exponent at least 2.
We now turn to $L^*$. For the purpose of the proof it is sufficient to
determine those vectors of norm less or equal to $2{N-1\over N+1}$ in $L^*$.
In fact, we will find that there are no such non-zero vectors.
Note that the techniques developed below can also be applied to determine  
the shortest non-zero vectors in the lattices $L^*$, which turn out to be of 
norm $2{N-1\over N}$. The calculations were carried out as a double check of 
the results, but will not be presented here.

\flexskip 

The projection $\pi_L$ onto the orthogonal complement of
$\Lambda^\sigma$ is identical to the composition of the following two
operations on vectors of $\RR^{24}$.\par
$\pi_1$: Kill the components which correspond to the $M$ fixed points of
$\sigma$. This is a map $\RR^{24} \to \RR^{24-M}$.
\par

$\pi_2$: For each of the $M$ $N$-cycles of $\sigma$, subtract a suitable
(possibly non-integer) multiple of $(0^{(N(M-1))},1^{(N)})$ (non-zero entries
in the positions of an $N$-cycle) such that the sum of the co-ordinates of that
cycle becomes 0. This maps $\RR^{24-M}$ into a $(24-2M)$-dimensional subspace.
\flexskip 

Recalling the description of the Leech lattice in chapter 4.1, formula
(4.2), we note that we can describe the elements of $\sqrt8L^*$ as all elements
$\pi_L(2c+4x)$ where $c$ is still in the Golay code, but $x$ now is arbitrary in
$\ZZ^{24}$. In particular, we only need to consider even entries.
The following argument will only depend on the allocation of the
co-ordinates to the individual cycles and will be independent of the particular
order of the co-ordinates. Hence we adopt the convention that vectors will be
represented cycle by cycle without specification of the order of the entries
within a cycle or of the order of the cycles.
A general element $\lambda$ of $\sqrt8\pi_1 \Lambda$ will thus be represented as
\begin{equation}
\lambda = \left((a_{1,1},...,a_{1,N}),...,
       (a_{M,1},...,a_{M,N}) \right) \in \ZZ^{24-M}.
\end{equation}
We introduce the average of a cycle $\bar a_m = {1\over N}\sum_{n=1}^N a_{mn}$
and proceed to calculate the norm of ${1\over\sqrt8}\pi_2(\lambda) \in L^*$.
\begin{equation}
\pi_2 \lambda = \left((a_{1,1} - \bar a_1,...),...\right) \in \RR^{24-M}.
\end{equation}

\begin{equation}
 \| \pi_2 {1\over \sqrt8}\lambda \| = (\pi_2 {1\over \sqrt8}\lambda)^2
     = {1\over8} \sum_m \sum_n (a_{mn}-\bar a_m)^2
     = {1\over8} \sum_{m=1}^M \left( \sum_{n=1}^N
		 a_{mn}^2 - N \bar a_m^2 \right).
\end{equation}

Note that, in the notation of formula (4.2), ${1\over\sqrt8}\lambda$ is the
actual element of $L^*$ and that the norm of an element is understood to be the
inner product with itself rather than the square root of this inner product.
Formula (4.8) proves the following

\begin{lemma} If $\lambda$ is as in formula (4.6) above then
\[
 \| \pi_2 \lambda \| \ge \| \pi_2 \left((|a_{1,1}|,...,|a_{1,N}|),...,
       (|a_{M,1}|,...,|a_{M,N}|) \right) \|. 
\]
The inequality will be strict if within a cycle there are strictly negative and
strictly positive entries.
\end{lemma} \theoremproven

Furthermore, to achieve minimal norm the entries must be spread out as evenly
as possible in the sense of 

\begin{lemma} Let $\lambda$ be as in formula (4.6) above. Suppose that
within a cycle co-ordinate entries differ by more than 2. Suppose further that
the entries of the cycle are not all equivalent modulo 4. Then there exists a
$\mu \in \sqrt8 \pi_1 \Lambda$ such that $\|\pi_2 \lambda \| > \| \pi_2 \mu \|
>0$. 
\end{lemma}
\begin{proof} We use the notation of formula (4.8).
Recall that we only need to consider $a_{mn} \in 2\ZZ$.
For fixed average $\bar a_m$ of a cycle the minimum norm is obviously attained
if the $a_{mn} \in 2\ZZ$ have no greater pairwise difference than two.
Note that any vector that does not satisfy the condition has a strictly greater
norm. Note further that the norm does not depend on the value of the
average because of $\pi_2$. The vectors $\left( (4,0^{(N-1)}),(0^{(N)}),...
\right)$ are elements of $\sqrt8 \pi_1 \Lambda$. Hence, if $\lambda
\in \sqrt8 \pi_1 \Lambda$ satisfies the assumptions of lemma 4.4 it is possible
to subtract suitable multiples of these
in order to obtain a vector $\mu$ which satisfies the claimed inequality.
$\pi_2 \mu =0$ is
equivalent to all the co-ordinate entries being equivalent modulo 4.
\end{proof}
\begin{lemma} Let 
\begin{equation}
\lambda_q =\left((4^{(q)},0^{(N-q)}),
(0^{(N)}),...\right).
\end{equation}
Then $\| \pi_2 \lambda_q\| =
2{q(N-q)\over N}$. There are no vectors with entries equivalent 0 modulo 4 which
are mapped by $\pi_2$ to non-zero vectors with norm less or equal $2{N-1\over
N+1}$.
\end{lemma}
\begin{proof} The first claim is a straightforward calculation. Obviously none
of the vectors $\lambda_q$ has norm less or equal $2{N-1\over N+1}$. For the
second claim it remains to consider vectors whose entries differ by at least 8.
We use the argument of the proof of lemma 4.4.
\end{proof}

The only $q$ with $||\pi_2\lambda_q||=2{N-1\over N}$ are
$q=1$ and $q=N-1$. As $\pi_2\lambda_1=-\pi_2\lambda_{N-1}$ there are $2(24-M)$
vectors of the type $\pi_2(\lambda_q)$ and norm $2{N-1\over N}$ in $L^*$.
\flexskip

Thus we only have to consider vectors of the form $\pi_L(2c)$ where $c$ is an
element of the 24-dimensional Golay code, understood as element of $\RR^{24}$.
We recall the following three properties of the Golay code, for
details see \cite{CS88}:\flexskip

\beginaxioms
\item[{\bf P1}] Given a Golay code element $d$ with $\pi_2 d = 0$ we conclude that a
Golay code element $c$ has (up to sign) the same image under $\pi_2$ as the
symmetric difference $c\triangle d$. Note that if $c$ and $d$ meet in, say, $a$
places then $c\triangle d$ and $d$ meet in $|d| - a$ places. This will reduce
the number of cases to be considered in the sequel. 
\item[{\bf P2}] An octad meets a dodecad in 2, 4, or 6, places. Any two distinct 
octads meet in 0, 2, or 4, places. (Otherwise their symmetric difference 
would have length not in $\{0,8,12,16,24\}$.)
\item[{\bf P3}] Any dodecad is the symmetric difference of two octads.
\endaxioms \flexskip

We furthermore introduce the following shorthand. For $N=11$,
\begin{equation}
(a_1,a_2):= \left( (0^{(11-a_1)},2^{(a_1)}), (0^{(11-a_2)},2^{(a_2)}) \right)
\in\sqrt8 \pi_1 \Lambda.
\end{equation}
Furthermore, given a
Golay code element $c$, we define its intersection numbers $(x_1,x_2)$ with
respect to the dodecads $d_1$ and $d_2$. An analogous notation will be used for
all $N$. \flexskip

{\bf Case N=23:} $\|\pi_2 (0^{(23-n)},2^{(n)}) \| = {1\over2}
{n(23-n)\over 23}.$ This is smaller than 2 only for $n<6$ or $n>17$.
The 24-dimensional Golay code is such that any element $c$ has 0,8,12,16 or 24
entries 1, the rest being zero. As the projection kills the first co-ordinate
the projected Golay codes can have 0, 7, 8, 11, 12, 15, 16 entries. This 
completes the analysis of any non-zero Golay code element. Hence there are 
no non-zero vectors of norm less or equal $2{N-1\over N+1}$. \flexskip

{\bf Case N=11:} Let $d_1, d_2$ denote the two dodecads which define the cycle
shape of the automorphism. Given an octad $o$ such that $\pi_1(o)=(a_1,a_2)$ in
the notation of (4.10) we obtain the norm as
\[
\| \pi_2 (a_1,a_2) \| = {1\over 22} \bigl( a_1(11-a_1) + a_2(11-a_2) \bigr).
\]
By (P2), up to permutation
there are the following possibilities for the intersection numbers $(x_j)$:
(6,2), or (4,4). Then $a_j=x_j$ or $a_j=x_j-1$ depending on the fact which
entries of $o$ are killed by $\pi_1$. Note further that $a_j=6$ and $a_j=5$ are
equivalent by (P1). Hence we need only consider the following cases:
$(a_j)=$ (6,2), (6,1), (4,4), (4,3), (3,3). The norms are easily obtained as
24/11, 20/11, 28/11, 26/11, 24/11 respectively. None of these
vectors has norm less or equal $2{N-1\over N+1}$. \flexskip

As dodecads cannot intersect in more than 8 points it is straightforward to see
that there cannot arise any vectors of norm less than 2 from
a dodecad. The vectors related to 16-point elements of the Golay code have been
accounted for by considering octads by (P1). \flexskip

Because of lemma 4.4 this concludes the proof that there are no non-zero
vectors of norm less or equal to $2{N-1\over N+1}$. \flexskip


{\bf Case N=7:} Let $o_1,o_2,o_3$ denote the three octads which define the
cycle shape of the automorphism. Given an octad $o$ such that
$\pi_1(o)=(a_1,a_2,a_3)$ in the notation of (4.10) we obtain the norm as
\[
\| \pi_2 (a_1,a_2,a_3) \| = {1\over14} \left( a_1(7-a_1) + a_2(7-a_2) +
a_3(7-a_3) \right)
\]
Given an octad $o\not= o_j$ of the Golay code it meets $o_j$ in
$x_j$ places where $x_j \in \{0,2,4\}$ because of (P2). Hence
up to permutation $(x_j)$ must be one of the following: (4,4,0), or (4,2,2).
Note that $a_j=x_j$ or $a_j=x_j-1$ ($\pi_1$ kills one of the places of the octad
$o_j$). Note further that $a_j=3$ and $a_j=4$ are
equivalent because of (P1). Hence we only need to consider the following cases
for $(a_j)$: (4,4,0), (4,2,2), (4,2,1), (4,1,1). These have norms 12/7, 16/7,
14/7, 12/7 respectively. There are no cases of norm less or equal to
$2{N-1\over N+1}$. \flexskip

Next we consider a dodecad $d$. If it intersects with $o_j$ in 6 places then
$d-o_j$ is an octad. Hence we have already accounted for it. Otherwise, by (P2)
$d$ meets each $o_j$ in 4 places. Hence the projection has norm 18/7. As
before, the vectors related to 16-point elements of the Golay code have been
accounted for by considering the octads. Because of lemma 4.4 this concludes the
list of vectors of norm less or equal to $2{N-1\over N+1}$. \flexskip

{\bf Case N=5:} Let $o_1,o_2,o_3,o_4$ denote the four octads which define the
cycle shape of the automorphism. Note that this time the defining octads
intersect, though only in places which are killed by $\pi_1$.
Given an octad $o$ such that $\pi_1(o)=(a_1,a_2,a_3,a_4)$ in the notation of
(4.10) we obtain the norm as
\[
 \| \pi_2 \bigl((a_j)_{j=1}^4\bigr)
 \| = {1\over10} \sum_{j=1}^4 a_j(5-a_j).
\]
Given an octad $o \not= o_j$ of the Golay code it meets the octads $o_j$ in
$x_j, j=1,...,4$ places respectively. Because of (P1) it is sufficient to
consider the cases where $a_j \in \{0,1,2\}$. The number $k$ of co-ordinates
killed by $\pi_1$ satisfies $0 \le k=8-\sum a_j\le 4$.
Hence up to permutation $(a_j)$ must be one of the following:
(2,2,2,2), (2,2,2,1), (2,2,2,0), (2,2,1,1), (2,2,1,0), (2,2,0,0), (2,1,1,1),
(2,1,1,0), (1,1,1,1). If $k=4$ the intersection numbers $x_j$ of $o$ with $o_j$
are just $x_j=a_j+3$. Because of (P2) the $x_j$ must be less or equal 4. This
excludes (2,2,0,0) and (2,1,1,0). The remaining cases - whether actually
occurring or not - have norm greater than $2{N-1\over N+1}$ and need
not be considered.  \flexskip

Next we consider a dodecad $d$. The sum of the intersection numbers $a_j$
with the 5-cycles must be at least 8. If $(a_j)=(2,2,2,2)$ the
projection has norm 12/5. Else there are intersection numbers of 3 or higher.
By (P1) the symmetric difference with a suitable $o_j$ will be a dodecad or
octad with smaller intersection numbers. Hence we have already accounted for it.
As before, the vectors related to 16-point elements of the Golay code have been
accounted for by considering the octads. Because of lemma 4.4 this concludes the
list of vectors of norm less or equal to $2{N-1\over N+1}$. \flexskip

{\bf Case N=3:} As established above, $L=K_{12}$. Hence
$L^*={1\over\sqrt3}K_{12}$ and we quote from \cite{CS88}, chapter 4.9 that
\[
 \theta_{L^*}(q) = 1 + 756 q^{2\over3} +\dots 
\]
{\bf Case N=2:} As established above, $L=\sqrt2 E_8$. Hence
$L^*={1\over\sqrt2}E_8$ and we quote from \cite{CS88}, chapter 4.8 that
\[
 \theta_{L^*}(q) = 1 + 240 q^{1\over2} +\dots
\]
\vspace{-0.2in}

\theoremproven

\flexskip
 
\section{The Conclusion of the Proof} 

So far we have established that $\theta_L$ is a modular form of $\Gamma_0(N)$
and have determined the leading coefficients of $\theta_L$ and $\theta_{L^*}$.
We now turn to the right hand side of formula (1.25).
For $N=2,3,5,7,11,23$, let $M={24\over N+1}$ and define for $r=0,\dots,N-1$
\par
\[
 \Theta(\tau) := \eta(\tau)^{MN} \left( \eta(N\tau)^{-M}
+\sum_{j=0}^{N-1}\psi_j(\tau)^{-M}\right), \]
\[
 \Theta_r(\tau) := \eta(\tau)^{MN} \left(
\sum_{j=0}^{N-1}e^{-j{-2r\over N}\pi i}\psi_j(\tau)^{-M}\right)
= \eta(\tau)^{MN} \left(
\sum_{j=0}^{N-1}e^{jr{2\pi i\over N}}\psi_j(\tau)^{-M}\right). 
\]
(Remember that we chose $\lambda_r^*$ such that $(\lambda_r^*)^2\equiv
{-2r\over N}\ \hbox{mod}\ 2\ZZ$.)

\begin{theorem} $\Theta$ is a modular form of $\Gamma_0(N)$ of weight
$12{N-1\over N+1}={NM-M\over2}$ and character $\chi_N$. For $N=2,3,5,11$,
$\chi_N$ is trivial, for $N=7,23$, $\chi_N \Bigl({ab \choose cd},*\Bigr) =
\left({d\over N}\right).$ (Here, as in theorems 3.4 and 3.5, we implicitly
understand $*$ to be the branch as obtained when the group element is
represented as a product of generators $(T,1)$ and $(V_k,*)$ of formula (2.1).)
\end{theorem}
\begin{proof} Consider $N=7,23$.
$\Gamma_0(N)$ is generated by $T$ and the set of
$V_k=ST^kST^{k'}S$ (formula (2.1)). Now from the results of chapter 2.2,
\begin{equation} \begin{split}
 \Theta &\vert_{\scriptstyle(T,1,{NM-M\over2})} (\tau) \\
  &= \eta(\tau)^{NM}\vert_{\scriptstyle(T,1,{NM\over2})}
     \left( \eta(N\tau)^{-M}\vert_{\scriptstyle(T,1,-{M\over2})} + \sum_{j=0}^{N-1}
     \psi_j(\tau)^{-M} \vert_{\scriptstyle(T,1,-{M\over2})} \right) \\
  &= e^{{\pi i\over 12}{NM}}\eta(\tau)^{NM} \left( e^{N{\pi i\over 12}(-M)}\eta(N\tau)^{-M}
     +\sum_{j=0}^{N-1}e^{-{\pi i\over 12}(-M)}
     \psi_j(\tau)^{-M} \right) \\
  &= \Theta (\tau). \\
\end{split}\end{equation}
Here we used $NM+M \equiv 0 \;(\hbox{mod}\; 24)$.
\begin{equation} \begin{split}
\Theta &\vert_{\scriptstyle(V_k,*,{NM-M\over2})} (\tau) \\
 &= \eta(\tau)^{\scriptscriptstyle NM}\vert_{\scriptstyle(V_k,*,{NM\over2})}
    \left( \eta(N\tau)^{\scriptscriptstyle -M}\vert_{\scriptstyle(V_k,*,-{M\over2})} +
    \sum_{j=0}^{N-1}
    \psi_j(\tau)^{\scriptscriptstyle -M} \vert_{\scriptstyle(V_k,*,-{M\over2})}
    \right) \\
\end{split}\end{equation}

Now, $(V_k,*,{NM\over2}) ={ \bigl((S,+)(T^k,1)(S,+)(T^{k'},1)(S,+),
{NM\over2}\bigr)} $ acts on $ \eta(\tau)^{NM}$ as multiplication by \par
\begin{equation}
\hbox{exp}\bigl( ({-3\pi i\over 4}+{(k+k')\pi i\over 12})NM \bigr)
\end{equation}

We now proceed to calculate the action of $(V_k,*,-{M\over2})$ in the second
factor of the right hand side of equation (4.12). This will have to be done
case by case for all the summands. We adopt the following shorthand: $\phi_1
{G\atop\mapsto}\phi_2$ means that the action of $G$ takes $\phi_1$ to
$\phi_2$. The individual actions applied in every step have been calculated in
chapter 2.2, in particular theorem 2.3. To simplify notation we work for $M$
rather than $-M$.

\begin{align*} \eta(N\tau)^M \mskip 30mu
  {{\scriptscriptstyle (S,+,{M\over2})} \atop \longmapsto}
    \mskip 30mu &{e^{-M{\pi i\over4}} \over {\sqrt N}^M} \psi_0^M \\
  {{\scriptscriptstyle (T^k,1,{M\over2})} \atop \longmapsto}
    \mskip 30mu &{e^{-M{\pi i\over4}} \over {\sqrt N}^M} e^{-kM{\pi i\over12}}
		       \psi_k^M \\
  {{\scriptscriptstyle(S,+,{M\over2})} \atop \longmapsto}
    \mskip 30mu &{e^{-4M{\pi i\over4}} \over {\sqrt N}^M}
		\left({-k\over N}\right) e^{-kM{\pi i\over12}} \psi_{-k'}^M \\
  {{\scriptscriptstyle (T^{k'},1,{M\over2})} \atop \longmapsto}
    \mskip 30mu &{e^{-4M{\pi i\over4}} \over {\sqrt N}^M}
	       \left({-k\over N}\right) e^{(-k-k')M{\pi i\over12}} \psi_0^M \\
  {{\scriptscriptstyle (S,+,{M\over2})} \atop \longmapsto}
    \mskip 30mu &e^{-5M{\pi i\over4}}\left({-k\over N}\right) e^{(-k-k')M
		{\pi i\over12}} \eta(N\tau)^M
\intertext{For $j\ne 0,k'$ we calculate}
 \psi_j^M  \mskip 30mu
  {{\scriptscriptstyle (S,+,{M\over2})} \atop \longmapsto}
    \mskip 30mu &e^{-3M{\pi i\over4}} {\scriptstyle \left({-j'\over N}\right)}
		       \psi_{-j'}^M \\
  {{\scriptscriptstyle (T^k,1,{M\over2})} \atop \longmapsto}
    \mskip 30mu &e^{-3M{\pi i\over4}} {\scriptstyle \left({-j'\over N}\right)}
		       e^{-kM{\pi i\over12}} \psi_{-j'+k}^M \\
  {{\scriptscriptstyle (S,+,{M\over2})} \atop \longmapsto}
    \mskip 30mu &e^{-6M{\pi i\over4}} {\scriptstyle \left({-j'\over N}\right)
       \left({-(-j'+k)'\over N}\right)} e^{-kM{\pi i\over12}}
			\psi_{-(-j'+k)'}^M \\
  {{\scriptscriptstyle (T^{k'},1,{M\over2})} \atop \longmapsto} \mskip 30mu
       &e^{-6M{\pi i\over4}} {\scriptstyle \left({-j'\over N}\right)
       \left({-(-j'+k)'\over N}\right)}
       e^{(-k-k')M{\pi i\over12}} \psi_{-(-j'+k)'+k'}^M \\
  {{\scriptscriptstyle (S,+,{M\over2})} \atop \longmapsto} \mskip 30mu
       &e^{-9M{\pi i\over4}} {\scriptstyle \left({-j'\over N}\right)
       \left({-(-j'+k)'\over N}\right)
       \left({-(-(-j'+k)'+k')'\over N}\right)}
       e^{{\scriptscriptstyle(-k-k')M}{\pi i\over12}} \times \\
       & \times \psi_{-(-(-j'+k)'+k')'}^M \\
       & = e^{-5M{\pi i\over4}} \left({-k\over N}\right)
       e^{(-k-k')M{\pi i\over12}} \psi_{-(-(-j'+k)'+k')'}^M .
\end{align*}

For the last equality we note that
$\bigl({k\over N}\bigr) = \bigl({k'\over N}\bigr) =
- \bigl({-k\over N}\bigr)$ and that the Legendre symbol is multiplicative. Thus
$\scriptstyle{
   \left({-j'\over N}\right) \left({-(-j'+k)'\over N}\right)
   \left({-(-(-j'+k)'+k')'\over N}\right) = }$ \\
$\scriptstyle{
   \left({-j\over N}\right) \left({(-j'+k)\over N}\right)
   \left({(-(-j'+k)'+k')\over N}\right)= 
   \left({-j\over N}\right) \left({-1+(-j'+k)k'\over N}\right)=
   \left({-j\over N}\right) \left({-j'k'\over N}\right)=
   \left({k\over N}\right).}$ 
\flexskip

Similar calculations yield
\begin{align*} \psi_0^M  \mskip 30mu
  &{{\scriptscriptstyle (V_k,*,{NM\over2})} \atop \longmapsto}
    \mskip 30mu e^{-5M{\pi i\over4}}\left({-k\over N}\right) e^{(-k-k')M
             {\pi i\over12}} \psi_{-k}^M, \\
   \psi_{k'}^M  \mskip 30mu
  &{{\scriptscriptstyle (V_k,*,{NM\over2})} \atop \longmapsto}
    \mskip 30mu e^{-5M{\pi i\over4}} \left({-k\over N}\right) e^{(-k-k')M
          {\pi i\over12}} \psi_0^M.
\end{align*}

We combine formula (4.13) with the above four calculations to obtain that
\[
(V_k,*,{NM-M\over2}) = {\bigl((S,+)(T^k,1)(S,+)(T^{k'},1)
(S,+)},{NM-M\over2}\bigr)
\]
acts on $\Theta$ as multiplication by
$ e^{(-3N+5)M\pi i/4} \left({-k\over N}\right) = \left({-k\over N}\right)$
(remember to revert to $-M$ for $M$ in the second factor!).
Combining this with equation (4.11) we have proven the claimed modularity
properties of $\Theta$ for $N=7,23$ for all generators and hence all elements
of $\Gamma_0(N)$. \flexskip

The proof for $N=3,5,11$ is a simplification of the above argument as we use
theorem 2.2 instead of theorem 2.3. For $N=2$, consider the transformation
behaviour under $(S,+)(T^2,1)(S,+)$ instead. This concludes the proof of
theorem 4.4.
\qed

A comparison with the results of theorem 4.1 then proves
\begin{lemma} $\theta_L$ and $\Theta$ are both modular forms of $\Gamma_0(N)$
of weight $12{N-1\over N+1}$ with identical character.
\end{lemma} \theoremproven

Hence all that remains to be checked is whether the first
coefficients of \\
$\theta_L\vert_{\scriptstyle (S,+,{NM-M\over2})}$ equal those
of $\Theta\vert_{\scriptstyle (S,+,{NM-M\over2})}$.
The transformation under $S$ follows from the results of chapter 2.2

\[
{{(-i {\sqrt N})^M}} \Theta\vert_{(S,+,
{NM-M\over2})}(q) = 
\]

\[
=\eta(q)^{NM} \left( \eta(q^N)^{-M} + { N^M} \psi_0(q)^{-M} +
{(i{\sqrt N})^M} \sum_{j=1}^{N-1} {\left({j\over N}
\right)}^M  \psi_j(q)^{-M} \right) 
\]

It is a straightforward evaluation of the above formula to determine that the
lowest positive power of $q$ with non-zero coefficient is ${N-1\over N}$ for
all $N$.

\begin{theorem} $\theta_L(q) = \Theta(q)$.
\end{theorem}
\begin{proof} We recall that the corollary to theorem 2.1 stated that a modular
form $\phi$ of $\Gamma_0(N)$ of weight $k$ vanishes identically if the
coefficients of the $q$-expansions of both $\phi$ and $\phi\vert_S$ vanish
to order $k\over12$. We will apply this theorem to $\phi(q)=\theta_L(q)-\Theta(q)$.
Note that it follows from lemma 4.6 above that this is a modular form of
$\Gamma_0(N)$ with character of order 2. Furthermore, we obtain
$k=12{N-1\over N+1}$. Thus we need to consider the exponents less than
${N-1\over N+1}$ in the $q$-expansions of $\phi$ and $\phi\vert_S$.
These are trivially zero for $\phi$ itself. In order to calculate the leading
coefficients in the $q$-expansions of $\phi\vert_{\scriptstyle(S,+,{M\over2})}=
\Theta\vert_{\scriptstyle(S,+,{M\over2})}-\theta_L\vert_{\scriptstyle(S,+,
{M\over2})}$ we
note that, by the corollary to theorem 3.3, $(-i)^M\sqrt{N}^M\theta_L
\vert_{\scriptstyle(S,+,{M\over2})}=\theta_{L^*}.$
The leading coefficients of $\theta_{L^*}$ up to and including exponent
${N\over N+1}$ have been determined in chapter 4.3.
From those results we conclude that the coefficients of the $q$-expansion of
$\phi|_S$ up to $N-1\over N+1$ are zero.
\end{proof}

To complete the proof of theorem 1.7 we have to verify formula (1.25) which
can be rewritten as 
\[
\theta_{L+\lambda_r^*} (q) = \Theta_r(q).
\]
We begin by considering the sums
\begin{multline*}
  \sum_{j=0}^{N-1} e^{jr{2\pi i\over N}} \psi_j(q)^{-M}=
  \sum_{j=0}^{N-1} e^{jr{2\pi i\over N}}
	\sum_m c(m)\left( e^{j{2\pi i\over N}}q^{1\over N}\right)^m \\
  =\sum_m \sum_{j=0}^{N-1} \left(e^{2\pi i\over N}\right)^{jr+jm}
	c(m)q^{m\over N}.
\end{multline*}

As $\sum_{j=0}^{N-1} e^{-jx{2\pi i\over N}} = 0$ unless
$x \equiv 0$ mod $N$ it follows that $\Theta_r$ has a $q$-expansion with
exponents in $\ZZ-{r\over N}$ only. \flexskip

Furthermore, we calculate for $N=7,23$, using $\rho(r)=\rho_M(r,N)$ of theorem
4.3, 
\begin{align*}
\Theta &+ \sum_{r=0}^{N-1} \rho(r) \Theta_r \\
 &= \eta(q)^{NM} \Bigl( \eta(q^N)^{-M} + \sum_{j=0}^{N-1}\psi_j(q)^{-M} +
    \sum_{r=0}^{N-1} \rho(r) \sum_{j=0}^{N-1} e^{jr{2\pi i\over N}} \psi_j(q)^{-M}
    \Bigr)  \\
 &= \eta(q)^{NM} \Bigl( \eta(q^N)^{-M} + \left( 1+ \sum_{r=0}^{N-1} \rho(r)
    \right) \psi_0(q)^{-M} \\
 &  \qquad\qquad + \sum_{j=1}^{N-1} \left( 1 + \sum_{r=0}^{N-1} \rho(r) 
    e^{jr{2\pi i\over N}} \right) \psi_j(q)^{-M} \Bigr)
\end{align*}
Now $1 + \sum_{r=0}^{N-1} \rho(r)$ is the number of all possibilities to
represent any residue modulo $N$ as a sum of $M$ squares (cf. equation (4.5)),
hence it equals $N^M$. And
\begin{align*}
 1 + &\sum_{r=0}^{N-1} \rho(r) e^{jr{2\pi i\over N}}  =
      1 + \left( N^{M-1} - 1 \right) + \sum_{r=1}^{N-1} \left(N^{M-1}+
      \left({r\over N}\right)N^{M-1\over2}\right) e^{jr{2\pi i\over N}}  \\
   &= N^{M-1} + N^{M-1} \sum_{r=1}^{N-1} e^{jr{2\pi i\over N}} + N^{M-1\over2}
      \sum_{r=1}^{N-1} \left({r\over N}\right) e^{jr{2\pi i\over N}} \\
   &= N^{M-1} + N^{M-1} (-1) + N^{M-1\over2} i\sqrt{N} \left({j\over N}\right)
      = i{\sqrt N}^M \left({j\over N}\right). 
\end{align*}
Hence we have shown that for $N=7$ and $N=23$ 
\[
\Theta + \sum_{r=0}^{N-1} \rho(r) \Theta_r 
 = (-i\sqrt N)^M \Theta\vert_{(S,+,{NM-M\over2})}
\]
An analogous calculation shows that this formula is correct for $N=2,3,5,11$ as
well. The proof is a simplification of the argument above as the numbers
$\rho(r)$ are equal for all non-zero $r$. \flexskip Furthermore,
\[
(-i\sqrt N)^M \Theta\vert_{\scriptstyle(S,+,{NM-M\over2})} - \Theta =
\sum_{r=0}^{N-1} \rho(r) \Theta_r
\]
From theorem 3.3 and chapter 4.2 (in particular theorem 4.2) we know that
\[
(-i\sqrt N)^M \theta_L \vert_{\scriptstyle(S,+,{NM-M\over2})} - \theta_L =
 \theta_{L^*} - \theta_L = \sum_{r=0}^{N-1} \rho(r) \theta_{L+\lambda_r^*}
\]
where $\theta_{L+\lambda_r^*}$ is the theta-function of the lattice $L$
translated by an element $\lambda_r^*\in L^*$ such that
${(\lambda_r^*)^2\over2} \equiv {-r\over N}\ \hbox{mod}\ \ZZ$.
Now, if $\lambda\in L$, ${1\over2}(\lambda_r^*+\lambda)^2 =
{(\lambda_r^*)^2\over2} + (\lambda_r^*,\lambda) +{\lambda^2\over2}
\in \ZZ-{r\over N}$. This
shows that $\theta_{L+\lambda_r^*}$ has a $q$-expansion with
exponents in $\ZZ-{r\over N}$ only, just as $\Theta_r$. This proves that
$\theta_{L+\lambda_r^*} = \Theta_r$ for all $r$ and concludes the proof of
theorem 1.7.
\end{proof}

\chapter{The Real Simple Roots} 

Let $\sigma$ be an automorphism of the Leech lattice $\Lambda$ of order $N$
and cycle shape $1^MN^M$, as before. In this chapter we determine the real
simple roots of the GKM ${\cal G}_N$ constructed from $\sigma$. As
the imaginary simple roots of ${\cal G}_N$ were explicitly described in
theorem 1.6, this will therefore complete the explicit identification
of all simple roots of ${\cal G}_N$, and enable us to calculate all entries
into the generalized Cartan matrix of ${\cal G}_N$.
We identify the set of real simple roots with a set ${\cal R}$
which consists of the fixed point lattice $\Lambda^\sigma$ and some
elements of its dual ${\Lambda^\sigma}^*$. We prove a number of results
concerning the decomposition of space induced by the set ${\cal R}$. This
programme is a generalization of work by Borcherds. In \cite{Bor92} the real
simple roots of the monster Lie algebra are identified with the
24-dimensional Leech lattice. The holes of the Leech lattice,
corresponding to the finite and affine subalgebras of the monster Lie
algebra, are enumerated in \cite{CS88}, chapter 25. Returning to the
GKMs ${\cal G}_N$, the set ${\cal R}$ does not form a lattice though
many of the properties are preserved. In particular, the decomposition of
space will only produce holes which have finite and affine diagrams. We
prove that they provide the complete classification of all finite and
affine subalgebras of the GKM ${\cal G}_N$. The results of chapter 5 will
be used in chapter 6 to carry out the explicit decomposition of space into
finite and affine holes. This in turn enables us to identify hyperbolic
subalgebras of the GKM ${\cal G}_N$ and thus to obtain upper bounds for
their root multiplicities. \flexskip

As identified in theorem 1.6, the real simple roots of the GKM ${\cal G}_N$
are the simple roots of its Weyl group. In section 5.1 we determine the set
of simple roots of the Lorentzian lattice explicitly. We identify this set
with the set ${\cal R}$. Section 5.2 will show how subsets of ${\cal R}$
translate to Dynkin diagrams and which Dynkin diagrams can be expected.
As a by-product, we will be able to give an explicit description of the
generalized Cartan matrix of ${\cal G}_N$.
Section 5.3 establishes a number of results on the volumes of the holes
corresponding to those Dynkin diagrams. Section 5.4 looks into the symmetries
and automorphisms of the holes of ${\cal R}$. These will be used to
establish relations between the hole decomposition of the Leech lattice
and the hole decompositions of interest in this work.

\section{The Set of Real Simple Roots}

This section identifies the real simple roots of the GKM ${\cal G}_N$
associated with an automorphism $\sigma$ of the Leech lattice $\Lambda$
of order $N$, as constructed in the theorem 1.6. We consider all relevant
orders $N=$ 2, 3, 5, 7, 11, 23. Let, as before, $\Lambda^\sigma$ denote
the fixed point lattice, $L=\Lambda^\sigma\oplus II_{1,1}$ the
corresponding Lorentzian lattice. Let
$\pi_{\Lambda^\sigma} : \Lambda \to \Lambda^\sigma$ and
$\pi_L : \Lambda \oplus II_{1,1} \to L$
be the respective projections. Formula (1.22) expressed the
multiplicities of roots in terms of the trace of $\sigma$ on the root
space $E_r$ (see section 1.5). Here $r$ is an element of $W^\sigma$,
the Weyl group of the fixed point lattice $\Lambda^\sigma \oplus II_{1,1}$,
which was characterised in lemma 1.1. \flexskip

By lemma 1.3 and 1.2a, the twisted Weyl group $W^\sigma$ is the full
reflection group of $L$. Lemma 1.2b and 1.2c give an explicit description
of the simple roots of $L$.

\begin{theorem} The (primitive) simple roots of the lattice $L$, and
thus the (primitive) real simple roots of the GKM ${\cal G}_N$, are the
following:
\begin{equation}
\bigl(\lambda,1,{\lambda^2\over2}-1\bigr) \mskip 30mu {\it for}  \mskip 30mu
\lambda\in \Lambda^\sigma
\end{equation}
and
\begin{equation}
\bigl(\lambda,N,{\lambda^2\over2N}-1\bigr)
\mskip 30mu {\it for} \mskip 30mu \lambda\in N{\Lambda^\sigma}^* \mskip 15mu 
{\it such}\ {\it that}\mskip 15mu N|\ \bigl({\lambda^2\over 2N}-1\bigr).
\end{equation}
The former simple roots have height $1$ and norm $2$ whereas the latter have
height $N$ and norm $2N$.
\end{theorem}
\begin{proof} We begin by determining all roots of $L$. As introduced in the
context of lemma 1.2, we only consider primitive roots. Suppose,
$r\in L$ is a (primitive) root of $L$. By definition, the norm of $r$ must
be greater than 0. Because $r$ is primitive, if $n\in \ZZ$, $\vert n\vert>1$,
then $r\over n$ is not in $L$. We recall that the dual lattice of $L$
satisfies $L^*= (\Lambda^\sigma)^*\oplus II_{1,1}$. For any $v\in L$, the
inner product $(r,v)$ is integer, hence the greatest common divisor $d=(r,L)$
of the absolute values of all inner products $(r,v)$, $v\in L$, is defined.
It follows that ${r\over d} \in L^*$. Lemma 4.1 implies that $NL^*
\subset L$. As $N$ is prime the only possible cases are $d=1$ or $d=N$. 
If $d=(r,L)=1$ then $r^2$ must necessarily be equal to 1 or 2 because
$2(r,v)/(r,r)$ must be integer for all $v\in L$. As $L$ is
even, only the case $r^2=2$ occurs. Hence, this case yields precisely
all $r\in L$ of norm 2 as roots. 
The remaining case, $d=(r,L)=N$, implies that $r^2 = N$ or $r^2 =2N$.
For $N\not=2$, the case $r^2 = N$ is impossible because $N$ is odd and
$L$ even. We conclude that $r\in NL^*$ will be a primitive root if and
only if $r^2=2N$. This determines the primitive roots of $L$.  \flexskip

We now turn to the identification of the simple roots of $L$, which, at
the same time, are the real simple roots of ${\cal G}_N$ following theorem
1.6. We use the description of the simple roots of $W^\sigma$ in lemma 1.2
to derive formulas (5.1) and (5.2). Theorem 1.6 established that the
Weyl vector $\rho$ has the coordinates $(0,0,1)$. Let $r=(\lambda,h,n)$ be a
(primitive) simple root of $W^\sigma$, so of norm 2 or norm $2N$ by
the first part of the proof.
If $N$ does not divide $h$ then $r\not\in NL^*$, hence $(r,L)=1$. As seen in
the first part of this proof this implies $r^2=2$, and hence lemma 1.2b
requires $(r,\rho)=-1$ which is $h=1$. From $r^2=2$ we conclude that $r$ is
of the type of formula (5.1). The root $r$ in formula (5.1) furthermore
satisfies the scaling condition of lemma 1.2c so that $r$ is the primitive
simple root identified in lemma 1.2. \flexskip

If, on the other hand $N$ does divide $h$, then $h=Nh'$ and lemma 1.2b implies
that $-(\rho,r)=Nh'$ divides $(r,v)$ for all $v\in L$; hence we conclude that
$r\in h'NL^*$, and ${r\over h'} \in NL^* \subset L$. Thus $r$ will only be
primitive if $h'=1$. As seen in the first part of this proof, this implies
$r^2=2N$. Hence $r$ is of the type of formula (5.2). The root $r$ in formula
(5.2) furthermore satisfies the scaling condition of lemma 1.2c so that $r$
is the primitive simple root identified in lemma 1.2.
\end{proof}

If $(\lambda,N,N*)$ is a real simple root of norm $2N$ then
$({\lambda\over N},1,*)$ is an element of $L^*$ that defines the same
reflection. In particular, $\lambda\over N$ is an element of
$(\Lambda^\sigma)^*$ of norm ${2\over N}$ mod $2\ZZ$. Conversely, theorem 5.1
shows that any such element of the dual represents a real simple root of the
GKM ${\cal G}_N$. Let us introduce the notation ${\cal R}_{dual}$ for those
elements of the dual lattice which represent real simple roots of norm $2N$.
Further, we will write ${\cal R}_{fix}$ for the fixed point lattice when it
is understood to represent the real simple roots of norm 2 of the GKM
${\cal G}_N$. The totality of all real simple roots of the GKM ${\cal G}_N$
will then be represented by the set ${\cal R}$ which is the disjoint union
of ${\cal R}_{fix}$ and ${\cal R}_{dual}$. \flexskip

We have the following result of \cite{Bor90a} about
the geometry of the fixed point lattice:

\begin{theorem} The fixed point lattice $\Lambda^\sigma$ can be covered
by spheres of radius $\sqrt 2$ centered at the elements of ${\cal R}_{fix}$
plus spheres of radius $\sqrt{2\over N}$ centered at the elements of 
${\cal R}_{dual}$. 
\end{theorem}
\begin{proof} This is a direct consequence of theorem 3.1
of \cite{Bor90a}. The affine space described there can be identified with the
fixed point lattices considered in this work. Then, the radius of the spheres
is generally identified as $\sqrt{\bigr( {r\over(r,\rho)} \bigl)^2}$.
\end{proof}

\begin{remark} We can understand this result as follows. We recall from chapter 4 that
the shortest non-zero vectors $\lambda^*=\pi_{(\Lambda^\sigma)^\perp}(\lambda)$
of ${{\Lambda^\sigma}^\perp}^*$ have norm $2-{2\over N}$. This implies
$\pi_{\Lambda^\sigma}(\lambda)$ is an element of ${\cal R}_{dual}$. Hence,
if we intersect a covering of the Leech lattice (spheres of radius $\sqrt2$) with
the span of $\Lambda^\sigma$ we obtain precisely the kind of covering as
described in theorem 5.2. 
\end{remark}

We furthermore observe that - even though the set of real simple roots
${\cal R}$ is not a lattice - the group of automorphisms of this set is the
group of automorphisms of the fixed point lattice ${\cal R}_{fix}$.

\section{Holes and Dynkin Diagrams}

Given a GKM with (symmetric) Cartan matrix $C$ and simple roots $\{r_i\}$,
we may define a diagonal matrix $D=(d_{ij})$ such that the diagonal
elements of the matrix $A = DC$ are either equal to 2, or less or equal
zero. We may now associate a Dynkin diagram with any subset of simple roots
which consists entirely of real simple roots. We use the conventions for 
Dynkin diagrams laid out in \cite{Wan91}, because in chapter 6 we will relate
the results of this work to the classifications in \cite{Wan91}. \flexskip

We will use the following conventions: The Dynkin diagram of a subalgebra
contains one node for every real simple root. Suppose the submatrix of the
matrix $A$, corresponding to the real simple roots $r_i$ and $r_j$ is of
form ${{\ 2\ a_{ij}} \choose{a_{ji}\ 2\ }}$ then their respective nodes in
the Dynkin diagram will be linked by max$(|a_{ij}|,|a_{ji}|)$ bonds with
additional arrows pointing towards $i$ if $|a_{ij}|>1$. Note that this
notation does contain all information about the Cartan matrix in the cases
of the finite, affine, and hyperbolic diagrams we will use it for. In
particular, ${{\ 2\ -2}\choose{-2\ 2\ }}$ corresponds to a double bond
with arrows pointing to both sides. ${{\ 2\ -1} \choose{-1\ 2\ }}$
corresponds to a single bond. Note further that the conventions imply that
one-sided arrows point toward the shorter root should there be one. \flexskip

The Dynkin diagram has so far been associated to the real simple roots of the
GKM ${\cal G}_N$. We will now establish how to associate it directly to
subsets of ${\cal R}={\cal R}_{fix} \cup {\cal R}_{dual}$.
There are three cases. \par

{\bf Case A:} Both real simple roots are of norm 2.
Then $r_1=\bigl(\lambda,1,{\lambda^2\over2}-1\bigr)$ and
$r_2=\bigl(\mu,1,{\mu^2\over2}-1\bigr)$. They correspond to
$\lambda,\mu\in{\cal R}_{fix}$ and
\begin{subequations}
\begin{equation}
r_1 \cdot r_2 = \lambda\cdot\mu - ({\lambda^2\over2}-1)-({\mu^2\over2}-1)
=2-{1\over2}(\lambda-\mu)^2
\end{equation}
Thus the two nodes will be unlinked if $(\lambda-\mu)^2=4$, they
will be linked by one bond if $(\lambda-\mu)^2=6$, they will be linked by a
double bond with arrows on both sides if $(\lambda-\mu)^2=8$.
Other types of bond will not arise in this work. \par

{\bf Case B:} Both real simple roots are of norm $2N$.
Then $r_1=\bigl(\lambda,N,{\lambda^2\over2N}-1\bigr)$ and
$r_2=\bigl(\mu,N,{\mu^2\over2N}-1\bigr)$. They correspond to ${\lambda\over N},
{\mu\over N}\in{\cal R}_{dual}$ and
\begin{multline}
r_1 \cdot r_2 = \lambda\cdot\mu - ({\lambda^2\over2}-N)-({\mu^2\over2}-N) \\
=2N-{1\over2}(\lambda-\mu)^2 ={N^2\over2} \Bigl[ {4\over N} -
\bigl({\lambda\over N}-{\mu\over N} \bigr)^2 \Bigr].
\end{multline}
Thus the two nodes will be unlinked if $({\lambda\over N}-
{\mu\over N})^2= {4\over N}$, they will be linked by one bond if
$({\lambda\over N}-{\mu\over N})^2={6\over N}$, they will be linked by a double
bond with arrows on both sides if $({\lambda\over N}-{\mu\over N})^2=
{8\over N}$. Other types of bond will not arise in this work. \par

{\bf Case C:} The real simple root $r_1$ is of norm 2, $r_2$ is of norm $2N$.
Then $r_1=\bigl(\lambda,1,{\lambda^2\over2}-1\bigr)$ and
$r_2=\bigl(\mu,N,{\mu^2\over2N}-1\bigr)$. They correspond to
$\lambda\in{\cal R}_{fix}$, ${\mu\over N}\in{\cal R}_{dual}$  and
\begin{equation}
r_1 \cdot r_2 = \lambda\cdot\mu - N({\lambda^2\over2}-1)-({\mu^2\over2N}-1)
=N+1-{N\over2}(\lambda-{\mu\over N})^2.
\end{equation}
\end{subequations}
Thus the two nodes will be unlinked if $(\lambda-{\mu\over N})^2=2+{2\over N}$,
that is, if they are of minimal distance; they will be linked by $N$ bonds with
an arrow pointing toward $r_1$ if
$(\lambda-{\mu\over N})^2=4+{2\over N}$. Note that for a fixed element
of the dual $\mu \over N$ this exactly describes the set of lattice elements
closest but one to $\mu \over N$.
Other types of bond will not arise in this work. \flexskip

The above calculations complete, as a corollary, the explicit specification
of the elements of the generalized Cartan matrices for the GKMs constructed
in this work.
\begin{theorem} The generalized Cartan matrix $C_N$ of the GKM
${\cal G}_N$ has been completely and explicitly determined.
\end{theorem}
\begin{proof} The set of imaginary simple roots has been identified in
theorem 1.6, the set of real simple roots in theorem 5.1. The products of
any two real simple roots have been determined in equations (5.3). The
remaining products involving imaginary simple roots, that is positive
multiples $n\rho, m\rho$ of the Weyl vector $\rho$ are as follows;
\begin{equation} \begin{split}
           (n\rho, m\rho)                                        &=0, \\
           \Bigl(n\rho, (\lambda, 1, {\lambda^2\over2-1})\Bigr)  &=-n,\\
           \Bigl(n\rho, (\lambda, N, {\lambda^2\over2N-1})\Bigr) &=-nN.
\end{split}\end{equation}
This completes the explicit specification of the generalized Cartan matrix.
\end{proof}

We need to recall some facts about Lie algebras from \cite{Kac90}. Given a system
of $n$ simple roots $\alpha_1$, ..., $\alpha_n$ of a Lie algebra with
Cartan matrix $C=(a_{ij})$, there is a system of
co-roots $\alpha_1^\vee$, ..., $\alpha_n^\vee$ defined by 
\begin{equation}
\langle\alpha_i^\vee,\alpha_j\rangle=a_{ij}\mskip 20mu (i,j=1,\dots,n)
\end{equation}
We define a scalar product in the space of roots $[\alpha]$ and the space of
co-roots, $[\alpha^\vee]$.
We can then find an isomorphism $\nu:[\alpha^\vee] \to [\alpha]$ such that
\begin{equation}
\nu(\alpha_i^\vee)={2\over(\alpha_i,\alpha_i)}\alpha_i.
\end{equation}
We may form a Cartan matrix from a subset of the real simple roots of the
GKM ${\cal G}_N$. In that situation, the above works out as follows:
$\alpha_i^2=2$ implies $(\alpha_i^\vee)^2=2$, and $\alpha_i^2=2N$ implies
$(\alpha_i^\vee)^2={2\over N}$. We further recall that for any simple
finite-dimensional Lie algebra there exists a Weyl vector. This is a vector
$\rho$ of the root space satisfying
\begin{subequations}
\begin{equation}
(\rho,\alpha_i)={(\alpha_i,\alpha_i)\over2}\mskip 20mu \hbox{for\
all\ } i.
\end{equation}
Unlike section 1.1.3, in this case, we use the standard sign convention for
the Weyl vector to simplify comparison with the statements of \cite{Kac90}.
If we define $\rho^\vee:=\nu^{-1}\rho$ then
\begin{equation}
(\rho^\vee,\alpha_i^\vee)=\left(\nu^{-1}\rho,
\nu^{-1}\bigl({2\over\alpha_i^2}\alpha_i\bigr)\right) =
  {2\over\alpha_i^2}(\rho,\alpha_i)=1.
\end{equation}
\end{subequations}
We can express $\rho^\vee$ as an integer linear combination of the simple
co-roots: $\rho^\vee=\sum_i n_i^\vee\alpha_i^\vee$. Then
\begin{equation}
\rho^2={\rho^\vee}^2=\sum_in_i^\vee(\rho^\vee,\alpha_i^\vee)=\sum_i
n_i^\vee. 
\end{equation}
For any simple Lie algebra of affine type there exists a vector $\delta$ such
that
\begin{subequations}
\begin{equation}
\delta^2=(\delta,\alpha_i)=0\mskip 20mu \hbox{for\ all\ } i.
\end{equation}
$\delta$ can be expressed as a linear combination of the simple roots:
$\delta=\sum_i n_i \alpha_i$. $\delta$ is then characterised by the condition
that the $n_i$ are integer with greatest common divisor 1. Recall from
section 1.1.2 that, in the framework of \cite{Jur96}, we consider roots in a root
space which is extended by degree derivations. The vector $\delta$ will then
have non-zero components in the degree derivations (which are part of the
null space of the bilinear form). Similarly, there exists a vector
$\delta^\vee=\sum_i n_i^\vee\alpha_i^\vee$ such that 
\begin{equation}
{\delta^\vee}^2
=(\delta^\vee,\alpha_i^\vee)=0\mskip 20mu \hbox{for\ all\ } i,
\end{equation}
\end{subequations}
uniquely determined by the condition that the $n_i^\vee$ are integer with
greatest common divisor 1. If $\delta=\sum n_i\alpha_i$ let $d$ denote the
greatest common divisor of the integers $n_i{\alpha_i^2\over2}$.
We find that 
\begin{equation}
\delta^\vee=\sum_i n_i^\vee\alpha_i^\vee={1\over
d}\nu^{-1}\delta,
\end{equation}
We finally recall the definition of the Coxeter number and the dual Coxeter
number of an affine Lie algebra:
\begin{equation}
h=\sum_i n_i,\mskip 60mu h^\vee=\sum_i n_i^\vee.
\end{equation}
Next, we need to turn to semisimple finite and affine Lie algebras, that is,
direct sums of simple Lie algebras. The roots of different components are
mutually orthogonal. Hence, we find that there exist $\rho$ (or $\delta$
respectively) as a direct sum of the $\rho$ (or $\delta$) of the components.
\flexskip

We are now in the position to generalize the notion of a hole in a lattice
to an analogue in the set ${\cal R}$ representing the real simple roots of
the GKM ${\cal G}_N$. For any point in the $2M$-dimensional space spanned by
the fixed point lattice we define a radius function as follows
\begin{equation}
\hbox{r}_{\cal R}(x) =\ \ \hbox{min}\{\sqrt{(x-v)^2},
\ \sqrt{(x-d)^2+2-{2\over N}}\  |
\ v\in {\cal R}_{fix}, \ d\in {\cal R}_{dual}\}.
\end{equation}
The proof of lemma 5.1 below will clarify the motivation behind this definition
for a generalized radius.
Because of theorem 5.2, the radius function is bounded above by $\sqrt2$. A
local maximum of the radius function can then be considered as a generalized
hole of the set of real simple roots ${\cal R}$. The value of the radius
function will be called the generalized radius of the hole. As the dimension
of space is $2M$ there will be at least $2M+1$ real simple roots closest
(in the sense of the radius function) to a generalized hole.
The real simple roots which are closest in the sense of the radius function,
shall be referred to as the generalized vertices of the generalized hole.
They will include both elements of ${\cal R}_{fix}$ and ${\cal R}_{dual}$.
They form a Dynkin diagram in the way described above, corresponding to a
subalgebra of the GKM ${\cal G}_N$. \flexskip

For the remainder of this work let us introduce some notational conventions
with respect to generalized holes. Let $H$ denote the set of generalized
vertices of such a hole, let $\langle H\rangle$ denote the convex hull of
the vertices, that is the `hole' in its spatial meaning, and let $\Delta(H)$
denote the Dynkin diagram associated to the hole. We will drop the brackets
$\langle\rangle$ occasionally when there is no need to distinguish between a
hole and its set of vertices. 

\begin{lemma} Suppose that we are given a subset $H$ of $\cal R$ such
that $\Delta(H)$ is of finite, or affine, type. Then there exists a
(generalized) centre $c$ of (generalized) radius $r_0$ from all elements of
$H$. $c$ lies within the convex hull $\langle H\rangle$. \par
b) Suppose that $H$ is as in a). Then there exists no element of ${\cal R}$
which has smaller (generalized) distance from $c$ than $r_0$. \par
c) Now suppose that $\Delta(H)$ is of finite type and that $H$ contains $n$
elements where $n\leq2M$ ($2M$ is the dimension of the fixed point lattice).
For any (generalized) radius $r$, such that $r_0<r\leq \sqrt2$, there exists a
$(2M-n)$-dimensional manifold of points which have (generalized) radius $r$ from
all elements of $H$.
\end{lemma}
\begin{proof} a) For any Cartan matrix of finite or affine type we can choose
a system of corresponding roots (or co-roots) in Euclidean space. Now in ${\cal
R}$ we consider vectors of norm 2 and norm $2\over N$. Thus we represent these
by the relevant co-roots $r_i^\vee=(\lambda_i,1,*)$, rather than the roots.
Let us first consider the affine case.
We begin by an explicit identification of the centre $c$.
Let $n_i^\vee,h^\vee$ be defined as in equations (5.10), (5.11).
From (5.9) we obtain $\delta^\vee=(\sum_i n_i^\vee\lambda_i, h^\vee, *)$.
Let 
\begin{subequations}
\begin{equation}
c=\sum_i {n_i^\vee\over h^\vee}\lambda_i
\end{equation}
Then 
\[
(c-\lambda_j)^2=(c-\lambda_j,0,*)^2=({\delta^\vee\over h^\vee}-
r_j^\vee)^2= {r_j^\vee}^2.
\]
Hence, $c$ is the generalized centre of a hole of generalized radius
$r_0=\sqrt2$. It is well known that the solutions $\delta^\vee$ have
coefficients all greater zero which by definition is equivalent to the claim
that $c$ lies within the convex hull. \par

In the case of finite Lie algebras
let $\rho,\ \rho^\vee$ be defined as in equation (5.7). The multiple
$\rho^\vee\over \rho^2$ is contained within the affine span of the co-roots
because of (5.8). Then
\[
({\rho^\vee\over\rho^2}-r_i^\vee)^2={\rho^2\over\rho^4}-
{2(\rho^\vee,r_i^\vee)\over\rho^2}+(r_i^\vee)^2= (r_i^\vee)^2-{1\over\rho^2}.
\]
Let us now define 
\begin{equation}
c={\sum_i n_i^\vee \lambda_i \over\rho^2}
\end{equation}
\end{subequations}
Then ${\rho^\vee\over\rho^2}=(c,1,*)$ and thus
$({\rho^\vee\over\rho^2}-r_j^\vee)^2=
(c-\lambda_j,0,*)^2=(c-\lambda_j)^2.$ Hence we have identified $c$ as the
generalized centre in terms of elements of ${\cal R}$ and justified the
definition of the generalized radius function. This also provides the
minimal radius $r_0=2-{1\over\rho^2}$. A case by case analysis shows that
the coefficients obtained for the centre again are all greater than zero and
thus $c$ is contained in the convex hull $\langle H\rangle$. \flexskip

b) We begin by considering the following special case:
Suppose that the finite diagram $H$ does contain a root of norm 2, say $r_1$.
We will prove that no element $p$ of ${\cal R}_{fix}$ lies closer to $c$
than $r_0$.\par

Without loss of generality we can assume that the element $p$ is the origin
in ${\cal R}$ as ${\cal R}$ is translation invariant. Thus we consider the
following setup: the finite-dimensional Lie algebra is represented by elements
of ${\cal R}_{fix}$ of norms 4, 6, 8,... and elements of 
${\cal R}_{dual}$ of norms $2+{2\over N}$, $4+{2\over N}$,... We 
need to show that $c^2 \geq (c-\lambda_1)^2.$
This is equivalent to $(2c-\lambda_1,\lambda_1) \geq 0.$
We return to the equation $\rho^\vee=\sum n_i^\vee r_i^\vee$ and recall that
$r_i^\vee=(\lambda_i,1,x_i)$. If $(r_i^\vee)^2=2$, then
$x_i={\lambda_i^2\over2}-1 \ge 1$. If $(r_i^\vee)^2={2\over N}$, then
$x_i={\lambda_i^2\over2}-{1\over N} \ge {{2+{2\over N}}\over 2}-{1\over N}= 1$.
Now $(\rho^\vee,r_1^\vee)=1={(r_1^\vee,r_1^\vee)\over2}$. This implies
\[
\bigl( 2{\rho^\vee\over\rho^2}-r_1^\vee,\ r_1^\vee \bigr) = 2({1\over\rho^2}-1).
\]
Rewriting this in terms of the elements $\lambda_i$ of ${\cal R}$ we obtain:
\begin{multline*}
{2\over\rho^2}-2=
\Bigl(\bigl(2c-\lambda_1,1,2{\sum_i n_i^\vee x_i\over \sum_i n_i^\vee}-
 ({\lambda_1^2\over2}-1)\bigr),
  \bigl(\lambda_1,1,{\lambda_1^2\over2}-1\bigr)\Bigr) \\
  = (2c-\lambda_1,\lambda_1)-2{\sum_i n_i^\vee x_i\over \sum_i n_i^\vee}.
\end{multline*}
\[
(2c-\lambda_1,\lambda_1)= {2\over\rho^2}-2+2{\sum_i n_i^\vee x_i\over \sum_i
n_i^\vee}  \geq {2\over\rho^2}>0,
\]
as every individual $x_i$ was greater or equal 1.
This concludes the argument for the special case. There remain the following
cases for the finite types: 1) $p\in{\cal R}_{fix}$ but all roots of the
finite diagram have norm $2N$. 2) $p\in{\cal R}_{dual}$. It suffices to
remark for these cases that the proof works along exactly the same lines of
argument. We only have to adjust the norms of the vectors involved.
For the affine cases we use the vector $\delta^\vee$ of equation (5.9) instead
of $\rho^\vee$.\flexskip

c) We now identify the manifold of centres for radius greater than the minimal
radius $r_0$. Suppose that $H=H_{f}\cup H_{d}$ with $H_{f} \subset
{\cal R}_{fix}$ and $H_{d} \subset {\cal R}_{dual}$. Suppose
further that $|H_{f}|=n_{f}$ and $|H_{d}|=n_{d}$. Because
the $n_{f}$ vectors are affinely independent there exists an affine centre
$c_{f}$ for them. We consider the space $F^\perp$ orthogonal
to the affine span of the $n_{f}$ vectors, passing through $c_{f}$.
$F^\perp$ has $2M-n_{f}+1$ dimensions. Equally, because the $n_{d}$
dual vectors are affinely independent there exists an affine centre $c_{d}$
for them. We consider the space $D^\perp$ orthogonal to the affine span of the
dual vectors, passing through $c_{d}$. $D^\perp$ has $2M-n_{d}+1$
dimensions. \par

The set $H$ was assumed to be affinely independent. Hence
the affine spaces $F^\perp$ and $D^\perp$ intersect in a $(2M-n+2)$-dimensional
affine space. Every point of this space has a certain radius $r_{f}$ which
is the distance to any element of $H_{f}$, and a certain radius
$r_{d}$ which is the distance to any element of $H_{d}$. In turn, all
points of this space of a fixed $r_{f}$ form a sphere in the affine space
centered at $\tilde c_{f}$, which is the projection of $c_{f}$. Note
that there is a minimal $r_{0f}<\sqrt2$ which corresponds to
the distance of $\tilde c_{f}$ from the elements of $H_{f}$. Similarly,
the points of the intersection space of distance $r_{d}$
from the elements of $H_{d}$ form spheres
centered at some $\tilde c_{d}$. Again, there exists a minimal
$r_{0d}<\sqrt{2\over N}$. Clearly,
$\tilde c_{d} \ne \tilde c_{f}$. \par

An element of the $2M$-dimensional space spanned by the fixed point lattice
will have distance $r$ from all elements of $H_{f}$ and distance
$\sqrt{r^2-2+{2\over N}}$ from all elements of $H_{d}$ if and only if it
lies in the intersection of the corresponding spheres
\[
  S\left(\tilde c_{f},\sqrt{r^2-r_{0f}^2}\right) \cap
  S\left(\tilde c_{d},\sqrt{r^2-2+{2\over N}-r_{0d}^2}\right).
\]

Two spheres in space of different radius and centre can intersect in the
following ways: \par

\beginaxioms
\item [$\alpha$)] the volumes enclosed by the spheres are disjoint. \par

\item [$ \beta$)] the volumes enclosed by the spheres intersect in one
                  point. \par

\item [$\gamma$)] the intersection is a proper subset of both volumes
                  enclosed, but does not consist of a single point. It then
                  is a $(2M-n)$-dimensional manifold.\par

\item [$\delta$)] the volumes enclosed lie within one another.
\endaxioms

Clearly, for the centre $c$ and radius $r_0$, situation $\beta$) holds. For
radius $\sqrt2$ situation $\gamma$) must hold because for instance the elements
of the dual lattice and those of the lattice itself have distance greater than
2. Hence we cannot already be in situation $\delta$). From this it follows that
situation $\gamma$) is attained for any radius between $r_0$ and $\sqrt2$. Thus
the intersection provides the manifold as claimed for any radius between the
two bounds given in the claim. 
\end{proof}

The question whether a hole has covering radius less or equal to $\sqrt2$ is
significant because of
\begin{theorem} The Dynkin diagram of a generalized hole is of finite
type if and only if the hole has generalized radius less than $\sqrt2$. The
Dynkin diagram is of affine type if and only if the hole has generalized radius
equal to $\sqrt2$.
\end{theorem}
\begin{proof} The argument is a straightforward generalization of the arguments
used in chapter 23 of \cite{CS88}. It is well known (see e.g. \cite{Kac90},
chapter 4) that an indecomposable Cartan matrix is of finite type if and only
if all its principal minors are positive definite. Also, an indecomposable
Cartan matrix is of affine type if and only if it has determinant 0 and all its
proper principal minors are positive definite. We can rephrase this as follows:
For any Cartan matrix of finite or affine type, and only for these, we can
choose a co-root system corresponding to the Cartan matrix, in Euclidean space.
Then the Cartan
matrix is of finite type if and only if the co-roots are linearly independent.
Then the origin will not lie within the affine span of the simple roots. The
origin is, however a point at (generalized) radius $\sqrt2$ from all simple
roots. Thus the (generalized) centre of the affine span has radius less than
$\sqrt2$. This proves the claim for the finite case.
The argument for the affine case is similar, using the existence of the vector
$\delta^\vee$ of equation (5.9). 
\end{proof}

We observe, as a corollary of theorem 5.4, that the hole diagrams only contain
bonds as described in the first part of section 5.2. In particular, for
$N\geq5$, the elements of ${\cal R}_{fix}$ and those of
${\cal R}_{dual}$ will form disjoint components of the diagram of a hole
because, if $\lambda\in\Lambda^\sigma$ and $\mu\over N$ in the dual are both
vertices of a generalized hole, then 
\[
\bigl(\lambda-{\mu\over N} \bigr)^2 \leq \bigl( \sqrt2+\sqrt{2\over N}\bigr)^2 
= 2 + {2\over N} + {4\over\sqrt N} < 4+{2\over N}.
\] 

We now use lemma 5.1 to prove that the set of finite and affine subalgebras
identified by theorem 5.4 is in fact the complete classification of all finite
and affine subalgebras contained in the GKM ${\cal G}_N$. Thus the task of
identifying subalgebras will be reduced to a classification of generalized
holes. 

\begin{theorem} Suppose there is given a subset $H$ of $\cal R$ such
that the Dynkin diagram $\Delta(H)$ is of either finite or affine type. Then
there exists a generalized hole $K$ in $\cal R$ such that $H\subset K$.
\end{theorem}
\begin{proof} Given a set $H$ of finite, or affine, type,
then by lemma 5.1 we identify its (generalized) centre $c$ and radius $r_0$.
We prove the theorem by induction on the dimension of the affine span of
$H$. Suppose the radius $r_0$ of $H$ equals $\sqrt2$. Then there are no points
in the set ${\cal R}$ closer to $c$ by lemma 5.1b. Let $H'$ be the
collection of all points at (generalized) distance $\sqrt2$ from $c$. Then the
affine span of $H'$ must be the total space because $\sqrt2$ is the generalized
covering radius of the set ${\cal R}$. Thus $c$ is the centre of a hole and
we have detected the diagram $\Delta(H)$ as part of its diagram. \flexskip

Now suppose $H$ has radius $r_0<\sqrt2$. If the dimension of the affine span of
$H$ is less than $2M$ we consider the manifold ${\cal M}(r)$ of points at
(generalized) distance $r\geq r_0$ from each point of $H$. Again, by lemma 5.1b
there are no elements of ${\cal R}$ closer to $c$ than distance $r_0$. Hence
there will be a minimal $r_1\geq r_0$ such that for some $c_1\in{\cal M}(r_1)$
there is a further point of ${\cal R}$ at that same distance from $c_1$. (That
such an $r_1$ exists follows from the fact that there exists a generalized
covering radius to the lattice. Hence $r_1$ is bounded from above by $\sqrt2$.)
We then collect all points at that generalized distance $r_1$ from $c_1$ to
form the set $H_1$. Now either $r_1=\sqrt2$. In this case we are done (see
above). Or, $r_1<\sqrt2$. This implies that the extended set of simple roots
still is of finite type. If the dimension of the affine span of $H_1$ equals
$2M$, we are done. Else, the assumptions of lemma 5.1 are satisfied. Thus the
theorem follows by induction.
\end{proof}

We conclude this section with the observation that any affine component
determines the generalized centre of a hole and hence the complete hole.
On the other hand, finite components can be part of several holes.

\section{The Volume Formula}

This section establishes that the generalized holes partition the vector spaces
spanned by the real simple roots $\cal R$ of the GKM ${\cal G}_N$
and calculates the volumes of all finite and affine holes which will occur in
the decomposition. In chapter 6 we will then carry out the explicit
decomposition of space for all relevant $N$. We will have a check of the
decomposition because the sum of the volumes of
the holes found in a fundamental domain must equal the total volume of this
fundamental domain. \flexskip

We begin by establishing that the holes as described in sections 5.1 and 5.2
actually partition space in the following sense.\flexskip

\begin{lemma} We consider generalized holes as defined in sections
5.1 and 5.2. If two such holes intersect they do so in a common face, the convex
hull of at most $2M$ points, where $2M$ is the dimension of space.
Every point of space lies either within such a face or within the interior of a
unique (generalized) hole.
\end{lemma}
\begin{proof} The first claim is straightforward. For the second, it remains to
show the following: Suppose we are given a generalized hole $H$ and a
$(2M-1)$-dimensional face $F$
of it. Claim: If we pass through the interior of $F$ we will enter the
interior of another generalized hole. \par

By lemma 5.1 we know that $F$ has a generalized centre $c$, generalized radius
$r_0$, and that there exists a
1-dimensional manifold ${\cal M}$ of points which have equal generalized
distance to all points spanning $F$. For the points $x\in{\cal M}$, let
$r_F(x)$ be the generalized distance $x$ to the points spanning $F$, and let
$r_{{\cal R}-F}(x)$ be the minimal generalized distance of $x$ to any point of
${\cal R}-F$. Here we always use the minimal distance function as defined in
formula (5.12) of section 5.2. Then, by lemma 5.1, the set $${\cal M}(r)=\{x\in
{\cal M}\mskip 10mu | \mskip 10mu r_F(x)=r\}$$ is non empty for $r_0\le
r\le\sqrt2$. The manifold ${\cal M}$ is symmetric with respect to
the hyperplane spanned by $F$. Hence the set ${\cal M'}(r) = {\cal M}(r) -
\langle H\rangle$ is non-empty for $r_0< r\le \sqrt2$. From lemma 5.1b we know
that $r_F(x)<r_{{\cal R}-F}(x)$ for $x\in {\cal M}(r_0)$. That means that no
point of ${\cal R}$ is closer or equally close to the generalized centre of
$F$ than those spanning $F$. On the other hand $r_{\cal R}(x)\le \sqrt2$ for
all $x$. Hence there certainly will be points $p\in {\cal R}$, $p \not\in F$
closer or equally close to any $x\in {\cal M}(\sqrt2)$. Hence there will be a
minimal $r_1$, $r_0<r_1\le\sqrt2$ such that there exists a point $x\in
{\cal M}'(r_1)$ with $r_{{\cal R}-F}(x) = r_F(x)= r_1$. That implies that $x$
is the centre of a generalized hole of generalized radius $r_1$ and concludes
the proof.
\end{proof}

We consider the fixed point lattice $\Lambda^\sigma$ among the roots
${\cal R}$. We define a fundamental region of $\Lambda^\sigma$ as in, e.g.,
\cite{CS88}. Let $\hbox{Aut}(\Lambda^\sigma)$ denote the group of automorphisms
of the fixed point lattice fixing the origin.
Let $\hbox{Aut}(H)$ denote the group of automorphisms fixing a generalized hole
$H$. Note that, for every $H$, $\hbox{Aut}(H)$ can be identified with a subgroup
of $\hbox{Aut}(\Lambda^\sigma)$ using a unique translation. \flexskip

As a corollary to lemma 5.2 we obtain a volume formula for the fundamental
volume $V$.

\begin{theorem} 
\[
  V=\sum_{\hbox{\rm holes } H} \hbox{\rm vol}(H) =
   |\hbox{\rm Aut}(\Lambda^\sigma)| \sum_{\hbox{\rm orbits of holes}}
  {\hbox{\rm vol}(H)\over{|\hbox{\rm Aut}(H)|}}
\]
\end{theorem}
\begin{proof} Consider all the holes $H$ within a fundamental region which are
equivalent under $\hbox{Aut}(\Lambda^\sigma)$. Their number is the order of
$\hbox{Aut}(\Lambda^\sigma)$ divided by the number of automorphisms that map $H$
to itself. Lemma 5.2 shows that the total volume $V$ will be accounted for.
\end{proof}
 
The theory which we developed so far enables us to derive a result
concerning the covering radii of the fixed point lattices which provides an
interesting insight into the relations between the fixed point lattices
$\Lambda^\sigma$ and the Leech lattice,
even though we will not use it in the remainder of this work.

\begin{theorem} The covering
radius of $\Lambda^\sigma$ is $\sqrt{2+{2\over N}}$.
\end{theorem}
\begin{proof} We begin with the cases $N=2$ and $N=3$. Here it is well known
that the covering radius is as claimed, see, e.g., \cite{CS88}, chapter 4. \flexskip

We can now restrict ourselves to the remaining cases $N=5,7,11,23$. It is
straightforward to see that the covering radius of the fixed point lattice must
be greater or equal to $\sqrt{2+{2\over N}}$ because the elements of
${\cal R}_{dual}$ have distance greater or equal $\sqrt{2+{2\over N}}$ from
every element of $\Lambda^\sigma$. (See, e.g. equations (5.3).)
It remains to show that space is covered by spheres of radius
$\sqrt{2+{2\over N}}$ centered at the elements of $\Lambda^\sigma$.
We have established in theorem 5.6 that space is decomposed into convex holes
which correspond to finite or affine subalgebras of ${\cal G}_N$. Hence it
suffices to prove that all holes, which occur in the decomposition, are
covered by the above spheres. \flexskip

Let $H$ be the set of vertices of a such a hole, of finite or affine type.
Let $H_{f}=H \cap {\cal R}_{fix}$, and $H_{d}=H \cap 
{\cal R}_{dual}$. As $\Delta(H)$ is of finite or affine type we can 
represent the elements of $H_{f}$ by vectors $v_i$ of norm 2, and the 
elements of $H_{d}$ by vectors $w_j$ of norm $2\over N$, in Euclidean 
space. Now in the cases $N=5,7,11,23$ there exist no finite or affine Lie 
algebras with bonds between the long and the short roots. This implies that 
all the $v_i$ are orthogonal to all the $w_j$. \flexskip

We now consider the convex hull $C$ of the points $0, v_i, w_j$. (If the
original $\Delta(H)$ was of finite type this is a simplex, one face of which
is the original hole $\langle H\rangle$. If $\Delta(H)$ was of affine type this is the
original hole $\langle H\rangle$.)
Let the span of the $v_i$ be denoted $V$. We show that $C$ is covered by spheres
$B(v_i,\sqrt{2+{2\over N}})$, centered at the $v_i$.
We begin by observing that every $w_j$ is contained in every sphere
$B(v_i,\sqrt{2+{2\over N}})$. Because $C$ is convex, no point of $C$ has a
distance from $V$ which is greater than $\sqrt{2\over N}$, which is the distance
from $V$ of the vertices $w_j$. Because of the mutual orthogonality of the
$v_i$ and the $w_j$ it follows that the projection into $V$ takes every point
of $C$ into $C\cap V$. The spheres $B(v_i,\sqrt{2+{2\over N}})$ contain
multi-dimensional cylinders as follows: a sphere of radius ${\sqrt 2}$ in the
intersection with $V$ times a sphere of radius $\sqrt{2\over N}$ in the orthogonal
dimensions. It follows that $C$ will certainly be covered by the spheres
$B(v_i,\sqrt{2+{2\over N}})$ if $C\cap V$ is covered by spheres
$B(v_i,\sqrt{2})$. This however is true because the $v_i$ represent a Lie
algebra of finite or affine type. 
\end{proof}
\begin{remark} The same argument can also be used to identify the covering
radius in the cases $N=2$ and $N=3$. The only additional requirement would be
to check the diagrams of type $b,c,f$ (finite and affine) case by case as to
whether the convex hull of its vertices is covered by spheres centered at the
roots of norm 2. Note that by the above theorem we have related the covering
radius of the Leech lattice $\Lambda$ to the radii of the fixed point lattices
$\Lambda^\sigma$. Note furthermore that we have established that the deep holes
of $\Lambda^\sigma$ are precisely the elements of ${\cal R}_{dual}$.
\end{remark}

The remainder of this section identifies the volumes of the various finite and
affine holes that form the partition of space. Some, though not all of the
formulas have been documented in \cite{CS88}, chapter 25. None of the results 
which will follow in the remainder of this section 5.3 are original,
however I do not know of any reference which states or proves them. \flexskip

We calculate the volumes of the convex hulls of the finite and affine Dynkin
diagrams obtained in the decomposition of the root-space ${\cal R}$ of the GKM
${\cal G}_N$. We consider a finite Dynkin diagram, $\Delta$. The construction
of the set ${\cal R}$ implies that $\Delta$ can be realized by co-roots
$r_i^\vee$ of lengths $2$ and ${2\over N}$ with the origin not contained in the
convex hull of the vectors $r_i^\vee$. We want to determine the
$(n-1)$-dimensional volume $v(\Delta)$ of the convex hull of the roots
$r_i^\vee$. The volume $V$ of the $n$-simplex of the points
$0,r_1^\vee,\dots,r_n^\vee$ can be calculated in two ways: First, $$V=\hbox{vol}
(0,r_1^\vee,\dots,r_n^\vee)={1\over n!} \hbox{det}(r_1^\vee,\dots, r_n^\vee)$$
Now, for any matrix $X$, $(\hbox{det}X)(\hbox{det}X^T) = \hbox{det} (XX^T)$.
If $X$ consists of columns $x_j$ then the entries of $XX^T$ will be the inner
products $(x_i,x_j)$. In the case at hand, this is the Cartan matrix of the
Lie algebra $\Delta$ symmetrized such that the diagonal elements are either
$2$ or $2\over N$. Let this matrix be denoted $DC$. Here, $D$ is a diagonal
matrix and $C$ is the Cartan matrix of $\Delta$. We conclude that
$V={1\over n!} \sqrt{\hbox{det} (DC)}$.  \flexskip

Second, $V$ can be calculated as the product of the $(n-1)$-dimensional 
volume $v(\Delta)$ and the height of the origin above the affine hyperplane 
spanned by the $r_i^\vee$, divided by the dimension $n$. The height vector 
$h$ stands orthogonal on the hyperplane, that is $(h,h-r_i^\vee)=0$ for all 
$i$. It follows that the direction of $h$ is characterised by the fact that 
$(h,r_i^\vee)$ is constant for all $i$. We recall that there exists a unique 
such direction, given by the Weyl vector $\rho^\vee$ of equation (5.7). 
Hence $h$ is a multiple of $\rho^\vee$. We further know that $h$ is part of 
the affine hyperplane of the $r_i$. Using equation (5.8) we conclude that 
$h={\rho^\vee \over \rho^2}$. It follows that
$V={1\over n} |h| v(\Delta)= {1\over n} {v(\Delta)\over\sqrt{\rho^2}}$. 
Hence 
\begin{equation}
 v(\Delta)= n \sqrt{\rho^2}\ V= 
{1\over (n-1)!}\sqrt{\rho^2}  \sqrt{\hbox{det} (DC)}.
\end{equation}
From this we can furthermore derive the following rules for the volumes of
diagrams made up from disjoint components $\Delta_1$ and $\Delta_2$. Disjoint
components imply that the affine spans are orthogonal. Hence the norms of the
Weyl vectors are additive and the determinants are multiplicative. \flexskip

We now consider an indecomposable affine diagram, $\Delta_n$. Again, it is
represented by co-roots $r_i^\vee$ of norm 2 and ${2\over N}$.
We consider the vector $\delta^\vee=\sum_i n_i^\vee r_i^\vee$ as constructed in
equation (5.9). We label the co-roots $r_i^\vee$, $i=0,\dots,n$ such that
$n_0^\vee=1$. (See \cite{Kac90}, chapter 4
for the fact that this is possible for all affine subalgebras.) \flexskip

The convex hull of $\Delta$ is a simplex as there are $n+1$ points spanning an
$n$-dimensional space. We understand the $r_i^\vee$ to be vectors written as
columns. Then the volume can be calculated as
\begin{equation} \begin{split}
n! v(\Delta)
 &=\hbox{det}(r_1^\vee-r_0^\vee,r_2^\vee-r_0^\vee,\dots,r_n^\vee-r_0^\vee) \\
 &=\hbox{det}
   \begin{pmatrix}
    1 &                 0 &                 0 & \dots &                0 \\
    0 & r_1^\vee-r_0^\vee & r_2^\vee-r_0^\vee & \dots & r_n^\vee-r_0^\vee
    \end{pmatrix} \\
 &=\hbox{det}
   \begin{pmatrix}
   1       &               0 &               0 &\dots&              0 \\
   r_0^\vee&r_1^\vee-r_0^\vee&r_2^\vee-r_0^\vee&\dots& r_n^\vee-r_0^\vee
   \end{pmatrix} \\
 &=\hbox{det}
   \begin{pmatrix}
   1        &        1 &        1 & \dots &        1 \\
   r_0^\vee & r_1^\vee & r_2^\vee & \dots & r_n^\vee
   \end{pmatrix}
  =\hbox{det}
   \begin{pmatrix}
   h^\vee      &        1 &        1 & \dots &       1 \\
   \delta^\vee & r_1^\vee & r_2^\vee & \dots & r_n^\vee
   \end{pmatrix} \\
 &=h^\vee \hbox{det} (r_1^\vee, r_2^\vee, \dots, r_n^\vee)
   =h^\vee \sqrt{\hbox{det} (DC)}
\end{split}\end{equation}

Here, the last step is precisely as in the case of finite diagrams, because
the co-roots $r_1^\vee, \dots, r_n^\vee$ correspond to a finite diagram.
Note that the choice of the co-root labelled 0 is not always unique, the only
condition we used was that $n_0^\vee=1$ in the construction of $\delta^\vee$.
Still, the result is the same for any such co-root.
For disconnected Dynkin diagrams we observe that orthogonality implies that the
volumes are multiplicative: $(n_1+n_2)! v(\Delta_1\Delta_2)=
n_1! v(\Delta_1) n_2! v(\Delta_2)$. \flexskip

We note that the simply laced diagrams do appear within the decomposition in
two varieties. Let $\Delta_n$ denote an indecomposable simply laced diagram of
rank $n$ such that all roots have norm 2. Here, $\Delta_n$ may be either finite
or affine. Then $N\Delta_n$ shall denote the diagram of type $\Delta_n$ such
that all roots have norm $2N$. For the symmetrized Cartan matrices of the
co-roots we obtain:
$DC(N\Delta_n)= {1\over N}DC(\Delta_n)$. This proves the following
formulas for the determinant and the norm of the Weyl vector:
\begin{equation}
\hbox{det}(N\Delta_n)={1\over N^n}\hbox{det}(\Delta_n)
\end{equation}
\begin{equation}
\rho^2(N\Delta_n)=N\rho^2(\Delta_n).
\end{equation} 

This concludes the results needed to calculate the volumes of the generalized
holes. The actual values of determinant and norm of the Weyl vector are well
documented in the literature. For convenience, we reproduce in table 5.1 the
values for those finite and affine Lie algebras that will occur in chapter 6.

We conclude this section with a remark which will provide a further useful check
for forthcoming calculations. If all the co-ordinates of the vertices of a hole
are rational so will be the volume. Thus the terms under the square root will
have to prove to be squares. Similarly, if $M$ co-ordinates are
rational, and $M$ co-ordinates are elements of $\sqrt N\QQ$ then the volume will
be element of $\sqrt N^M\QQ$.

\begin{center}
\begin{tabular}[t]{|c|c|c||c|c|c|}
\multicolumn{1}{|c|} {diagram} &
\multicolumn{1}{|c|} {$\rho^2$} &
\multicolumn{1}{|c||} {det} &
\multicolumn{1}{|c|} {diagram} &
\multicolumn{1}{|c|} {$h^\vee$} &
\multicolumn{1}{|c|} {det} \\ 
\hline
  & & & & & \\
\vspace{-0.1in}
 $a_n$ & ${1\over12}n(n+1)(n+2)$ & $n+1$ & $A_n$ & $n+1$ & $n+1$ \\
  & & & & & \\
\vspace{-0.1in}
 $b_n$ & ${1\over6}n(2n-1)(2n+1)$ & $2^{2-n}$ & $B_n$ & $2n-1$ & $2^{2-n}$ \\
  & & & & & \\
\vspace{-0.1in}
 $c_n$ & ${1\over6}n(n+1)(2n+1)$ & 1 & $C_n$ & $n+1$ & 1 \\
  & & & & & \\
\vspace{-0.1in}
 $d_n$ & ${1\over6}(n-1)n(2n-1)$ & 4 & $D_n$ & $2n-2$ & 4 \\
  & & & & & \\
\vspace{-0.1in}
 $e_6$ & 78 & 3 & $E_6$ & 12 & 3 \\
  & & & & & \\
\vspace{-0.1in}
 $e_7$ & ${399\over2}$ & 2 & $E_7$ & 18 & 2 \\
  & & & & & \\
\vspace{-0.1in}
 $e_8$ & 620 & 1 & $E_8$ & 30 & 1 \\
  & & & & & \\
\vspace{-0.1in}
 $f_4$ & 78 & ${1\over4}$ & $F_4$ & 9 & ${1\over4}$ \\
  & & & & & \\
\vspace{-0.1in}
 $g_2$ & 14 & ${1\over3}$ & $G_2$ & 4 & ${1\over3}$ \\
  & & & & & \\
\vspace{-0.1in}
       & & & $A^{(2)}_{2n-1}$ & $2n$ & 1 \\
  & & & & & \\
\vspace{-0.1in}
       & & & $D^{(2)}_{n+1}$ & $2n$ & $2^{2-n}$ \\
  & & & & & \\
\vspace{-0.1in}
       & & & $E^{(2)}_6$ & 12 & ${1\over4}$ \\
  & & & & & \\
\vspace{-0.1in}
       & & & $D^{(3)}_4$ & 6 & ${1\over3}$ \\
  & & & & & \\
\hline
\end{tabular}
\end{center}
\vspace{-0.1in}
 
\begin{table}[ht]
\caption{Parameters of finite and affine Lie algebras}
\end{table}

\section{The Automorphism Groups}

The aim of chapter 6 will be to provide a complete account of the types of
generalized holes contained in the root system ${\cal R}$. Two holes will only
be accepted as equivalent if there exists an automorphism of the fixed point
lattice taking one to the other. Now it turns out that in higher dimensions the
Dynkin diagram formed by the vertices of a hole does not necessarily determine
the equivalence class of the hole uniquely. There are examples of such holes
in the Leech lattice (see \cite{CS88}, chapter 25). As we will see in chapter 6,
the set ${\cal R}$ for $N=3$ does provide further examples. Hence, there is no
hope to prove a general result which lifts holes in the set ${\cal R}$ to holes
in $\Lambda$ accounting for equivalence. Therefore,
this section will describe a number of results and techniques which
will enable us to assert the equivalence of holes on a case by case basis in
chapter 6. \flexskip

We will need the following description of the faces of generalized holes: \par
\begin{lemma} Let $H$ be a subset of ${\cal R}$ such that $\Delta(H)$ is
of finite type. Then $\langle H\rangle$  has $2M+1$ faces, obtained by removing
any one of its vertices. Let $H$ be a subset of ${\cal R}$ such that
$\Delta(H)$ is of affine type. Any face of $\langle H\rangle$ is the
convex hull of a subset of its vertices, chosen to the following rule: Remove
one vertex (simple root) of every component of the Dynkin diagram.
\end{lemma}
\begin{remark} Note that the vertices of a face of a hole of affine type form a
Dynkin diagram of finite type.
\end{remark}

\begin{proof} For the finite type, we remember that the generalized hole is a
simplex. For the affine type we observe that we must delete at least one vertex
of every component because otherwise the centre of the whole generalized hole
would be contained in the span. On the other hand, that leaves at most $2M$
vertices. They will all be needed to form the face in a $2M$-dimensional space.
\end{proof}

The purpose of this section is to identify the relations between holes in the
Leech lattice (as listed in \cite{CS88}, chapter 25) and generalized holes in
the sets ${\cal R}$. We consider the set of vertices $H$ of a generalized hole
$\langle H\rangle$ in the decomposition of ${\cal R}$. As before, write $H$ as
the disjoint
union $H_{f}\cup H_{d}$ with $H_{f}\subset {\cal R}_{fix}$ and
$H_{d}\subset {\cal R}_{dual}$. Now consider $\lambda^*\in H_{d}$.
Theorem 3.1 shows that the elements of the dual lattice are images under the
projection $\pi_{\Lambda^\sigma}$ onto the space spanned by the fixed point
lattice. We consider the set of preimages $\pi_{\Lambda^\sigma}^-(\lambda^*)$.
If $\lambda \in \Lambda$ such that $\pi_{\Lambda^\sigma}\lambda=\lambda^*$ then
$(\lambda-\pi_{\Lambda^\sigma}\lambda)^2\ge 2-{2\over N}$ because
$(\lambda^*)^2\equiv {2\over N} \ \hbox{mod}\ 2\ZZ$.
We consider a dual vector of type $\lambda^*=
{1\over\sqrt8}(4,0^{(M-1)},{4\over N}^{(N)},0^{(N)},...,0^{(N)})$ in the
notation of chapter 4.1. It has norm $2+{2\over N}$. Any $\lambda\in\Lambda$
such that 
\begin{equation}
\pi_{\Lambda^\sigma}(\lambda)= \lambda^* \ \hbox{and\ }
(\lambda-\pi_{\Lambda^\sigma}\lambda)^2=2-{2\over N}
\end{equation}
thus has norm 4. Obviously there are precisely $N$ such, namely
\begin{equation}
\lambda_n=(4,0^{(M-1)},0^{n-1},4,0^{(N-n)},0^{(N)}, \dots, 0^{(N)}),
\end{equation}
for $n=1,\dots,N$.
For a general dual vector $\lambda^*$ we recall from chapter 4.4 (see proof of
theorem 4.2) that any two vectors of
norm ${2\over N}$ modulo $2\ZZ$ are equivalent modulo $\Lambda^\sigma$.
Hence the preimages at distance $2-{2\over N}$ of any element of
${\cal R}_{dual}$ form a cycle of $N$ points, and
$\sum_{n=0}^{N-1} {1\over N}\sigma^n(\lambda)=\lambda^*$. Note that this does
not imply that all roots of norm $2+{2\over N}$ are equivalent under the
automorphism group of $\Lambda^\sigma$. This is, in fact, not true in the cases
$N=11,7,5$. \flexskip

For the generalized hole $\langle H\rangle$ with set of vertices $H$ and
generalized
centre $c$ we now consider the following set:
\begin{equation}
H_\Lambda=\{\lambda\in\Lambda \ | \ \lambda\in H \ \hbox{or} \
\lambda \ \hbox{satisfies} \ (5.18)\ \hbox{for}\ \lambda^*\in H \ \}
\end{equation}
From the construction of $c$ in formula (5.13) it is immediate that
$c$ has distance less or equal to 2 from all points of $H_\Lambda$.
Theorem 5.4 holds for the Leech lattice (see
\cite{CS88}, chapter 25). Hence the points of $H_\Lambda$ are part
of either a finite, or an affine, diagram of the Leech lattice. \flexskip

Suppose that $\Delta(H)$ was of finite type in $\Lambda^\sigma$ such that
$|H_{f}|=n$ and (consequently) $|H_{d}|=2M+1-n$.
First, assume that $n\leq M$. The set $H_\Lambda$ then
will contain $n+N\times(2M+1-n)=MN+n+MN-nN+N=
24-M+n + (M-n)N+N=24+(N-1)(M-n)+N\ge 24+N$ elements.
Hence we have found a hole of $\Lambda$ of more than 25 vertices and a radius
strictly smaller than $\sqrt2$. This contradicts the fact that the roots of
finite diagrams are affinely independent. Next, assume that $n=M+1$. Then
$H_\Lambda$ contains $M+1+NM=M(N+1)+1=25$ vertices. Thus this describes a
complete finite hole $H_\Lambda$ of the Leech lattice which thus is unique.
Further, it follows that $\sigma$ must preserve the hole ${H}_\Lambda$,
or equivalently that $\sigma\in\hbox{Aut}(H_\Lambda)$. \par

Now suppose that $\Delta(H)$ is of finite type, and that $n>M+1$. Then
$H_\Lambda$ has less than 25 elements and hence forms part of faces of
a number of finite and affine holes. There is no uniqueness, and $\sigma$ acts
by permuting cycles of these holes which share the fixed vertices. \par

Finally, suppose that $\Delta(H)$ was of affine type in $\Lambda^\sigma$. Then
the distance of the generalized centre $c$ was exactly $\sqrt2$ from all points
of $H_\Lambda$. It then follows from the analogue of theorem 5.4 that
there exists an affine hole $H_\Lambda$ of the Leech lattice with the
same centre $c$. This hole is unique in so far as $c$ lies within its interior.
This implies that $\sigma$ must preserve the hole $H_\Lambda$,
or equivalently that $\sigma\in\hbox{Aut}(H_\Lambda)$. This, in turn,
implies that all vertices in $H_\Lambda$ will be preimages of
vertices of $H$. This proves
\begin{theorem} The generalized holes of ${\cal R}$ are of the following
three types: \par

(I) The interior of each affine hole of ${\cal R}$ lies within the
interior of a unique affine hole ${H}_\Lambda$ of the Leech lattice
such that $\sigma\in\hbox{Aut}({H}_\Lambda)$. \par

(II) The interior of each finite hole with $M+1$ elements of
${\cal R}_{fix}$ among its vertices
lies within the interior of a unique finite hole ${H}_\Lambda$ of the Leech
lattice such that $\sigma\in\hbox{Aut}({H}_\Lambda)$. \par

(III) Finite holes with more than $M+1$ elements of ${\cal R}_{fix}$
among their vertices lie
within the faces of some affine and finite holes of the Leech lattice.
\end{theorem} \theoremproven

\begin{remark} The faces of affine holes correspond to finite diagrams, hence there
are no `affine faces'. That there exist holes of each of the three types will
become apparent in the explicit decomposition, as provided in appendix A. 
\end{remark}

It is also worth noting that the preimages of an element of ${\cal R}_{dual}$
form a diagram of type $a_1^N$, as can be seen in the explicit case
considered above (equation (5.19)). Now consider a hole $H_\Lambda$ of the
Leech lattice such that $\sigma\in\hbox{Aut} (H_\Lambda)$. The action of
$\sigma$ is to permute a number of vectors which belong to the same cycle, and
to fix the rest. Thus $\sigma$ induces a map taking $\Delta(H_\Lambda)$ to
the Dynkin diagram of some generalized hole of ${\cal R}$. We want to describe
this action directly as an action on the diagrams. It is obvious that
$\sigma$ identifies some sets of nodes, and fixes the rest. We need to
understand the action on the bonds. First, consider the case when
$\lambda\in\Lambda$ is fixed, and $\mu_i$, $i=1,\dots,N$ are cyclically
permuted by $\sigma$. It follows that all $(\lambda-\mu_i)^2$, $i=1,\dots,N$
must be equal. From the construction it follows immediately that the vectors
$\lambda$, $\mu_j$, and $\pi(\mu_j)={1\over N}\sum_i\mu_i$ form an rectangular
triangle. Hence, $(\lambda-
{1\over N}\sum_i\mu_i)^2=(\lambda-\mu_j)^2-(2-{2\over N})$. Thus, if $\lambda$
and $\mu_j$ are unlinked in $\Lambda$, so are $\lambda$ and ${1\over N} \sum_i
\mu_i$ in ${\cal R}$. And, if $\lambda$ and $\mu_j$ are linked in $\Lambda$ by
a single bond, then $\lambda$ and ${1\over N}\sum_i\mu_i$ will be linked by an
arrow with $N$ bonds. Now consider the case of two
cycles. Suppose, $\lambda_i$, $i=1,\dots,N$, and $\mu_i$, $i=1,\dots,N$ are
such. Name the corresponding simple roots of the monster Lie algebra $r_i=
(\lambda_i,1,{\lambda_i^2\over2}-1)$ and $s_i=(\mu_i,1,{\mu_i^2\over2}-1)$.
As observed above, each cycle forms an $a_1^N$ diagram. Hence
$(\lambda_i-\lambda_j)^2=4$ for all $i,j=1,\dots,N$. Hence,
$$(\lambda_i,\lambda_j)={\lambda_i^2+\lambda_j^2-(\lambda_i-\lambda_j)^2\over2}
={\lambda_i^2+\lambda_j^2\over2}-2$$
$$(\sum_i\lambda_i)^2=\sum_{i,j}(\lambda_i,\lambda_j)=
\sum_i\sum_{j\not=i}({\lambda_i^2+\lambda_j^2\over2}-2)+\sum_i\lambda_i^2=
N\sum_i\lambda_i^2-2N(N-1).$$ Hence, we obtain
$${(\sum_i\lambda_i)^2\over2N}-1=\sum_i{\lambda_i^2\over2}-N.$$
We recall that the root $r$ of norm $2N$ had the form
$r=(\sum_i\lambda_i,N,{(\sum_i\lambda_i)^2\over2N}-1)$. Thus we can write
\[
r=(\sum_i\lambda_i,N,\sum_i{\lambda_i^2\over2}-N)=
\sum_i(\lambda_i,1,{\lambda_i^2\over2}-1)=\sum_ir_i
\]
Analogously we decompose
\[
s=(\sum\mu_i,N,{(\sum\mu_i)^2\over2N}-1)=\sum(\mu_i,1,{\mu_i^2\over2}-1)=
\sum s_i.
\]
Now, $(r,s)=\sum_{i,j}(r_i,s_j)$. We will only describe the most important
cases:
a) the $r_i,s_j$ form a diagram of type $a_1^{2N}$. Then all products $r_i,s_j$
are zero and hence $(r,s)=0$, leaving $r$ and $s$ unlinked, forming an
${\mathbf N}a_1^2$ (here the bold ${\mathbf N}$ is introduced to indicate long roots).
b) the $r_i,s_j$ form a diagram of the following type: for every $i$ there is
a unique $j$ such that $r_i$ and $s_j$ are joined by a single bond.
(The resulting diagram is of type $a_2^N$.) Adding up the inner products
$r_i,s_j$ shows that $r$ and $s$ will be joined by
a single bond, thus forming ${\mathbf N}a_2$. It is now obvious how to identify
the action of $\sigma$ for some other constellations that may be needed.
This completes the description of the induced action of $\sigma$ on
Dynkin diagrams. \flexskip

The above results describe the relation of holes in the Leech lattice to those
in ${\cal R}$. We will now proceed to develop some techniques which
will help to identify the automorphism groups of generalized holes in ${\cal
R}$. Throughout chapter 6, we will work with the normalizer of
$\sigma$ in $Co_0$ as the group of known automorphisms of the fixed point
lattices. The normalizer (as in chapter 4.1, above) is the subgroup of $Co_0$
of automorphisms which normalize the cyclic group $\langle\sigma\rangle$.
The normalizer of $\sigma$ fixes $\Lambda^\sigma$. There may, or may not, exist
further automorphisms of $\Lambda^\sigma$. The full automorphism group
$\hbox{Aut}(\Lambda^\sigma)$ will be identified case by case in the course of
the calculations of chapter 6.
\begin{theorem} Suppose $\langle H\rangle$ is a generalized hole of
${\cal R}$ of type (I) or (II), as defined in theorem 5.8. Let $\langle
H_\Lambda\rangle$ denote the unique hole of the Leech lattice associated with
$\langle H\rangle$ in theorem 5.8. Suppose further that in the Leech lattice
all holes $\langle K_\Lambda\rangle$ with $\Delta(K_\Lambda)=\Delta(H_\Lambda)$
are equivalent to $\langle H_\Lambda\rangle$ under the automorphism group
$Co_0$ of $\Lambda$. Suppose that the order of $\hbox{Aut} (H_\Lambda)$ is
divisible by $N$ but not by $N^2$. Then all holes $\langle K\rangle$ in
${\cal R}$ with $\Delta(K)=\Delta(H)$ are equivalent to $\langle H\rangle$
under the automorphism group of $\Lambda^\sigma$.
\end{theorem}
\begin{proof} Suppose $\langle K\rangle$ is another hole of ${\cal R}$ with
$\Delta(K)=\Delta(H)$. Then there exists a unique corresponding
hole $\langle K_\Lambda\rangle$ of the Leech lattice. $\Delta(K_\Lambda)$
must be equal to $\Delta(H_\Lambda)$ because, in cases (I) and (II), all
vertices of the Leech lattice holes are preimages of those
vertices in ${\cal R}$. (The considerations above showed how to fully
reconstruct the Leech lattice hole.) Hence, there is a $\phi\in Co_0$ such that
$K_\Lambda=\phi(H_\Lambda)$. $\sigma$ acts on $K_\Lambda$.
Consider $\phi^{-1}\sigma\phi$ and $\sigma$. Both are automorphisms of order
$N$, acting on ${H}_\Lambda$.
The cyclic groups generated by $\phi^{-1}\sigma\phi$ and
$\sigma$ must be conjugate in $\hbox{Aut}({H}_\Lambda)$ because, by assumption,
the order of this group is divisible by $N$, and not divisible by $N^2$
(see Sylow's theorem, e.g., \cite{Jac85}, p.80), say by some $\psi\in\hbox{Aut}
(H_\Lambda)$. Then $\sigma^n=\psi^{-1}\phi^{-1}\sigma\phi\psi$ for some
integer $n$, which implies that $\phi\psi$ is an element of the normalizer of
$\sigma$. This shows that $\phi\psi$ is element of
$\hbox{Aut}(\Lambda^\sigma)$. 
\end{proof}

The above theorem will cover many of the occurring diagrams. However, the rest
will have to be treated on a case by case basis, using the following technique:
\begin{theorem} Suppose we consider two types of diagram,
$\Delta_1$, and $\Delta_2$. They may be of any of the three types (I), (II), or
(III) specified in theorem 5.8. Suppose there are $n_1$ generalized holes of
type $\Delta_1$, and $n_2$ generalized holes of type $\Delta_2$, within a
fundamental region of the fixed point lattice. Suppose it is known that all
holes of type $\Delta_1$ are equivalent under $\hbox{Aut}(\Lambda^\sigma)$.
Suppose a hole $\langle H_1\rangle$ such that $\Delta(H_1)=\Delta_1$ has $f_1$
faces ${F}$ which border on holes $\langle H_2\rangle$ of type $\Delta(H_2)=
\Delta_2$ such that there exist diagram automorphisms in $\hbox{Aut}(H_1)$
acting transitively on the $f_1$ faces ${F}$. Suppose the neighbouring holes
$\langle H_2\rangle$ of type $\Delta_2$ have $f_2$ faces of type ${F}$
bordering on holes of type $\Delta_1$. Suppose that $n_1f_1=n_2f_2$. Then all
diagrams of type $\Delta_2$ are equivalent under $\hbox{Aut}(\Lambda^\sigma)$.
\end{theorem}
\begin{proof} We use the automorphisms identifying the holes of type $\Delta_1$,
and those of $\hbox{Aut}(H_1)$ to identify all holes of type $\Delta_2$
which border on holes of type $\Delta_1$. The assumption $n_1f_1=n_2f_2$
then implies that we have accounted for all diagrams of type
$\Delta_2$. 
\end{proof}

We conclude the chapter with a number of remarks. The check of the
assumptions of theorem 5.10 will be simplified if $f_1=1$.
We know from lemma 5.3 about the faces of holes. We understand the action of
$\sigma$ on Dynkin diagrams. Hence we can use theorem 5.9 to
provide plenty of starting points so that some careful planning of the correct
route through the remaining types of generalized holes will almost always
succeed in selecting such neighbours.  Finally, in chapter 6 we will first
develop a strategy to count the numbers of holes of any type within a
fundamental region of the fixed point lattice. The application of theorem 5.10
will provide a double check of the consistency of these numbers. This is all
the more valuable as the complexity of the search means that it will only be
possible to present the results, not however all individual calculations.

\chapter{Hyperbolic Lie Algebras}
  
A Lie algebra is called hyperbolic if its Cartan matrix satisfies the
condition that every proper principal minor is of finite or affine type.
\cite{Wan91} gave a complete list of all 238 hyperbolic Lie algebras with more
than 2 simple roots. The significant difference to finite and affine Lie
algebras is that hyperbolic Lie algebras possess roots of negative norm.
Some of their multiplicities will be greater than 1. The purpose of this
chapter is to identify the hyperbolic subalgebras of the GKMs ${\cal G}_N$
and thus to find upper bounds for the root multiplicities of these
hyperbolic Lie algebras.
\flexskip
 
The project of finding hyperbolic subalgebras does not present any
theoretical problems as we simply have to check all subsets (of correct size)
of the set of those real simple roots in the GKM ${\cal G}_N$ such that the
representing vectors in ${\cal R}$ have norm less or equal to 8 (as the root
system ${\cal R}$ is translation invariant). However, in the case $N=2$ this
would require checking subsets of up to 10 roots of a set of several hundred
thousand roots. No computer is fast enough to carry out such an undirected
search. Hence we will devise a strategy for a more organized search.
\flexskip
 
As before, let $\sigma$ be an automorphism of the Leech lattice of order
$N$, and $N=$ 23, 11, 7, 5, 3, 2. We consider the GKM ${\cal G}_N$
corresponding to $\sigma$ as constructed in chapter 1. Section 6.1 recalls
the classification of \cite{Wan91} and identifies which hyperbolic Lie algebras
are candidates to be subalgebras of one or any of the GKMs ${\cal G}_N$.
Section 6.2 indicates how to carry out the explicit complete enumeration of
all finite, affine, and hyperbolic, subalgebras of the GKMs ${\cal G}_N$ for
all relevant $N$. The 2-dimensional case ($N=23$, section 6.2.1) is almost
trivial but very useful for visualizing the problem. The 4-dimensional case
($N=11$, section 6.2.2) already shows all problems of the cases in higher
dimensions. As the calculations become very tedious for the remaining four
cases we will restrict ourselves in section 6.2.3.1 to an account
of respective bases and symmetries used in the actual calculations.
Section 6.2.3.2 provides a list of all those hyperbolic Lie algebras which can
be identified as subalgebras of one of the GKMs ${\cal G}_N$  in such a way
that the resulting upper bounds are sharp at least for some roots and such that
they represent an improvement on existing bounds. Again, we can only give the
results of the trivial but tedious calculations. For each such hyperbolic Lie
algebra we give the Dynkin diagram and the root multiplicities for some
imaginary roots of small height. I used Peterson's recursion formula (\cite{Kac90},
p.210) to carry out the numerical calculations of these root multiplicities. We
contrast the multiplicities with the resulting upper bounds.
\flexskip
 
The complete classifications of finite and affine diagrams for all relevant
$N$ are provided in appendix A. 
Given that computer calculations were required to complete the proof, the
function of sections 6.2.1 and 6.2.2 is not so much to carry
out some trivial calculations but to show that the output of a computer program
actually constitutes a mathematical proof. The partition of space as listed in
appendix A is not merely auxiliary but a result in its own right as well and it
extends the results obtained by Borcherds, Conway, and Queen for the Leech
lattice (see \cite{CS88}, chapter 25).
\flexskip
 
The concluding section 6.3 tries to put into perspective the results obtained in
the previous sections of chapter 6. In section 6.2 we observed that the root
multiplicities of ${\cal G}_N$ provided good upper bounds for some hyperbolic
Lie algebras while for others they did not bear any resemblance to the
correct multiplicities. We use some well known results about the multiplicities
of particularly the norm 0 vectors of hyperbolic Lie algebras to identify some
conditions which are necessary to obtain sharp upper bounds. This is not simply
a descriptive exercise but will provide us with a strategy which could help to
determine sharp upper bounds for more hyperbolic Lie algebras. We will
describe this strategy and relate it to some explicit numerical calculations
of root multiplicities. We will further relate it
to some results and conjectures of Borcherds. We will conclude with a number
of observations which present open problems related to the results of this work.
 
\section{Wan's classification}  
 
\cite{Wan91} classified all hyperbolic Lie algebras. The notation introduced in
chapter 5.2 was based on Wan's notation. We will continue using it here.
To achieve the aims of this chapter we need to identify all hyperbolic Lie
algebras which are subalgebras of the GKMs ${\cal G}_N$ constructed in
chapter 1. Wan's classification provides the pool of our search. However,
we cannot expect to find all hyperbolic Lie algebras as such subalgebras. For
a start, all subalgebras of the GKMs ${\cal G}_N$ constructed in chapter 1
will be symmetrizable. Further, among the ${\cal G}_N$ we do not have any Lie
algebra whose real simple roots have
ratio of norms 4:1. Finally, from chapter 5.1 it follows that in any
subalgebra of the GKMs ${\cal G}_N$ there are real simple roots of at most
two different norms, as we consider prime $N$ only. Hence we may restrict our
attention to symmetrizable hyperbolic Lie algebras with simple roots whose
norms have ratio 2:1, 3:1 only. This leaves us with a list of altogether 96
hyperbolic Lie algebras which are potential candidates.
 
\section{Finite, Affine, and Hyperbolic Subalgebras}
 
We now proceed to identify the finite, affine, and hyperbolic Lie subalgebras
of the GKMs. We will present the cases $N=23$ (section 6.2.1) and $N=11$
(section 6.2.2) explicitly. The case $N=23$ can be visualized easily as it is
2-dimensional. The case $N=11$ already shows all essential features of the
general case. Section 6.2.3 summarizes the results for the remaining $N$, that
is $N=7,5,3,$ and 2.

\subsection{N=23}
 
\subsubsection{Finite and Affine Subalgebras of N=23}
 
Let ${\cal G}_{23}$  denote the GKM constructed in
chapter 1 from an automorphism $\sigma$ of cycle shape $1^1 23^1$. Let
$\Lambda^\sigma$ be the 2-dimensional fixed point lattice and
$L=\Lambda^\sigma \oplus II_{1,1}$ be the corresponding Lorentzian lattice.
Its real simple roots $r$ have been calculated in theorem 5.1 to be the norm
2 roots $\bigl(\lambda,1, {\lambda^2\over2}-1\bigr)$ where
$\lambda\in\Lambda^\sigma$ and the norm 46 roots
$\bigl(\lambda,23,{\lambda^2\over46}-1\bigr)$ where $\lambda\in
23(\Lambda^\sigma)^*$. A basis of $\Lambda^\sigma$ is easily obtained as
\begin{subequations}
\begin{equation}
  \lambda_1={1\over\sqrt8}(-3,\sqrt{23}),\mskip 50mu
  \lambda_2={1\over\sqrt8}(5,\sqrt{23}).
\end{equation}
We can therefore choose our fundamental region as the convex hull of the
vectors $0,$ $\lambda_1$, $\lambda_2$, and $\lambda_1+\lambda_2$. We recall
from chapter 4.1 (a remark following lemma 4.1) that the additive group
$(\Lambda^\sigma)^*/\Lambda^\sigma$ has just one generator, which is of order
23. Hence, as a basis of the dual lattice we may choose
\begin{equation}
{1\over23}(9\lambda_2-8\lambda_1)={1\over\sqrt8}(3,{1\over23}\sqrt{23}),
\mskip 50mu  \lambda_1-\lambda_2={1\over\sqrt8}(8,0).
\end{equation}
\end{subequations}
There are precisely 2 elements of ${\cal R}_{dual}$ within the
fundamental region:
$\lambda_1^*={1\over\sqrt8}(1,{19\over23}\sqrt{23})$ (of norm $2+{2\over23})$ and
 $\lambda_2^*={1\over\sqrt8}(1,{27\over23}\sqrt{23})$ (of norm $4+{2\over23})$.
A diagram of the root system ${\cal R}$ for the fundamental
region is shown in figure 6.1.
 
\begin{figure}
\begin{picture}(325,300)(-125,-30)
\put(0,0){\line(1,1){125}}
\put(0,0){\line(-3,5){75}}
\put(50,250){\line(-1,-1){125}}
\put(50,250){\line(3,-5){75}}
\put(0,0){\line(1,4){25}}
\put(-75,125){\line(4,1){100}}
\put(-75,125){\line(4,-1){100}}
\put(125,125){\line(-4,1){100}}
\put(125,125){\line(-4,-1){100}}
\put(50,250){\line(-1,-4){25}}
\put(0,0){\circle*{5}}
\put(15,0){0}
\put(-75,125){\circle*{5}}
\put(-100,125){$\lambda_1$}
\put(125,125){\circle*{5}}
\put(140,125){$\lambda_2$}
\put(50,250){\circle*{5}}
\put(60,250){$\lambda_1+\lambda_2$}
\put(25,100){\circle*{5}}
\put(30,110){$\lambda_1^*$}
\put(25,150){\circle*{5}}
\put(30,155){$\lambda_2^*$}
\put(25,125){\circle*{2}}
\put(8,88){\circle*{2}}
\put(45,162){\circle*{2}}
\put(43,82){\circle*{2}}
\put(8,163){\circle*{2}}
\end{picture}
\caption{The root system of the GKM ${\cal G}_{23}$}
\end{figure}
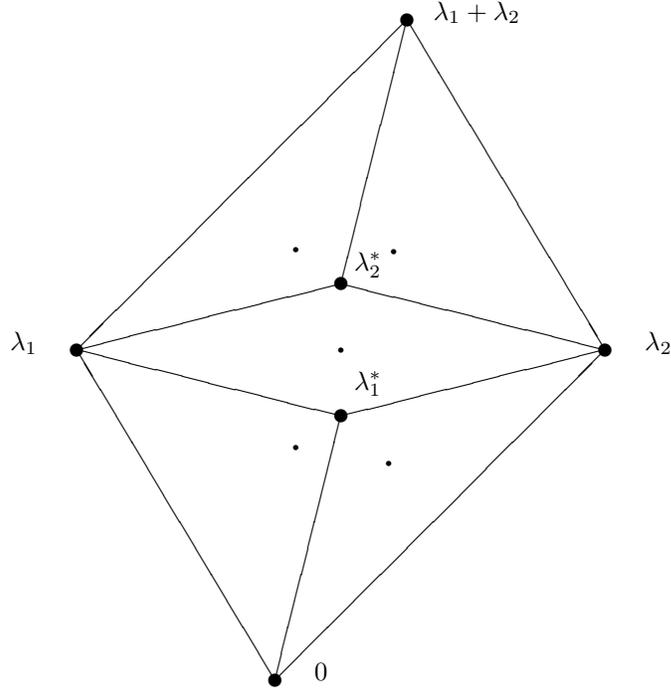
 
\begin{remark} The smaller dots indicate the positions of the generalized centres
of the generalized holes in the partition. We further observe that the balls
of radius $\sqrt2$ and centres at the fixed points ${\cal R}_{fix}$ do not
cover the space. Correspondingly, the simple roots of norm 2 do
not constitute the full set of real simple roots of the GKM ${\cal G}_{23}$.
\end{remark}

The volume of the fundamental region is the determinant of the
matrix formed by the vectors $\lambda_1$ and $\lambda_2$:
\begin{equation}
 V(23)=\hbox{det}\begin{pmatrix}
         {1\over\sqrt8}5&{1\over\sqrt8}\sqrt{23}\\
        -{1\over\sqrt8}3&{1\over\sqrt8}\sqrt{23} \end{pmatrix} = \sqrt{23}
\end{equation}

Within the fundamental region we identify 5 generalized holes. There is one
affine hole of type $A_1\ {\bf 23}A_1$, two finite ones of type $a_1^2\
{\bf 23}a_1$, and two finite ones of type $a_2\ {\bf 23}a_1$.
Using the automorphism $-\hbox{id}$ of the fixed point
lattice (and suitable translations by lattice elements) we observe that the
finite holes of equal type are in fact equivalent. It is easily recognized
that, in fact, the full automorphism group fixing 0 (see chapter 5.4) is of
order 2.
Formulae (5.14) to (5.17) of chapter 5.3 in conjunction with table 5.1 provide
the volumes of these holes:
\begin{subequations}
\begin{align}
  V(A_1\ {\bf 23}A_1) 
  =& {1\over2!} V(A_1)V({\bf 23}A_1)= {1\over2!}
     (2\sqrt2)(2{\sqrt2\over\sqrt{23}}) = {4\over\sqrt{23}} \\
\begin{split}
  V(a_1^2\ {\bf 23}a_1)
  =& {1\over2!} \sqrt{\Bigl(\rho^2(a_1)+\rho^2(a_1)+\rho^2({\bf 23}a_1)\Bigr)} \\
   &\times \sqrt{\Bigl(\hbox{det}(a_1)\hbox{det}(a_1)\hbox{det}({\bf 23}a_1)\Bigr)} \\
  =& {1\over2}\sqrt{\Bigl({1\over2}+{1\over2}+{23\over2}\Bigr)
       \Bigl(2\times2\times{2\over23}\Bigr)} ={5\over\sqrt{23}}
\end{split} \\
\begin{split}
  V(a_2\ {\bf 23}a_1)
  =& {1\over2!}\sqrt{\Bigl(\rho^2(a_2)+\rho^2({\bf 23}a_1)\Bigr)
       \Bigl(\hbox{det}(a_2)\hbox{det}({\bf 23}a_1)\Bigr)} \\
  =& {1\over2}\sqrt{\Bigl(2+{23\over2}\Bigr)
       \Bigl(3\times{2\over23}\Bigr)} ={4.5\over\sqrt{23}} 
\end{split}
\end{align}\end{subequations}

The sum of volumes of the individual holes within any fundamental region hence
works out to
\begin{equation}
1\times{4\over\sqrt{23}} + 2\times{5\over\sqrt{23}} +
2\times{4.5\over\sqrt{23}} = {23\over\sqrt{23}}
\end{equation}
as expected from
formula (6.2). This provides a check that the decomposition into generalized
holes was both complete and error free.
 
\subsubsection{Hyperbolic Subalgebras for N=23}
 
We now search for hyperbolic Lie algebras contained in the root system of
${\cal G}_{23}$. As we
do not have any affine, or finite, subalgebra containing more than 2 roots we
can restrict ourselves to hyperbolic Lie algebras of rank 3. Furthermore,
we will not obtain pointed arrows as the norm of the longer root is 46.
We recall from chapter 5.2 that in the notation of \cite{Wan91} a double arrow
corresponds to a Cartan matrix ${{2\ \ -2}\choose{-2\ \ 2}}$, and that a
single bond corresponds to a Cartan matrix ${{2\ \ -1}\choose{-1\ \ 2}}$.
Checking \cite{Wan91}, we are left with the five potential candidates only which
are shown in figure 6.2.
All of these contain an affine subalgebra $A_1$. The decomposition of the
volume has proved that up to automorphism there is a unique representative of
$A_1$ within ${\cal R}$, which without loss of generality can be taken to be
$\lambda_1,\lambda_2$. (Recall chapter 5.2, particularly formula (5.3), for
details on how to translate Dynkin-diagrams into
distances in ${\cal R}$.) We fix a double arrow within the first of the above
diagrams. The rest of the diagram then translates into a pair of distances
(4,6) of the remaining third root from $\lambda_1$ and $\lambda_2$. Clearly
the 0 vector does satisfy these conditions. Hence we have found that the
corresponding real simple roots of ${\cal G}_{23}$,
\begin{equation}
\alpha_1=(\lambda_1,1,1), \mskip 30mu
\alpha_2=(\lambda_2,1,2)), \mskip 30mu \alpha_3=(0,1,-1)
\end{equation}
do represent the Dynkin diagram. This algebra is named $H^{(3)}_{96}$ in Wan's
classification and $AE_3$ in \cite{Kac90} (chapter 4).
It is straightforward to check the remaining four candidates for hyperbolic
subalgebras. They, however, cannot be realized within $\cal R$.
 
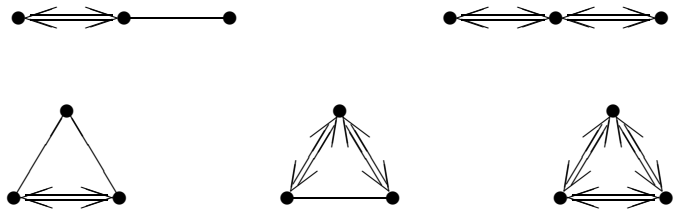
\begin{figure}
\begin{center} \begin{picture}(170,-35)(0,35)
\put(15,0){\circle*{5}}
\put(55,0){\circle*{5}}
\put(95,0){\circle*{5}}
\put(55,0){\line(1,0){40}}
\put(19.5,0.9) {\line(1,0){31}}
\put(19.5,-0.9) {\line(1,0){31}}
\put(17.5,0) {\line(3,1){12}}
\put(17.5,0) {\line(3,-1){12}}
\put(52.5,0) {\line(-3,1){12}}
\put(52.5,0) {\line(-3,-1){12}}
\end{picture}
\begin{picture}(100,-35)(0,35)
\put(5,0){\circle*{5}}
\put(45,0){\circle*{5}}
\put(85,0){\circle*{5}}
\put(9.5,0.9) {\line(1,0){31}}
\put(9.5,-0.9) {\line(1,0){31}}
\put(7.5,0) {\line(3,1){12}}
\put(7.5,0) {\line(3,-1){12}}
\put(42.5,0) {\line(-3,1){12}}
\put(42.5,0) {\line(-3,-1){12}}
\put(49.5,0.9) {\line(1,0){31}}
\put(49.5,-0.9) {\line(1,0){31}}
\put(47.5,0) {\line(3,1){12}}
\put(47.5,0) {\line(3,-1){12}}
\put(82.5,0) {\line(-3,1){12}}
\put(82.5,0) {\line(-3,-1){12}}
\end{picture} \end{center}
 
\begin{center} \begin{picture}(100,100)(0,0)
\put(30,0){\circle*{5}}
\put(30,0) {\line(3,5){20}}
\put(50,33){\circle*{5}}
\put(70,0) {\line(-3,5){20}}
\put(70,0){\circle*{5}}
\put(34.5,0.9) {\line(1,0){31}}
\put(34.5,-0.9) {\line(1,0){31}}
\put(32.5,0) {\line(3,1){12}}
\put(32.5,0) {\line(3,-1){12}}
\put(67.5,0) {\line(-3,1){12}}
\put(67.5,0) {\line(-3,-1){12}}
\end{picture}
\begin{picture}(100,100)(0,0)             
\put(30,0){\circle*{5}}
\put(31.5,2.5) {\line(1,6){1.9}}
\put(31.5,2.5) {\line(4,3){10}}
\put(48.75,30.5) {\line(-1,-6){1.9}}
\put(48.75,30.5) {\line(-4,-3){10}}
\put(32.25,5.25) {\line(3,5){13.5}}
\put(34.25,4.65) {\line(3,5){13.5}}
\put(50,33){\circle*{5}}
\put(67.75,5.25) {\line(-3,5){13.5}}
\put(65.75,4.65) {\line(-3,5){13.5}}
\put(68.5,2.5) {\line(-1,6){1.9}}
\put(68.5,2.5) {\line(-4,3){10}}
\put(51.25,30.5) {\line(1,-6){1.9}}
\put(51.25,30.5) {\line(4,-3){10}}
\put(70,0){\circle*{5}}
\put(30,0) {\line(1,0){40}}
\end{picture}
\begin{picture}(100,100)(0,0)             
\put(30,0){\circle*{5}}
\put(31.5,2.5) {\line(1,6){1.9}}
\put(31.5,2.5) {\line(4,3){10}}
\put(48.75,30.5) {\line(-1,-6){1.9}}
\put(48.75,30.5) {\line(-4,-3){10}}
\put(32.25,5.25) {\line(3,5){13.5}}
\put(34.25,4.65) {\line(3,5){13.5}}
\put(50,33){\circle*{5}}
\put(67.75,5.25) {\line(-3,5){13.5}}
\put(65.75,4.65) {\line(-3,5){13.5}}
\put(68.5,2.5) {\line(-1,6){1.9}}
\put(68.5,2.5) {\line(-4,3){10}}
\put(51.25,30.5) {\line(1,-6){1.9}}
\put(51.25,30.5) {\line(4,-3){10}}
\put(70,0){\circle*{5}}
\put(34.5,0.9) {\line(1,0){31}}
\put(34.5,-0.9) {\line(1,0){31}}
\put(32.5,0) {\line(3,1){12}}
\put(32.5,0) {\line(3,-1){12}}
\put(67.5,0) {\line(-3,1){12}}
\put(67.5,0) {\line(-3,-1){12}}
\end{picture} \end{center}
\caption{Candidates in the case $N=23$}
\end{figure}
  
Let us return to $AE_3$. The corollary to theorem 1.7 now shows that
\begin{subequations}
\begin{equation}
\hbox{mult}(r) \le p_\sigma(1-{r^2\over2}) \mskip 30mu \hbox{if} \mskip
20mu r\not\in 23L^*
\end{equation}
\begin{equation}
\hbox{mult}(r) \le p_\sigma(1-{r^2\over2})+p_\sigma(1-{r^2\over46})
\mskip 30mu \hbox{if} \mskip 20mu r\in 23L^*.
\end{equation}
\end{subequations}
Here $p_\sigma(1+n)$ is defined to
be the coefficient of $q^{1+n}$ in the $q$-expansion of
$q\eta^{-1}(q)$ $\eta^{-1}(q^{23})$
as in formula (1.24). We will want to compare these coefficients with some
obtained in \cite{FF83}. For convenience, we rewrite the definition of the
$p_\sigma$ as follows: 
\begin{equation}
\sum_n p_\sigma(n)q^n=
\prod_n(1-q^n)^{-1}\bigl(1+q^{23}+2q^{46}+\dots\bigr).
\end{equation}
 
The imaginary roots can be partially ordered by norm and height and are
described as a linear combination of the simple roots. In the table 6.1 below,
`2,2,1' represents the imaginary root $2\alpha_1+2\alpha_2+1\alpha_3$ with roots
as labelled in the Dynkin diagram. We list a number of imaginary roots of small
height and give their multiplicities. The multiplicities can be found in
\cite{Kac90}, p.215. They are reproduced here for convenience. Obviously, it is
sufficient to register the roots of a Weyl chamber. The column `bound' shows
the upper bounds as obtained in formula (6.6).
 
\vfil\eject
 
\begin{center}
\begin{picture}(170,-35)(-110,35)
\put(-120,0) {$H^{(3)}_{96}= AE_3 \subset {\cal G}_{23}$ }
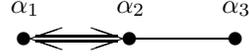
\put(10,10){$\alpha_1$}
\put(50,10){$\alpha_2$}
\put(90,10){$\alpha_3$}
\put(15,0){\circle*{5}}
\put(55,0){\circle*{5}}
\put(95,0){\circle*{5}}
\put(55,0){\line(1,0){40}}
\put(19.5,0.9) {\line(1,0){31}}
\put(19.5,-0.9) {\line(1,0){31}}
\put(17.5,0) {\line(3,1){12}}
\put(17.5,0) {\line(3,-1){12}}
\put(52.5,0) {\line(-3,1){12}}
\put(52.5,0) {\line(-3,-1){12}}
\end{picture}
\end{center}
 
\vspace{0.8in}
 
\begin{center}
\begin{minipage}[t]{2.3in}
\begin{tabular}[t]{c|c|c|c|}
 \multicolumn{1}{|c|} {coefficients} & norm & mult & bound \\ \hline
 \multicolumn{1}{|c|} {1, 1, 0} &  0 & 1 & 1\vspace{-0.015in}\\
 \multicolumn{1}{|c|} {2, 2, 1} & -2 & 2 & 2\vspace{-0.015in}\\
 \multicolumn{1}{|c|} {3, 3, 1} & -4 & 3 & 3\vspace{-0.015in}\\
 \multicolumn{1}{|c|} {3, 4, 2} & -6 & 5 & 5\vspace{-0.015in}\\
 \multicolumn{1}{|c|} {4, 4, 1} & -6 & 5 & 5\vspace{-0.015in}\\
 \multicolumn{1}{|c|} {4, 4, 2} & -8 & 7 & 7\vspace{-0.015in}\\
 \multicolumn{1}{|c|} {5, 5, 1} & -8 & 7 & 7\vspace{-0.015in}\\
 \multicolumn{1}{|c|} {4, 5, 2} & -10 & 11 & 11\vspace{-0.015in}\\
 \multicolumn{1}{|c|} {6, 6, 1} & -10 & 11 & 11\vspace{-0.015in}\\
 \multicolumn{1}{|c|} {5, 5, 2} & -12 & 15 & 15\vspace{-0.015in}\\
 \multicolumn{1}{|c|} {7, 7, 1} & -12 & 15 & 15\vspace{-0.015in}\\
 \multicolumn{1}{|c|} {5, 6, 2} & -14 & 22 & 22\vspace{-0.015in}\\
 \multicolumn{1}{|c|} {8, 8, 1} & -14 & 22 & 22\vspace{-0.015in}\\
 \multicolumn{1}{|c|} {5, 6, 3} & -16 & 30 & 30\vspace{-0.015in}\\
 \multicolumn{1}{|c|} {6, 6, 2} & -16 & 30 & 30\vspace{-0.015in}\\
 \multicolumn{1}{|c|} {9, 9, 1} & -16 & 30 & 30\vspace{-0.015in}\\
 \multicolumn{1}{|c|} {6, 6, 3} & -18 & 42 & 42\vspace{-0.015in}\\
 \multicolumn{1}{|c|} {6, 7, 2} & -18 & 42 & 42\vspace{-0.015in}\\
 \multicolumn{1}{|c|} {10, 10, 1} & -18 & 42 & 42\vspace{-0.015in}\\
 \multicolumn{1}{|c|} {7, 7, 2} & -20 & 56 & 56\vspace{-0.015in}\\
 \multicolumn{1}{|c|} {11, 11, 1} & -20 & 56 & 56\vspace{-0.015in}\\
 \multicolumn{1}{|c|} {6, 7, 3} & -22 & 77 & 77\vspace{-0.015in}\\
 \multicolumn{1}{|c|} {7, 8, 2} & -22 & 77 & 77\vspace{-0.015in}\\
 \multicolumn{1}{|c|} {12, 12, 1} & -22 & 77 & 77\vspace{-0.015in}\\
 \multicolumn{1}{|c|} {6, 8, 4} & -24 & 101 & 101\vspace{-0.015in}\\
 \multicolumn{1}{|c|} {7, 7, 3} & -24 & 101 & 101\vspace{-0.015in}\\
 \multicolumn{1}{|c|} {8, 8, 2} & -24 & 101 & 101\vspace{-0.015in}\\
 \multicolumn{1}{|c|} {13, 13, 1} & -24 & 101 & 101\vspace{-0.015in}\\
 \multicolumn{1}{|c|} {8, 9, 2} & -26 & 135 & 135\vspace{-0.015in}\\
 \multicolumn{1}{|c|} {14, 14, 1} & -26 & 135 & 135\vspace{-0.015in}\\
 \multicolumn{1}{|c|} {7, 8, 3} & -28 & 176 & 176\vspace{-0.015in}\\
 \multicolumn{1}{|c|} {9, 9, 2} & -28 & 176 & 176\vspace{-0.015in}\\
 \multicolumn{1}{|c|} {15, 15, 1} & -28 & 176 & 176\vspace{-0.015in}\\
 \multicolumn{1}{|c|} {7, 8, 4} & -30 & 231 & 231\vspace{-0.015in}\\
 \multicolumn{1}{|c|} {8, 8, 3} & -30 & 231 & 231\vspace{-0.015in}\\
 \multicolumn{1}{|c|} {9, 10, 2} & -30 & 231 & 231\vspace{-0.015in}\\
 \multicolumn{1}{|c|} {16, 16, 1} & -30 & 231 & 231\vspace{-0.015in}\\
 \multicolumn{1}{|c|} {7, 9, 4} & -32 & 297 & 297\vspace{-0.015in}\\
 \multicolumn{1}{|c|} {8, 8, 4} & -32 & 297 & 297 \\ \hline
\end{tabular}
\end{minipage}
$\;\;\;$
\begin{minipage}[t]{2.3in}
\begin{tabular}[t]{c|c|c|c|}
 \multicolumn{1}{|c|} {coefficients} & norm & mult & bound \\ \hline
 \multicolumn{1}{|c|} {10, 10, 2} & -32 & 297 & 297\vspace{-0.015in}\\
 \multicolumn{1}{|c|} {17, 17, 1} & -32 & 297 & 297\vspace{-0.015in}\\
 \multicolumn{1}{|c|} {8, 9, 3} & -34 & 385 & 385\vspace{-0.015in}\\
 \multicolumn{1}{|c|} {10, 11, 2} & -34 & 385 & 385\vspace{-0.015in}\\
 \multicolumn{1}{|c|} {18, 18, 1} & -34 & 385 & 385\vspace{-0.015in}\\
 \multicolumn{1}{|c|} {9, 9, 3} & -36 & 490 & 490\vspace{-0.015in}\\
 \multicolumn{1}{|c|} {11, 11, 2} & -36 & 490 & 490\vspace{-0.015in}\\
 \multicolumn{1}{|c|} {19, 19, 1} & -36 & 490 & 490\vspace{-0.015in}\\
 \multicolumn{1}{|c|} {8, 9, 4} & -38 & 627 & 627\vspace{-0.015in}\\
 \multicolumn{1}{|c|} {11, 12, 2} & -38 & 626 & 627\vspace{-0.015in}\\
 \multicolumn{1}{|c|} {20, 20, 1} & -38 & 627 & 627\vspace{-0.015in}\\
 \multicolumn{1}{|c|} {8, 10, 4} & -40 & 792 & 792\vspace{-0.015in}\\
 \multicolumn{1}{|c|} {9, 9, 4} & -40 & 792 & 792\vspace{-0.015in}\\
 \multicolumn{1}{|c|} {9, 10, 3} & -40 & 792 & 792\vspace{-0.015in}\\
 \multicolumn{1}{|c|} {12, 12, 2} & -40 & 791 & 792\vspace{-0.015in}\\
 \multicolumn{1}{|c|} {21, 21, 1} & -40 & 792 & 792\vspace{-0.015in}\\
 \multicolumn{1}{|c|} {8, 10, 5} & -42 & 1002 & 1002\vspace{-0.015in}\\
 \multicolumn{1}{|c|} {10, 10, 3} & -42 & 1002 & 1002\vspace{-0.015in}\\
 \multicolumn{1}{|c|} {12, 13, 2} & -42 & 1001 & 1002\vspace{-0.015in}\\
 \multicolumn{1}{|c|} {13, 13, 2} & -44 & 1253 & 1256\vspace{-0.015in}\\
 \multicolumn{1}{|c|} {9, 10, 4} & -46 & 1574 & 1576\vspace{-0.015in}\\
 \multicolumn{1}{|c|} {10, 11, 3} & -46 & 1574 & 1576\vspace{-0.015in}\\
 \multicolumn{1}{|c|} {13, 14, 2} & -46 & 1571 & 1576\vspace{-0.015in}\\
 \multicolumn{1}{|c|} {9, 10, 5} & -48 & 1957 & 1960\vspace{-0.015in}\\
 \multicolumn{1}{|c|} {9, 11, 4} & -48 & 1957 & 1960\vspace{-0.015in}\\
 \multicolumn{1}{|c|} {10, 10, 4} & -48 & 1957 & 1960\vspace{-0.015in}\\
 \multicolumn{1}{|c|} {11, 11, 3} & -48 & 1956 & 1960\vspace{-0.015in}\\
 \multicolumn{1}{|c|} {14, 14, 2} & -48 & 1953 & 1960\vspace{-0.015in}\\
 \multicolumn{1}{|c|} {10, 10, 5} & -50 & 2434 & 2439\vspace{-0.015in}\\
 \multicolumn{1}{|c|} {14, 15, 2} & -50 & 2429 & 2439\vspace{-0.015in}\\
 \multicolumn{1}{|c|} {9, 11, 5} & -52 & 3007 & 3015\vspace{-0.015in}\\
 \multicolumn{1}{|c|} {11, 12, 3} & -52 & 3005 & 3015\vspace{-0.015in}\\
 \multicolumn{1}{|c|} {15, 15, 2} & -52 & 3000 & 3015\vspace{-0.015in}\\
 \multicolumn{1}{|c|} {9, 12, 6} & -54 & 3712 & 3725\vspace{-0.015in}\\
 \multicolumn{1}{|c|} {10, 11, 4} & -54 & 3713 & 3725\vspace{-0.015in}\\
 \multicolumn{1}{|c|} {12, 12, 3} & -54 & 3710 & 3725\vspace{-0.015in}\\
 \multicolumn{1}{|c|} {15, 16, 2} & -54 & 3702 & 3725\vspace{-0.015in}\\
 \multicolumn{1}{|c|} {10, 12, 4} & -56 & 4557 & 4576\vspace{-0.015in}\\
 \multicolumn{1}{|c|} {11, 11, 4} & -56 & 4557 & 4576 \\ \hline
\end{tabular}
\end{minipage}
\end{center}

\vspace{-0.1in}

\begin{table}[ht]
\caption{root multiplicities of $AE_3$}
\end{table}
 
In the context of hyperbolic Lie algebras we require the definition of
the level of a root of a hyperbolic Lie algebra.
Suppose we are given a hyperbolic Lie algebra of, say, $n$ simple roots.
Suppose that its Dynkin diagram contains a unique affine subalgebra, which
then necessarily has $n-1$ simple roots, determining a unique simple root
$\alpha_0$ not in the affine diagram. If we express any root $r$ as a positive
linear combination of the simple roots $r=\sum n_i\alpha_i$ the level of $r$
is defined to be the coefficient $n_0$ of $\alpha_0$.
Hence in our labelling of the simple roots the level of a root of $AE_3$ is
simply the coefficient of $\alpha_3$.
It follows from proposition 10.10 and exercise 11.7 of \cite{Kac90} that for simply
laced hyperbolic Lie algebras all roots of level 1 and rank $n$ satisfy
\begin{equation}
\hbox{mult}(r)=p_{n-2}(1-{r^2\over2}) 
\end{equation}
Here, $p_n(x)$ is the number of partitions of $x$ into parts of $n$ colours.
\flexskip
 
Let us return to $AE_3$. This algebra has been studied extensively in \cite{FF83}.
\cite{FF83} define the
following series $p(n)$ and $p'(n)$: 
\begin{subequations}
\begin{equation}
\sum_n p(n)q^n = \prod_n (1-q^n)^{-1}
\end{equation}
Thus $p(n)$ is the number of partitions of $n$.
\begin{equation}
\sum_n p'(n)q^n = \bigl(\prod_n (1-q^n)^{-1}\bigr)\bigl(1-q^{20}+q^{22}\pm
\dots \bigr)
\end{equation}
\end{subequations}
\cite{FF83} then restate the root multiplicities for level 1 explicitly as follows:
\begin{subequations}
\begin{equation}
\hbox{mult}(r)=p(1-{r^2\over2})
\end{equation}
\cite{FF83} is mainly concerned with level 2. For roots of level 2 they obtain:
\begin{equation}
\hbox{mult}(r)=p'(1-{r^2\over2})
\end{equation}
\end{subequations}
 
Let us consider a root of level $r_3$, say
\[
 r=r_1\alpha_1+r_2\alpha_2+r_3\alpha_3=
 \bigl(r_1\lambda_1+r_2\lambda_2,r_1+r_2+r_3,r_1+2r_2-r_3\bigr).
\]
This will be an element of $23L^*$ if and only if
\begin{subequations}
\begin{align}
r_1\lambda_1+r_2\lambda_2  &\in 23 (\Lambda^\sigma)^* \\
r_1+r_2+r_3                &\equiv 0(23)              \\
r_1+2r_2-r_3               &\equiv 0(23)
\end{align}
\end{subequations}
 
Without loss of generality we can consider coefficients $r_1,r_2,r_3$ modulo
23. Comparing with the basis (formula (6.1b)) of $(\Lambda^\sigma)^*$,
we conclude that (6.11a) will only be satisfied
if $r_1\equiv -8n$, $r_2\equiv 9n$, for some integer $n$. We feed this into
(6.11b) and obtain $r_3\equiv -n$. Now, (6.11c) reads $11n\equiv 0$. Thus, in
summary, $r$ will be an element of $23L^*$ if and only if all three
coefficients $r_1,r_2,r_3$ are divisible by 23.
 
In particular, we observe that, for $r$ of level 1 or 2, $r\not\in 23L^*$.
Hence, for these $r$ the upper bounds of (6.6a) apply. We can now compare our
upper bound (6.7) to the exact values (6.9a) and (6.9b). We find that our
results are very close to the exact multiplicities where these are known. At
the same time, our results apply to all roots of $AE_3$.
 
Furthermore, we observe that the five real simple roots $0,\lambda_1,
\lambda_2, \lambda_1^*,\lambda_2^*$ form a symmetrized Cartan matrix as
follows:
\[
 \tilde C=\begin{pmatrix} 2&-2&0&0&0\\ -2&2&-1&0&0\\ 0&-1&2&-23&0\\
0&0&-23&46&-46\\ 0&0&0&-46&46\end{pmatrix}
\]
In particular, the equation $23(\lambda_1+\lambda_2)=\lambda_1^*+\lambda_2^*$
shows that for every root in $23L^*$ there are additional roots in
${\cal G}_{23}$ which are not in the hyperbolic Lie algebra $AE_3$.
This provides some idea why the upper bounds are not exact for all roots.
 
\subsection{N=11}
\subsubsection{Finite and Affine Subalgebras}
 
Let ${\cal G}_{11}$ denote the GKM constructed in chapter 1 from an
automorphism $\sigma$ of cycle shape $1^2 11^2$. Let $\Lambda^\sigma$ be the
4-dimensional fixed point lattice and $L=\Lambda^\sigma \oplus II_{1,1}$ be
the corresponding Lorentzian lattice. The set of its real simple roots
${\cal R}$ has been calculated in theorem 5.1. A basis of $\Lambda^\sigma$
is easily obtained as
\[
  (-3,1,\sqrt{11},\sqrt{11}),\mskip 20mu(2,0,2\sqrt{11},0),
  \mskip 20mu  (4,4,0,0), \mskip 20mu(8,0,0,0).
\]
Here, as in the rest of this section, we operate in basic units of
$1\over\sqrt8$, that is, we suppress this factor, which is common to all vectors
considered. The basis vectors form a diagonal matrix and we can read off the
fundamental volume to be $\sqrt{11}^2$.
 
The fixed point lattice is a sublattice of the Leech lattice. Hence the
shortest vectors will have norm 4. From the considerations of chapter 5.2,
formula (5.3) it follows that any finite or affine diagram has bonds which
correspond to distances between the roots either 4, or 6, or 8. Without
loss of generality we can choose the 0 vector as one of the roots. Then we can
restrict our attention to vectors of norm  4,6, and 8.
A complete list of these is easily obtained from the explicit basis and given
in table 6.2:
 
\begin{table}[ht]
\renewcommand\arraystretch{1.5}
\noindent\[
\begin{array}{|c|c|c|}
\text{norm 4 vectors} & \text{norm 6 vectors} & \text{norm 8 vectors} \\
\hline
(4,4,0,0) & (2,0,2\sqrt{11},0) & (8,0,0,0) \\
(4,-4,0,0)& (0,2,0,2\sqrt{11}) & (0,8,0,0) \\
(3,-1,-\sqrt{11},-\sqrt{11})&(5,1,\sqrt{11},\sqrt{11})
                                  &(4,-2,0,2\sqrt{11}) \\
(-1,3,-\sqrt{11},-\sqrt{11})&(1,5,\sqrt{11},\sqrt{11})
                                  &(-4,-2,0,2\sqrt{11}) \\
(3,1,-\sqrt{11},\sqrt{11})&(5,-1,\sqrt{11},-\sqrt{11})
                                  &(-2,4,2\sqrt{11},0) \\
(1,3,-\sqrt{11},\sqrt{11})&(-1,5,-\sqrt{11},\sqrt{11})
                                  &(-2,-4,2\sqrt{11},0) \\
\hline
\end{array}
\]
\caption{Short vectors of $\Lambda^\sigma$}
\end{table}
  
We now turn to the elements of ${\cal R}_{dual}$. As established in
chapter 5.2 (remark following theorem 5.4), there cannot be any bonds between
elements of ${\cal R}_{fix}$ and ${\cal R}_{dual}$ in a diagram of
finite or affine type. If we require 0 to be one of the roots we can restrict
our attention to those elements of the dual which have norm exactly
$2+{2\over11}$. To identify those, we begin by recalling the remark following
lemma 4.1 that the additive group $(\Lambda^\sigma)^*/\Lambda^\sigma$ has two
generators, each of order 11. We can choose them as
\[
  [\lambda_1^*]=(3,-1,{\sqrt{11}\over11},{\sqrt{11}\over11}), \mskip 40mu
  [\lambda_2^*]=(3,1,{\sqrt{11}\over11},{-\sqrt{11}\over11}). 
\]
To find elements of the dual lattice of norm $2+{2\over11}$ we require
$(a[\lambda_1^*]+b[\lambda_2^*])^2 \equiv {2\over11}\ \hbox{mod}\ 2\ZZ$, where
$a,b$ are integers modulo 11. This
is equivalent to the condition 
\begin{equation}
a^2+b^2 \equiv 8\ \hbox{mod}\ 11.
\end{equation}
As established in lemma 4.2, this identity has 12 solutions: They are
$\{(\pm2,\pm2),$ $(\pm4,\pm5),$ $(\pm5,\pm4)\}$.
We consider the solution $(2,2)$. The corresponding equivalence class is
represented by the vector
$(4,0,{4\over11}\sqrt{11},0)$. A check of all short vectors establishes that
there are precisely 6 representatives of norm $2+{2\over11}$, namely
$(\pm4,0,{4\over11}\sqrt{11},0)$, $(0,\pm4,{4\over11}\sqrt{11},0)$,
 $(-1,\pm1,{-7\over11}\sqrt{11},{\pm11\over11}\sqrt{11})$.
 
We recall the claim, proved and used in the deduction of theorem 
4.2, that $GO_2(11)$ acts transitively on the 12 solutions of (6.12). It 
furthermore lifts to automorphisms of the fixed point lattice. Hence we 
have found that there are $12\times 6=72$ elements of norm $2+{2\over11}$ 
within ${\cal R}_{dual}$. From here, it is straightforward to calculate 
their co-ordinates. We will, however, omit a complete list.
 
The search for generalized holes in ${\cal R}$ will be simplified if we can
make use of the automorphisms of the lattice. There are a number of obvious
symmetries of order 2:
\begin{equation}
  \phi_1=-\hbox{id} \mskip 40mu
  \phi_2=\begin{pmatrix}1&0&0&0\\ 0&-1&0&0\\ 0&0&1&0\\ 0&0&0&-1 \end{pmatrix} \mskip 40mu
  \phi_3=\begin{pmatrix}0&1&0&0\\ 1&0&0&0\\ 0&0&0&1\\ 0&0&1&0 \end{pmatrix} 
\end{equation}
Referring to the complete list of short vectors we observe that there are just
two types of norm 8 vectors which cannot be identified by the above 
automorphisms.
As the centralizer of $\sigma$ in $Co_0$ has order 66 we may now conjecture
that there exists an automorphism of order 3 identifying the two types. (It
must be stressed that it can only be a pious hope to expect that all norm 8
vectors will be equivalent, even though it is true for
Leech lattice vectors. It is, for example, false in the case $N=2$, where
there exist two types of diagrams, $A_1^{16}$ and $A_1^8 \ {\bf 2}A_1^8$.)
Let us try to construct an automorphism $\phi$ mapping $(8,0,0,0)^T$ to
$(4,-2,0,2\sqrt{11})^T$. Thus the first column of $\phi$ will be $({1\over2},
-{1\over4},0,{1\over4}\sqrt{11})^T$. $(0,8,0,0)^T$ must now be mapped to some
other norm 8 vector which still will be orthogonal to $\phi(8,0,0,0)^T$.
Checking the list of norm 8 vectors there are precisely 2 possible choices:
$\pm(2,4,-2\sqrt{11},0)^T$. (These two choices are equivalent because of
the automorphism $\phi_2$, formula (6.13)). Continuing in this way we obtain
\begin{equation}
\phi=\begin{pmatrix} {1\over2} & {1\over4} & 0 & -{1\over4}\sqrt{11} \\
		 -{1\over4} & {1\over2} & {1\over4}\sqrt{11} & 0 \\
		 0 & -{1\over4}\sqrt{11} & {1\over2} & -{1\over4}  \\
	     {1\over4}\sqrt{11} & 0 & {1\over4} & {1\over2}  \\ \end{pmatrix}
\end{equation}
We observe that $\phi$ does not only identify the two types of norm 8 vectors
which were not equivalent under the automorphism group generated by $\phi_1,
\phi_2, \phi_3$ but it does so for the norm 4 and norm 6 vectors as well.
\flexskip 
 
We begin the decomposition of space into generalized holes by
searching for $A_1$ diagrams, the only type requiring norm 8 vectors. As all
norm 8 vectors are equivalent under the automorphism group we choose the vectors
0 and $(8,0,0,0)$ without loss of generality. This is a complete affine
component, hence all remaining elements of this hole must have minimal distance
from both. This obviously leaves the vectors $(4,\pm4,0,0)$, and the dual
vectors $(4,0,\pm{4\over11}\sqrt{11},0)$, and $(4,0,\pm{4\over11}\sqrt{11},0)$.
Thus we have identified an $A_1^2\ {\bf 11}A_1^2$. Now there are 12 norm 8
vectors which each form one such diagram with 0. On the other hand, within each
diagram there are 4 ways of placing the 0 and norm 8 vectors. Hence, within any
fundamental volume, there are ${12\over4}=3$ such diagrams.
\flexskip
 
Next, let us search for diagrams of the type $a_1^n \ {\bf 11}\Delta$ where
$\Delta$ is any (not necessarily undecomposable) simply laced diagram of finite
type. Theorem 5.8 asserts that among the vertices of any generalized hole there
are at least $M+1$ elements of ${\cal R}_{fix}$.
Hence $n$ will be greater or equal 3. Without loss of generality, choose
0 and any of the norm 4 vectors, say $(4,4,0,0)$. There are only 2 vectors left
which have distance 4 to both the above: $(3,1,-\sqrt{11},\sqrt{11})$ and
$(1,3,\sqrt{11},-\sqrt{11})$. Automorphism $\phi_3$ shows that the 2 choices
for a third vector are equivalent. Furthermore, we observe that there cannot
be an $a_1^4$ diagram. We choose,
say, $(3,1,-\sqrt{11},\sqrt{11})$. Of the 72 dual vectors there are precisely 6
left which have minimal distance from the 3 chosen vectors. We refer to them
as vectors $a$ to $f$, respectively.
We calculate their distance table, where a `$*$' stands in the place of any
distance greater $6\over11$ as such a distance cannot occur within a finite
diagram. The distances are given in units of $1\over11$:
 
\begin{center}\begin{tabular}{c|cccccc}
  & $a$ & $b$ & $c$ & $d$ & $e$ & $f$ \\
\hline
$a$ & & 4 & $*$ & $*$ & $*$ & $*$\vspace{-0.015in}\\
$b$ & & & $*$ & $*$ & $*$ & $*$\vspace{-0.015in}\\
$c$ & & & & $*$ & $*$ & 4\vspace{-0.015in}\\
$d$ & & & & & 4 & $*$\vspace{-0.015in}\\
$c$ & & & & & & $*$ \\
\end{tabular} \end{center}
 
Hence, we identify three diagrams of type $a_1^3\ {\bf 11}a_1^2$. We have had
12 choices of the first norm 4 vector, another 2 choices for the second.
On the other hand, there are 3! ways of choosing which of the vectors
corresponds to which of the $a_1^3$. Hence, within any fundamental volume,
there are ${12\times2\times3\over3!}=12$ holes of type $a_1^3\ {\bf 11}a_1^2$.
\flexskip
 
Any generalized hole which is not of either of the above two types will contain
a single bond, that is a pair of vectors of ${\cal R}_{fix}$ at distance 6.
We choose without loss of generality 0 and $(5,1,\sqrt{11},\sqrt{11})$. We
first complete the search for affine diagrams. There are only 2
affine holes in the Leech lattice which admit an automorphism of order 11. They
are $A_1^{24}$ and $A_2^{12}$. Theorem 5.8 (I) asserts that all generalized
affine holes of ${\cal R}$ lift to affine holes of the Leech lattice whose
automorphism group has order divisible by $N$. Hence we can
restrict the affine search to $A_1$ and $A_2$. There are exactly 2 vectors of
norm 6 which have distance 6 from $(5,1,\sqrt{11},\sqrt{11})$. They are
$(5,-1,\sqrt{11},-\sqrt{11})$ and $(2,0,2\sqrt{11},0)$. Obviously these two
choices to form an $A_2$ are equivalent as the respective triplets of vectors
are translations of one another. Choosing either, and collecting all dual
vectors that have minimal distance to the three vectors we obtain a unique
diagram $A_2\ {\bf 11}A_2$. Taking into consideration the choices made we found
${12\times2\over3\times2}=4$ such diagrams within the fundamental volume.
 
We now turn to finite diagrams. Any such contains roots which are joined only
to one other root. Without loss of generality we can assume that the 0 root is
chosen to be such.
Then, apart from the fixed choice of one norm 6 vector, we will only have to
consider norm 4 vectors which have distance less or equal to 6 from the chosen
vector of norm 6. (This may not seem important, however, it does make a
difference if there are 4320 vectors of norm 4 and 61440 of norm 6, as in the
case $N=2$.) In the case at hand, there are 10 elements of
${\cal R}_{dual}$ of minimal
distance to both 0 and $(5,1,\sqrt{11},\sqrt{11})$. We will refer to them as
$a$ to $j$, respectively.
The distance table (in units of $1\over11$, as above) is as follows:
 
\begin{center}\begin{tabular}{c|cccccccccc}
  & $a$ & $b$ & $c$ & $d$ & $e$ & $f$ & $g$ & $h$ & $i$ & $j$ \\
\hline
$a$ & & 4 & $*$ & $*$ & 6 & 6 & $*$ & $*$ & $*$ & $*$\vspace{-0.015in}\\
$b$ & & & $*$ & $*$ & $*$ & $*$ & $*$ & $*$ & 6 & $*$\vspace{-0.015in}\\
$c$ & & & & 4 & $*$ & $*$ & 6 & $*$ & $*$ & 6\vspace{-0.015in}\\
$d$ & & & & & $*$ & $*$ & $*$ & 6 & $*$ & $*$\vspace{-0.015in}\\
$e$ & & & & & & 6 & $*$ & $*$ & $*$ & $*$\vspace{-0.015in}\\
$f$ & & & & & & & $*$ & 4 & $*$ & $*$\vspace{-0.015in}\\
$g$ & & & & & & & & $*$ & $*$ & 6\vspace{-0.015in}\\
$h$ & & & & & & & & & $*$ & $*$\vspace{-0.015in}\\
$i$ & & & & & & & & & & 4 \\
\end{tabular} \end{center}
 
We list the norm 4 vectors which satisfy all conditions imposed so far. Let
`dist' denote the distance to $(5,1,\sqrt{11},\sqrt{11})$. `duals' lists
those elements of ${\cal R}_{dual}$
which have minimal distance from the respective norm 4 vector:
 
\begin{center}\begin{tabular}{lll}
$A=(4,4,0,0)$ & $\hbox{dist}=4$ & $\hbox{duals}\ a,b,f,g,h$ \\
$B=(4,-4,0,0)$ & $\hbox{dist}=6$ & $\hbox{duals}\ a,b,e,i$ \\
$C=(1,-3,\sqrt{11},\sqrt{11})$ & $\hbox{dist}=4$ & $\hbox{duals}\ c,d,e,i,j$ \\
$D=(3,1,-\sqrt{11},\sqrt{11})$ & $\hbox{dist}=6$ & $\hbox{duals}\ b,g,i,j$ \\
\end{tabular} \end{center}
 
Every generalized hole of finite type in the case $N=11$ has 5 vertices. We
begin by searching for diagrams that contain the minimal number of elements of
${\cal R}_{fix}$, which is 3 (theorem 5.8). We then need 2 elements of
${\cal R}_{dual}$ to complete the diagram.
Choose $A$. $ab$, and $fh$, yield $a_2a_1\ {\bf 11}a_1^2$. $af$ yields $a_2a_1
\ {\bf 11}a_2$.
Choose $C$. $cd$, and $ij$, yield $a_2a_1\ {\bf 11}a_1^2$. $cj$ yields $a_2a_1
\ {\bf 11}a_2$.
Choose $B$. $ab$ yields $a_3\ {\bf 11}a_1^2$. $ae$, and $bi$, yield
	    $a_3\ {\bf 11}a_2$.
Choose $D$. $ij$ yields $a_3\ {\bf 11}a_1^2$. $bi$, and $gj$, yield
	    $a_3\ {\bf 11}a_2$.
Next, let us search for diagrams containing 4 elements of ${\cal R}_{fix}$
 and 1 element of ${\cal R}_{dual}$.
Choose $AD$. $b$, and $g$, yield $a_3a_1\ {\bf 11}a_1$.
Choose $BC$. $e$, and $i$, yield $a_3a_1\ {\bf 11}a_1$.
There are no further diagrams containing 4 elements of ${\cal R}_{fix}$
because the distances $AB$, $AC$, $BD$ are greater than 6. For the same reason,
there can be no diagrams containing 5 elements of ${\cal R}_{fix}$.
In conclusion, we have obtained:
 
\begin{center}\begin{tabular}{lrl}
type $a_3\ {\bf 11}a_1^2$:   & ${12\times(1+1)\over2}$ & $=12$ copies. \\
type $a_3\ {\bf 11}a_2$:     & ${12\times(2+2)\over2}$ & $=24$ copies. \\
type $a_2a_1\ {\bf 11}a_1^2$:& ${12\times(2+2)\over2}$ & $=24$ copies. \\
type $a_2a_1\ {\bf 11}a_2$:  & ${12\times(1+1)\over2}$ & $=12$ copies. \\
type $a_3a_1\ {\bf 11}a_1$:  & ${12\times(2+2)\over2}$ & $=24$ copies. \\
\end{tabular} \end{center}
 
We can now proceed to calculate the volumes of the individual holes and carry
out the volume check, just as in section 6.2.1 for the case $N=23$. The details
for the case $N=11$ can be found in appendix A.
 
So far we have established the numbers of holes of the various types. We now
turn to the question whether all holes of same type are equivalent under the
automorphism group of the fixed point lattice. We know that $11^2$ does not
divide the order of $Co_0$, thus certainly not the order of any automorphism
group of a hole in the Leech lattice. Hence, theorem 5.9 covers all affine
diagrams and all finite diagrams which contain precisely 3 elements of ${\cal
R}_{fix}$. This only leaves the diagram $a_3a_1\ {\bf 11}a_1$ to be
investigated. We consider a diagram of type $a_3\ {\bf 11}a_2$. This has two
faces of type $a_3\ {\bf 11}a_1$. If we consider a representative of this type
we will discover that the adjacent holes are $a_3\ {\bf 11}a_1^2$ and $a_3a_1\
{\bf 11}a_1$, respectively. Hence, the conditions of theorem 5.10 are satisfied
and we conclude that there are 24 holes of type $a_3a_1\ {\bf 11}a_1$
equivalent under the automorphism group. This concludes the argument as there
are precisely 24 such copies within any fundamental region.
 
We now turn to the size of the automorphism group.
Again, we consider the diagram $a_3\ {\bf 11}a_2$. We established
already that its two faces of type $a_3\ {\bf 11}a_1$ are not symmetric. But so
are the two faces of type $a_2\ {\bf 11}a_2$, the adjacent diagrams being
$A_2\ {\bf 11}A_2$ and $a_2a_1\ {\bf 11}a_2$. Hence, the hole $a_3\ {\bf 11}
a_2$ will only be preserved by the identity automorphism. As there exist 24
copies of it within a fundamental volume, the total automorphism group must
have order 24. Thus we have verified the complete decomposition of the
fundamental region as given in appendix A.
 
\subsubsection{Hyperbolic Subalgebras}
 
We can now carry out the search for hyperbolic subalgebras. We observe that
there are no finite, or affine, subalgebras of more than 3 roots. Hence, we can
restrict ourselves to hyperbolic Lie algebras of rank 3 and 4, whose roots, as
in the case $N=23$, must all be of equal length. That leaves us with 5
hyperbolic Lie algebras of rank 3 and 3 of rank 4 from the list of \cite{Wan91}.
They contain an $A_1$ or $A_2$ subdiagram, respectively. From the decomposition
of space in section 6.2.2.1 we know that up to isomorphism there are two unique
representatives of these in ${\cal R}$, one consisting of roots of norm 2, the
other consisting of roots of norm $2N$. All that remains now is to check the
elements of ${\cal R}$ close to the representing diagrams $A_1$, $A_2$, ${\bf
11}A_1$, and ${\bf 11}A_2$. We find that there exist exactly 3 hyperbolic Lie
subalgebras of ${\cal G}_{11}$, shown in figure 6.3.
 
\begin{figure}
\begin{center}
\begin{picture}(100,50)(0,0)
\put(15,0){\circle*{5}}
\put(55,0){\circle*{5}}
\put(95,0){\circle*{5}}
\put(55,0){\line(1,0){40}}
\put(19.5,0.9) {\line(1,0){31}}
\put(19.5,-0.9) {\line(1,0){31}}
\put(17.5,0) {\line(3,1){12}}
\put(17.5,0) {\line(3,-1){12}}
\put(52.5,0) {\line(-3,1){12}}
\put(52.5,0) {\line(-3,-1){12}}
\end{picture}
\begin{picture}(100,50)(0,0)             
\put(30,0){\circle*{5}}
\put(31.5,2.5) {\line(1,6){1.9}}
\put(31.5,2.5) {\line(4,3){10}}
\put(48.75,30.5) {\line(-1,-6){1.9}}
\put(48.75,30.5) {\line(-4,-3){10}}
\put(32.25,5.25) {\line(3,5){13.5}}
\put(34.25,4.65) {\line(3,5){13.5}}
\put(50,33){\circle*{5}}
\put(67.75,5.25) {\line(-3,5){13.5}}
\put(65.75,4.65) {\line(-3,5){13.5}}
\put(68.5,2.5) {\line(-1,6){1.9}}
\put(68.5,2.5) {\line(-4,3){10}}
\put(51.25,30.5) {\line(1,-6){1.9}}
\put(51.25,30.5) {\line(4,-3){10}}
\put(70,0){\circle*{5}}
\put(34.5,0.9) {\line(1,0){31}}
\put(34.5,-0.9) {\line(1,0){31}}
\put(32.5,0) {\line(3,1){12}}
\put(32.5,0) {\line(3,-1){12}}
\put(67.5,0) {\line(-3,1){12}}
\put(67.5,0) {\line(-3,-1){12}}
\end{picture}
\begin{picture}(100,50)(0,0)
\put(10,0){\circle*{5}}
\put(50,0){\circle*{5}}
\put(90,0){\circle*{5}}
\put(30,33.3){\circle*{5}}
\put(10,0){\line(1,0){40}}
\put(10,0){\line(3,5){20}}
\put(30,33.3){\line(3,-5){20}}
\put(50,0){\line(1,0){40}}
\end{picture} \end{center}
\caption{Hyperbolic Subalgebras of ${\cal G}_{11}$}
\end{figure}
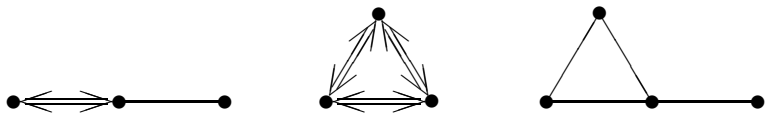
 
The first is the algebra $AE_3$  which we identified as a
subalgebra of ${\cal G}_{23}$. As the upper bounds arising from $N=11$ are
worse than those arising from $N=23$ we can discard this case. The second can
be realized by the vectors 0, $(8,0,0,0)$, and $(4,2,0,-2\sqrt{11})$. The
third is represented through the vectors 0, $(5,1,\sqrt{11},\sqrt{11})$,
$(5,-1,\sqrt{11},-\sqrt{11})$, and $(4,4,0,0)$.
We do not find any additional hyperbolic subalgebras consisting of long roots.
 
\flexskip
 
As in the case $N=23$, we obtain upper bounds for the root multiplicities of
the above hyperbolic Lie algebras. We proceed to compare these bounds with the
exact values for some imaginary roots, ordered by norm and height. The notation
is as before (see section 6.2.1.2), and I calculated the exact root
multiplicities by a program based on the Peterson recursion formula for root
multiplicities (see \cite{Kac90}, p.210).

\vfil\eject
{\parindent=0pt
\begin{minipage}[t]{2.2in}
\begin{picture}(100,-60)(0,60)
\put(0,60){$H^{(3)}_{71}\subset {\cal G}_{11}$}
\put(10,0){$\alpha_1$}
\put(30,0){\circle*{5}}
\put(31.5,2.5) {\line(1,6){1.9}}
\put(31.5,2.5) {\line(4,3){10}}
\put(48.75,30.5) {\line(-1,-6){1.9}}
\put(48.75,30.5) {\line(-4,-3){10}}
\put(32.25,5.25) {\line(3,5){13.5}}
\put(34.25,4.65) {\line(3,5){13.5}}
\put(50,40) {$\alpha_2$}
\put(50,33){\circle*{5}}
\put(67.75,5.25) {\line(-3,5){13.5}}
\put(65.75,4.65) {\line(-3,5){13.5}}
\put(68.5,2.5) {\line(-1,6){1.9}}
\put(68.5,2.5) {\line(-4,3){10}}
\put(51.25,30.5) {\line(1,-6){1.9}}
\put(51.25,30.5) {\line(4,-3){10}}
\put(70,0){\circle*{5}}
\put(75,0) {$\alpha_3$}
\put(34.5,0.9) {\line(1,0){31}}
\put(34.5,-0.9) {\line(1,0){31}}
\put(32.5,0) {\line(3,1){12}}
\put(32.5,0) {\line(3,-1){12}}
\put(67.5,0) {\line(-3,1){12}}
\put(67.5,0) {\line(-3,-1){12}}
\end{picture}
\end{minipage}
\begin{minipage}[t]{2.7in}
\begin{tabular}[t]{p{1.15in}c|c|c|}
 & \multicolumn{2}{c}{} & \multicolumn{1}{c}{\bf N=11} \\
 \multicolumn{1}{|c|} {root-coefficients} & norm & mult & bound\\ \hline
 \multicolumn{1}{|c|} {1, 1, 0} & 0 & 1 & 2\vspace{-0.015in}\\
 \multicolumn{1}{|c|} {1, 0, 1} & 0 & 1 & 2\vspace{-0.015in}\\
 \multicolumn{1}{|c|} {0, 1, 1} & 0 & 1 & 2\vspace{-0.015in}\\
 \multicolumn{1}{|c|} {1, 1, 1} & -6 & 2 & 20\vspace{-0.015in}\\
 \multicolumn{1}{|c|} {2, 1, 1} & -8 & 3 & 36\vspace{-0.015in}\\
 \multicolumn{1}{|c|} {1, 2, 1} & -8 & 3 & 36\vspace{-0.015in}\\
 \multicolumn{1}{|c|} {1, 1, 2} & -8 & 3 & 36\vspace{-0.015in}\\
 \multicolumn{1}{|c|} {2, 2, 1} & -14 & 6 & 185 \\
\end{tabular}
\end{minipage}
 
\rule{4.95in}{0.01in}
\vspace{0.1 in} 

\begin{minipage}[t]{2.2in}
\begin{picture}(100,-60)(0,60)
\put(0,60){$H^{(4)}_{3}= AE_4\subset {\cal G}_{11}$}
\put(10,-10){$\alpha_1$}
\put(30,40){$\alpha_2$}
\put(50,-10){$\alpha_3$}
\put(90,-10){$\alpha_4$}
\put(10,0){\circle*{5}}
\put(50,0){\circle*{5}}
\put(90,0){\circle*{5}}
\put(30,33.3){\circle*{5}}
\put(10,0){\line(1,0){40}}
\put(10,0){\line(3,5){20}}
\put(30,33.3){\line(3,-5){20}}
\put(50,0){\line(1,0){40}}
\end{picture}
\end{minipage}
\begin{minipage}[t]{2.7in}
\begin{tabular}[t]{p{1.15in}c|c|c|}
 & \multicolumn{2}{c}{} & \multicolumn{1}{c}{\bf N=11} \\
 \multicolumn{1}{|c|} {root-coefficients} & norm & mult & bound\\ \hline
 \multicolumn{1}{|c|} {1, 1, 1, 0} & 0 & 2 & 2\vspace{-0.015in}\\
 \multicolumn{1}{|c|} {2, 2, 2, 1} & -2 & 5 & 5\vspace{-0.015in}\\
 \multicolumn{1}{|c|} {3, 3, 3, 1} & -4 & 10 & 10\vspace{-0.015in}\\
 \multicolumn{1}{|c|} {3, 3, 4, 2} & -6 & 20 & 20\vspace{-0.015in}\\
 \multicolumn{1}{|c|} {4, 4, 4, 1} & -6 & 20 & 20\vspace{-0.015in}\\
 \multicolumn{1}{|c|} {4, 4, 4, 2} & -8 & 36 & 36\vspace{-0.015in}\\
 \multicolumn{1}{|c|} {5, 5, 5, 1} & -8 & 36 & 36\vspace{-0.015in}\\
 \multicolumn{1}{|c|} {4, 4, 5, 2} & -10 & 65 & 65\vspace{-0.015in}\\
 \multicolumn{1}{|c|} {6, 6, 6, 1} & -10 & 65 & 65\vspace{-0.015in}\\
 \multicolumn{1}{|c|} {5, 5, 5, 2} & -12 & 110 & 110\vspace{-0.015in}\\
 \multicolumn{1}{|c|} {5, 4, 6, 3} & -12 & 110 & 110\vspace{-0.015in}\\
 \multicolumn{1}{|c|} {4, 5, 6, 3} & -12 & 110 & 110\vspace{-0.015in}\\
 \multicolumn{1}{|c|} {7, 7, 7, 1} & -12 & 110 & 110\vspace{-0.015in}\\
 \multicolumn{1}{|c|} {5, 5, 6, 2} & -14 & 185 & 185\vspace{-0.015in}\\
 \multicolumn{1}{|c|} {8, 8, 8, 1} & -14 & 185 & 185\vspace{-0.015in}\\
 \multicolumn{1}{|c|} {5, 5, 6, 3} & -16 & 300 & 300\vspace{-0.015in}\\
 \multicolumn{1}{|c|} {6, 6, 6, 2} & -16 & 300 & 300\vspace{-0.015in}\\
 \multicolumn{1}{|c|} {6, 6, 7, 2} & -18 & 481 & 481\vspace{-0.015in}\\
 \multicolumn{1}{|c|} {6, 6, 6, 3} & -18 & 481 & 481\vspace{-0.015in}\\
 \multicolumn{1}{|c|} {6, 5, 7, 3} & -18 & 481 & 481\vspace{-0.015in}\\
 \multicolumn{1}{|c|} {5, 6, 7, 3} & -18 & 481 & 481\vspace{-0.015in}\\
 \multicolumn{1}{|c|} {7, 7, 7, 2} & -20 & 752 & 754\vspace{-0.015in}\\
 \multicolumn{1}{|c|} {6, 6, 7, 3} & -22 & 1165 & 1169\vspace{-0.015in}\\
 \multicolumn{1}{|c|} {7, 7, 8, 2} & -22 & 1164 & 1169\vspace{-0.015in}\\
 \multicolumn{1}{|c|} {7, 7, 7, 3} & -24 & 1770 & 1780\vspace{-0.015in}\\
 \multicolumn{1}{|c|} {7, 6, 8, 3} & -24 & 1769 & 1780\vspace{-0.015in}\\
 \multicolumn{1}{|c|} {6, 7, 8, 3} & -24 & 1769 & 1780\vspace{-0.015in}\\
 \multicolumn{1}{|c|} {6, 6, 8, 4} & -24 & 1769 & 1780\vspace{-0.015in}\\
 \multicolumn{1}{|c|} {8, 8, 8, 2} & -24 & 1767 & 1780\vspace{-0.015in}\\
 \multicolumn{1}{|c|} {7, 6, 8, 4} & -26 & 2663 & 2685\vspace{-0.015in}\\
 \multicolumn{1}{|c|} {6, 7, 8, 4} & -26 & 2663 & 2685\vspace{-0.015in}\\
 \multicolumn{1}{|c|} {7, 7, 8, 3} & -28 & 3950 & 3996\vspace{-0.015in}\\
 \multicolumn{1}{|c|} {7, 7, 8, 4} & -30 & 5812 & 5894 \\
\end{tabular}
\end{minipage}
 
\rule{4.95in}{0.01in}
\vspace{+0.05in}

\begin{table}[ht]
\caption{Root Multiplicities of Subalgebras of ${\cal G}_{11}$}
\end{table}
} 

We observe that the upper bounds for $AE_4$ are again very useful, though for
roots of larger height there do occur discrepancies between the correct values
and the upper bounds, similarly to the results for $AE_3$. However, the results
for $H^{(3)}_{71}$ are so poor that they must be regarded as useless. In fact,
the general bounds provided by theorem 1.4b are actually far closer than our
result. We will return to this in section 6.3.
 
\subsection{The remaining N}
 
\subsubsection{Finite and Hyperbolic Subalgebras}
 
The previous section showed in detail how to analyse the root system ${\cal R}$
of ${\cal G}_{11}$. I carried out the analysis of ${\cal G}_N$  for $N=7,5,3,2$
using the same principal line of the argument in all these cases. The numbers of
vectors of norms 4 and 6 do, however, increase to levels which are better
treated by computers. Again, it was necessary to make extensive use of the
known symmetries of the fixed point lattice in order to reduce the number of
individual representatives of each type that had to be identified and counted.
I determined the symmetries by explicit inspection. Given an understanding of
the symmetries, I was then in the position to carry out a computerized search
for generalized holes, using preselections similar to those of $A_1$, $a_1^n$,
$a_2$ in the case $N=11$ demonstrated in section 6.2.2.1, above. Later, I
transformed the numbers of representatives found into numbers
per fundamental region and carried out the volume check manually. I
subsequently analysed the automorphism groups of the individual holes by
inspection (with computerized search of representing holes). This provided an
important double check for the results. The checks of the automorphism groups
yielded the following result: In the cases $N=2,5,7,11,23$ it is true
that, whenever two generalized holes have the same Dynkin diagram, they are
equivalent under the automorphism group of $\Lambda^\sigma$.
However, in the case $N=3$, there do exist two pairs of generalized holes
which represent the same Dynkin diagram but are inequivalent under the
automorphism group of $\Lambda^\sigma$. One of these pairs concerns the Dynkin
diagram $a_8\ {\bf 3}a_5$ (of type III of theorem 5.8), the other the diagram
$a_5a_2\ {\bf 3}a_2^3$ (of type II of theorem 5.8).
This result implies that the theoretical analysis of the holes as carried out
in section 5.4 cannot be carried much further. In particular, one cannot
infer the equivalence of two holes in $\Lambda^\sigma$ from the equivalence
of their lifts in $\Lambda$, unless the automorphism group of the hole in
$\Lambda$ satisfies the assumptions of theorem 5.9.
\flexskip
 
In this section we restrict ourselves to listing the bases of the fixed point
lattices from which we can read off the volume of a fundamental region. The
complete decompositions of the fundamental regions are given in appendix A.
\flexskip
 
For $N=7$ we choose as basis of the fixed point lattice
(in units of ${1\over\sqrt8}$): 

\[
      \begin{pmatrix}-3 \\ 1 \\ 1 \\ \sqrt7\\ \sqrt7\\\sqrt7\end{pmatrix} ,
\;\;  \begin{pmatrix} 0 \\ 2 \\ 0 \\ 0     \\2\sqrt7\\ 0    \end{pmatrix} ,
\;\;  \begin{pmatrix} 2 \\ 0 \\ 0 \\2\sqrt7\\ 0     \\ 0    \end{pmatrix} ,
\;\;  \begin{pmatrix} 4 \\ 0 \\ 4 \\ 0     \\ 0     \\ 0    \end{pmatrix} ,
\;\;  \begin{pmatrix} 4 \\ 4 \\ 0 \\ 0     \\ 0     \\ 0    \end{pmatrix} ,
\;\;  \begin{pmatrix} 8 \\ 0 \\ 0 \\ 0     \\ 0     \\ 0    \end{pmatrix} 
\]
 
The fundamental volume thus is $\sqrt7^3$. If we check the table of results we
can confirm the size of the total automorphism group easily: The automorphism
group of $A_6$ has size at most 14, as this already accounts for all
possible moves of the (isolated) diagram. Hence the total group cannot be
bigger than 1176. On the other hand, there are holes of automorphism group 1.
\flexskip
 
For $N=5$ we choose as basis of the fixed point lattice
 (in units of ${1\over\sqrt8}$):
\[
      \begin{pmatrix}-3 \\ 1 \\ 1 \\ 1 \\ \sqrt5\\ \sqrt5\\ \sqrt5\\\sqrt5\end{pmatrix} ,
\;\;  \begin{pmatrix} 2 \\ 2 \\ 0 \\ 2 \\ 0     \\     0 \\2\sqrt5\\    0 \end{pmatrix} ,
\;\;  \begin{pmatrix} 2 \\ 0 \\ 2 \\ 2 \\ 0     \\2\sqrt5\\     0 \\    0 \end{pmatrix} ,
\;\;  \begin{pmatrix} 0 \\ 2 \\ 2 \\ 2 \\2\sqrt5\\     0 \\     0 \\    0 \end{pmatrix} ,
\;\;  \begin{pmatrix} 4 \\ 0 \\ 0 \\ 4 \\ 0     \\     0 \\     0 \\    0 \end{pmatrix} ,
\;\;  \begin{pmatrix} 4 \\ 0 \\ 4 \\ 0 \\ 0     \\     0 \\     0 \\    0 \end{pmatrix} ,
\;\;  \begin{pmatrix} 4 \\ 4 \\ 0 \\ 0 \\ 0     \\     0 \\     0 \\    0 \end{pmatrix} ,
\;\;  \begin{pmatrix} 8 \\ 0 \\ 0 \\ 0 \\ 0     \\     0 \\     0 \\    0 \end{pmatrix} 
 \]
 
The fundamental volume thus is $\sqrt5^4$. If we check the table of results we
can confirm the size of the total automorphism group easily: The automorphism
group of $A_4^2$ has size at most 200, as this already accounts for all
possible moves of the (isolated) diagram. Hence the total group cannot be
bigger than 14400. On the other hand, there are holes of automorphism group 1.
\flexskip 
 
The fixed point lattices in the cases $N=3$ and $N=2$ are well known to be
$K_{12}$ and $\Lambda_{16}$, respectively. Hence there is no need to give an
explicit basis here. The fundamental volumes and total automorphism groups can
be quoted from \cite{CS88}, chapter 4. The volumes are $\sqrt3^6$, and $\sqrt2^8$
respectively. (Observe that the determinant referred to in \cite{CS88} is
the square of the volume of interest here.)
  
\subsubsection{Hyperbolic Subalgebras}
  
Before we continue to list examples of hyperbolic Lie algebras we need to
reflect on the results obtained so far. We observed that the bounds are sharp
only for some cases. In the remainder of this section we will restrict
ourselves to those hyperbolic Lie algebras for which the upper bounds bear
some resemblance to the true multiplicities and at the same time represent an
improvement on known upper bounds. In section 6.3 we will return in more detail
to the question which cases are successful.
\flexskip
  
To completely analyse all the GKMs ${\cal G}_N$ and to determine all
hyperbolic subalgebras again requires tedious calculations such that we can
only report the results. However, we do not print a classification of all 
such hyperbolic subalgebras as, on the one hand, there is no obvious use 
for such a classification, on the other hand, given the classification of 
all finite and affine subalgebras as in appendix A, the classification of 
all hyperbolic subalgebras is a straightforward exercise. We restrict 
ourselves to recording the following experimental fact (though I cannot 
offer any explanation): in no ${\cal G}_N$ there exist hyperbolic subalgebras 
whose rank oversteps the value given by the following dimension formula,
$$ \hbox{rank} = {\hbox{dim}\Lambda^\sigma\over 2} +2$$
even though there do exist affine subalgebras of higher rank (such as $A_8$ 
in $N=3$, $A_6$ in $N=7$, and $D_6$ in $N=5$).
\flexskip

This means that hyperbolic Lie algebras of rank 7 or 8 have only been
identified, if at all, as subalgebras of ${\cal G}_3$ and ${\cal G}_2$,
hyperbolic algebras of ranks 9 and 10 have been found, if at all, as
subalgebras of ${\cal G}_2$. Let us now compare the bounds as obtained by
identification as subalgebra with the global bounds of theorem 1.4b, given
for rank 7-10. Note that we use the values for
roots not in $NL^*$ in the case of the subalgebras as these are the smaller
values:

\begin {table}[ht]
\renewcommand\arraystretch{1.5}
\noindent\[ 
\begin{array}{|r||r|r|r||r|r|r|}
\multicolumn{1}{|c||} {r^2} &
\multicolumn{1}{|c|} {\hbox{Rk 7}} &
\multicolumn{1}{|c|} {\hbox{Rk 8}} &
\multicolumn{1}{|c||} {{\cal G}_3} &
\multicolumn{1}{|c|} {\hbox{Rk 9}} &
\multicolumn{1}{|c|} {\hbox{Rk 10}} &
\multicolumn{1}{|c|} {{\cal G}_2} \\ \hline
 0  & 5 & 6 & 6 & 7 & 8 & 8 \\
-2  & 21 & 28   &  27   &  36 & 45   &  52   \\
-4  & 71 & 105   &  104   &  148 & 201   &  256   \\
-6  & 217 & 350   &  351   &  534 & 780   &  1122  \\
-8  & 603 & 1057   &  1080   & 1738 & 2718   &  4352   \\
-10  & 1574 & 2975   &  3107 & 5240 & 8730   &  15640   \\
-12  & 3880 & 7883   &  8424 & 14824 & 26226   &  52224   \\
-14  & 9153 & 19900 & 21762 & 39809 & 74556 & 165087 \\
-16  & 20755 & 48160 & 53976 & 102223 & 202180 & 495872   \\
\hline
\end{array}
\]
\caption{Comparison of quality of upper bounds, rank 7 to 10}
\end{table}

We see that we only obtain minimal improvements in some rare cases by
${\cal G}_3$, and worse bounds throughout by ${\cal G}_2$. Therefore we will
not list any hyperbolic Lie algebras of rank greater or equal to 7 in this
chapter, as we do not obtain improved upper bounds. Nevertheless, the
analysis of their root lattices remains valuable as it provides the key to
answering questions such as to any kind of subalgebra. Similarly, the
decomposition of the Leech lattice has been used in a number of problems since
it was established (see \cite{CS88}, chapter 25).
\flexskip
 
Let us now proceed to give an account of those hyperbolic Lie algebras which
are subalgebras of ${\cal G}_7$ and ${\cal G}_5$ such that the upper bounds
obtained from theorem 1.7 are sharp for some roots of small norm and height.
We compare these bounds with the exact values for some imaginary roots,
ordered by norm and height. The notation is as before (see section 6.2.1.2),
and I calculated the exact root multiplicities by a program based on the
Peterson recursion formula for root multiplicities (see \cite{Kac90}, p.210).
\vfil\eject

{\parindent=0pt
\centerline{\bf Hyperbolic Lie algebras of Rank 5}
\vspace{0.15in}

\begin{minipage}[t]{2.2in}
\begin{picture}(100,-60)(0,60)
\put(0,60){$H^{(5)}_1= AE_5\subset {\cal G}_7$}
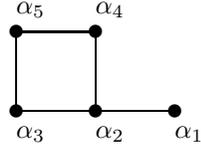
\put(80,-10){$\alpha_1$}
\put(50,-10){$\alpha_2$}
\put(20,-10){$\alpha_3$}
\put(50,37){$\alpha_4$}
\put(20,37){$\alpha_5$}
\put(20,0){\circle*{5}}
\put(50,0){\circle*{5}}
\put(80,0){\circle*{5}}
\put(20,30){\circle*{5}}
\put(50,30){\circle*{5}}
\put(20,0) {\line(1,0){30}}
\put(50,0) {\line(1,0){30}}
\put(20,30) {\line(1,0){30}}
\put(20,0) {\line(0,1){30}}
\put(50,0) {\line(0,1){30}}
\end{picture}
\end{minipage}
\begin{minipage}[t]{2.7in}
\begin{tabular}[t]{p{1.15in}c|c|c|}
 & \multicolumn{2}{c}{} & \multicolumn{1}{c}{\bf N=7} \\
 \multicolumn{1}{|c|} {root-coefficients} & norm & mult & bound\\ \hline
 \multicolumn{1}{|c|} {0, 1, 1, 1, 1} & 0 & 3 & 3\vspace{-0.015in}\\
 \multicolumn{1}{|c|} {1, 2, 2, 2, 2} & -2 & 9 & 9\vspace{-0.015in}\\
 \multicolumn{1}{|c|} {1, 3, 3, 3, 3} & -4 & 22 & 22\vspace{-0.015in}\\
 \multicolumn{1}{|c|} {2, 4, 3, 3, 2} & -4 & 22 & 22\vspace{-0.015in}\\
 \multicolumn{1}{|c|} {2, 4, 3, 3, 3} & -6 & 51 & 51\vspace{-0.015in}\\
 \multicolumn{1}{|c|} {1, 4, 4, 4, 4} & -6 & 51 & 51\vspace{-0.015in}\\
 \multicolumn{1}{|c|} {2, 5, 4, 4, 3} & -8 & 108 & 108\vspace{-0.015in}\\
 \multicolumn{1}{|c|} {2, 4, 4, 4, 4} & -8 & 108 & 108\vspace{-0.015in}\\
 \multicolumn{1}{|c|} {1, 5, 5, 5, 5} & -8 & 108 & 108\vspace{-0.015in}\\
 \multicolumn{1}{|c|} {2, 5, 4, 4, 4} & -10 & 221 & 221\vspace{-0.015in}\\
 \multicolumn{1}{|c|} {3, 6, 5, 4, 4} & -12 & 429 & 432 \\
\end{tabular}
\end{minipage}
 
\rule{4.95in}{0.01in}

\vspace{0.15in}
\centerline{\bf Hyperbolic Lie algebras of Rank 6}
\vspace{0.15in}

\begin{minipage}[t]{2.2in}
\begin{picture}(158,-80)(0,80)
\put(0,80){$H^{(6)}_1= AE_6\subset {\cal G}_5$}
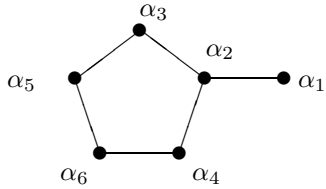
\put(55,52){$\alpha_3$}
\put(25,-10){$\alpha_6$}
\put(75,-10){$\alpha_4$}
\put(5,25){$\alpha_5$}
\put(80,37){$\alpha_2$}
\put(115,25){$\alpha_1$}
\put(40,0){\circle*{5}}
\put(40,0){\line(1,0){30}}
\put(70,0){\circle*{5}}
\put(70,0){\line(1,3){9.5}}
\put(79.5,28.5){\circle*{5}}
\put(79.5,28.5){\line(1,0){30}}
\put(109.5,28.5){\circle*{5}}
\put(79.5,28.5) {\line(-4,3){24.5}}
\put(55,46.6){\circle*{5}}
\put(55.4,46.6) {\line(-4,-3){24.5}}
\put(30.5,28.5){\circle*{5}}
\put(30.5,28.5) {\line(1,-3){9.5}}
\end{picture}
\end{minipage}
\begin{minipage}[t]{2.7in}
\begin{tabular}[t]{p{1.15in}c|c|c|}
 & \multicolumn{2}{c}{} & \multicolumn{1}{c}{\bf N=5} \\
 \multicolumn{1}{|c|} {root-coefficients} & norm & mult & bound\\ \hline
 \multicolumn{1}{|c|} {0, 1, 1, 1, 1, 1} & 0 & 4 & 4\vspace{-0.015in}\\
 \multicolumn{1}{|c|} {1, 2, 2, 2, 2, 2} & -2 & 14 & 14\vspace{-0.015in}\\
 \multicolumn{1}{|c|} {2, 4, 3, 3, 2, 2} & -4 & 40 & 40\vspace{-0.015in}\\
 \multicolumn{1}{|c|} {1, 3, 3, 3, 3, 3} & -4 & 40 & 40\vspace{-0.015in}\\
 \multicolumn{1}{|c|} {2, 4, 3, 3, 3, 3} & -6 & 105 & 105\vspace{-0.015in}\\
 \multicolumn{1}{|c|} {1, 4, 4, 4, 4, 4} & -6 & 105 & 105\vspace{-0.015in}\\
 \multicolumn{1}{|c|} {2, 5, 4, 4, 3, 3} & -8 & 252 & 256\vspace{-0.015in}\\
 \multicolumn{1}{|c|} {2, 4, 4, 4, 4, 4} & -8 & 251 & 256\vspace{-0.015in}\\
 \multicolumn{1}{|c|} {2, 5, 4, 4, 4, 4} & -10 & 572 & 590\vspace{-0.015in}\\
 \multicolumn{1}{|c|} {3, 6, 5, 4, 4, 3} & -10 & 574 & 590\vspace{-0.015in}\\
 \multicolumn{1}{|c|} {3, 6, 5, 4, 3, 4} & -10 & 574 & 590 \\
\end{tabular}
\end{minipage}
 
\rule{4.95in}{0.01in}
\vspace{0.25in}
 
\begin{minipage}[t]{2.2in}
\begin{picture}(158,-50)(0,50)
\put(0,50){$H^{(6)}_6=DE_6\subset {\cal G}_5$}
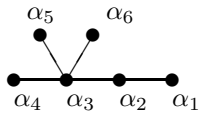
\put(70,-10){$\alpha_1$}
\put(50,-10){$\alpha_2$}
\put(30,-10){$\alpha_3$}
\put(10,-10){$\alpha_4$}
\put(15,23){$\alpha_5$}
\put(45,23){$\alpha_6$}
\put(10,0){\circle*{5}}
\put(30,0){\circle*{5}}
\put(50,0){\circle*{5}}
\put(70,0){\circle*{5}}
\put(20,16.7){\circle*{5}}
\put(40,16.7){\circle*{5}}
\put(10,0){\line(1,0){60}}
\put(30,0){\line(3,5){10}}
\put(30,0){\line(-3,5){10}}
\end{picture}
\end{minipage}
\begin{minipage}[t]{2.7in}
\begin{tabular}[t]{p{1.15in}c|c|c|}
 & \multicolumn{2}{c}{} & \multicolumn{1}{c}{\bf N=5} \\
 \multicolumn{1}{|c|} {root-coefficients} & norm & mult & bound\\ \hline
 \multicolumn{1}{|c|} {0, 1, 2, 1, 1, 1} & 0 & 4 & 4\vspace{-0.015in}\\
 \multicolumn{1}{|c|} {1, 2, 4, 2, 2, 2} & -2 & 14 & 14\vspace{-0.015in}\\
 \multicolumn{1}{|c|} {1, 3, 6, 3, 3, 3} & -4 & 40 & 40\vspace{-0.015in}\\
 \multicolumn{1}{|c|} {2, 4, 6, 3, 3, 2} & -4 & 40 & 40\vspace{-0.015in}\\
 \multicolumn{1}{|c|} {2, 4, 6, 3, 2, 3} & -4 & 40 & 40\vspace{-0.015in}\\
 \multicolumn{1}{|c|} {2, 4, 6, 2, 3, 3} & -4 & 40 & 40\vspace{-0.015in}\\
 \multicolumn{1}{|c|} {2, 4, 6, 3, 3, 3} & -6 & 105 & 105\vspace{-0.015in}\\
 \multicolumn{1}{|c|} {1, 4, 8, 4, 4, 4} & -6 & 105 & 105\vspace{-0.015in}\\
 \multicolumn{1}{|c|} {2, 5, 8, 4, 4, 3} & -8 & 252 & 256\vspace{-0.015in}\\
 \multicolumn{1}{|c|} {2, 5, 8, 4, 3, 4} & -8 & 252 & 256\vspace{-0.015in}\\
 \multicolumn{1}{|c|} {2, 5, 8, 3, 4, 4} & -8 & 252 & 256 \\
\end{tabular}
\end{minipage}

\rule{4.95in}{0.01in}
\vspace{-0.05in}

\begin{table}[ht]
\caption{Root multiplicities of some hyperbolic Lie algebras of rank 5 and 6}
\end{table}
} 
 
\section{Conclusions}
 
This final section tries to put into perspective the results obtained in
the previous sections of chapter 6. In section 6.2 we observed that the root
multiplicities of ${\cal G}_N$ provided good upper bounds for some hyperbolic
Lie algebras while for others they did not bear any resemblance to the
correct multiplicities.
\flexskip
 
We recall that in the cases $N=5,7,11,23$ there existed examples where the
upper bounds obtained in this work represented
a significant improvement on the global bounds of theorem 1.4b. Let us
therefore turn to these cases. This automatically restricts us to simply
laced algebras. Still, we did not obtain useful, sharp upper bounds for each
hyperbolic subalgebra contained in one of the ${\cal G}_N$. If we recall
table 6.3 of $N=11$ given in section 6.2.2.2 we found that only the rank 4
algebra provided useful bounds whereas the bounds for the rank 3 algebra were
far off the true values.  We recall from section 6.2.3 that the situation is
similar for the remaining $N$.
\flexskip
 
One necessary condition for useful upper bounds can be found
experimentally. If we consider the numerical evidence for those hyperbolic Lie
algebras which contain more than one affine subalgebra we will find that, even
for roots of the same small negative norm, the root multiplicities vary
considerably. Thus we will never be able to provide sharp bounds which
depend on norm only. We therefore now concentrate on the following case:
 
We consider a simply laced hyperbolic Lie algebra ${\cal A}$ with a unique
affine subalgebra ${\cal A}_0$, thus allowing the definition of a level, as
in section 6.2.1.2, above.
If we consider numerical examples of successful and useless upper bounds we
observe another necessary condition for the construction of useful bounds:
They must be sharp for at least some norm 0 vectors.
Therefore, we will now determine the exact root multiplicity of norm 0 vectors
in hyperbolic Lie algebras and then ask when will this number be equal to the
multiplicity of those norm 0 vectors in ${\cal G}_N$ which are not elements
of $NL^*$ (for the notation, see the corollary to theorem 1.7).
 
We recall from formula (5.9) that for every simple affine Lie algebra there
exists a unique norm 0 vector $\delta=\sum_1^n n_i\alpha_i$ where the $n_i$ are
positive integers whose greatest common divisor be 1. For the remainder
of this section let $\delta({\cal A}_0)$ denote this norm 0 vector of the
affine algebra ${\cal A}_0$.
\flexskip
 
We quote three results of \cite{Kac90}:

\begin{proposition} [\cite{Kac90}, Prop. 5.10c] Let ${\cal A}$ be a Lie
algebra of finite, affine, or hyperbolic type. Let $Q$ denote the lattice
spanned by the simple roots. Then the set of all imaginary roots is
\[
\{ \; \; \alpha\in Q - \{0\} \; \; \; | \; \; \; \alpha^2 \leq 0 \; \; \}.
\]
\end{proposition}
 
\begin{remark} This is, in fact, an equivalent characterization of the 
three types of Lie algebras among ordinary Kac-Moody algebras. For every 
other Kac-Moody algebra, there exist vectors of the root lattice of negative 
norm which are not roots. Note, however, that the series of GKMs ${\cal G}_N$ 
constructed in this work provides examples of further Lie algebras where all 
imaginary root lattice vectors are roots.
\end{remark}
 
\begin{proposition} [\cite{Kac90}, Prop. 5.7] Let ${\cal A}$ be
symmetrizable and let it have Cartan matrix $C$. Let
the Weyl group be denoted $W$. A root $\alpha$ is isotropic (i.e. $\alpha^2=0$)
if and only if it is $W$-equivalent to an imaginary root $\beta$ such that
${\rm supp}\beta$ is a subdiagram of affine type of the Dynkin diagram of
${\cal A}$. Let the corresponding affine subalgebra be called ${\cal A}_0$.
Then $\beta=k\delta({\cal A}_0)$.
\end{proposition}
  
\begin{proposition} [\cite{Kac90}, Cor. 7.4] Let ${\cal A}_0$ be a simply
laced affine Lie algebra of $n$ simple roots. Then the multiplicity of every
imaginary root of ${\cal A}_0$ is $n-1$.
\end{proposition}
 
We rephrase this for the cases of interest to us.
 
\begin{theorem} Suppose ${\cal A}$ is a simply laced hyperbolic Lie
algebra of rank $n$ with a unique affine subalgebra ${\cal A}_0$. Then the
multiplicity of every isotropic $\alpha\in Q-\{0\}$ is $n-2$.
\end{theorem}
 
\begin{proof} The Dynkin diagram of the affine subalgebra ${\cal A}_0$ must be
the Dynkin diagram of ${\cal A}$ with one simple root removed, hence it will
have $n-1$ simple roots. Now consider $\alpha\ne 0$, such that $\alpha^2=0$. By
proposition 6.1 it is a root of ${\cal A}$. By proposition 6.2 it
is $W$-equivalent to some $\beta$ which corresponds to some affine subalgebra.
Now there is only one such subalgebra, ${\cal A}_0$, of $n-1$ simple roots.
Hence, by proposition 6.3 the multiplicity of $\beta$ as a root of
${\cal A}_0$ is $n-2$. Hence its multiplicity as a root of ${\cal A}$ is $n-2$.
This in turn implies that the multiplicity of $\alpha$ is $n-2$.
\end{proof}
 
Theorem 1.1 of \cite{Bor90a} shows that, for any root $r$ of ${\cal G}_N$ 
with $r^2>0$, there exists a unique representation as a sum of positive 
simple roots. As we will see below, no such uniqueness holds for roots $r$ 
such that $r^2 \le 0$. Instead, we have the following situation. 
Let ${\cal A}$ be a hyperbolic Lie algebra with unique affine subalgebra
${\cal A}_0$ such that ${\cal A}$ is subalgebra of one of the GKMs
${\cal G}_N$. As before, let ${\cal R}=\{r_i\}$ denote the set of real simple
roots of ${\cal G}_N$. The vector $\delta=\delta({\cal A}_0)$ can be
represented as $\delta=\sum_j n_{i_j} r_{i_j}$ where the $r_{i_j}$,
$j=1,\dots,n$ are the simple roots of ${\cal A}_0$. The vector $\delta$ lies
within the (closure of the) Weyl chamber of the affine Lie algebra
${\cal A}_0$. It also lies within the (closure of the) Weyl chamber of the
GKM ${\cal G}_N$. Suppose, it did not. Then there would exist a root $r$ of
${\cal G}_N$ such that the reflection $\phi_r$ takes $\delta$ to a root of
smaller height. As $\delta$ lies within the Weyl chamber of ${\cal A}_0$ $r$
must be distinct from any of the simple roots $r_{i_j}$, $j=1,\dots,n$. Then
$\phi_r(\delta)$ is a root which is neither positive nor negative, a
contradiction. We conclude that if $\delta=\sum_j n_{i_j} r_{i_j} =
\sum_k n_{i_k} r_{i_k}$ then by proposition 6.2 not only the collection of
$r_{i_j}$ but also the collection of $r_{i_k}$ corresponds to an affine
subalgebra of ${\cal G}_N$.

The multiplicity of $\delta$ as a root of
${\cal G}_N$ and the multiplicity of $\delta$ as a root of the hyperbolic Lie
algebra will be equal if and only if the representation $\delta= \sum n_{i_j}
r_{i_j}$ is unique in ${\cal G}_N$.
We recall formula (5.13) for the generalized centre $c$ of an affine hole
and formula (5.10) for the relation $\delta^\vee= {1\over d}\nu^{-1}(\delta)$,
using all notation as in chapter 5. Then
$$(c,1,*) = {\delta^\vee\over h^\vee} = {\nu^{-1}(\delta)\over dh^\vee}.$$
Thus we can identify the centre $c$
from $\delta$ alone, using the fact that $(c,1,*)$ must have height 1.
This can be reformulated as follows:
If $\delta=\sum_j n_{i_j} r_{i_j} = \sum_k n_{i_k} r_{i_k}$ within the GKM
${\cal G}_N$, then the representatives of $r_{i_j}$ and those of $r_{i_k}$
are associated to the same affine hole of ${\cal R}$. To confirm that $\delta$
can only be represented in a unique way it therefore suffices to check the
the products $dh^\vee$ of the
components of the affine hole in ${\cal R}$ determined by $\delta$.
 
\begin{theorem} Suppose that a simply laced hyperbolic Lie algebra
${\cal A}$ with unique affine subalgebra ${\cal A}_0$ is contained in one of
the GKMs ${\cal G}_N$ constructed in theorem 1.6. Consider the affine hole
$H$ of ${\cal R}$ containing the elements of ${\cal R}$ representing
$\delta({\cal A}_0)$. If the product $d({\cal A}_0)h^\vee({\cal A}_0)$ is
smaller than the product $dh^\vee$ of any other
component of $H$ then the upper bound for the multiplicity of $\delta$ in
${\cal A}$ is sharp. 
\end{theorem} \theoremproven

\begin{remark} By theorem 6.1, the multiplicity of any norm 0 vector in ${\cal A}$
equals that of $\delta$. Suppose we obtain a sharp upper bound for the
multiplicity of $\delta$. If the upper bounds depend on norm only they will be
sharp for every norm 0 vector. However, in the cases
${\cal G}_N$ at hand, we obtain discrepancies for norm 0 vectors which
are elements of $NL^*$.
\end{remark}
 
\begin{example} It is now straightforward to understand the quality of the upper
bounds obtained in section 6.2. For example, $A_1$ in $A_1\ {\bf 23}A_1$
satisfies the conditions of theorem 6.2, $A_1$ in $A_1^2\ {\bf 11}A_1^2$ does
not. Correspondingly, the bounds for $AE_3$ as a subset of ${\cal G}_{23}$ are
useful, those for $AE_3$ as a subset of ${\cal G}_{11}$ are not.
 
Similarly, $A_2$ in $A_2\ {\bf 11}A_2$ satisfies the condition. Hence the bounds
for $AE_4$ obtained from $N=11$ are useful. As a last example, the bounds for
$H_{71}^{(3)}$ are useless, again because $A_1^2\ {\bf 11}A_1^2$ does not
satisfy the conditions of the theorem.
\end{example} 
 
Let us look at the results of this work from the point of view that it
provides a strategy to calculate sharp upper bounds for the root multiplicities
of some hyperbolic Lie algebras. We are now in the position to conjecture
how far it may be possible to generalize this strategy. Suppose we consider
a hyperbolic Lie algebra ${\cal A}$ of small rank, that is, the global upper 
bounds as provided by theorem 1.4b are substantially greater than the true
values. We will have to find a suitable automorphism $\sigma$ of the Leech
lattice such that the root system ${\cal R}$ contains ${\cal A}$ as a
subalgebra and such that at the same time the conditions of theorem 6.2 are
satisfied. Recall that we can search for such $\sigma$ easily as we understand
the action of $\sigma$ on Dynkin diagrams directly (chapter 5.4). Note that
theorem 6.2 provides some idea of the quality of the bounds before we begin to
calculate denominator formulas.
\flexskip
 
The twisted denominator formula will then describe a Lie superalgebra,
as indicated in the final remark of section 1.6.
We will be able to obtain upper bounds for the root multiplicities of
${\cal A}$ if the twisted denominator formula describes in fact either a
GKM or a Lie superalgebra with strictly alternating multiplicities.
Borcherds states some conjectures about the properties of such automorphisms
in chapter 6 of \cite{Bor90b}. 
\flexskip
 
For all hyperbolic Lie algebras of rank 7 to 10, the root multiplicities for 
some roots of small negative norm are collected in appendix B, and compared 
with both theorem 1.4b and the results of this work. 
One remarkable result of the calculations concerns the algebra $T_{4,3,3}$
which contains roots of multiplicities both greater and smaller than $p_6$,
disproving a number of conjectures. The data shows furthermore
that for simply laced hyperbolic Lie algebras of 
high rank the upper bound of theorem 1.4b already proves useful, which is
one of the reasons why the new upper bounds could not improve on them.
However, in the cases of those $BE_n$ and $CE_n$ which we identified as
subalgebras of the ${\cal G}_N$, the new bounds
are substantially greater than the true values, and theorem 1.4b does not 
apply.
We note that the bounds for $CE_n$ are not even sharp for the some of the 
norm 0 vectors. A look on the numerical results
for the $BE_n$ shows that, even though the bounds for norm 0 are sharp, those
for any other norm are useless. This provides a note of caution: Constructing
sharp bounds for norm 0 vectors is necessary if we want to obtain useful upper
bounds but not sufficient. The question of useful upper bounds remains wide
open in these cases.
\flexskip\flexskip

{\bf Acknowledgements:} I wish to thank Richard Borcherds, my research 
supervisor, for suggesting many of the problems discussed in this work
and for the advice and encouragement he provided throughout my research.
Furthermore, I would like to thank Simon Norton for his patient
explanations of the ATLAS, and Elizabeth Jurisich for clarifications
regarding the nature of specializations. Finally, I would like to thank
the referee for valuable comments on the first version of this paper.

\backmatter

\chapter*{Appendix A}

This appendix contains the complete tables of affine and finite diagrams for
the sets ${\cal R}$ as described in chapter 5. The conventions
denoting Lie algebras are standard with capital letters denoting affine Lie
algebras and small letters denoting finite ones. The presence of a bold-faced
integer ${\mathbf N}$ in front of a component indicates that this component
consists of long roots (the ratio of the norms of being N to 1). \flexskip

Each table lists first all affine and then all finite diagrams.
The types of diagrams are in the first instance ordered by their indices.
They are ordered by natural order of partitions. This convention is adopted
even in those cases (like $\Delta^{(2)}$ or $\Delta^{(3)}$) where the index does
not coincide with the rank. For a fixed set of indices, the order is
alphabetic in the components. Note that the alphabet is extended such that any
component ${\mathbf N}\Delta$ immediately follows $\Delta$. \flexskip

The second column of each table contains the size of the automorphism group of
the respective diagrams. The total number of representatives of any diagram
within a fundamental volume of the respective fixed point lattice is the size
of the automorphism group of the lattice divided by the automorphism group of
the diagram. The third column contains the unit volume of any one
representative of a diagram. The necessary calculations are described in
chapter 5.3. The total volume is the product of the
unit volume by the total number of representatives of a diagram. The volume
formula (theorem 5.6) states that the total volumes of the complete list of
types must add up to the volume of the fundamental region. This equals
$\sqrt N^M$
as follows from the explicit bases provided in chapter 6.2. Here, $N=$ 23, 11,
7, 5, 3, 2, and $M=24/(N+1)$. The unit volume and total volume are given in
units indicated at the top of each table. These units are chosen purely for
convenience.

{\parindent=0pt
\section*{N=23}

Order of the automorphism group of the fixed point lattice: 2.

There are 1 type of affine diagram and 2 types of finite diagrams.\flexskip

\begin{center}
\begin{tabular}[t]{lrrrl}
\\type of diagram $\mskip 40mu$ 
   & $\mskip 30mu$ group 
   & $\mskip 30mu$ unit volume 
   & $\mskip 20mu$ total volume 
   & $\mskip 10mu$ \\ 
& & $\scriptstyle{[u=2!^{-1}\sqrt{23}^{-1}]}$ 
   & $\scriptstyle{[u=\sqrt{23}^{-1}]}$ \\
\hline
$A_1 \ \ {\bf 23}A_1$ & 2 & 8 & 4\\
\vspace{-0.15in} \\
$a_2 \ \ {\bf 23} a_1$ & 1 & 9 & 9\\
${a_1}^2 \ \ {\bf 23}a_1$ & 1 & 10 & 10\\
\hline
&&& 23\\
\\
\end{tabular}
\end{center}
Total volume of 3 types (units of 1) = $23/\sqrt{23}= \sqrt{23}^1$.
\vfil\eject

\section*{N=11}

Order of the automorphism group of the fixed point lattice: $2^3\ 3^1$ = 24.

There are 2 types of affine diagrams and 6 types of finite diagrams.\flexskip

\begin{center}
\begin{tabular}[t]{lrrrl}
\\type of diagram $\mskip 40mu$ 
   & $\mskip 30mu$ group 
   & $\mskip 30mu$ unit volume 
   & $\mskip 20mu$ total volume 
   & $\mskip 10mu$ \\ 
& & $\scriptstyle{[u=(11*4!)^{-1}]}$ 
   & $\scriptstyle{[u=11^{-1}]}$ \\
\hline
\vspace{-0.1in} \\
$A_2 \ \ {\bf 11}A_2$ &  6& 27& 4&$\negthickspace\negthickspace\negthickspace.5$\\
${A_1}^2\ \ {\bf 11}{A_1}^2$ & 8& 64& 8&\\
\vspace{-0.1in} \\
$a_3\ \ {\bf 11}a_2$ & 1& 18& 18&\\
$a_3\ \ a_1 \ \ {\bf 11}a_1$ & 1& 44& 44&\\
$a_3 \ \ {\bf 11}{a_1}^2$ & 2& 16& 8&\\
$a_2\ \ {\bf 11}a_2\ \ a_1$ & 2& 21& 10&$\negthickspace\negthickspace\negthickspace.5$\\
$a_2\ \ a_1 \ \ {\bf 11}{a_1}^2$ & 1& 18& 18&\\
${a_1}^3\ \ {\bf 11}{a_1}^2$ & 2& 20& 10&\\
\hline
&&& 121\\
\\
\end{tabular}
\end{center}
Total volume of 8 types (units of 1) $= 121/11=\sqrt{11}^2$.
\flexskip\flexskip

\section*{N=7}

Order of the automorphism group of the fixed point lattice: $2^5\ 3^1\ 7^1$ =
672.

There 3 types of affine diagrams and 11 types of finite diagrams. \flexskip

\begin{center}
\begin{tabular}[t]{lrrrl}
\\type of diagram $\mskip 40mu$ 
   & $\mskip 30mu$ group 
   & $\mskip 30mu$ unit volume 
   & $\mskip 20mu$ total volume 
   & $\mskip 10mu$ \\ 
  & & $\scriptstyle{[u=(7^2*6!)^{-1}\sqrt7]}$ & $\scriptstyle{[u={1\over315}\sqrt7]}$ \\
\hline
\vspace{-0.1in} \\
$A_6$ & 14& 343& 147\\
$A_3 \ \ {\bf 7}A_3$ & 8& 64& 48\\
${A_1} ^3 \ \ {\bf 7} {A_1} ^3$ & 48& 512& 64\\
\vspace{-0.1in} \\
$a_6\ \ {\bf 7}a_1$ & 2& 147& 441\\
$a_5\ \ {\bf 7}a_2$ & 2& 63& 189\\
$a_4\ \ {\bf 7}a_3$ & 1& 30& 180\\
$d_4\ \ {\bf 7}a_3$ & 2& 28& 84\\
$a_4\ \ a_1\ \ {\bf 7}{a_1}^2$ & 1& 70& 420\\
$d_4\ \ {\bf 7}{a_1}^3$ & 6& 28& 28\\
$a_3\ \ {\bf 7}a_3\ \ a_1$ & 2& 36& 108\\
$a_3\ \ a_1\ \ {\bf 7}{a_1}^3$ & 2& 32& 96\\
$a_2\ \ {a_1}^3\ \ {\bf 7}{a_1}^2$ & 2& 84& 252\\
$a_2\ \ {a_1}^2\ \ {\bf 7}{a_1}^3$ & 2& 36& 108\\
${a_1}^4\ \ {\bf 7}{a_1}^3$ & 6& 40& 40\\
\hline
&&& 2205\\
\\
\end{tabular}
\end{center}
Total volume of 14 types (units of 1) $= 2205 {\sqrt7\over315}=\sqrt{7}^3$.

\vfil\eject
\section*{N=5}

Order of the automorphism group of the fixed point lattice:
$2^7\ 3^2\ 5^2$ = 28800.

There are 6 types of affine diagrams and 24 types of finite diagrams.

\flexskip

\begin{center}
\begin{tabular}[t]{lrrrl}
\\type of diagram $\mskip 40mu$ 
   & $\mskip 30mu$ group 
   & $\mskip 30mu$ unit volume 
   & $\mskip 20mu$ total volume 
   & $\mskip 10mu$ \\ 
& & $\scriptstyle{[u=(25*8!)^{-1}]}$ 
   & $\scriptstyle{[u=1/210]}$ \\
\hline
\vspace{-0.1in} \\
$D_6\ \ {\bf 5}{A_1}^2$ & 8& 800& 600\\
${A_4}^2 $ & 200& 3125& 93&$\negthickspace\negthickspace\negthickspace.75$\\
${A_4}\ \ {\bf 5} {A_4} $ & 20& 125& 37&$\negthickspace\negthickspace\negthickspace.5$\\
$D_4\ \ {\bf 5}D_4$ & 24& 144& 36\\
\end{tabular}
\end{center}
Total volume of 30 types (units of 1) $=5250 /210=25=\sqrt{5}^4$.
\vfil\eject

\section*{N=3}

Order of the automorphism group of the fixed point lattice:
$2^{10} \ 3^7 \ 5^1\ 7^1$ = 78382080.

There are 15 types of affine diagrams and 93 types of finite diagrams.
\flexskip

\begin{center}
\begin{tabular}[t]{lrrrl}
\\type of diagram $\mskip 40mu$ 
   & $\mskip 30mu$ group 
   & $\mskip 30mu$ unit volume 
   & $\mskip 20mu$ total volume 
   & $\mskip 10mu$ \\ 
& & $\scriptstyle{[u=(27*12!)^{-1}]}$ & $\scriptstyle{[u=1/118800]}$ \\
\hline
\vspace{-0.1in} \\
${A_8}\ \ {\bf 3}{A_2}^{2}\ \ $
 & 36&          2187&        43740\\
${E_7}\ \ {\bf 3}{A_5}\ \ $
 & 4&           648&       116640\\
${D_7}\ \ {\bf 3}{A_3}\ \ {G_2}\ \ $
 & 8&          2304&       207360\\
${A_6}\ \ {\bf 3}{A_6}\ \ $
 & 14&           343&        17640\\
\end{tabular}
\vfil\eject

\begin{tabular}[t]{lrrrl}
\\type of diagram $\mskip 40mu$ 
   & $\mskip 30mu$ group 
   & $\mskip 30mu$ unit volume 
   & $\mskip 20mu$ total volume 
   & $\mskip 10mu$ \\ 
& & $\scriptstyle{[u=(27*12!)^{-1}]}$ & $\scriptstyle{[u=1/118800]}$ \\
\hline
\vspace{-0.1in} \\
${d_6}\ \ {\bf 3}{a_5}\ \ {g_2}\ \ $
 & 2&            54&        19440\\
${e_6}\ \ {\bf 3}{a_5}\ \ {g_2}\ \ $
 & 4&            51&         9180\\
${\bf 3}{e_6}\ \ {a_5}\ \ {a_2}\ \ $
 & 4&           117&        21060\\
${\bf 3}{d_6}\ \ {d_5}\ \ {a_1}^{2}\ \ $
 & 2&           112&        40320\\
\end{tabular}
\vfil\eject

\begin{tabular}[t]{lrrrl}
\\type of diagram $\mskip 40mu$ 
   & $\mskip 30mu$ group 
   & $\mskip 30mu$ unit volume 
   & $\mskip 20mu$ total volume 
   & $\mskip 10mu$ \\ 
  & & $\scriptstyle{[u=(27*12!)^{-1}]}$ & $\scriptstyle{[u=1/118800]}$ \\
\hline
\vspace{-0.1in} \\
${d_4}\ \ {a_3}\ \ {\bf 3}{a_3}\ \ {a_1}\ \ {\bf 3}{a_1}^{2}\ \ $
 & 4&           240&        43200\\
${\bf 3}{d_4}\ \ {a_3}\ \ {g_2}\ \ {a_1}^{3}\ \ {\bf 3}{a_1}\ \ $
 & 6&           128&        15360\\
${d_4}\ \ {a_3}\ \ {a_1}\ \ {\bf 3}{a_1}^{5}\ \ $
 & 12&           288&        17280\\
${a_4}\ \ {g_2}^{4}\ \ {\bf 3}{a_1}\ \ $
 & 24&            45&         1350\\
\end{tabular}
\end{center}
Total volume of 108 types (units of 1) $= 3207600/118800=27=\sqrt{3}^6$.
\vfil \eject

\section*{N=2}

Order of the automorphism group of the fixed point lattice:
$2^{21}\ 3^5\ 5^2\ 7^1$ = 89181388800.

There are 40 types of affine diagrams and 435 types of finite diagrams.
\flexskip

\begin{center}
\begin{tabular}[t]{lrrrl}
\\type of diagram $\mskip 40mu$ 
   & $\mskip 30mu$ group 
   & $\mskip 30mu$ unit volume 
   & $\mskip 20mu$ total volume 
   & $\mskip 10mu$ \\ 
& & $\scriptstyle{[u=(16*16!)^{-1}]}$ & $\scriptstyle{[u=1891890^{-1}]}$ \\
\hline
\vspace{-0.1in} \\
${A_{17}^{(2)}}\ \ {\bf 2}{E_{7}}\ \ $
 & 2&           648&       163296\\
${A_{15}^{(2)}}\ \ {D_{9}^{(2)}}\ \ $
 & 2&           512&       129024\\
${A_{11}^{(2)}}\ \ {D_{7}^{(2)}}\ \ {E_6^{(2)}}\ \ $
 & 2&          3456&       870912\\
${D_{10}^{(2)}}\ \ {E_{7}}\ \ $
 & 2&           648&       163296\\
\end{tabular}
\vfil\eject

\begin{tabular}[t]{lrrrl}
\\type of diagram $\mskip 40mu$ 
   & $\mskip 30mu$ group 
   & $\mskip 30mu$ unit volume 
   & $\mskip 20mu$ total volume 
   & $\mskip 10mu$ \\ 
& & $\scriptstyle{[u=(16*16!)^{-1}]}$ & $\scriptstyle{[u=1891890^{-1}]}$ \\
\hline
\vspace{-0.1in} \\
${A_{3}}^{4}\ \ {\bf 2}{A_{1}}^{4}\ \ $
 & 24576&       1048576&        21504\\
${A_{2}}^{4}\ \ {\bf 2}{A_{2}}^{4}\ \ $
 & 15552&        531441&        17222&$\negthickspace\negthickspace\negthickspace.625$\\
${A_{1}}^{16}\ \ $
 & 660602880&     268435456&          204&$\negthickspace\negthickspace\negthickspace.8$\\
${A_{1}}^{8}\ \ {\bf 2}{A_{1}}^{8}\ \ $
 & 2752512&      16777216&         3072\\
\end{tabular}
\vfil\eject

\begin{tabular}[t]{lrrrl}
\\type of diagram $\mskip 40mu$ 
   & $\mskip 30mu$ group 
   & $\mskip 30mu$ unit volume 
   & $\mskip 20mu$ total volume 
   & $\mskip 10mu$ \\ 
& & $\scriptstyle{[u=(16*16!)^{-1}]}$ & $\scriptstyle{[u=1891890^{-1}]}$ \\
\hline
\vspace{-0.1in} \\
${d_{8}}\ \ {\bf 2}{d_{8}}\ \ {a_{1}}\ \ $
 & 1&           116&        58464\\
${e_{8}}\ \ {\bf 2}{e_{8}}\ \ {a_{1}}\ \ $
 & 1&            61&        30744\\
${b_{8}}\ \ {a_{7}}\ \ {a_{2}}\ \ $
 & 1&           192&        96768\\
${b_{8}}\ \ {e_{7}}\ \ {a_{2}}\ \ $
 & 1&           114&        57456\\
\end{tabular}
\vfil\eject

\begin{tabular}[t]{lrrrl}
\\type of diagram $\mskip 40mu$ 
   & $\mskip 30mu$ group 
   & $\mskip 30mu$ unit volume 
   & $\mskip 20mu$ total volume 
   & $\mskip 10mu$ \\ 
& & $\scriptstyle{[u=(16*16!)^{-1}]}$ & $\scriptstyle{[u=1891890^{-1}]}$ \\
\hline
\vspace{-0.1in} \\
${c_{8}}\ \ {f_{4}}^{2}\ \ {\bf 2}{a_{1}}\ \ $
 & 2&            76&        19152\\
${d_{8}}\ \ {\bf 2}{a_{4}}\ \ {b_{4}}\ \ {a_{1}}\ \ $
 & 1&           180&        90720\\
${d_{8}}\ \ {b_{4}}^{2}\ \ {\bf 2}{a_{1}}\ \ $
 & 2&           120&        30240\\
${c_{8}}\ \ {\bf 2}{a_{4}}\ \ {c_{3}}\ \ {a_{2}}\ \ $
 & 1&           240&       120960\\
\end{tabular}
\vfil\eject

\begin{tabular}[t]{lrrrl}
\\type of diagram $\mskip 40mu$ 
   & $\mskip 30mu$ group 
   & $\mskip 30mu$ unit volume 
   & $\mskip 20mu$ total volume 
   & $\mskip 10mu$ \\ 
  & & $\scriptstyle{[u=(16*16!)^{-1}]}$ & $\scriptstyle{[u=1891890^{-1}]}$ \\
\hline
\vspace{-0.1in} \\
${a_{7}}\ \ {\bf 2}{d_{5}}\ \ {f_{4}}\ \ {a_{1}}\ \ $
 & 2&           152&        38304\\
${\bf 2}{a_{7}}\ \ {d_{5}}\ \ {a_{4}}\ \ {\bf 2}{a_{1}}\ \ $
 & 2&           200&        50400\\
${\bf 2}{a_{7}}\ \ {d_{5}}\ \ {c_{4}}\ \ {a_{1}}\ \ $
 & 2&           136&        34272\\
${c_{7}}\ \ {\bf 2}{a_{5}}\ \ {f_{4}}\ \ {a_{1}}\ \ $
 & 1&            78&        39312\\
\end{tabular}
\vfil\eject

\begin{tabular}[t]{lrrrl}
\\type of diagram $\mskip 40mu$ 
   & $\mskip 30mu$ group 
   & $\mskip 30mu$ unit volume 
   & $\mskip 20mu$ total volume 
   & $\mskip 10mu$ \\ 
  & & $\scriptstyle{[u=(16*16!)^{-1}]}$ & $\scriptstyle{[u=1891890^{-1}]}$ \\
\hline
\vspace{-0.1in} \\
${c_{6}}^{2}\ \ {\bf 2}{d_{5}}\ \ $
 & 2&            88&        22176\\
${c_{6}}\ \ {\bf 2}{d_{6}}\ \ {c_{5}}\ \ $
 & 1&            64&        32256\\
${c_{6}}\ \ {\bf 2}{e_{6}}\ \ {a_{5}}\ \ $
 & 2&           138&        34776\\
${\bf 2}{d_{6}}\ \ {e_{6}}\ \ {c_{5}}\ \ $
 & 1&           108&        54432\\
\end{tabular}
\vfil\eject

\begin{tabular}[t]{lrrrl}
\\type of diagram $\mskip 40mu$ 
   & $\mskip 30mu$ group 
   & $\mskip 30mu$ unit volume 
   & $\mskip 20mu$ total volume 
   & $\mskip 10mu$ \\ 
& & $\scriptstyle{[u=(16*16!)^{-1}]}$ & $\scriptstyle{[u=1891890^{-1}]}$ \\
\hline
\vspace{-0.1in} \\
${d_{6}}\ \ {\bf 2}{a_{5}}\ \ {c_{4}}\ \ {a_{1}}\ \ {\bf 2}{a_{1}}\ \ $
 & 1&           216&       108864\\
${d_{6}}\ \ {b_{5}}\ \ {b_{4}}\ \ {a_{1}}^{2}\ \ $
 & 1&           152&        76608\\
${d_{6}}\ \ {b_{5}}\ \ {c_{4}}\ \ {a_{1}}\ \ {\bf 2}{a_{1}}\ \ $
 & 1&           208&       104832\\
${a_{6}}\ \ {\bf 2}{a_{5}}\ \ {b_{3}}\ \ {c_{3}}\ \ $
 & 1&           126&        63504\\
\end{tabular}
\vfil\eject

\begin{tabular}[t]{lrrrl}
\\type of diagram $\mskip 40mu$ 
   & $\mskip 30mu$ group 
   & $\mskip 30mu$ unit volume 
   & $\mskip 20mu$ total volume 
   & $\mskip 10mu$ \\ 
& & $\scriptstyle{[u=(16*16!)^{-1}]}$ & $\scriptstyle{[u=1891890^{-1}]}$ \\
\hline
\vspace{-0.1in} \\
${d_{6}}\ \ {c_{4}}\ \ {\bf 2}{a_{3}}\ \ {a_{1}}^{2}\ \ {\bf 2}{a_{1}}^{2}\ \ $
 & 2&           448&       112896\\
${d_{6}}\ \ {\bf 2}{d_{4}}\ \ {c_{3}}\ \ {a_{1}}^{2}\ \ {\bf 2}{a_{1}}^{2}\ \ $
 & 2&           320&        80640\\
${c_{6}}\ \ {a_{3}}\ \ {\bf 2}{a_{3}}^{2}\ \ {b_{2}}\ \ $
 & 2&           176&        44352\\
${d_{6}}\ \ {b_{3}}^{3}\ \ {b_{2}}\ \ $
 & 2&           120&        30240\\
\end{tabular}
\vfil\eject

\begin{tabular}[t]{lrrrl}
\\type of diagram $\mskip 40mu$ 
   & $\mskip 30mu$ group 
   & $\mskip 30mu$ unit volume 
   & $\mskip 20mu$ total volume 
   & $\mskip 10mu$ \\ 
& & $\scriptstyle{[u=(16*16!)^{-1}]}$ & $\scriptstyle{[u=1891890^{-1}]}$ \\
\hline
\vspace{-0.1in} \\
${d_{5}}\ \ {b_{4}}^{2}\ \ {c_{4}}\ \ $
 & 2&            96&        24192\\
${d_{5}}\ \ {d_{4}}\ \ {\bf 2}{d_{4}}^{2}\ \ $
 & 4&           160&        20160\\
${a_{5}}\ \ {b_{4}}\ \ {c_{4}}\ \ {b_{3}}\ \ {\bf 2}{a_{1}}\ \ $
 & 1&           144&        72576\\
${a_{5}}\ \ {f_{4}}^{2}\ \ {b_{3}}\ \ {\bf 2}{a_{1}}\ \ $
 & 2&            96&        24192\\
\end{tabular}
\vfil\eject

\begin{tabular}[t]{lrrrl}
\\type of diagram $\mskip 40mu$ 
   & $\mskip 30mu$ group 
   & $\mskip 30mu$ unit volume 
   & $\mskip 20mu$ total volume 
   & $\mskip 10mu$ \\ 
& & $\scriptstyle{[u=(16*16!)^{-1}]}$ & $\scriptstyle{[u=1891890^{-1}]}$ \\
\hline
\vspace{-0.1in} \\
${f_{4}}^{3}\ \ {b_{3}}\ \ {a_{2}}\ \ $
 & 6&            39&         3276\\
${b_{4}}^{2}\ \ {c_{4}}\ \ {a_{2}}\ \ {b_{2}}\ \ {a_{1}}\ \ $
 & 2&           108&        27216\\
${b_{4}}\ \ {c_{4}}^{2}\ \ {b_{2}}^{2}\ \ {a_{1}}\ \ $
 & 2&           120&        30240\\
${b_{4}}\ \ {d_{4}}^{2}\ \ {b_{2}}^{2}\ \ {\bf 2}{a_{1}}\ \ $
 & 4&           288&        36288\\
\end{tabular}
\vfil\eject

\begin{tabular}[t]{lrrrl}
\\type of diagram $\mskip 40mu$ 
   & $\mskip 30mu$ group 
   & $\mskip 30mu$ unit volume 
   & $\mskip 20mu$ total volume 
   & $\mskip 10mu$ \\ 
& & $\scriptstyle{[u=(16*16!)^{-1}]}$ & $\scriptstyle{[u=1891890^{-1}]}$ \\
\hline
\vspace{-0.1in} \\
${\bf 2}{d_{4}}\ \ {a_{3}}\ \ {c_{3}}^{2}\ \ {a_{1}}^{2}\ \ {\bf 2}{a_{1}}^{2}\ \ $
 & 4&           256&        32256\\
${d_{4}}\ \ {b_{3}}^{2}\ \ {b_{2}}^{2}\ \ {a_{1}}^{3}\ \ $
 & 4&           352&        44352\\
${d_{4}}\ \ {b_{3}}\ \ {b_{2}}^{3}\ \ {a_{1}}^{3}\ \ {\bf 2}{a_{1}}\ \ $
 & 6&           448&        37632\\
${d_{4}}\ \ {c_{3}}\ \ {a_{1}}^{4}\ \ {\bf 2}{a_{1}}^{6}\ \ $
 & 24&           768&        16128\\
\end{tabular}
\end{center}
Total volume of 475 types (units of 1) $= 30270240 /1891890=16=\sqrt2^8$.
\vfil\eject}


\chapter*{Appendix B}
 
There are 18 hyperbolic Lie algebras of rank greater or equal to 7, falling
into 5 families, $AE_n$, $BE_n$, $CE_n$, $DE_n$, and $T_{r,s,t}$. We list all
of them and give multiplicities for a limited number of roots of small height.
For simply laced algebras, we compare to the global upper bound given in 
theorem 1.4. (Note that this theorem does not apply to the families $BE_n$ 
and $CE_n$.) Even though we can only give very few values they still provide 
some indication of the quality of the upper bounds. It may be observed that 
the multiplicities for roots of level 1 all conform to formula (6.8). The 
values $$p_{{\rm rank}-2}(1-{r^2\over2})$$ provide material for many 
intriguing conjectures in this context. They are printed for this purpose 
only. From the point of disproving any existing conjectures the table of 
$T_{4,3,3}$ may be the most interesting. It shows that $T_{4,3,3}$ contains 
roots of multiplicities both larger and smaller than $p_6$. All root 
multiplicities in this appendix were calculated by a program based
on the Peterson recursion formula (\cite{Kac90}, p.210).
 
The algebra $E_{10}=T_{7,3,2}$ has been studied extensively in \cite{KMW88} where
the root multiplicities for level 0,1, and 2 were determined explicitly.
The table for $E_{10}$ contains, for some  small negative norms $r^2$, the
explicit multiplicities as determined by \cite{KMW88} for levels 1 and 2, and the
global upper bound as provided by theorem 1.4. We see that this bound
provides a reasonably good approximation.
\vfil \eject

{\parindent=0pt \vsmall

\begin{minipage}[t]{2.2in}

\end{minipage}
\begin{minipage}[t]{2.7in}
\begin{tabular}[t]{|c|c|c|c|c|}
   $\mskip 150mu$
 & $\mskip 25mu$
 & $\mskip 25mu$
 & $\mskip 25mu$
 & $\mskip 25mu$ \vspace{-0.2in}\\
 & & & thm & \vspace{-0.08in}\\
  {root-coefficients} & norm & mult & 1.4 & $p_5$ \\ \hline
      {0, 1, 1, 1, 1, 1, 1} & 0 & 5 & 5 & 5\vspace{-0.01in}\\
     {1, 2, 2, 2, 2, 2, 2} & -2 & 20 & 21 & 20\vspace{-0.01in}\\
     {2, 4, 3, 2, 1, 2, 3} & -2 & 20 & 21 & 20\vspace{-0.01in}\\
     {2, 4, 3, 2, 2, 2, 3} & -4 & 65 & 71 & 65\vspace{-0.01in}\\
     {1, 3, 3, 3, 3, 3, 3} & -4 & 65 & 71 & 65\vspace{-0.01in}\\
     {2, 4, 3, 3, 3, 3, 3} & -6 & 189 & 217 & 190\vspace{-0.01in}\\
     {2, 5, 4, 3, 2, 3, 4} & -6 & 190 & 217 & 190\vspace{-0.01in}\\
     {1, 4, 4, 4, 4, 4, 4} & -6 & 190 & 217 & 190\vspace{-0.01in}\\
     {3, 6, 5, 4, 3, 2, 4} & -6 & 190 & 217 & 190\vspace{-0.01in}\\
     {2, 5, 4, 3, 3, 3, 4} & -8 & 502 & 603 & 506\vspace{-0.01in}\\
     {2, 4, 4, 4, 4, 4, 4} & -8 & 500 & 603 & 506 \\
\end{tabular}
\end{minipage}
\begin{minipage}[t]{2.7in}
\begin{tabular}[t]{|c|c|c|c|c|}
   $\mskip 150mu$
 & $\mskip 25mu$
 & $\mskip 25mu$
 & $\mskip 25mu$
 & $\mskip 25mu$ \vspace{-0.2in}\\
 & & & thm & \vspace{-0.08in}\\
     {root-coefficients} & norm & mult & 1.4 & $p_6$ \\ \hline
    {0, 1, 1, 1, 1, 1, 1, 1} & 0 & 6 & 6 & 6\vspace{-0.01in}\\
    {1, 2, 2, 2, 2, 2, 2, 2} & -2 & 27 & 28 & 27\vspace{-0.01in}\\
    {2, 4, 3, 2, 1, 1, 2, 3} & -2 & 27 & 28 & 27\vspace{-0.01in}\\
    {2, 4, 3, 2, 2, 2, 2, 3} & -4 & 97 & 105 & 98\vspace{-0.01in}\\
    {1, 3, 3, 3, 3, 3, 3, 3} & -4 & 98 & 105 & 98\vspace{-0.01in}\\
    {2, 4, 3, 3, 3, 3, 3, 3} & -6 & 310 & 350 & 315\vspace{-0.01in}\\
    {2, 5, 4, 3, 2, 2, 3, 4} & -6 & 309 & 350 & 315\vspace{-0.01in}\\
    {3, 6, 5, 4, 3, 2, 2, 4} & -6 & 309 & 350 & 315\vspace{-0.01in}\\
    {3, 6, 4, 2, 2, 3, 4, 5} & -6 & 309 & 350 & 315\vspace{-0.01in}\\
    {1, 4, 4, 4, 4, 4, 4, 4} & -6 & 315 & 350 & 315\vspace{-0.01in}\\
    {2, 5, 4, 3, 3, 3, 3, 4} & -8 & 894 & 1057 & 918 \\
\end{tabular}
\end{minipage}
\begin{minipage}[t]{2.7in}
\begin{tabular}[t]{|c|c|c|c|c|}
   $\mskip 150mu$
 & $\mskip 25mu$
 & $\mskip 25mu$
 & $\mskip 25mu$
 & $\mskip 25mu$ \vspace{-0.2in}\\
 & & & thm & \vspace{-0.08in}\\
     {root-coefficients} & norm & mult & 1.4 & $p_7$ \\ \hline
    {0, 1, 1, 1, 1, 1, 1, 1, 1} & 0 & 7 & 7 & 7\vspace{-0.01in}\\
    {1, 2, 2, 2, 2, 2, 2, 2, 2} & -2 & 35 & 36 & 35\vspace{-0.01in}\\
    {2, 4, 3, 2, 1, 1, 1, 2, 3} & -2 & 34 & 36 & 35\vspace{-0.01in}\\
    {2, 4, 3, 2, 2, 2, 2, 2, 3} & -4 & 136 & 148 & 140\vspace{-0.01in}\\
    {1, 3, 3, 3, 3, 3, 3, 3, 3} & -4 & 140 & 148 & 140\vspace{-0.01in}\\
    {2, 5, 4, 3, 2, 1, 2, 3, 4} & -4 & 132 & 148 & 140\vspace{-0.01in}\\
    {3, 6, 5, 4, 3, 2, 1, 2, 4} & -4 & 133 & 148 & 140\vspace{-0.01in}\\
    {3, 6, 4, 2, 1, 2, 3, 4, 5} & -4 & 133 & 148 & 140\vspace{-0.01in}\\
    {2, 5, 4, 3, 2, 2, 2, 3, 4} & -6 & 464 & 534 & 490\vspace{-0.01in}\\
    {2, 4, 3, 3, 3, 3, 3, 3, 3} & -6 & 475 & 534 & 490\vspace{-0.01in}\\
    {2, 5, 4, 3, 3, 3, 3, 3, 4} & -8 & 1464 & 1738 & 1547 \\
\end{tabular}
\end{minipage}
\begin{minipage}[t]{2.2in}
\begin{tabular}[t]{|c|c|c|}
   $\mskip 150mu$
 & $\mskip 25mu$
 & $\mskip 25mu$ \vspace{-0.2in}\\
     {root-coefficients} & norm & mult \\ \hline
{2, 2, 2, 2, 1, 1, 0} & 0 & 5\vspace{-0.01in}\\
{3, 3, 3, 3, 1, 2, 1} & -2 & 6\vspace{-0.01in}\\
{4, 4, 4, 4, 2, 2, 1} & -4 & 21\vspace{-0.01in}\\
{5, 5, 5, 5, 2, 3, 1} & -6 & 26\vspace{-0.01in}\\
{3, 4, 5, 6, 3, 4, 2} & -6 & 26\vspace{-0.01in}\\
{4, 4, 5, 6, 3, 4, 2} & -8 & 71\vspace{-0.01in}\\
{6, 6, 6, 6, 3, 3, 1} & -8 & 71\vspace{-0.01in}\\
{6, 6, 6, 6, 2, 4, 2} & -8 & 71\vspace{-0.01in}\\
{5, 5, 5, 6, 3, 4, 2} & -10 & 91\vspace{-0.01in}\\
{7, 7, 7, 7, 3, 4, 1} & -10 & 91\vspace{-0.01in}\\
{6, 6, 6, 6, 3, 4, 2} & -12 & 217 \\
\end{tabular}
\end{minipage}
\begin{minipage}[t]{2.2in}
\begin{tabular}[t]{|c|c|c|}
   $\mskip 150mu$
 & $\mskip 25mu$
 & $\mskip 25mu$ \vspace{-0.2in}\\
     {root-coefficients} & norm & mult \\ \hline
{2, 2, 2, 2, 2, 1, 1, 0} & 0 & 6\vspace{-0.01in}\\
{3, 3, 3, 3, 3, 1, 2, 1} & -2 & 7\vspace{-0.01in}\\
{4, 4, 4, 4, 4, 2, 2, 1} & -4 & 28\vspace{-0.01in}\\
{2, 3, 4, 5, 6, 3, 4, 2} & -4 & 28\vspace{-0.01in}\\
{3, 3, 4, 5, 6, 3, 4, 2} & -6 & 34\vspace{-0.01in}\\
{5, 5, 5, 5, 5, 2, 3, 1} & -6 & 34\vspace{-0.01in}\\
{4, 4, 4, 5, 6, 3, 4, 2} & -8 & 105\vspace{-0.01in}\\
{6, 6, 6, 6, 6, 3, 3, 1} & -8 & 105\vspace{-0.01in}\\
{6, 6, 6, 6, 6, 2, 4, 2} & -8 & 105\vspace{-0.01in}\\
{5, 5, 5, 5, 6, 3, 4, 2} & -10 & 131\vspace{-0.01in}\\
{6, 6, 6, 6, 6, 3, 4, 2} & -12 & 349 \\
\end{tabular}
\end{minipage}
\begin{minipage}[t]{2.2in}
\begin{tabular}[t]{|c|c|c|}
   $\mskip 150mu$
 & $\mskip 25mu$
 & $\mskip 25mu$ \vspace{-0.2in}\\
     {root-coefficients} & norm & mult \\ \hline
{2, 2, 2, 2, 2, 2, 1, 1, 0} & 0 & 7\vspace{-0.01in}\\
{3, 3, 3, 3, 3, 3, 1, 2, 1} & -2 & 8\vspace{-0.01in}\\
{1, 2, 3, 4, 5, 6, 3, 4, 2} & -2 & 8\vspace{-0.01in}\\
{4, 4, 4, 4, 4, 4, 2, 2, 1} & -4 & 36\vspace{-0.01in}\\
{2, 2, 3, 4, 5, 6, 3, 4, 2} & -4 & 36\vspace{-0.01in}\\
{3, 3, 3, 4, 5, 6, 3, 4, 2} & -6 & 42\vspace{-0.01in}\\
{5, 5, 5, 5, 5, 5, 2, 3, 1} & -6 & 43\vspace{-0.01in}\\
{4, 4, 4, 4, 5, 6, 3, 4, 2} & -8 & 147\vspace{-0.01in}\\
{6, 6, 6, 6, 6, 6, 3, 3, 1} & -8 & 148\vspace{-0.01in}\\
{6, 6, 6, 6, 6, 6, 2, 4, 2} & -8 & 148\vspace{-0.01in}\\
{5, 5, 5, 5, 5, 6, 3, 4, 2} & -10 & 179\\
\end{tabular}
\end{minipage}
\begin{minipage}[t]{2.2in}
\begin{tabular}[t]{|c|c|c|}
   $\mskip 150mu$
 & $\mskip 25mu$
 & $\mskip 25mu$ \vspace{-0.2in}\\
     {root-coefficients} & norm & mult \\ \hline
{2, 2, 2, 2, 2, 2, 2, 1, 1, 0} & 0 & 8\vspace{-0.01in}\\
{3, 3, 3, 3, 3, 3, 3, 1, 2, 1} & -2 & 9\vspace{-0.01in}\\
{1, 1, 2, 3, 4, 5, 6, 3, 4, 2} & -2 & 8\vspace{-0.01in}\\
{4, 4, 4, 4, 4, 4, 4, 2, 2, 1} & -4 & 45\vspace{-0.01in}\\
{2, 2, 2, 3, 4, 5, 6, 3, 4, 2} & -4 & 44 \\
\end{tabular}
\end{minipage}
\begin{minipage}[t]{2.2in}
\begin{tabular}[t]{|c|c|c|}
   $\mskip 150mu$
 & $\mskip 25mu$
 & $\mskip 25mu$ \vspace{-0.2in}\\
     {root-coefficients} & norm & mult \\ \hline
{0, 1, 2, 1, 2, 2, 1} & 0 & 4\vspace{-0.01in}\\
{0, 2, 4, 2, 4, 4, 2} & 0 & 5\vspace{-0.01in}\\
{1, 2, 4, 2, 4, 4, 2} & -2 & 15\vspace{-0.01in}\\
{2, 4, 6, 3, 5, 4, 2} & -4 & 44\vspace{-0.01in}\\
{1, 3, 6, 3, 6, 6, 3} & -4 & 44\vspace{-0.01in}\\
{2, 4, 6, 2, 6, 6, 3} & -4 & 50\vspace{-0.01in}\\
{2, 4, 6, 3, 6, 6, 3} & -6 & 122\vspace{-0.01in}\\
{1, 4, 8, 4, 8, 8, 4} & -6 & 121\vspace{-0.01in}\\
{2, 5, 8, 4, 7, 6, 3} & -8 & 304\vspace{-0.01in}\\
{2, 5, 8, 3, 8, 8, 4} & -8 & 304\vspace{-0.01in}\\
{2, 4, 8, 4, 8, 8, 4} & -8 & 311\vspace{-0.01in}\\
{3, 6, 9, 4, 8, 7, 3} & -10 & 725 \\
\end{tabular}
\end{minipage}
\begin{minipage}[t]{2.2in}
\begin{tabular}[t]{|c|c|c|}
   $\mskip 150mu$
 & $\mskip 25mu$
 & $\mskip 25mu$ \vspace{-0.2in}\\
     {root-coefficients} & norm & mult \\ \hline
{0, 1, 2, 1, 2, 2, 2, 1} & 0 & 5\vspace{-0.01in}\\
{0, 2, 4, 2, 4, 4, 4, 2} & 0 & 6\vspace{-0.01in}\\
{1, 2, 4, 2, 4, 4, 4, 2} & -2 & 21\vspace{-0.01in}\\
{2, 4, 6, 3, 5, 4, 3, 1} & -2 & 21\vspace{-0.01in}\\
{2, 4, 6, 3, 5, 4, 4, 2} & -4 & 71\vspace{-0.01in}\\
{1, 3, 6, 3, 6, 6, 6, 3} & -4 & 70\vspace{-0.01in}\\
{2, 4, 6, 2, 6, 6, 6, 3} & -4 & 77\vspace{-0.01in}\\
{2, 4, 6, 3, 6, 6, 6, 3} & -6 & 215\vspace{-0.01in}\\
{2, 5, 8, 4, 7, 6, 5, 2} & -6 & 218\vspace{-0.01in}\\
{1, 4, 8, 4, 8, 8, 8, 4} & -6 & 212\vspace{-0.01in}\\
{2, 5, 8, 4, 7, 6, 6, 3} & -8 & 596\vspace{-0.01in}\\
{3, 6, 9, 4, 8, 7, 6, 3} & -10 & 1555 \\
\end{tabular}
\end{minipage}
\begin{minipage}[t]{2.2in}
\begin{tabular}[t]{|c|c|c|}
   $\mskip 150mu$
 & $\mskip 25mu$
 & $\mskip 25mu$ \vspace{-0.2in}\\
     {root-coefficients} & norm & mult \\ \hline
{0, 1, 2, 1, 2, 2, 2, 2, 1} & 0 & 6\vspace{-0.01in}\\
{0, 2, 4, 2, 4, 4, 4, 4, 2} & 0 & 7\vspace{-0.01in}\\
{1, 2, 4, 2, 4, 4, 4, 4, 2} & -2 & 28\vspace{-0.01in}\\
{2, 4, 6, 3, 5, 4, 3, 2, 1} & -2 & 29\vspace{-0.01in}\\
{2, 4, 6, 3, 5, 4, 4, 4, 2} & -4 & 106\vspace{-0.01in}\\
{1, 3, 6, 3, 6, 6, 6, 6, 3} & -4 & 104\vspace{-0.01in}\\
{2, 4, 6, 2, 6, 6, 6, 6, 3} & -4 & 112\vspace{-0.01in}\\
{2, 4, 6, 3, 6, 6, 6, 6, 3} & -6 & 349\vspace{-0.01in}\\
{2, 5, 8, 4, 7, 6, 5, 4, 2} & -6 & 357\vspace{-0.01in}\\
{2, 5, 8, 4, 7, 6, 6, 6, 3} & -8 & 1058 \\
\end{tabular}
\end{minipage}
\begin{minipage}[t]{2.2in}
\vspace{.1in}
\begin{tabular}[t]{|c|c|c|}
   $\mskip 150mu$
 & $\mskip 25mu$
 & $\mskip 25mu$ \vspace{-0.2in}\\
     {root-coefficients} & norm & mult \\ \hline
    {0, 1, 2, 1, 2, 2, 2, 2, 2, 1} & 0 & 7\vspace{-0.01in}\\
    {0, 2, 4, 2, 4, 4, 4, 4, 4, 2} & 0 & 8\vspace{-0.01in}\\
    {1, 2, 4, 2, 4, 4, 4, 4, 4, 2} & -2 & 36\vspace{-0.01in}\\
    {2, 4, 6, 3, 5, 4, 3, 2, 2, 1} & -2 & 37 \\
\end{tabular}
\end{minipage}
\begin{minipage}[t]{2.7in}
\begin{tabular}[t]{|c|c|c|c|c|}
   $\mskip 150mu$
 & $\mskip 25mu$
 & $\mskip 25mu$
 & $\mskip 25mu$
 & $\mskip 25mu$ \vspace{-0.2in}\\
 & & & thm & \vspace{-0.08in}\\
     {root-coefficients} & norm & mult & 1.4 & $p_5$ \\ \hline
    {0, 1, 2, 1, 2, 1, 1} & 0 & 5 & 5 & 5\vspace{-0.01in}\\
    {1, 2, 4, 2, 4, 2, 2} & -2 & 20 & 21 & 20\vspace{-0.01in}\\
    {2, 4, 6, 3, 5, 2, 2} & -4 & 65 & 71 & 65\vspace{-0.01in}\\
    {1, 3, 6, 3, 6, 3, 3} & -4 & 65 & 71 & 65\vspace{-0.01in}\\
    {2, 4, 6, 2, 6, 3, 3} & -4 & 65 & 71 & 65\vspace{-0.01in}\\
    {2, 4, 6, 3, 6, 3, 3} & -6 & 190 & 217 & 190\vspace{-0.01in}\\
    {1, 4, 8, 4, 8, 4, 4} & -6 & 190 & 217 & 190\vspace{-0.01in}\\
    {2, 5, 8, 4, 7, 3, 3} & -8 & 506 & 603 & 506\vspace{-0.01in}\\
    {2, 5, 8, 3, 8, 4, 4} & -8 & 505 & 603 & 506\vspace{-0.01in}\\
    {2, 4, 8, 4, 8, 4, 4} & -8 & 505 & 603 & 506\vspace{-0.01in}\\
    {2, 5, 8, 4, 8, 4, 4} & -10 & 1263 & 1574 & 1265 \\
\end{tabular}
\end{minipage}
\begin{minipage}[t]{2.7in}
\begin{tabular}[t]{|c|c|c|c|c|}
   $\mskip 150mu$
 & $\mskip 25mu$
 & $\mskip 25mu$
 & $\mskip 25mu$
 & $\mskip 25mu$ \vspace{-0.2in}\\
 & & & thm & \vspace{-0.08in}\\
     {root-coefficients} & norm & mult & 1.4 & $p_6$ \\ \hline
    {0, 1, 2, 1, 2, 2, 1, 1} & 0 & 6 & 6 & 6\vspace{-0.01in}\\
    {1, 2, 4, 2, 4, 4, 2, 2} & -2 & 27 & 28 & 27\vspace{-0.01in}\\
    {2, 4, 6, 3, 5, 4, 2, 1} & -2 & 27 & 28 & 27\vspace{-0.01in}\\
    {2, 4, 6, 3, 5, 4, 1, 2} & -2 & 27 & 28 & 27\vspace{-0.01in}\\
    {2, 4, 6, 3, 5, 4, 2, 2} & -4 & 98 & 105 & 98\vspace{-0.01in}\\
    {1, 3, 6, 3, 6, 6, 3, 3} & -4 & 98 & 105 & 98\vspace{-0.01in}\\
    {2, 4, 6, 2, 6, 6, 3, 3} & -4 & 98 & 105 & 98\vspace{-0.01in}\\
    {2, 4, 6, 3, 6, 6, 3, 3} & -6 & 314 & 350 & 315\vspace{-0.01in}\\
    {2, 5, 8, 4, 7, 6, 3, 2} & -6 & 315 & 350 & 315\vspace{-0.01in}\\
    {2, 5, 8, 4, 7, 6, 2, 3} & -6 & 315 & 350 & 315\vspace{-0.01in}\\
    {2, 5, 8, 4, 7, 6, 3, 3} & -8 & 914 & 1057 & 918 \\
\end{tabular}
\end{minipage}
\begin{minipage}[t]{2.7in}
\begin{tabular}[t]{|c|c|c|c|c|}
   $\mskip 150mu$
 & $\mskip 25mu$
 & $\mskip 25mu$
 & $\mskip 25mu$
 & $\mskip 25mu$ \vspace{-0.2in}\\
 & & & thm & \vspace{-0.08in}\\
     {root-coefficients} & norm & mult & 1.4 & $p_7$ \\ \hline
    {0, 1, 2, 1, 2, 2, 2, 1, 1} & 0 & 7 & 7 & 7\vspace{-0.01in}\\
    {1, 2, 4, 2, 4, 4, 4, 2, 2} & -2 & 35 & 36 & 35\vspace{-0.01in}\\
    {2, 4, 6, 3, 5, 4, 3, 1, 1} & -2 & 35 & 36 & 35\vspace{-0.01in}\\
    {2, 4, 6, 3, 5, 4, 4, 2, 2} & -4 & 139 & 148 & 140\vspace{-0.01in}\\
    {1, 3, 6, 3, 6, 6, 6, 3, 3} & -4 & 140 & 148 & 140\vspace{-0.01in}\\
    {2, 4, 6, 2, 6, 6, 6, 3, 3} & -4 & 140 & 148 & 140\vspace{-0.01in}\\
    {2, 4, 6, 3, 6, 6, 6, 3, 3} & -6 & 485 & 534 & 490\vspace{-0.01in}\\
    {2, 5, 8, 4, 7, 6, 5, 2, 2} & -6 & 484 & 534 & 490\vspace{-0.01in}\\
    {2, 5, 8, 4, 7, 6, 6, 3, 3} & -8 & 1522 & 1738 & 1547 \\
\end{tabular}
\end{minipage}
\begin{minipage}[t]{2.7in}
\begin{tabular}[t]{|c|c|c|c|c|}
   $\mskip 150mu$
 & $\mskip 25mu$
 & $\mskip 25mu$
 & $\mskip 25mu$
 & $\mskip 25mu$ \vspace{-0.2in}\\
 & & & thm & \vspace{-0.08in}\\
     {root-coefficients} & norm & mult & 1.4 & $p_8$ \\ \hline
   {0, 1, 2, 1, 2, 2, 2, 2, 1, 1} & 0 & 8 & 8 & 8\vspace{-0.01in}\\
   {1, 2, 4, 2, 4, 4, 4, 4, 2, 2} & -2 & 44 & 45 & 44\vspace{-0.01in}\\
   {2, 4, 6, 3, 5, 4, 3, 2, 1, 1} & -2 & 43 & 45 & 44\vspace{-0.01in}\\
   {1, 3, 6, 3, 6, 6, 6, 6, 3, 3} & -4 & 192 & 201 & 192 \\
\end{tabular}
\end{minipage}
\begin{minipage}[t]{2.7in}
\begin{tabular}[t]{|c|c|c|c|c|}
   $\mskip 132mu$
 & $\mskip 25mu$
 & $\mskip 32mu$
 & $\mskip 32mu$
 & $\mskip 32mu$ \vspace{-0.2in}\\
 & & & thm & \vspace{-0.08in}\\
     {root-coefficients} & norm & mult & 1.4 & $p_6$ \\ \hline
    {0, 1, 2, 3, 2, 1, 2, 1} & 0 & 5 & 5 & 5\vspace{-0.01in}\\
    {1, 2, 4, 6, 4, 2, 4, 2} & -2 & 27 & 28 & 27\vspace{-0.01in}\\
    {2, 4, 6, 8, 5, 2, 5, 2} & -4 & 98 & 105 & 98\vspace{-0.01in}\\
    {1, 3, 6, 9, 6, 3, 6, 3} & -4 & 98 & 105 & 98\vspace{-0.01in}\\
    {2, 4, 6, 9, 6, 3, 6, 3} & -6 & 315 & 350 & 315\vspace{-0.01in}\\
    {1, 4, 8, 12, 8, 4, 8, 4} & -6 & 315 & 350 & 315\vspace{-0.01in}\\
    {2, 5, 8, 11, 7, 3, 7, 3} & -8 & 918 & 1057 & 918\vspace{-0.01in}\\
    {2, 4, 8, 12, 8, 4, 8, 4} & -8 & 917 & 1057 & 918\vspace{-0.01in}\\
    {3, 6, 9, 12, 8, 4, 7, 2} & -6 & 316 & 350 & 315\vspace{-0.01in}\\
    {3, 6, 9, 12, 7, 2, 8, 4} & -6 & 316 & 350 & 315\vspace{-0.01in}\\
    {2, 5, 8, 12, 8, 4, 8, 4} & -10 & 2491 & 2975 & 2492\vspace{-0.01in}\\
    {3, 6, 9, 12, 8, 4, 7, 3} & -10 & 2493 & 2975 & 2492\vspace{-0.01in}\\
    {3, 6, 9, 12, 7, 3, 8, 4} & -10 & 2493 & 2975 & 2492\vspace{-0.01in}\\
    {3, 6, 9, 12, 8, 4, 8, 4} & -12 & 6372 & 7883 & 6372\vspace{-0.01in}\\
    {2, 6, 10, 14, 9, 4, 9, 4} & -12 & 6368 & 7883 & 6372\vspace{-0.01in}\\
    {3, 6, 10, 14, 9, 4, 9, 4} & -14 & 15524 & 19900 & 15525\\
\end{tabular}
\end{minipage}
\begin{minipage}[t]{2.7in}
\begin{tabular}[t]{|c|c|c|c|c|}
   $\mskip 132mu$
 & $\mskip 25mu$
 & $\mskip 32mu$
 & $\mskip 32mu$
 & $\mskip 32mu$ \vspace{-0.2in}\\
 & & & thm & \vspace{-0.08in}\\
     {root-coefficients} & norm & mult & 1.4 & $p_7$ \\ \hline
   {0, 1, 2, 3, 4, 2, 3, 2, 1} & 0 & 7 & 7 & 7\vspace{-0.01in}\\
   {1, 2, 4, 6, 8, 4, 6, 4, 2} & -2 & 35 & 36 & 35\vspace{-0.01in}\\
   {2, 4, 6, 8, 10, 5, 7, 4, 1} & -2 & 35 & 36 & 35\vspace{-0.01in}\\
   {2, 4, 6, 8, 10, 5, 7, 4, 2} & -4 & 140 & 148 & 140\vspace{-0.01in}\\
   {1, 3, 6, 9, 12, 6, 9, 6, 3} & -4 & 140 & 148 & 140 \\
\end{tabular}
\end{minipage}
\begin{minipage}[t]{2.1in}
\begin{tabular}[t]{|r|r|r|r|}
   $\mskip 44mu$
 & $\mskip 44mu$
 & $\mskip 44mu$
 & $\mskip 44mu$ \vspace{-0.2in}\\
\multicolumn{1}{|c|} {$\mskip 44mu$} &
\multicolumn{1}{|c|} {level 0,1} &
\multicolumn{1}{|c|} {level 2} &
\vspace{-0.08in} \\
\multicolumn{1}{|c|} {norm} &
\multicolumn{1}{|c|} {[KMW88]} &
\multicolumn{1}{|c|} {[KMW88]} &
\multicolumn{1}{|c|} {thm 1.4} \\ \hline
 0  &   8   &  -   &  8   \\
-2  &  44   &  44   &  45   \\
-4  &  192   &  192   &  201   \\
-6  &  726   &  727   &  780  \\
-8  &  2464   &  2472   &  2718   \\
-10  &  7704   &  7747   &  8730   \\
-12  &  22528   &  22712   &  26226   \\
-14  &  62337   &  63020   &  74556  \\
-16  &  164560   &  166840   &  202180  \\ \hline
\end{tabular}
\end{minipage}
} 


\bibliographystyle{amsalpha}

\begin{thebibliography}{A}

\bibitem [Apo76]{Apo76}
           Apostol, T.M.: \textit{ Modular Functions and Dirichlet Series in Number
           Theory}, Springer Verlag (Graduate Texts in Mathematics 41), New York,
           Berlin, Heidelberg, 1976 
\bibitem [BM79]{BM79}
           Berman, S., Moody, R.V.: \textit {Multiplicities in Lie Algebras}, Proc.
           Amer. Math. Soc. 76 (1979), 223--228
\bibitem [Bor85]{Bor85}
           Borcherds, R.E.: \textit {The Leech Lattice}, Proc. R. Soc. Lond. A 398
           (1985), 365--376
\bibitem [Bor86]{Bor86}
           Borcherds, R.E.: \textit {Vertex Algebras, Kac-Moody Algebras, and the
           Monster}, Proc. Natl. Acad. Sci. USA 83 (1986), 3068-3071 
\bibitem [Bor88]{Bor88}
           Borcherds, R.E.: \textit {Generalized Kac-Moody algebras},
           J. Algebra 115 (1988), 501-512
\bibitem [Bor90a]{Bor90a}
           Borcherds, R.E.: \textit {Lattices like the Leech Lattice}, J. Algebra 130
           (1990), 219-234               
\bibitem [Bor90b]{Bor90b}
           Borcherds, R.E.: \textit {The Monster Lie algebra}, Adv. Math. 83 (1990),
           30-47                         
\bibitem [Bor91]{Bor91}
           Borcherds, R.E.: \textit {Central extensions of generalized Kac-Moody 
           algebras}, J. Algebra 140 (1991), 330-335
\bibitem [Bor92]{Bor92}
           Borcherds, R.E.: \textit {Monstrous Moonshine and Monstrous Lie
           Superalgebras}, Invent. Math. 109 (1992), 405-444 
\bibitem [CE56]{CE56}
           Cartan, H., Eilenberg, S.: \textit {Homological Algebra}, Princeton:
           Princeton University Press 1956
\bibitem [Con85]{Con85}
           Conway, J.H., Curtis, R.T., Norton, S.P., Parker, R.A., Wilson,
           R.A.: \textit {ATLAS of Finite Groups}, Oxford Univ. Press, 1985
\bibitem [CS83]{CS83}
           Conway, J.H., Sloane, N.J.A.: \textit {The Coxeter-Todd Lattice, the Mitchell
           Group, and Related Sphere Packings}, Proc. Camb. Phil. Soc. 93
           (1983), 421-440                                
\bibitem [CS88]{CS88}
           Conway, J.H., Sloane, N.J.A.: \textit {Sphere Packings, Lattices and Groups},
           Springer Verlag (Grundlehren d. math. Wiss. 290), New York, 
           Berlin, Heidelberg, 1988                       
\bibitem [FF83]{FF83}
           Feingold, A.J., Frenkel, I.B.: \textit {A Hyperbolic Kac-Moody Algebra
           and the Theory of Siegel Modular Forms of Genus 2}, Math. Ann. 263
           (1983), 87-144                                 
\bibitem [FLM88]{FLM88}
           Frenkel, I.B., Lepowsky, J., and Meurman, A.: \textit {Vertex Operator
           Algebras and the Monster}, MA Academic Press, Boston, 1988
\bibitem [GL76]{GL76}
           Garland, H., Lepowsky, J.: \textit {Lie algebra homology and the
           Macdonald-Kac formulas}, Invent. Math. 34 (1976), 37-76
\bibitem [GN97]{GN97}
           Gebert, R.W., Nicolai, H.: \textit {On the imaginary simple roots of the
           Borcherds algebra {\fraktur g}$_{II_{9,1}}$}, preprint 
           IASSNS-HEP-97-53, AEI-937, 1997 
\bibitem [Jac85]{Jac85}
           Jacobson, N.: \textit {Basic Algebra I}, 2nd ed., Freeman, New York, 1985
\bibitem [Jan95]{Jan95}
           Jansen, C., Lux, K., Parker, R., Wilson, R.:
           \textit {An Atlas of Brauer Characters}, Oxford University Press, 1995
\bibitem [Jur96]{Jur96}
           Jurisich, E.: \textit {An Exposition of Generalized Kac-Moody algebras},
           Contemporary Maths 194 (1996) 
\bibitem [Jur98]{Jur98}
           Jurisich, E.: \textit {Generalized Kac-Moody Lie algebras, free Lie
           algebras and the structure of the Monster Lie algebra},
           J. Pure Appl. Algebra 126 (1998), 233-266              
\bibitem [Kac90]{Kac90}
           Kac, V.G.: \textit {Infinite Dimensional Lie Algebras}, 3rd ed., Cambridge
           University Press, 1990                                 
\bibitem [KP83]{KP83}
           Kac, V.G., Peterson, D.H.: \textit {Regular Functions on Certain
           Infinite-dimensional Groups}, in: Arithmetic and Geometry, 141-166,
           Progress in Math. 36, Birkh\"auser, Boston, 1983
\bibitem [KMW88]{KMW88}
           Kac, V.G., Moody, R.V., and Wakimoto, M.: \textit {On $E_{10}$}, Differential
           Methods in Theoretical Physics, Bleuler, K., Werner, M. (eds.),
           Kluwer Academic Publishers, 109-128             
\bibitem [Kob84]{Kob84}
           Koblitz, N.: \textit {Introduction to Elliptic Curves and Modular Forms},
           Springer Verlag (Graduate Texts in Mathematics), New York, Berlin,
           Heidelberg, 1984                                
\bibitem [Rad29]{Rad29}
           Rademacher, H.: \textit {\"Uber die Erzeugenden von Kongruenzuntergruppen der
           Modulgruppe}, Abh. Math. Sem. Hamburg 7 (1929), 134-148
\bibitem [Shi71]{Shi71}
           Shimura, G.: \textit {Introduction to the Arithmetic Theory of Automorphic
           Functions}, Publ. Math. Soc. Japan 11, 1971            
\bibitem [Wan91]{Wan91}
           Wan Zhe-xian: \textit {Introduction to Kac-Moody Algebra}, World Scientific
           Publishing, Singapore, 1991          

\end{thebibliography}

\chapter*{Notation}

The following list comprises all non-standard notations used in this work. It
refers to the chapter where the notation is introduced and also gives a brief
indication of its meaning.
Note that we use the standard set in \cite{Kac90} for finite and affine Lie
algebras. They are therefore not listed below. \flexskip

{\parindent=0pt
\begin{tabular} {lll} 
\vspace{-0.1in} \\
 $AE_n$& B& hyperbolic Lie algebra \\
 $(A,j)$&2.1 &element of the metaplectic group\\
 Aut$(\Lambda^\sigma)$& 5.3& group of automorphisms of $\Lambda^\sigma$ that
                              fix the origin\\
 Aut$(H)$& 5.3& group of automorphisms fixing a hole $H$ of ${\cal R}$\\
 $BE_n$& B& hyperbolic Lie algebra \\
 $C$& 1.1.1, 5.2& Cartan matrix\\
 ${\cal C}$& 4.1& 24-dimensional Golay-code\\
 $CE_n$& B& hyperbolic Lie algebra \\
 $D$& 5.2& diagonal matrix\\
 $DE_n$& B& hyperbolic Lie algebra\\
 $d$& 5.2& factor relating $\delta^\vee={1\over d}\nu^{-1}\delta$\\
 $E_r$& 1.5& part of $M_\Lambda$ of grade $r\in II_{25,1}$ \\
 $E_{10}$& B& hyperbolic Lie algebra\\
 $e^r$& 1.3& element of the central extension of the Leech lattice\\
 $F$& 2.1& the Fricke involution\\
 $G(C)$ & 1.1.1& GKM generated from generalized Cartan matrix $C$ \\
 $G^e(C)$ & 1.1.2& GKM centrally extended by degree derivations\\
 ${\cal G}_N$& 1.6 & the GKM constructed in theorem 1.6\\
 $GO^\epsilon_M(N)$& 4.1& the general orthogonal group over
                  $(\ZZ_N)^M$ of Witt type $\epsilon$\\
 $H$& 5.2& subset of $\cal R$ that constitutes the vertices of a hole\\
 $\langle H \rangle$& 5.2& hole in its spatial meaning\\
 ${\cal H}$& 2.1& the upper half complex plane\\
 $II_{n,1}$& 1.4& $(n+1)$-dimensional even unimodular Lorentzian lattice\\
 $K_{12}$& 4.0& the Coxeter-Todd lattice\\
 $L$& 4.1& ${\Lambda^\sigma}^\perp$\\
 $L$& 1.6, 5.1& $\Lambda^\sigma\oplus II_{1,1}$\\
 $L^*$& 3.1& the dual lattice of $L$\\
 $L^\perp$& 3.1& the orthogonal complement of $L$ with respect to $\Lambda$\\
 ${L^*}^+$& 1.5& the positive roots in $L^*$\\
 $L_n$& 1.3& operators forming the Virasoro algebra\\        
 $M$& 1.6& =$24\over N+1$\\
 $Mp$& 2.1& the metaplectic group, double cover of $\Gamma$\\
 $M_\Lambda$& 1.4& the fake monster Lie algebra\\
 $N$& 1.6& order of the automorphism $\sigma$, that is 2, 3, 5, 7, 11, or 23\\
 $N\Delta$& 5.3& Dynkin diagram entirely of long roots \\
\end{tabular} 

\begin{tabular} {lll}
 $O^\epsilon_M(N)$& 4.1& the `generically simple' orthogonal group over
                  $(\ZZ_N)^M$ of \\
 & & Witt type $\epsilon$\\
 $P^n$& 1.3& physical space as subspace of a vertex algebra\\        
 $p_\sigma$& 1.7& coefficients of the $q$ expansion of $1/\eta_\sigma$\\
 $\QQ$& & the rational numbers\\
 $Q$& 1.3& vertex operator\\        
 $q$& 2.1& element of the unit disc, $q=e^{2\pi i\tau}$\\
 $\RR$&& the real numbers\\
 ${\cal R}$& 5.1& elements of $\Lambda^\sigma$ and ${\Lambda^\sigma}^*$
                   representing the simple roots of \\
 &&${\cal G}_N$\\
 ${\cal R}_{dual}$& 5.1& $={\cal R}\cap {\Lambda^\sigma}^*$ \\
 ${\cal R}_{fix}$& 5.1& $={\cal R}\cap \Lambda^\sigma$ \\
 $r_{\cal R}$& 5.2& radius function\\
 $(r,L)$& 1.5, 5.1& greatest common divisor of $(r,v), v\in L$\\
 $S$& 2.1& one of the generators of $\Gamma$\\
 $S_n$& 1.4& the part of the symmetric algebra $S$ of $\ZZ$ grading $n$\\
 $SO^\epsilon_M(N)$& 4.1& the special orthogonal group over
                  $(\ZZ_N)^M$ of Witt type $\epsilon$\\
 $T_{r,s,t}$& B& hyperbolic Lie algebra\\
 $U(C)$ & 1.1.4& universal GKM generated from generalized Cartan \\
 &&matrix $C$\\
 $V(L)$& 1.3& the vertex algebra of the lattice $L$\\        
 $(V_k,*)$& 2.1& $V_k$ is the typical generator of $\Gamma_0(N)$, \\
          &    & $*$ stands for the unspecified branch of $j$ (see $(A,j)$)\\
 $W^\sigma$& 1.5& the Weyl group of the twisted algebra\\
 $\ZZ$&& the integers\\
 $\ZZ_N$& 4.1& the finite field $\ZZ/N\ZZ$\\
 $\Gamma$& 2.1& the modular group\\
 $\Gamma_0(N)$& 2.1& subgroup of $\Gamma$\\
 $\Delta$& 5.3& arbitrary Dynkin diagram\\
 $\Delta(H)$& 5.2& Dynkin diagram associated to hole $H$\\
 $\delta,\delta^\vee$& 5.2& unique norm 0 vector of affine Lie algebra\\
 $\eta$& 2.2& the Dedekind eta-function\\
 $\Theta,\Theta_r$& 4.4& modular forms, right hand side of eqn. (1.25)\\
 $\theta$& 3.1& the theta-function of a lattice\\
 $\Lambda$& 1.4& the Leech lattice\\
 $\bigwedge$& 1.2& the wedge product\\
 $\Lambda_{16}$& 4.0& the Barnes-Wall lattice\\
 $(\lambda,m,n)$& 1.4& typical element of $II_{25,1}$\\
 $\nu$& 5.3& the isomorphism between Cartan subalgebra and its dual \\
 $\pi_L,\pi_{\Lambda^\sigma}$& 4.1& projections to the span of the relevant
lattice\\
 $\pi_1,\pi_2$& 4.3& $\pi_L=\pi_1 \circ \pi_2$\\
 $\rho,\rho^\vee$& 1.1.3, 5.2& the Weyl vector of a Lie algebra\\
 $\rho_M,\tilde\rho_M$& 4.2& number of representations as sum of squares\\
 $\sigma$& 1.5& automorphism of the Leech lattice, of prime order\\
 $\sigma$& 1.6& automorphism of the Leech lattice, cycle shape $1^MN^M$\\
 $\phi_{\sigma,V}$& 1.7& $=\sum {\rm Tr}(\sigma,V_n)q^n$ \\
 $\psi_j$& 2.2& modular form, $\psi_j(\tau)=\eta({\tau+j\over N}+j)$\\
 $[ {^.}]$& & space spanned by any basis\\
 $( {^.}, {^.})$ & 1.1.3& bilinear form\\
\end{tabular}}

\end{document}